\documentclass[12pt,a4paper]{book}

\def\P{{\mathbf P}}
\usepackage{index}
\usepackage{theorem}
\input{amssym.def}

\newindex{not}{adx}{and}{Notations}
\newindex{sub}{idx}{ind}{Index}

\newtheorem{theorem}{Theorem}[section]
\newtheorem{lemma}[theorem]{Lemma}
\newtheorem{proposition}[theorem]{Proposition}
\newtheorem{corollary}[theorem]{Corollary}
\newtheorem{definition}[theorem]{Definition}
\newtheorem{remarkk}[theorem]{Remark}
\newtheorem{example}[theorem]{Example}
\newcommand{\be}{\begin{equation}}
\newcommand{\ee}{\end{equation}}
\newcommand{\bea}{\begin{eqnarray}}
\newcommand{\eea}{\end{eqnarray}}
\newcommand{\beq}[1]{\begin{equation}\label{#1}}
\newcommand{\eeq}{\end{equation}}
\newcommand{\beqn}[1]{\begin{eqnarray}\label{#1}}
\newcommand{\eeqn}{\end{eqnarray}}
\newcommand{\beaa}{\begin{eqnarray*}}
\newcommand{\eeaa}{\end{eqnarray*}}

\newcommand{\La}{\Lambda}
\newcommand{\la}{\lambda}
\newcommand{\eps}{\varepsilon}

\newcommand{\Om}{\Omega}
\def\t{\mbox{trace}}
\newcommand{\veps}{\varepsilon}

\newcommand{\B}{{\mathcal B}}
\newcommand{\F}{{\mathcal F}}
\newcommand{\E}{{\mathcal E}}
\newcommand{\K}{{\mathcal K}}
\renewcommand{\L}{{\mathcal L}}
\renewcommand{\H}{{\mathcal H}}
\newcommand{\Ad}{\rm {Ad}}
\newcommand{\ad}{\rm {ad}}
\newcommand{\tad}{\widetilde{\ad}}
\newcommand{\tA}{\widetilde{\Ad}}
\newcommand{\G}{{\mathcal G}}
\newcommand{\ka}{\kappa}

\newcommand{\reals}{{\rm I\!R}}

\newcommand{\Z}{{\Bbb Z}}
\newcommand{\NN}{{\rm I\!N}}
\newcommand{\DD}{{\rm I\!D}}
\newcommand{\Dll}{\DD^{loc}}
\newcommand{\ga}{\gamma}
\newcommand{\LL}{{\rm I\!L}}

\def\complex{\mathop{\raise .45ex\hbox{${\bf\scriptstyle{|}}$}
     \kern -0.40em {\rm \textstyle{C}}}\nolimits}
\def\hilbert{\mathop{\raise .21ex\hbox{$\bigcirc$}}\kern -1.005em {\rm\textstyle{H}}} 

\newcommand{\calB}{{\cal B}}
\newcommand{\calC}{{\cal C}}

\newcommand{\calF}{{\cal F}}

\newcommand{\calH}{{\cal H}}
\newcommand{\calL}{{\cal L}}
\newcommand{\calM}{{\cal M}}

\newcommand{\calS}{{\cal S}}

\newcommand{\calY}{{\cal Y}}

\newcommand{\al}{{\alpha}}

\newcommand{\R}{{\rm I\!R}}

\newcommand{\N}{{\rm I\!N}}

\renewcommand{\varrho}{{\ell}}
\newcommand{\dett}{{\textstyle{\det_2}}}

\newcommand{\op}{{\mbox{\footnotesize{\rm op}}}}

\newcommand{\won}{{\mbox{\bf 1}}}
\newcommand{\Nl}{\nabla^{loc}}
\newcommand{\dl}{\delta^{loc}}

\newcommand{\trace}{{\,\,\rm trace\,\,}}
\newcommand{\Dom}{{\rm Dom}}

\newcommand{\proof}{\noindent {\bf Proof:\ }}
\newcommand{\remark}{\noindent {\bf Remark:\ }}
\newcommand{\remarks}{\noindent {\bf Remarks:\ }}

\newcommand{\ra}{ \rightarrow }

\newcommand{\half}{\frac{1}{2}\:}

\def\ge{\geq}

\def\squarebox#1{\hbox to #1{\hfill\vbox to #1{\vfill}}}
\newcommand{\qed}{\hspace*{\fill}
           \vbox{\hrule\hbox{\vrule\squarebox{.667em}\vrule}\hrule}\bigskip}
\begin{document}

\pagestyle{myheadings}

\thispagestyle{empty}
\begin{titlepage}
\title{\bf\Huge  Analysis on Wiener Space and Applications}
\author{A. S. \"Ust\"unel}
\date{}
\end{titlepage}

\maketitle

\newpage
\noindent{\huge\bf Introduction\\}

 The aim of this book is  to give a   rigorous introduction for the
 graduate students  to  Analysis  on  Wiener space, a subject which
 has grown up very quickly these recent 
years under the new  impulse of the Stochastic Calculus of Variations of Paul
Malliavin (cf. \cite{Mall}). A  portion  of the  material exposed
is our own research, in particular, with Moshe Zakai and Denis
Feyel for the rest we have used  the works listed  in  the  bibliography.

The origin of this  book goes back to a series of  seminars  that I
had given in Bilkent University of Ankara  in the summer of 1987 and
also  during the spring and some portion  of the summer of 1993 at 
the Mathematics Institute of Oslo University and a graduate course
dispensed at the University of Paris VI.  An initial and rather naive  version of
these notes  has  been published in Lecture Notes in Mathematics
series of Springer at 1995. Since then we have assisted to  a very quick development
and progress of  the subject  in several directions. In particular,
its use has been remarked by mathematical economists. Consequently I
have decided to write a more complete text with additional
contemporary  applications to illustrate the  strength and the
applicability of the subject. Several new results like  the logarithmic
Sobolev inequalities, applications of the capacity theory to the
local and global differentiability of Wiener functionals,
probabilistic notions of the convexity and log-concavity, the Monge
and the 
Monge-Kantorovitch measure transportation problems in the infinite
dimensional setting  and the analysis on the path space of a compact Lie group are added.

Although some concepts are given  in the first chapter, I assumed
that the students had already acquired the notions of stochastic
calculus with semimartingales, Brownian motion and some rudiments
of the theory of Markov processes.

 The second  chapter deals with the definition of the  (so-called)
Gross-Sobolev derivative and the Ornstein-Uhlenbeck operator which
are indispensable  tools of the analysis on Wiener space. In  the
third chapter we begin the proof of the Meyer inequalities, for
which the hypercontractivity property of the Ornstein-Uhlenbeck
semi-group is needed. We expose this last topic in the fourth
chapter and give the classical proof of the logarithmic Sobolev
inequality of L. Gross  for the Wiener measure. In chapter V, we  
 complete the proof of Meyer inequalities and study the distribution
 spaces which are defined via the Ornstein-Uhlenbeck operator. In
 particular we show that the derivative and divergence operators
 extend continuously to distribution spaces. In the appendix we
 indicate how one can transfer all these results to arbitrary abstract
 Wiener spaces using the notion of time associated to a continuous
 resolution of identity of the underlying Cameron-Martin space.

Chapter VI begins with an extension of Clark's formula to the
distributions defined in the preceding chapter. This formula  is
applied to  prove the classical $0-1$-law and as an application of the
latter, we prove the positivity improving property of the
Ornstein-Uhlenbeck semigroup. We then show that the functional
composition of  a non-degenerate
Wiener functional with values in $\R^n$,  (in the sense of Malliavin)
with a real-valued smooth function on $\R^n$ can be extended when the
latter is a tempered distribution if we look at to the result as a
distribution on the Wiener space. This result contains the fact that
the probability density of a non-degenerate functional is not only
$C^\infty$ but also it is rapidly decreasing. This observation is then
applied to prove the regularity of the solutions of Zakai equation of
the nonlinear filtering and to an extension of the Ito formula to the
space of tempered distributions with non-degenerate Ito processes.
We complete this chapter with two non-standart applications of Clark's
formula, the first concerns the equivalence between the independence
of two measurable sets and the orthogonality of the corresponding
kernels of their Ito-Clark representation and the latter is another
proof of the logarithmic Sobolev inequality via Clark's formula.

Chapter VII begins with the characterization of positive (Meyer)  distributions
as Radon measures and an application of this result to local
times. Using capacities defined with respect to the Ornstein-Uhlenbeck
process, we prove also a stronger version of
the $0-1$-law alraedy exposed in Chapter VI: it  says that any
$H$-invariant subset of the Wiener space or its complement has zero 
$C_{r,1}$-capacity. This result is then used that the $H$- gauge
functionals of measurable sets are finite quasi-everywhere instead of
almost everywhere. We define also there the local Sobolev spaces,
which is  a useful  notion when we study the problems where the
integrability is not a priori obvious. We show how to patch them
together to obtain global functionals. Finally we give a short section
about the distribution spaces defined with the second quantization of
a general ``elliptic'' operator, and as an example show that the
action of a shift define a distribution in this sense.

In chapter eight  we study the independence of some
Wiener functionals with the previously developed tools. 

The ninth
chapter  is devoted to some series of moment inequalities which
are important in   applications like large
 deviations, stochastic differential equations, etc. In  the tenth
 chapter we expose the contractive version of  Ramer's theorem as
 another example of the applications of moment inequalities developed
 in the preceding chapter and as an application we show  the validity
 of the logarithmic Sobolev inequality under this perturbated
 measures. Chapter XI deals with a rather new notion of convexity and
concavity which is quite appropriate for the equivalence classes of
Wiener functionals.  We believe that it will have important
applications in the field of convex analysis and financial
mathematics.
Chapter XII can be regarded as  an immediate application of Chapter
XI, where we study the problem of G. Monge and its generalization,
called the Monge-Kantorovitch\footnote{Another spelling is
  ''Kantorovich''.} measure transportation problem {\bf {for general
measures}} with a singular quadratic cost function, namely the square of
the Cameron-Martin norm. Later we study in detail when the initial
measure is the Wiener measure.

 The last chapter is devoted to construct a similar Sobolev
 analysis on the path space over a compact Lie group, which is the
 simplest non-linear situation. This problem has been studied in the
 more general case of compact Riemannian manifolds
 (cf. \cite{Mall-1}, \cite{M-M1}), however, I think that the case of Lie
 groups, as an intermediate step to clarify the ideas, is quite useful.

\hspace*{9cm}Ali S\"uleyman  \"Ust\"unel

\tableofcontents
\newpage
\pagenumbering{arabic}
\setcounter{page}{1}
\chapter{Introduction to Stochastic Analysis}
\label{ch.intro}
\markboth{Brownian Motion}{Wiener Measure}

This chapter is devoted to the basic results about the Wiener measure,
 Brownian motion, construction of the Ito stochastic integral,
 Cameron-Martin and Girsanov theorems,
 representation of the Wiener functionals with stochastic integrals  and the
 Ito-Wiener chaos decomposition which results from it. The proofs are
 rather sketchy whenever they are given; for complete treatment of
 these results  we refer the reader to the excellent references given
 in the bibliography.

\section{The Brownian Motion and the Wiener Measure}

Let $W=C_0([0,1])$, define $W_t$ as to be the
coordinate functional, i.e., for $w\in W$ and  $t\in[0,1]$, let  $W_t(w)=w(t)$
\index[not]{w@$W$}. If we note by
${\calB}_t=\sigma\{W_s;s\leq t\}$,
then, the following  theorem is well-known (cf. for instance \cite{S-V}):
\begin{theorem}
\label{wiener-thm}
There is one and only one  measure $\mu$
\index[sub]{Wiener  measure} on $W$ which satisfies the following
properties:  \index[not]{m@$\mu$}
\index[not]{b@${\calB}_t$}
\begin{itemize}
\item[i)]
$\mu\left\{w\in W:\,W_0(w)=0\right\}=1$,
\item[ii)]
For any  $ f\in{\calC}_b^\infty(\R)$, the stochastic process  process
\index[not]{w@$W_t$}
$$
(t,w)\mapsto f(W_t(w))-\frac{1}{2}\int_0^t
\Delta f(W_s(w))ds
$$
is a $({\calB}_t,\mu)$-martingale, where $\Delta$ denotes the Laplace
operator. $\mu$ is called the (standard)  Wiener measure.
\end{itemize}
\end{theorem}

\noindent
From Theorem \ref{wiener-thm}, it follows  that, for $t>s$,
$$
E_\mu\left[e^{i\alpha(W_t-W_s)}|{\calB}_s\right]=
\exp\left\{-\frac{1}{2}\alpha^2(t-s)\right\}\,,
$$
 hence $(t,w)\mapsto W_t(w)$ is  a continuous additive process (i.e.,a
 process  with independent increments) and  $(W_t;t\in [0,1])$ is
 also a continuous martingale.

\section{Stochastic Integration}
\markboth{Brownian Motion}{Stochastic Integration}
\index[sub]{stochastic integration}
The stochastic integration with respect to the Brownian motion is
first defined on the adapted step processes and then extended to their
completion by isometry. A mapping  $K: [0,1]\times W\to\R$ is called  a
step process if it can be represented in  the following form:
$$
K_t(w)=\sum_{i=1}^n a_i(w)\cdot1_{[t_i,t_{i+1})}(t),\qquad
a_i(w)\in L^2({\calB}_{t_i})\,.
$$
\noindent
For such a step process, we define \index[not]{i@$I(K)$} its
stochastic integral with respect to the Brownian motion, which is
denoted by
$$
I(K)=\int_0^1K_s dW_s(w)
$$
as to be
$$
\sum_{i=1}^n
a_i(w)\,\left(W_{t_{i+1}}(w)-W_{t_i}(w)\right).
$$
Using the the independence of the increments of $(W_t,t\in [0,1])$, it
is easy to see that
$$
E\left[\left| \int_0^1 K_sdW_s\right|^2\right]=E\int_0^1
\left|K_s\right|^2ds\,, 
$$
i.e.,  $I$ is an isometry from the adapted step processes into
$L^2(\mu)$, hence it has a unique extension as an isometry from
$$
L^2([0,1]\times W,{\mathcal A},dt\times
d\mu)\stackrel{I}{\longrightarrow} L^2(\mu)\,
$$
where ${\mathcal A}$ denotes the sigma algebra on $[0,1]\times W$
generated by the adapted, left (or right) continuous processes.
\noindent
The extension of $I(K)$ is called the \index[sub]{ stochastic
  integral}stochastic integral  of $K$ and it is  denoted as
$\int_0^1K_sdW_s$. If we define
$$
I_t(K)=\int_0^t K_s dW_s
$$
as
$$
\int_0^1 \won_{[0,t]}(s)K_s dW_s,
$$
it follows from the Doob inequality  that the stochastic process
$t\mapsto I_t(K)$ is a continuous, square integrable martingale.
 With some localization techniques using stopping times,
$I$ can be extended to any adapted process $K$
such that  $\int_0^1K_s^2(w)ds<\infty$ a.s. In this case the process
 $t\mapsto I_t(K)$ becomes a local martingale, i.e., there exists a sequence
of stopping times increasing to one, say $(T_n,n\in \N)$ such that the
process $t\mapsto I_{t\wedge T_n}(K)$ is a (square integrable)
martingale. Vector (i.e. $\reals^n$)- valued Brownian motion is
defined as a process whose components are independent, real-valued
 Brownian motions. A stochastic process $(X_t,t\geq 0)$ with values in
 a finite dimensional  Euclidean space is called an Ito process if it
 can be represented as 
$$
X_t=X_0+\int_0^t a_sdW_s+\int_0^tb_sds\,,
$$
where $(W_t,t\geq 0)$ is a vector valued Brownian motion and  $a$ and
$b$ are respectively matrix and vector valued, adapted, 
measurable processes with $\int_0^t(|a_s|^2+|b_s|)ds<\infty$ almost
surely for any $t\geq 0$. In the sequel the notation 
$\int_0^tH_sdX_s$ will mean $\int_0^tH_sa_sdW_s+\int_0^tH_sb_sds$, we
shall also denote by $([X,X]_t,t\geq 0)$ the Doob-Meyer process
defined as 
$$
[X,X]_t=\int_0^t\trace(a_sa^*_s)ds\,.
$$
This is the unique increasing process such that
$(|X_t|^2-[X,X]_t,t\geq 0)$ is a (continuous) local martingale. It can
be calculated as the limit of the sums 
\begin{eqnarray*}
\lim\sum_{i=1}^n\left(|X_{t_{i+1}}|^2-|X_{t_i}|^2\right)&=&
\lim\sum_{i=1}^n\left(|X_{t_{i+1}}|-|X_{t_i}|\right)^2\\
&=&\lim\sum_{i=1}^nE\left[|X_{t_{i+1}}|^2-|X_{t_i}|^2|\calF_{t_i}\right]\,,
\end{eqnarray*}
where the limit is taken as the length of the partition
$\{t_1,\ldots,t_{n+1}\}$ of $[0,t]$, defined by $\sup_i|t_{i+1}-t_i|$,
tends to zero.
\section{Ito formula\index[sub]{Ito formula}}
\markboth{Brownian Motion}{Ito Formula}
The  following result is one of the most  important  applications of
the stochastic integration:
\begin{theorem}
Let  $f\in C^2(\R)$ and let $(X_t,t\in [0,1])$ be an Ito process,
i.e.,
$$
X_t=X_0+\int_0^t K_r dW_r+\int_0^t U_rdr
$$
where $X_0$ is $\calB_0$-measurable,  $K$ and $U$ are adapted processes with
\be
\label{int-cond}
\int_0^1\left[|K_r|^2+|U_r|\right]dr<\infty
\ee
almost surely.
Then
\beaa
f(X_t)&=&f(X_0)+\int_0^t f'(X_s)K_s dW_s+\frac{1}{2} \int_0^t
f''(X_s)K_s^2ds\\
&&+\int_0^tf'(X_r)U_rdr\,.
\eeaa
\end{theorem}
\begin{remarkk}{\rm
This formula is also valid in the several dimensional case.  In fact,
if $K$ is and
$U$ are adapted processes with values in $\R^n\otimes \R^m$ and $\R^m$
respectively whose components are satisfying the condition
(\ref{int-cond}), then we have, for any $f\in C^2(\R^m)$,
\beaa
f(X_t)&=&f(X_0)+\int_0^t\partial_i f(X_r)K_{ij}(r)dW^j_r
+\int_0^t\partial_i f(X_r)U_i(r)dr\\
&&+\frac{1}{2}\int_0^t\partial^2_{ij}f(X_r)(K_rK_r^*)_{ij}dr
\eeaa
almost surely.}
\end{remarkk}

To prove the Ito formula we shall proceed by
\begin{lemma}
\label{ipp}
Let $X=(X_t,t\geq 0)$ and  $Y=(Y_t,t\geq 0)$ be two Ito real-valued processes,
then 
\begin{equation}
\label{eq1}
X_tY_t=X_0Y_0+\int_0^tX_sdY_s+\int_0^tY_sdX_s+[X,Y]_t\,,
\end{equation}
almost surely, where $[X,Y]$ denotes the Doob-Meyer process. In
particular 
\begin{equation}
\label{eq2}
X_t^2=X_0^2+2\int_0^tX_sdX_s+[X,X]_t\,.
\end{equation}
\end{lemma}
\proof
 Evidently it suffices to prove the relation (\ref{eq2}), since we
can obtain (\ref{eq1}) via a polarization argument. 
Since $X$ has almost surely continuous trajectories, using a
stopping time argument we can assume that $X$ is almost surely
bounded. Assume now that $\{t_1,\ldots,t_n\}$ is a partition of
$[0,t]$ and  denote by $(M_t,t\geq 0)$ the local martingale part and by
$(A_t,t\geq 0)$ the finite variaton part of $X$. We have 
\begin{eqnarray}
X_t^2-X_0^2&=&\sum_{k=1}^n(X_{t_k}-X_{t_{k-1}})^2+2\sum_{k=1}^nX_{t_{k-1}}
(X_{t_k}-X_{t_{k-1}})\label{part1}\\
&=&\sum_{k=1}^n(X_{t_k}^2-X_{t_{k-1}}^2)+2\sum_{k=1}^nX_{t_{k-1}}
(M_{t_k}-M_{t_{k-1}})\label{part2}\\
&&\hspace{3cm} +2\sum_{k=1}^nX_{t_{k-1}}(A_{t_k}-A_{t_{k-1}})\label{part3}\,.
\end{eqnarray}
Now, when $\sup_k|t_k-t_{k-1}|\to 0$, then  the first  term at the
right hand side of (\ref{part2}) converges to $[X,X]_t$ and the sum 
of the second term with  (\ref{part3}) converges to $2\int_0^tX_sdX_s$ in
probability.
\qed

\noindent
{\bf{Proof of the Ito formula:}}\\
Using a stopping argument we can assume that $X$ takes its values in a
bounded interval, say $[-K,K]$. The interest of this argument resides
in the fact that we can approach a $C^2$ function, as well as its first two
derivatives  uniformly by the polynomials on any compact interval. On
the other hand, using Lemma \ref{ipp}, we see that the formula is
valid for the polynomials. Let us denote by $(\Gamma f)_t$ the random
variable 
$$
f(X_t)-f(X_0)-\int_0^tf'(X_s)dX_s-\frac{1}{2}\int_0^tf''(X_s)d[X,X]_s\,.
$$
Assume moreover that $(p_n,n\geq 1)$ is a sequence of polynomials such
that $(p_n,\geq 1)$,  $(p'_n,\geq 1)$ and   $(p''_n,\geq 1)$ converge
uniformly on $[-K,K]$ to $f,f'$ and to $f''$ respectively. Choosing a
subsequence, if necessary, we may assume that 
$$
\sup_{x\in
  [-K,K]}\left(|f(x)-p_n(x)|+|f'(x)-p'_n(x)|+|f''(x)-p''_n(x)|\right)\leq 1/n\,.
$$
Using the Doob and the Chebytchev inequalities, it is easy to see that
$(\Gamma f)_t-(\Gamma p_n)_t$ converges to zero in probability, since
$(\Gamma p_n)_t=0$ almost surely, $(\Gamma f)_t$ should be also zero
almost surely and this completes the proof of the Ito formula.
\qed

As an immediate corollary of the Ito formula  we have
\begin{corollary}
For any $h\in L^2([0,1])$, the process defined by
$$
{\mathcal E}_t(I(h))=\exp\left(\int_0^t h_sdW_s-\frac{1}{2}\int_0^t
  h_s^2 ds\right)\,
$$
is a martingale.
\end{corollary}
\proof
Let us denote ${\mathcal E}_t(I(h))$ by $M_t$, then from the Ito
formula we have
$$
M_t=1+\int_0^tM_sh_sdW_s\,,
$$
hence $(M_t,t\in [0,1])$ is a local martingale, moreover, since $I(h)$
is Gaussian, $M_1$ is in all the $L^p$-spaces, hence $(M_t,t\in
[0,1])$ is a square integrable martingale.
\qed

\section{Alternative constructions of the Wiener measure}

\paragraph{A)}
Let us state first the celebrated theorem of Ito-Nisio about the
convergence of the random series of independent, Banach space valued
random variables (cf. \cite{I-N}):
\begin{theorem}[Ito-Nisio Theorem]
\label{Ito-N-thm}
Assume that $(X_n,n\in \NN)$ is a sequence of independent random
variables with values in a separable Banach space $B$ whose continuous
dual is denoted by $B^\star$. The sequence 
 $(S_n,n\in \NN)$ defined  as 
$$
S_n=\sum_{i=1}^nX_i\,,
$$
converges almost surely in the norm topology of $B$ if and only if
there exists a probability measure $\nu$ on $B$ such that
$$
\lim_n E\left[e^{i<\xi,S_n>}\right]=\int_B e^{i<\xi,y>}\nu(dy)
$$
for any $\xi\in B^\star$.
\end{theorem}
We can give another construction of the Brownian motion using Theorem
\ref{Ito-N-thm} as follows: Let $(\gamma_i;i\in\N)$ be an independent
sequence of $N_1(0,1)$-Gaussian  random variables. Let $(g_i)$ be a
complete, orthonormal basis  of $L^2([0,1])$. Then $W_t$ defined by
$$
W_t(w)=\sum_{i=1}^\infty \gamma_i(w)\cdot\int_0^t g_i(s)ds
$$
converges almost surely uniformly with respect to $t\in [0,1]$  and
$(W_t, t\in [0,1])$ is a Brownian motion. In fact to see this it
suffices to apply Theorem \ref{Ito-N-thm} to the sequence $(X_n,n\in
\NN)$ defined by
$$
X_n(w)=\gamma_n(w)\int_0^\cdot g_n(s)ds\,.
$$

\begin{remarkk}{\rm
In the sequel we shall denote by $H$ the so-called
Cameron-Martin space $H([0,1],\R^n)$ (in case $n=1$ we shall again
write simply $H$ or, in case of necessity $H([0,1])$) i. e., the
isometric image of 
$L^2([0,1],\R^n)$ under the mapping 
$$
g\to\int_0^\cdot g(\tau)d\tau\,.
$$ 
Hence for any  complete, orthonormal
basis $(g_i,i\in\N)$  of $L^2([0,1],\R^n)$,
$(\int_0^\cdot g_i(s)ds,\,i\in \N )$ is a complete orthonormal basis
of $H([0,1],\R^n)$. The use of the generic notation $H$ will be
preferred as long as the results are dimension independent.}
\end{remarkk}

\paragraph{B)}
Let $(\Omega,{\mathcal F},\P)$ be any abstract probability space and
let $H$ be any separable Hilbert space. If
$L:H\!\to\!L^2(\Omega,{\mathcal F},\P)$ is a linear operator such that
for any $h\in H$, $E[\exp iL(h)]=\exp-\frac{1}{2}|h|^2_H$, then there
exists a Banach space with dense injection
$$
H\stackrel{\scriptstyle \jmath}{\hookrightarrow}W
$$
dense, hence
$$
W^\ast\stackrel{\scriptstyle \jmath^\ast}{\hookrightarrow}H
$$
is also dense  and there exists  a probability measure $\mu$ on $W$ such that
$$
\int_W\exp\langle w^\ast,w\rangle
d\mu(w)=\exp-{\textstyle\frac{1}{2}} \mid j^\ast(w^\ast)\mid_H^2
$$

\noindent
and
$$
L(j^\ast(w^\ast))(w)=\langle w^\ast,w\rangle
$$
almost surely. $(W,H,\mu)$ is
called an Abstract Wiener space\index[sub]{abstract Wiener space} and
$\mu$ is the Wiener measure (cf. \cite{LG-1}). In the case $H$ is
chosen to be
$$
H([0,1])=\left\{h:h(t)=\int_0^t\dot{h}(s)ds,|h|_H=|\dot{h}|_{L^2([0,1])}\right\}
$$
then $\mu$ is the classical Wiener measure and $W$ can be taken
as $C_0([0,1])$.

\begin{remarkk}{\rm
In the case of the classical Wiener space, any element $\lambda$ of
$W^\ast$ is a signed measure on $[0,1]$, and its image in $H=H([0,1])$
can be represented as
$j^\ast(\lambda)(t)=\int_0^t\lambda([s,1])ds$. In fact, we have for
any $h\in H$
\begin{eqnarray*}
(j^{\ast}(\lambda),h)&=&<\lambda,j(h)>\\
  &=&\int_0^1 h(s)\lambda(ds)\\
&=&h(1)\lambda([0,1])-\int_0^1 \lambda([0,s])\dot{h}(s)ds\\
&=&\int_0^1(\lambda([0,1])-\lambda([0,s])\dot{h}(s)ds\\
&=&\int_0^1 \lambda([s,1])\dot{h}(s)ds.
\end{eqnarray*}}
\end{remarkk}

\section{Cameron-Martin and Girsanov  Theorems}
In the sequel we shall often need approximation of the Wiener
functional with cylindrical smooth functions on the Wiener space. This
kind of properties hold in every Wiener space since  this is due to the
analyticity of the characteristic function of the Wiener
measure. However, they are very easy to explain in the case of
classical Wiener space, that is why we have chosen to work in this
frame. In particular the Cameron-Martin theorem which is explained in
this section is absolutely indispensable for the development of the
next chapters.

\begin{lemma}
\label{app1-lem}
The set of random variables
$$
\Bigl\{f(W_{t_1},\ldots,W_{t_n}); t_i\in[0,1],
f\in{\calS}(\R^n); n\in\N\Bigr\}
$$
is dense in $L^2(\mu)$, where ${\calS}(\R^n)
$\index[not]{s@${{\calS}(\R^n)}$}
 denotes the space of infinitely differentiable, rapidly decreasing
 functions on $\R^n$.
\end{lemma}
\proof It follows from the martingale convergence theorem and the monotone
class theorem.
\qed
\begin{lemma}
\label{app2-lem}
The linear span of the set
$$
\Theta=\left\{\exp\left[\int_0^1 h_sdW_s-\frac{1}{2}\int_0^1
    h_s^2ds\right]:\,h\in L^2([0,1])\right\}
$$
is dense in $L^2(\mu)$.
\end{lemma}
\proof
It follows from Lemma  \ref{app1-lem}, via the Fourier transform.
\qed

\remark
Although the set $\Theta$ separates the points of $L^2(\mu)$, it does
not give any indication about the positivity.

\begin{lemma}
\label{app3-lem}
The polynomials are dense in  $L^2(\mu)$.
\end{lemma}
\proof
The proof follows by  the analyticity of the characteristic function
of the Wiener measure, in fact, due to this property,   the elements
of the set in Lemma \ref{app2-lem} can be approached by the
polynomials.
\qed

\begin{theorem}[Cameron-Martin Theorem]
For any bounded Borel measurable function $F$ on $C_0([0,1])$ and  $h\in
L^2([0,1])$, we have
$$
E_\mu\left[F\left(w+\int_0^\cdot h_sds\right)\,\exp\left\{ -\int_0^1
    h_sdW_s-\frac{1}{2} \int_0^1 h_s^2ds\right\}\right]=E_\mu[F]\,.
$$
This assertion  implies in particular  that the process
$(t,w)\to W_t(w)+\int_0^t h_sds$ is again a Brownian motion under the new
probability measure
$$
\exp\left\{-\int_0^1h_s dW_s-\frac{1}{2}\int_0^1h_s^2ds\right\} d\mu.
$$
\end{theorem}
\proof
It is sufficient to show that the new probability has the same characteristic
function  as $\mu$: if $x^\ast\in W^\ast$, then $x^\ast$ is a measure
on $[0,1]$ and
\begin{eqnarray*}
_{W^\ast}\langle x^\ast,w\rangle_W & = &
     \int_0^1W_s(w)x^\ast(ds)\\
&=&W_t(w)\cdot x^\ast([0,t]) \Big|_0^1-
                         \int_0^1 x^\ast([0,t])dW_t(w) \\
& = & W_1 x^\ast([0,1])-\int_0^1 x^\ast([0,t]).dW_t \\
& = & \int_0^1 x^\ast((t,1])dW_t\,.
\end{eqnarray*}
Consequently
\begin{eqnarray*}
\lefteqn{ E\left[\left\{\exp i \int_0^1
      x^\ast([t,1])dW_t\right\}\left(w+\int_0^\cdot  h_sds\right)\,\,
     {\mathcal E}(-I(h))\right]} \\
 = &&\!\!E\left[\exp\left\{ i \int_0^1 \!\!x^\ast([t,1])dW_t+i\int_0^1\!\! x^\ast([t,1]) h_tdt\! -
     \! \int_0^1 \!h_tdW_t-\!\frac{1}{2}\int_0^1\! h_t^2dt\right\}\right]\\
 = &&\!\! E\left[\exp\left\{i \int_0^1(ix^\ast([t,1])-h_t)dW_t\right\}\right.\\
&&\hspace{2cm}\left.\exp\left\{ i \int_0^1 x^\ast([t,1])
       h_tdt - \frac{1}{2}\int_0^1 h_t^2dt\right\}\right] \\
 = &&\!\!  \exp\left\{ \frac{1}{2}\int_0^1\!(ix^\ast([t,1])-h_t)^2dt +
   i\int_0^1 x^{\ast}([t,1])h_tdt - \frac{1}{2}\int_0^1\!h_t^2dt\right\} \\
 =&& \!\!  \exp - \frac{1}{2}\int_0^1(x^\ast([t,1]))^2dt \\
 =&& \!\!\exp-\frac{1}{2}|j(x^\ast)|_{H}^2\, ,
\end{eqnarray*}
and this achieves the proof.
\qed

\noindent
The following corollary is one of the most important results of the
modern probability theory:

\begin{corollary}[Paul L\'evy's Theorem] \index[sub]{Paul L\'evy's Theorem}
Suppose that $(M_t,\,t\in [0,1])$ is a continuous martingale with
$M_0=0$ and that  $(M_t^2-t,\,t\in [0,1])$ is again a martingale. Then
$(M_t,\,t\in [0,1])$ is a Brownian motion.
\end{corollary}
\proof
From  the Ito formula
$$
f(M_t)=f(0)+\int_0^t f'(M_s)\cdot dM_s+ \frac{1}{2}\int_0^t\Delta  f(M_s)\,\,
ds\,.
$$
Hence the law of $(M_t:t\in[0,1])$ is $\mu$.
\qed

As an application of Paul L\'evy's  theorem we can prove easily  the following
result known as the Girsanov theorem which generalizes the
Cameron-Martin theorem. This theorem is basic in several applications
like  the filtering of the
random signals corrupted by a Brownian motion, or the problem of
optimal control of Ito processes.
\begin{theorem}[Girsanov Theorem] \index[sub]{Girsanov Theorem}
Assume that $u:[0,1]\times W\to \reals^n$ is a measurable process
adapted to the Brownian filtration satisfying
$$
\int_0^1|u_s|^2ds<\infty
$$
$\mu$-almost surely. Let 
$$
\Lambda_t=\exp\left\{-\int_0^t(u_s,dW_s)-1/2\int_0^t|u_s|^2ds\right\}\,.
$$
Assume that 
\begin{equation}
\label{u.i-cond}
E\left[\Lambda_1\right]=1\,.
\end{equation}
Then the process $(t,w)\to W_t(w)+\int_0^tu_s(w)ds$ is a Brownian
motion under the probability $\Lambda_1d\mu$.
\end{theorem}
\begin{remarkk}{\rm
The condition (\ref{u.i-cond}) is satisfied in particular if
we have
$$
E\left[\exp\frac{1}{2}\int_0^1|u_s|^2ds\right]<\infty\,.
$$
This is called the Novikov condition (cf. \cite{Nov,UZ7}).
\index[sub]{Novikov condition} There is another, slightly more
general sufficient condition due to Kazamaki \cite{Kaz}, which is
$$
E\left[\exp\frac{1}{2}\int_0^1u_sdW_s\right]<\infty\,.
$$

\noindent
Note that the difference between the Cameron--Martin theorem and
the Girsanov theorem is that in the former the mapping $w\to
w+\int_0^\cdot h(s)ds$ is an invertible transformation of the
Wiener space $W$ and in the latter the corresponding map $w\to
w+\int_0^\cdot u_s(w)ds$ is not necessarily invertible.}
\end{remarkk}
\section{The Ito Representation Theorem}
The following result is known as the Ito representation formula:
\begin{theorem}
\label{Ito-rep-thm}
Any $\varphi\in L^2(\mu)$ can be represented as
$$
\varphi=E[\varphi]+\int_0^1 K_sdW_s
$$
where $K\in L^2([0,1]\times W)$ and it is  adapted.
\index[not]{ex@${\mathcal E}(I(h))$}
\end{theorem}
\proof
Since the \index[sub]{w@Wick exponentials} Wick exponentials
 $$
{\mathcal E}(I(h))=\exp\left\{\int_0^1h_sdW_s-1/2\int_0^1h_s^2ds\right\}
$$
can be represented  as claimed and since their finite linear
combinations are dense in $L^2(\mu)$,  the proof follows.
\qed
\begin{remarkk}{\rm
Let $\phi$ be an integrable real  random variable on the Wiener
space. We say that it belongs to the class $H^1$ if the martingale
$M=(M_t,t\in [0,1])$ satisfies the property that
$$
E[<M,M>^{1/2}_1]<\infty\,.
$$
The Ito  representation theorem extends via stopping techniques to the
random variables of class $H^1$.}
\end{remarkk}
\section{Ito-Wiener chaos representation \index[sub]{Ito-Wiener
    chaos representation}}

For any  $h\in L^2([0,1])$, define  $K_t=\int_0^th_sdW_s$, $t\in
[0,1]$.  Then, from the Ito formula, we can write
\begin{eqnarray*}
\lefteqn{K_1^p= p\int_0^1K_s^{p-1}h_sdW_s+\frac{p(p-1)}{2}
          \int_0^1K_s^{p-2}h_s^2ds}\\
&& = p\int_0^1\Big[ (p-1)\int_0^{t_1}K_{t_2}^{p-2}h_{t_2}dW_{t_2}+
          \frac{(p-1)(p-1)}{2}\int_0^{t_1}
          K_{t_2}^{p-3}h_{t_2}^2dt_2\Big] dW_{t_1} \\
&&\,\,\,\, + \cdots
\end{eqnarray*}
where $p$ is a positive integer.
\noindent
Iterating this procedure we   see that
$K_1^p$ can be written as the linear combination
of  the multiple integrals of deterministic integrands of the type
$$
J_p=\int_{0 <t_p< t_{p-1}<\cdots< t_1<1}\hspace*{-3em}
h_{t_1}h_{t_2}\ldots h_{t_p}\,dW_{t_1}^{i_1}\ldots dW_{t_p}^{i_p},
$$
$i_j=0$ or $1$ with  $dW_t^0=dt$ and  $dW_t^1=dW_t$. Hence we can
express the polynomials as multiple Wiener-Ito integrals. Let us now
combine this observation with the Ito representation:

Assume that  $\varphi\in L^2(\mu)$, then  from the Ito representation
theorem :
$$
\varphi =  E[\varphi]+\int_0^1K_sdW_s\,.
$$
Iterating the same procedure for the integrand of the above stochastic
integral:
\begin{eqnarray*}
\varphi & = & E[\varphi]+\int_0^1E[K_s]dW_s+\int_0^1\!\int_0^{t_1}
    E[K_{t_1,t_2}^{1,2}] dW_{t_2}dW_{t_1} \\
& &\hspace{0.8cm} +\int_0^1\!\int_0^{t_1}\!\int_0^{t_2}K_{t_1t_2t_3}^{1,2,3}
dW_{t_3}dW_{t_2}dW_{t_1}\,.
\end{eqnarray*}

\noindent
After $n$ iterations  we end up with
$$
\varphi=\sum_{p=0}^n J_p(K_p)+\varphi_{n+1}
$$
and each element of the sum is orthogonal to the other one. Hence $(\varphi_n;
n\in\N)$ is bounded in the Hilbert space  $L^2(\mu)$ and this means
that it is weakly relatively compact.  Let $(\varphi_{n_k})$ be a
weakly convergent
subsequence and $\varphi_\infty=\lim\limits_{k\to\infty}\varphi_{n_k}$. Then it
is easy from  the first part that $\varphi_\infty$ is orthogonal to the
polynomials, therefore $\varphi_\infty=0$ and the weak limit
$$
w-\lim_{n\to\infty}\sum_{p=0}^n J_p(K_p)
$$
exists and it is equal to $\varphi$ almost surely. Let 
$$
S_n=\sum_{p=0}^n J_p(K_p)\,,
$$
then, from the weak convergence, we have 
$$
\lim_nE[|S_n|^2]=\lim_n E[S_n\,\varphi]=E[|\varphi|^2]\,,
$$
hence $(S_n,n\geq 1)$ converges weakly to $\varphi$ and its $L^2$-norm
converges to the $L^2$-norm of $\varphi$ and this implies that  the series
$$
\sum_{p=1}^\infty J_p(K_p)
$$
converges to $\varphi$  in the strong topology of  $L^2(\mu)$. Let now
$\widehat{K}_p$ be an element of $\widehat{L}^2[0,1]^p$ (i.e.\ symmetric),
defined as $\widehat{K}_p=K_p$ on $C_p=\{t_1<\cdots<t_p\}$. We define
$I_p(\widehat{K}_p)=p!J_p(K_p)$ in such a way that \index[not]{im@$I_p$}
$$
E[|I_p(\widehat{K}_p)|^2]=(p!)^2\int_{C_p}K^2_p dt_1\ldots dt_p=p!
\int_{[0,1]^p}|\widehat{K}_p|^2 dt_1\ldots dt_p\,.
$$

\noindent
Let $\varphi_p=\frac{\widehat{K}_p}{p!}\,$, then we have proven
\begin{theorem}
\label{chaos-thm}
Any element $\varphi$ of $L^2(\mu)$, can be decomposed as an orthogonal
sum of multiple Wiener-Ito integrals
$$
\varphi=E[\varphi]+\sum_{p=1}^\infty I_p(\varphi_p)
$$
where $\varphi_p$ is a symmetric element of $ L^2[0,1]^p$. Moreover,
this decomposition is unique.
\end{theorem}
\remark
In the following chapters we shall give an explicit representation of
the kernels $\varphi_p$ using the Gross-Sobolev derivative.

\section*{Notes and suggested reading}
{\footnotesize{The basic references for the stochastic calculus
  are  the books of Dellacherie-Meyer \cite{D-M} and of
  Stroock-Varadhan \cite{S-V}. Especially in the former, the theory is
  established  for the general semimartingales with jumps. For the
  construction of  the Wiener measure on Banach spaces we refer the reader to
  \cite{LG-1} and especially to \cite{KUO}.
}}


\chapter{Sobolev Derivative, Divergence and Ornstein-Uhlenbeck
  Operators}
\label{ch.deriv}
\markboth{Derivative, Divergence}{Ornstein-Uhlenbeck Operators}

\setcounter{section}{0}

\section{Introduction}

Let $W=C_0([0,1],\R^d)$ be the classical Wiener space equipped  with $\mu$
the Wiener
measure. We want to construct on $W$ a Sobolev type analysis in such a way that
we can apply it to the random variables that we encounter in the applications.
Mainly we want to construct a differentiation operator and to be able to apply
it to practical examples. The Fr\'echet derivative is not satisfactory. In fact
the most frequently encountered Wiener functionals, as the multiple
(or single) Wiener integrals or the solutions of stochastic
differential equations with smooth coefficients are not even
continuous with respect to the Fr\'echet norm of the Wiener
space. Therefore, what we need is in fact to define a derivative on
the
$L^p(\mu)$-spaces of random variables, but in general, to be able to do this,
we need the following  property which is essential: if
$F,G\in L^p(\mu)$, and if we want
to define their directional derivative, in the direction, say $\tilde{w}\in W$,
we write $\frac{d}{dt}F(w+t\tilde{w})|_{t=0}$ and
$\frac{d}{dt}G(w+t\tilde{w})|_{t=0}$. If $F=G$ $\mu$-a.s., it is
natural to ask that
their derivatives are also equal a.s. For this, the only way is to choose
$\tilde{w}$ in some specific subspace of $W$, namely, the Cameron-Martin space
$H$:
$$
H=\left\{h:[0,1]\to\R^d/h(t)= \int_0^t\dot{h}(s)ds,\,\,
|h|_{H}^2=\int_0^1|\dot{h}(s)|^2 ds\right\}.
$$

\noindent
In fact, the theorem of Cameron-Martin says that for any $F\in L^p(\mu)$,
$p>1$, $h\in H$
$$
E_\mu\left[F(w+h)\exp\left\{-\int_0^1\dot{h}(s)\cdot
    dW_s-{\textstyle\frac{1}{2}} |h|_{H}^2\right\}\right]=E_\mu[F]\,,
$$
or equivalently
$$
E_\mu[F(w+h)]=E\left[F(w)\cdot\exp\left\{\int_0^1 \dot{h}_s\cdot dW_s-
{\textstyle\frac{1}{2}}|h|_{H}^2\right\}\right]\,.
$$

\noindent
That is to say, if $F=G$ a.s., then $F(\cdot+h)=G(\cdot+h)$ a.s. for all
$ h\in H$.\index[not]{h@$H$}

\section{The Construction of $\nabla$ and its properties}

If $F:W\to\R$ is a function of the following type (called\index[sub]{
cylindrical} cylindrical ):
$$
F(w)=f(W_{t_1}(w),\ldots,W_{t_n}(w)),\qquad f\in{\calS}(\R^n),
$$
\index[not]{sa@${\calS}$}
we define, for $h\in H$,
$$
\nabla_h F(w)=\frac{d}{d\lambda}F(w+\lambda h)|_{\lambda=0}\,.
$$
Noting that $W_t(w+h)=W_t(w)+h(t)$, we obtain
$$
\nabla_h F(w)=\sum_{i=1}^n \partial_i f(W_{t_1}(w),\ldots,W_{t_n}(w)) h(t_i),
$$

\noindent
in particular \index[not]{n@$\nabla$}
$$\nabla_hW_t(w)=h(t)=\int_0^t\dot{h}(s)ds=\int_0^1
1_{[0,t]}(s)\ \dot{h}(s)ds.
$$
If we denote by
$U_t$ the element of $H$ defined as
$U_t(s)=\int_0^s 1_{[0,t]}(r)dr$, we have
$\nabla_h W_t(w)=(U_t,h)_{H}$. Looking at the linear map
$h\mapsto\nabla_hF(w)$ we see that it defines a random element with
values in $H^\star$, since we have identified $H$ with $H^\star$,
$\nabla F$ is an $H$-valued  random variable.
Now we can prove:

\begin{proposition}
\label{clos-prop}
$\nabla$ is a closable operator\index[sub]{closable} on any $L^p(\mu)$ $(p>1)$.
\end{proposition}
\proof
Closable  means that if $(F_n:n\in\N)$ are cylindrical functions on
$W$, such that  $F_n\to 0$ in $L^p(\mu)$ and if $(\nabla F_n;n\in\N)$
is Cauchy in
$L^p(\mu,H)$, then its limit is zero. Hence
suppose that $\nabla F_n\to \xi$ in $L^p(\mu;H)$.
In order to prove $\xi=0$ $\mu$-a.s.,  we use the Cameron-Martin theorem: Let
$\varphi$ be any cylindrical
function. Since such $\varphi$'s are dense in $L^p(\mu)$, it is sufficient to
prove that $E[(\xi,h)_{H}\cdot\varphi]=0$  for any
$h\in H$. This follows from
\begin{eqnarray*}
& & E[(\nabla F_n,h)\varphi]=\frac{d}{d\lambda}E[ F_n(w+\lambda
h)\cdot\varphi]|_{\lambda=0} \\
& & = \frac{d}{d\lambda}E\left.\left[F_n(w)\varphi(w-\lambda h)\exp
   \left( \lambda \int_0^1\dot{h}(s) dW_s-\frac{\lambda^2}{2}
   \int_0^1|\dot{h}_s|^2ds\right)\right]\right|_{\lambda=0} \\
& & = E\left[F_n(w)\left(-\nabla_h\varphi(w)+\varphi(w)\int_0^1\dot{h}(s) dW_s\right)\right]
\,\longrightarrow\raisebox{-1.2ex}{\hspace*{-2em}$\scriptstyle n\to\infty$}\,0
\end{eqnarray*}
since  $(F_n,n\in \NN)$ converges to zero  in $L^p(\mu)$.
\qed

\noindent
Proposition \ref{clos-prop}  tells us that the operator $\nabla$ can
be extended to larger classes of Wiener functionals than the
cylindrical ones. In fact we  define first the extended  $L^p$-domain
of $\nabla$, denoted  by
${\Dom}_p(\nabla)$\index[not]{dom@${\Dom}_p$} as
\begin{definition}
$F\in{\Dom}_p(\nabla)$ if and only if there exists a sequence
$(F_n;n\in \N)$ of  cylindrical functions  such that
 $F_n\to F$ in $L^p(\mu)$ and $(\nabla F_n)$ is Cauchy in $L^p(\mu,H)$.
Then, for any $F\in {\Dom}_p(\nabla)$,  we define
$$
\nabla F=\lim_{n\to\infty}\nabla F_n .
$$
The extended operator $\nabla$ is called
{\bf {Gross-Sobolev derivative}}
\index[sub]{Gross-Sobolev derivative}.
\end{definition}
\begin{remarkk} Proposition \ref{clos-prop} implies that the definition of
$\nabla F$ is independent of the choice of the approximating sequence
$(F_n)$.
\end{remarkk}
Now we are ready to define
\begin{definition}
We will denote by $\DD_{p,1}$ the linear space ${\Dom}_p(\nabla)$ equipped
with the norm $\|F\|_{p,1}=\|F\|_p+\|\nabla F\|_{L^p(\mu,H)}$.
\index[not]{dd@$\DD_{p,1}$}
\end{definition}
\begin{remarkk}{\rm
\begin{enumerate}
\item  If ${\Xi}$ is a separable Hilbert space we can define
  $\DD_{p,1}({\Xi})$ exactly in the same way as before, the
  only difference is that we take ${\calS}_{\Xi}$ instead of
  ${\calS}$, i.e., the rapidly decreasing functions with
values in $\Xi$. Then we leave to the reader to prove that  the
same closability result holds.
\item Hence we can define $\DD_{p,k}$ by iteration:
\begin{itemize}
\item[i)]
We say that $F\in \DD_{p,2}$ if $\nabla F\in \DD_{p,1}(H)$, then write
$\nabla^2F=\nabla(\nabla F)$.\index[not]{n@$\nabla^k$}
\item[ii)]
$F\in \DD_{p,k}$ if $\nabla^{k-1}F\in \DD_{p,1}(H^{\otimes(k-1)})$.
\end{itemize}

\item  Note that, for $F\in \DD_{p,k}$, $\nabla^kF$ is in fact with values
$H^{\widehat{\otimes}k}$ (i.e.\ symmetric tensor product).
\index[not]{d@$\DD_{p,k}(X)$}
\item From the proof we have that if $F\in \DD_{p,1}$, $h\in H$ and $\varphi$ is
cylindrical, we have
$$
E[\nabla_h F\cdot\varphi]=-E[F\cdot\nabla_h\varphi]
  +E[I(h)\cdot F\cdot\varphi]\,,
$$
where $I(h)$ is the first order Wiener integral of the (Lebesgue) density
of $h$.
\noindent
If $\varphi\in \DD_{q,1}$ $(q^{-1}+p^{-1}=1)$, by a limiting argument, the same
relation holds  again. Let us note that this limiting procedure shows
in fact that if $\nabla F\in L^p(\mu,H)$ then $F.I(h)\in L^p(\mu)$,
i.e., $F$ is more than $p$-integrable. This observation gives rise to
the logarithmic Sobolev inequality.\index[sub]{lo@{logarithmic Sobolev
    inequality}}
\end{enumerate}
}
\end{remarkk}

\section{Derivative of  the Ito  integral}

Let $\varphi=f(W_{t_1},\ldots,W_{t_n})$, $t_i\leq t$, $f$ smooth. Then we have
$$\nabla_h\varphi(w)=\sum_{i=1}^n \partial_if(W_{t_1},\ldots,W_{t_n})h(t_i)\,,$$

\noindent
hence $\nabla\varphi$ is again a random variable which is ${\mathcal
B}_t$-measurable. In fact this property is satisfied by a larger
class of Wiener functionals:

\begin{proposition}
\label{mes-supp-prop}
Let $\varphi\in \DD_{p,1}$, $p>1$ and suppose that $\varphi$ is ${\mathcal
B}_t$-measurable for a given $t\geq0$. Then $\nabla\varphi$ is also ${\mathcal
B}_t$-measurable and furthermore, for any $h\in H$, whose support is in
$[t,1]$, $\nabla_h\varphi=(\nabla\varphi,h)_H=0$ a.s.
\end{proposition}
\proof
Let $(\varphi_n)$ be a sequence of cylindrical random variable
converging to $\varphi$ in
$\DD_{p,1}$. If $\varphi_n$ is of the form $f(W_{t_1},\ldots,W_{t_k})$, it is
easy to see that, even if $\varphi_n$ is not ${\calB}_t$-measurable,
$E[\varphi_n|{\calB}_t]$ is another cylindrical random variable, say
$\theta_n(W_{t_1\wedge t},\ldots,W_{t_k\wedge t})$. In fact, suppose that
$t_k>t$ and $t_1,\ldots,t_{k-1}\leq t$. We have
\begin{eqnarray*}
& & E[f(W_{t_1},\ldots,W_{t_k})|{\calB}_t]=
    E[f(W_{t_1}\ldots,W_{t_{k-1}},W_{t_k}-W_t+W_t)|{\calB}_t] \\
& & = \int_\R f(W_{t_1},\ldots, W_{t_{k-1}},W_t+x)p_{t_k-t}(x)dx \\
& & = \theta(W_{t_1},\ldots,W_{t_{k-1}},W_t)\,,
\end{eqnarray*}

\noindent
and $\theta\in{\calS}$ if $f\in{\calS}(\R^k)$, where $p_t$ denotes the heat
kernel. Hence we can choose a sequence $(\varphi_n)$ converging to $\varphi$ in
$\DD_{p,1}$ such that  $\nabla\varphi_n$ is ${\calB}_t$-measurable for
each $n\in\N$. Hence $\nabla\varphi$ is also ${\calB}_t$-measurable.
If $h\in H$ has its support in $[t,1]$, then, for each $n$, we have
$\nabla_h\varphi_n=0$ a.s., because $\nabla\varphi_n$ has its support in
$[0,t]$ as one can see from the explicit calculation for $\nabla\varphi_n$.
Taking an a.s.\ convergent subsequence, we see that $\nabla_h\varphi=0$ a.s.\
also.
\qed.

Let now $K$ be an adapted simple  process:
$$
K_t(w)=\sum_{i=1}^n a_i(w)1_{(t_i,t_{i+1}]}(t)
$$
where $a_i\in \DD_{p,1}$ and ${\calB}_{t_i}$-measurable for any $i$. Then
we have
$$
\int_0^1K_s dW_s=\sum_{i=1}^n a_i(W_{t_{i+1}}-W_{t_i})
$$
and
\begin{eqnarray*}
 \nabla_h\int_0^1 K_s dW_s&=&\sum_{i=1}^n \nabla_h
 a_i(W_{t_{i+1}}-W_{t_i})\\
&& +\sum_{i=1}^n a_i(h(t_{i+1})-h(t_i))\\
 &=& \int_0^1\nabla_h K_s dW_s+\int_0^1 K_s\dot{h}(s)ds\,.
\end{eqnarray*}
Hence
$$
\left| \nabla\int_0^1 K_s dW_s\right|_H^2\leq 2\Big\{\Big| \int_0^1 \nabla K_s
dW_s\Big|_H^2 + \int_0^1|K_s|^2ds\Big\}
$$
and
$$
E\left[\left(\Big|\nabla\int_0^1 K_s dW_s\Big|_H^2\right)^{p/2}\right]\leq
   2^p E\left[\left(\Big|\int_0^1\nabla K_s dW_s\Big|_H^p \right.\right.
$$
$$+ \left.\left.\int_0^1 |K_s|^2ds\right)^{p/2}\right]\,.
$$
Using the Burkholder-Davis-Gundy  inequality for the Hilbert space valued
martingales, the above quantity is majorized by
\beaa
\lefteqn{ 2c_p E\left\{\left(\int_0^1|\nabla K_s|_H^2 ds\right)^{p/2}
+ \left(\int_0^1 |K_s|^2ds\right)^{p/2}\right\}}\\
&&=\tilde{c}_p\|\nabla\tilde{K}\|^p_{L^p(\mu,H\otimes H)}+
\|\tilde{K}\|_{L^p(\mu,H)}\,,
\eeaa
where
$$
\tilde{K}.=\int_0^\cdot K_rdr\,.
$$
Thanks to this majoration, we have proved:

\begin{proposition}
Let $\tilde{K}\in \DD_{p,1}(H)$ such that  $K_t=\frac{d\tilde{K}(t)}{dt}$ be
${\calB}_t$-measurable for almost all $t$. Then we have
\be
\label{com-rel}
\nabla \int_0^1 K_s dW_s=\int_0^1\nabla_\cdot K_s dW_s+\tilde{K}
\ee
almost surely.
\end{proposition}
\begin{remarkk}
{\rm The relation \ref{com-rel} means that, for any $h\in H$, we have
$$
\nabla_h\int_0^1 K_s dW_s=\int_0^1\nabla_h
K_sdW_s+\int_0^1K_s{\dot{h}}(s)ds\,.
$$
}
\end{remarkk}
\begin{corollary}
If $\varphi=I_n(f_n)$, $f_n\in\hat{L}^2([0,1]^n)$, then we have, for $h\in H$,
$$
\nabla_h I_n(f_n)=n \int_{[0,1]^n}
f(t_1,\ldots,t_n)\dot{h}(t_n) dW_{t_1},\ldots,dW_{t_{n-1}}\,dt_n\,.
$$
\end{corollary}
\proof
Apply the above proposition $n$-times to the case in which, first $f_n$ is
$C^\infty([0,1]^n)$, then pass to the limit in $L^2(\mu)$.
\qed

\noindent
The following result will be extended in the sequel to much larger classes of
random variables:

\begin{corollary}
\label{I-W-decomp}
Let $\varphi:W\to\R$ be analytic in $H$-direction. Then we have
$$
\varphi=E[\varphi]+\sum_{n=1}^\infty {\tilde{I}}_n \left(
\frac{E[\nabla^n\varphi]}{n!}\right)\,,
$$
where $\tilde{I}_n(g)$, for a symmetric  $g\in H^{\otimes n}$, denotes
the multiple Wiener integral of
$$
\frac{\partial^n g}{\partial t_1\ldots \partial t_n}(t_1,\ldots,t_n)\,.
$$
In  other words  the kernel $\varphi_n\in\hat{L}^2[0,1]^n$ of the Wiener chaos
decomposition of $\varphi$ is equal to
$$
\frac{\partial^n}{\partial t_1\ldots \partial t_n}\frac{E[\nabla^n\varphi]}{n!}\,.
$$
\end{corollary}
\proof
We have, on one hand, for any $h\in H$,
$$
E[\varphi(w+h)]=E\Big[ \varphi\,\exp\int_0^1 \dot{h}_s dW_s-
{\textstyle\frac{1}{2}}\int_0^1\dot{h}_s^2ds\Big] =
     E[\varphi\,{\mathcal E}({\tilde{I}}(h))]\,.
$$
On the other hand, from Taylor's formula:
\begin{eqnarray*}
E[\varphi(w+h)] & = & E[\varphi]+\sum_1^\infty E\left[
    \frac{(\nabla^n\varphi(w),h^{\otimes n})}{n!}\right] \\
& = & E[\varphi]+\sum_1^\infty
     \frac{1}{n!}(E[\nabla^n\varphi],h^{\otimes n})_{ H^{\otimes n}} \\
& = & E[\varphi]+\sum_1^\infty \frac{1}{n!}
    \frac{E[{\tilde{I}}_n(E[\nabla^n\varphi])\,\,{\tilde{I}}_n(h^{\otimes n})]}{n!} \\
& = & E[\varphi]+\sum_1^\infty E\left[ \frac{{\tilde{I}}_n(E[\nabla^n\varphi])}{n!}
\frac{{\tilde{I}}_n(h^{\otimes n})}{n!}\right]
\end{eqnarray*}
hence, from the symmetry, we have
$${\tilde{I}}_n(\varphi_n)={\textstyle\frac{1}{n!}}{\tilde{I}}_n(E[\nabla^n\varphi])\,,$$

\noindent
where we have used the notation ${\tilde{I}}_1(h)={\tilde{I}}(h)=\int_0^1\dot{h}_s dW_s$ and
$$
{\tilde{I}}_n(\varphi_n)=\int_{[0,1]^n}
\frac{\partial^n\varphi_n}{\partial t_1\ldots\partial t_n}
(t_1,\ldots,t_n)dW_{t_1}\ldots dW_{t_n}\,.
$$
\qed

\section{The divergence operator}
The divergence operator, which is the adjoint of the Sobolev
derivative with respect to the Wiener measure, is one of the most
important tools of the Stochastic Analysis. We begin with  its formal
definition:
\begin{definition}
Let $\xi:W\to H$ be a random variable.\index[not]{de@$\delta$}
We say that $\xi\in\Dom_p(\delta)$, if
for any $\varphi\in \DD_{q,1}$ $(q^{-1}+p^{-1}=1)$, we have
$$
E[(\nabla\varphi,\xi)_H]\leq c_{p,q}(\xi).\|\varphi\|_q\,,
$$
and in this case we define $\delta\xi$ by
$$
E[\{\delta\xi\}\,\varphi]=E[(\xi,\nabla\varphi)_H]\,,
$$
i.e., $\delta\xi=\nabla^\ast\xi$, where $\nabla^\ast$ denotes the
adjoint of $\nabla$  with respect  to the Wiener  measure $\mu$, it is
called the \index[sub]{divergence operator} divergence operator.
\end{definition}
\remark
For the emergence of this operator cf. \cite{PK},  \cite{G-T} and  the
references there.

Let us give some properties of $\delta$:
\bigskip
\begin{enumerate}
\item[1.)]
Let $a:W\to\R$ be ``smooth", $\xi\in \Dom_p(\delta)$.
\index[not]{dz@$\Dom_p(\delta)$} Then we have, for any
$\varphi\in \DD_{q,1}\,$,
\begin{eqnarray*}
E\left[\delta(a\xi)\varphi\right] & = & E\left[(a\xi,\nabla\varphi)_H\right] \\
& = & E[(\xi,a\nabla\varphi)_H] \\
& = & E[(\xi,\nabla(a\varphi)-\varphi\,\,\nabla a)_H] \\
& = & E[(\delta\xi)\,\,a\varphi-\varphi\,\,(\nabla a,\xi)_H]\,,
\end{eqnarray*}

\noindent
hence
\begin{equation}
\label{identity-1}
\delta(a\xi)=a\delta\xi-(\nabla a,\xi)_H .
\end{equation}
\item[2.)]
Let $h\in H$, then we pretend that
$$
\delta h=\int_0^1\dot{h}(s)dW_s.
$$
To see this, it is sufficient to test this relation on the exponential
martingales: if $k\in H$, we have
\begin{eqnarray*}
\lefteqn{E\left[\delta h\,\,\exp\left\{\int_0^1\dot{k}_sdW_s-
{\textstyle\frac{1}{2}}\int_0^1\dot{k}_s^2ds\right\}\right]} \\
&=& E[(h,\nabla{\mathcal E}(I(k))_{H})] \\
& =& E[(h,k)_H {\mathcal E}(I(k))] \\
&=& (h,k)_{H}\,.
\end{eqnarray*}

\noindent
On the other hand, supposing first $h\in W^\ast$,
\begin{eqnarray*}
E[I(h)\,\,{\mathcal E}(I(k))] & = & E[I(h)(w+k)] \\
& = & E[I(h)]+(h,k)_{H} \\
& = & (h,k)_{H}\,.
\end{eqnarray*}

\noindent
Hence in particular, if we denote by $\tilde{1}_{[s,t]}$ the element of $H$
such that  $\tilde{1}_{[s,t]}(r)=\int_0^r 1_{[s,t]}(u)du$, we have that
\begin{equation}
\label{identity-2}
\displaystyle \delta(\tilde{1}_{[s,t]})=W_t-W_s\,.
\end{equation}

\item[3.)]
Let now $K$ be an adapted, simple  process
$$
K_t(w)=\sum_1^n a_i(w).1_{[t_i,t_{i+1}[}(t)\,,
$$
where $a_i\in \DD_{p,1}$ and ${\calB}_{t_i}$-measurable for each $i$. Let
$\tilde{K}$ be $\int_0^\cdot K_sds$. Then from  the identity
(\ref{identity-2}), we have
$$
\delta\tilde{K}=
\delta\Big( \sum_1^n a_i.\tilde{1}_{[t_i,t_{i+1}[}\Big)=
\sum_1^n\Big\{ a_i\delta(\tilde{1}_{[t_i,t_{i+1}[})-
(\nabla a_i,\tilde{1}_{[t_i,t_{i+1}[})\Big\}\,.
$$
From the relation (\ref{identity-2}), we have
$\delta(\tilde{1}_{[t_i,t_{i+1[}})=W_{t_{i+1}}-W_{t_i}\,$,
furthermore, from the Proposition \ref{mes-supp-prop}, the support of
$\nabla a_i$ is in $[0,t_i]$, consequently, we obtain
$$
\delta\tilde{K}=\sum_{i=1}^n a_i (W_{t_{i+1}}-W_{t_i})
    =\int_0^1 K_s dW_s\,.
$$
\end{enumerate}
\noindent
Hence we have the important result which says
that
\begin{theorem}
$\Dom_p(\delta)$ ($p>1$) contains the set consisting of  the
primitives of adapted stochastic  processes satisfying
$$
E\Big[\Big(\int_0^1 K_s^2ds\Big)^{p/2}\Big]<\infty\,.
$$
Moreover one has
$$
\delta \left\{\int_0^\cdot K_sds\right\}=\int_0^1K_s dW_s\,.
$$
\end{theorem}

\section{Local characters  of $\nabla$ and $\delta$}
\label{locality-section}
Before proceeding further, we shall
prove the locality of the Gross-Sobolev derivative and the divergence
operators in this section:
\begin{lemma}
\label{der-loc}
Let $\phi \in \DD_{p,1}$ for some $p>1$, then we have, for any
constant $c\in \reals$,
$$
\nabla\phi =0 \,{\mbox{ on}}\,\{\phi=c\}\,,
$$
almost surely.
\end{lemma}
\proof
Replacing $\phi$ by $\phi-c$, we may assume that $c=0$.
Let now  $f$ be a positive, smooth function of compact support on $\reals$
such that $f(0)=1$. Let $f_\eps(t)=f(t/\eps)$ and let $F_\eps$ be its
primitive. For any smooth, cylindrical, $H$-valued random variable
$u$, we have
\beaa
E[F_\eps(\phi)\,\delta u]&=&E[(\nabla F_\eps(\phi),u)_H]\\
  &=&E[f_\eps(\phi)(\nabla \phi,u)_H]\\
&\to&E[\won_{\{\phi=0\}}(\nabla \phi,u)_H]
\eeaa
as $\eps \to 0$. On the other hand $|F_\eps(\phi)|\leq
\eps\,\|f\|_{L^1(\reals, dt)}$, hence it follows that
$$
E[\won_{\{\phi=0\}}(\nabla \phi,u)_H]=0\,,
$$
since such $u$'s are dense in $L^q(\mu,H)$, the proof follows.
\qed

The divergence operator has an analogous property:
\begin{lemma}
\label{div-loc}
Assume that $u\in \Dom_p(\delta)$, $p>1$,  and  that the operator norm of
$\nabla u$, denoted by $\|\nabla u\|_{\op}$ is in $L^p(\mu)$. Then
$$
\delta u=0\,\,{\mbox{ a.s. on }} \{w\in W:\,u(w)=0\}\,.
$$
\end{lemma}
\proof
Let $f_\eps$ be as in the proof of Lemma \ref{der-loc}, then for any
cylindrical $\phi$, using the integration by parts formula:
\begin{eqnarray}
\label{ibp}
E\left[f_\eps\left(|u|_H^2\right)\delta
  u\,\phi\right]&=&E\left[f'_\eps\left(|u|_H^2\right)\left(u,\nabla
    |u|_H^2\right)_H\,\phi\right]\nonumber\\
&&+E\left[f_\eps\left(|u|_H^2\right)(u,\nabla \phi)_H\right]\,.
\end{eqnarray}
Note that
\beaa
\left|f'_\eps\left(|u|_H^2\right)\left(u,\nabla
    |u|_H^2\right)_H\right|&\leq&
|u|_H^2\left|f'_\eps\left(|u|_H^2\right)\right|\,\|\nabla u\|_\op\\
&\leq&\eps \sup_\reals|xf'(x)|\|\nabla u\|_\op\,.
\eeaa
Hence from the dominated convergence theorem, the first term at the
right of (\ref{ibp}) tends to zero with $\eps$. Evidently the second
one also converges to zero and this completes the proof.
\qed

\begin{remarkk}
{\rm Using the local character of the Sobolev derivative one can
  define the local Sobolev spaces as we shall see later.}
\end{remarkk}

\section{The Ornstein-Uhlenbeck Operator}
For a nice function $f$ on $W$, $t\geq 0$, we define
\begin{equation}
\label{mehler-formula}
P_tf(x)=\int_W f\left(e^{-t} x+\sqrt{1-e^{-2t}}\,y\right)\mu(dy)\,,
\end{equation}
this expression for $P_t$ is called Mehler's formula.
\index[sub]{me@Mehler's formula}
\noindent
Since $\mu(dx)\mu(dy)$ is invariant under the rotations of $W\times W$, i.e.,
$(\mu\times\mu)(dx,dy)$ is invariant under the transformation
$$
T_{t}(x,y)=\left(xe^{-t}+y(1-e^{-2t})^{1/2},x(1-e^{-2t})^{1/2}-ye^{-t}\right),
$$
we have obviously
\begin{eqnarray*}
\|P_tf(x)\|_{L^p(\mu)}^p & \leq & \int\!\!\int
      |(f\otimes1)(T_t(x,y))|^p \mu(dx)\mu(dy) \\
& = & \int\!\!\int|(f\otimes1)(x,y)|^p\mu(dx)\mu(dy) \\
& = & \int|f(x)|^p\mu(dx)\,,
\end{eqnarray*}
for any $p\geq1$, $\|P_tf\|_{L^p}\leq\|f\|_{L^p}\,$; hence also for $p=\infty$
by duality. A straightforward calculation gives that, for any $h\in H\cap
W^\ast\,(=W^\ast)$,
\begin{eqnarray*}
P_t({\mathcal E}(I(h)) & = & {\mathcal E}(e^{-t} I(h)) \\
& = & \sum_{n=0}^\infty e^{-nt}\frac{I_n(h^{\otimes n})}{n!}\,.
\end{eqnarray*}
Hence, by homogeneity, we have
$$P_t(I_n(h^{\otimes n}))=e^{-nt}I_n(h^{\otimes n})$$
and by a  density argument, we obtain
$$
P_tI_n(f_n)=e^{-nt}I_n(f_n)\,,
$$
for any $f_n\in\hat{L}^2([0,1]^n)$. Consequently $P_s\circ P_t=P_{s+t}\,$,
i.e., $(P_t)$\index[not]{p@$P_t$} is a measure preserving Markov
semi-group. Its infinitesimal
generator is denoted by  $-{\mathcal L}$ and  is ${\mathcal L}$ is  called the
Ornstein-Uhlenbeck \index[sub]{o@Ornstein-Uhlenbeck operator}
 or the number  operator\index[sub]{n@number operator}. Evidently,  we
have
\begin{equation}
\label{no-relation}
{\mathcal L}I_n(f_n)=n I_n(f_n)
\end{equation}
and this relation means that  the Wiener chaos are  its eigenspaces.
From the definition, it follows directly that (for $a_i$ being
${\mathcal F}_{t_i}$-measurable)
$$
 P_t\left(\sum a_i(W_{t_i+1}-W_{t_i})\right)=e^{-t}\sum(P_ta_i)(W_{t_i+1}-W_{t_i}),
$$
that is to say
$$
 P_t\int_0^1 H_sdW_s=e^{-t}\int_0^1 P_tH_sdW_s,
$$
and by differentiation
\begin{equation}
\label{com-rel1}
{\mathcal L}\int_0^1 H_sdW_s=\int_0^1(I+{\mathcal L})H_sdW_s\,.
\end{equation}
Also we have
\begin{equation}
\label{com-rel2}
\nabla P_t\varphi=e^{-t}P_t\nabla\varphi\,.
\end{equation}

The following lemma is a consequence of the relation (\ref{no-relation}):
\begin{lemma}
Assume that  $\phi\in L^2(\mu)$ with  the Wiener chaos representation
$$
\phi=\sum_{n=0}^\infty I_n(\phi_n)
$$
satisfying
$$
\sum_{n=1}^\infty n\,(n!)\|\phi_n\|_{H^{\circ n}}^2<\infty\,.
$$
Then
$$
\delta\circ\nabla\phi={\mathcal L}\phi\,,
$$
 where $\delta$ is the divergence operator \footnote{Sometimes, in the
   classical case,  it is also called \index[sub]{Hitsuda-Ramer-Skorohod
    integral} Hitsuda-Ramer-Skorohod integral.}.
\end{lemma}
\proof
It is sufficient to prove for  $\varphi={\mathcal E}(I(h))$. In this
case  from the identity (\ref{identity-1})
\begin{eqnarray*}
(\delta\circ\nabla)\varphi & = & \delta(h\,\,{\mathcal E}(I(h))) \\
& = &\left[I(h)-|h|_H^2\right]{\mathcal E}(I(h))\\
&=&\L{\mathcal{E}}(I(h))\,.
\end{eqnarray*}
\vspace*{-4ex}

\qed
\bigskip
\begin{remarkk}{\rm
Let us define for the smooth functions $\varphi$, a semi-norm
\index[not]{l@${\mathcal L}$}
$$
|\!|\!|\varphi|\!|\!|_{p,k}=\|(I+{\mathcal L})^{k/2}\varphi\|_{L^p(\mu)}\,.
$$
At  first glance, these semi-norms (in fact norms), seem different from the
one defined  by
$\|\varphi\|_{p,k}=\sum_0^k\|\nabla^j\varphi\|_{L^p(\mu,H^{\otimes
j})}\,$. We will show in the next chapters that they are
equivalent. }
\end{remarkk}



\qed
\section{Exercises}{\footnotesize{
These exercises are aimed to give some useful formulas about the
iterated divergence operator and related commutation properties.
\begin{enumerate}
\item Prove that
\begin{equation}
\label{com-rel-1}
\nabla P_t\phi=e^{-t}P_t\nabla \phi
\end{equation}
and
\begin{equation}
\label{com-rel-2}
P_t\delta u=e^{-t}\delta P_tu
\end{equation}
for any $\phi\in \DD_{p,1}$ and $u\in \DD_{p,1}(H)$.
\item Assume that $u:W\to H$ is a cylindrical random variable. Prove
  that
$$
\delta u=\sum_{i=1}^\infty \left\{(u,e_i)_H\delta e_i-\nabla_{e_i}(u,e_i)_H\right\}\,,
$$
for any complete, orthonormal basis $(e_i,i\in \NN)$ of $H$. In
particular, in the finite dimensional case we can write
$$
\delta u(w)=<u(w),w>-\trace \nabla u(w)\,,
$$
although in infinite dimensional case  such an expression is meaningless in
general. In case the $\trace \nabla u$ exists, the remaining part is
called the Stratonovitch integral. \index[sub]{Stratonovitch integral}

\item Assume that $u:W\to H$ is a cylindrical random variable. Prove
  that
$$
E[(\delta u)^2]=E[|u|_H^2]+E[\trace(\nabla u\,\nabla u)].
$$
\item Let $u$ be as above, prove the identity
$$
\delta^2 u^{\otimes 2}=(\delta u)^2-|u|_H^2-\trace(\nabla u\,\nabla
u)-2\delta (\nabla_uu)\,,
$$
where $\delta^2u^{\otimes 2}$ is defined by the integration by parts
formula as
$$
E[\delta^2u^{\otimes 2}\,\phi]=E[(\nabla^2\phi,u^{\otimes 2})_2]\,,
$$
for any test function $\phi$
and  $(\cdot,\cdot)_2$ denotes the inner product of the space of
Hilbert-Schmidt operators on $H$. Prove that more generally one has
\beaa
\delta \alpha\,\delta
\beta&=&\delta^2(\alpha\otimes\beta)+\trace(\nabla\alpha
\nabla\beta)\\
&&+\delta(\nabla_\alpha\beta+\nabla_\beta\alpha)+(\alpha,\beta)_H\,,
\eeaa
where $\alpha$ and $\beta$ are two $H$-valued, cylindrical random variables.
\item With the same hypothesis as above, show that one has
$$
\delta^{p+1}u^{\otimes p+1}=\delta u\,\delta^pu^{\otimes p}-\nabla_u
(\delta^pu^{\otimes p})-\delta(\nabla^p_{u^{\otimes p}}u)\,.
$$
\item For a $u:W\to H$ as above, prove that
\beaa
(\delta u)^p&=&\delta\left(u\,(\delta u)^{p-1}\right)\\
&&+(\delta
u)^{p-2}\left[(p-1)|u|_H^2+(p-2)\left(\delta(\nabla_uu)+
\trace(\nabla u\,\nabla u)\right)\right]
\eeaa
for any $p\in \NN$.
\end{enumerate}}}

\section*{Notes and suggested reading}
{\footnotesize{The notion of derivation in the setting of a Gaussian
  measure on an infinite dimensional setting can be found in the books
  of Quantum Field Theory, cf. \cite{Si}  also \cite{PK} and the references
  there. It has also been studied in a little bit more restricted case
  under the name $H$-derivative by L. Gross, cf. also \cite{KUO},
  \cite{PK}. However
  the full  use  of the quasi-invariance
  with respect to the translations from the Cameron-Martin space
  combined with the $L^p$-closure  of it in the sense of Sobolev has
  become popular with the advent of the stochastic calculus of
  variations of Paul Malliavin: cf. \cite{PAM}, \cite{Sh}, \cite{Mall-1}.}
}


\chapter{Meyer Inequalities}
\label{ch.meyer}
\markboth{Meyer Inequalities}{Riesz Transform}
\section*{Meyer Inequalities and Distributions}

Meyer inequalities are essential to control the Sobolev norms defined with the
Sobolev derivative with the norms defined via the Ornstein-Uhlenbeck operator.
They can be summarized as the equivalence of the two  norms defined on
the (real-valued) Wiener functionals as
$$
|||\phi|||_{p,k}=\sum_{i=0}^{k}\|\nabla^i\phi\|_{L^p(\mu,H^{\otimes i})},
$$
and
$$
\|\phi\|_{p,k}=\|(I+{\mathcal L})^{k/2}\phi\|_{L^p(\mu)},
$$
for any $p>1$ and $k\in \N$.\index[sub]{Meyer inequalities}
The key point is the continuity property of the Riesz transform on
$L^p([0,2\pi],dx)$, i.e., from a totally analytic origin, although the original
proof of P.~A.~Meyer was probabilistic (cf. \cite{PAM}). Here we
develop  the proof suggested by \cite{F}.

\section{Some Preparations}
\index[not]{p@$p.v.$}
Let $f$ be a function on $[0,2\pi]$, extended  to the whole $\R$ by
periodicity. We denote by $\tilde{f}(x)$ the function defined by

\begin{equation}
\label{pv}
\tilde{f}(x)=\frac{1}{\pi}{\rm{p.v.}}\int_0^\pi \frac{f(x+t)-f(x-t)}{2\tan
t/2}dt\,,
\end{equation}\index[not]{p@$p.v.$}
where p.v. denotes the the principal value of the integral in
(\ref{pv}). The famous theorem of M.~Riesz, cf. \cite{Z},  asserts
that, for any 
$f\in L^p[0,2\pi]$,
$\tilde{f}\in L^p([0,2\pi])$, for $1<p<\infty$ with
$$
\|\tilde{f}\|_p\leq A_p\|f\|_p\,,
$$
where $A_p$ is a constant depending only on $p$. Most of the classical
functional analysis of the 20-th century has been devoted to extend this result
to the case where the function $f$ was taking its values in more abstract
spaces than the real line. We will show that our problem also can be reduced to
this one.

In fact, the main result that we are going to show will be that
$$
\|\nabla(I+{\mathcal L})^{-1/2}\varphi\|_p\approx\|\varphi\|_p
$$

\noindent
by rewriting $\nabla(I+{\mathcal L})^{-1/2}$ as an $L^p(\mu,H)$-valued Riesz
transform. For this we need first, the following elementary

\begin{lemma}
\label{perturb-lem}
Let $K$ be any function on $[0,2\pi]$ such that
$$
K(\theta)-{\textstyle\frac{1}{2}}\cot{\textstyle\frac{\theta}{2}}\in
L^\infty([0,\pi])\,,
$$
then the operator $f\to T_Kf$ defined by
$$
T_Kf(x)=\frac{1}{\pi}p.v.\int_0^{\pi}(f(x+t)-f(x-t))K(t)dt
$$
is again a bounded operator on $L^p([0,2\pi])$ with
$$
\left\|T_Kf\right\|_p\leq B_p\|f\|_p\,\quad \mbox{for any $p\in(1,\infty)$}
$$
where $B_p$ depends only on $p$.
\end{lemma}
\proof
In fact we have
\begin{eqnarray*}
\left|T_Kf-\tilde{f}\right|(x) & \leq &
  \frac{1}{\pi}\int_0^\pi\left|f(x+t)-f(x-t)\right|\,\left|K(t)
      -{\textstyle\frac{1}{2}}\cot{\textstyle\frac{t}{2}}\right|dt \\
& \leq & c\left\|f\right\|_{L^p}\,\left\|K
       -{\textstyle\frac{1}{2}}\cot{\textstyle\frac{\theta}{2}}\right\|_{L^\infty}\,.
\end{eqnarray*}

\noindent
Hence
$$
\left\|T_Kf\right\|_p\leq
\left(c\left\|K
-{\textstyle\frac{1}{2}}\cot{\textstyle\frac{\theta}{2}}\right\|_{L^\infty}+A_p\right) \|f\|_p.
$$
\qed

\begin{remarkk}
\label{rem-1}
{\rm
If for some $a\not=0$, $\,aK(\theta)-\frac{1}{2}\cot\frac{\theta}{2}\in
L^\infty([0,2\pi])$, then we have
\begin{eqnarray*}
\left\|T_Kf\right\|_p & = & \frac{1}{|a|}\left\|aT_Kf\right\|_p\\
&\leq&\frac{1}{|a|}
        \left\{ \left\|aT_Kf-\tilde{f}\right\|_p+
\left\|\tilde{f}\right\|_p\right\} \\
& \leq & \frac{1}{|a|}\left\{\left\|aK-
{\textstyle\frac{1}{2}\cot\frac{\theta}{2}}\right\|_{L^\infty}\,
\|f\|_p+A_p\|f\|_p\right\} \\
& \leq & c_p\|f\|_p
\end{eqnarray*}
with another constant $c_p$.}
\end{remarkk}

\begin{corollary}
\label{cor-K}
Let $K$ be a function on $[0,\pi]$ such that $K=0$ on
$\Big[\frac{\pi}{2},\pi\Big]$ and $K-\frac{1}{2}\cot\frac{\theta}{2}\in
L^\infty\Big(\Big[0,\frac{\pi}{2}\Big]\Big)$. Then $T_K$ defined by
$$
T_Kf(x)=\int_0^{\pi/2}\left[f(x+t)-f(x-t)\right]K(t)dt
$$
is continuous from $L^p([0,2\pi])$ into itself for any $p\in[1,\infty\,)\,$.
\end{corollary}
\proof
We have
$$cK(\theta)1_{[0,\frac{\pi}{2}]}-\frac{1}{2}\cot
\frac{\theta}{2}\in L^\infty([0,\pi])
$$
since on the interval
$\Big[\frac{\pi}{2},\pi\Big]$, $\sin\frac{\theta}{2}\in
\Big[\frac{\sqrt{2}}{2},1\Big]$, then the result follows from the
Lemma \ref{perturb-lem}.
\qed

\section{$\nabla(I+{\mathcal L})^{-1/2}$ as  the Riesz Transform
  \index[sub]{Riesz transformation}}

Let us denote by $R_\theta(x,y)$ the rotation on $W\times W$ defined by
$$
R_\theta(x,y)=\Bigl(x\cos\theta+y\sin\theta,-x\sin\theta+y\cos\theta\Bigr)\,.
$$
Note that $R_\theta\circ R_\phi=R_{\phi+\theta}\,$. We have also, putting
$e^{-t}=\cos\theta$,
\begin{eqnarray*}
P_tf(x) & = & \int_W f(e^{-t}x+\sqrt{1-e^{-2t}}\,y)\mu(dy) \\
& = & \int_W(f\otimes1)(R_\theta(x,y))\mu(dy) \\
&=&P_{-\log\cos\theta}f(x)\,.
\end{eqnarray*}

\noindent
Let us now calculate $(I+{\mathcal L})^{-1/2}\varphi$ using this
transformation:
\begin{eqnarray*}
\lefteqn{(I+{\mathcal L})^{-1/2}\varphi(x) =
   \int_0^\infty t^{-1/2}e^{-t} P_t\varphi(x)dt}\\
& & = \int_0^{\pi/2}(-\log\cos\theta)^{-1/2}\cos\theta\cdot
\int_W(\varphi\otimes1)(R_\theta(x,y))\mu(dy)\tan\theta d\theta \\
& & = \int_W\mu(dy)\left[\int_0^{\pi/2}(-\log\cos\theta)^{-1/2}
     \sin\theta(\varphi\otimes1)(R_\theta(x,y))d\theta\right]\,.
\end{eqnarray*}

\noindent
On the other hand, we have, for $h\in H$
\begin{eqnarray*}
\lefteqn{\nabla_hP_t\varphi(x)} \\
& = &\!
     \frac{d}{d\lambda}P_t\varphi(x+\lambda h)|_{\lambda=0} \\
& = &\! \frac{d}{d\lambda}\int\varphi\left(e^{-t}(x+\lambda h)+
               \sqrt{1-e^{-2t}}\,y\right)\mu(dy)|_{\lambda=0} \\
& = &\! \frac{d}{d\lambda}\int\varphi\left(e^{-t} x+\sqrt{1-e^{-2t}}
   \Big(y+\frac{\lambda e^{-t}}{\sqrt{1-e^{-2t}}}h\Big)\right)
       \mu(dy)|_{\lambda=0} \\
& = &\! \frac{d}{d\lambda}\int\varphi\left(e^{-t}x+\sqrt{1-e^{-2t}}\,y\right)
   {\mathcal E}\Big(\frac{\lambda e^{-t}}{\sqrt{1-e^{-2t}}}
        I(h)\Big)(y)\mu(dy)|_{\lambda=0} \\
& = &\! \frac{e^{-t}}{\sqrt{1-e^{-2t}}} \int_W\varphi
\left(e^{-t}x+\sqrt{1-e^{-2t}}\,y\right)\delta h(y)\,\mu(dy)\,.
\end{eqnarray*}
Therefore
\begin{eqnarray*}
\lefteqn{\nabla_h(I+{\mathcal L})^{-1/2}\varphi(x)}\\
& = &\!\!\int_0^\infty t^{-1/2}e^{-t}\nabla_h P_t\varphi(x)dt \\
& = &\!\!  \int_0^\infty t^{-1/2}\frac{e^{-2t}}{\sqrt{1-e^{-2t}}}
    \int_W \delta h(y)\varphi\left(e^{-t} x+
           \sqrt{1-e^{-2t}}\,y\right)\mu(dy)dt \\
& = &\!\!  \int_0^{\pi/2}(-\log\cos\theta)^{-1/2}
    \frac{\cos^2\theta}{\sin\theta}\tan\theta \!\int \delta
           h(y)\left(\varphi\otimes1\right)\left(R_\theta(x,y)\right)\mu(dy)d\theta \\
& = &\!\!  \int_0^{\pi/2}(-\log\cos\theta)^{-1/2}\cos\theta \int_W
   \delta h(y)\cdot(\varphi\otimes1)(R_\theta(x,y))\mu(dy)d\theta
\end{eqnarray*}

\noindent
Since $\mu(dy)$ is invariant under the transformation $y\mapsto -y$, we have
$$
\int\delta h(y)(\varphi\otimes1)(R_\theta(x,y))\mu(dy)=-\int \delta
h(y)(\varphi\otimes1)(R_{-\theta}(x,y))\mu(dy),
$$
therefore:
\begin{eqnarray*}
\lefteqn{\nabla_h(I+{\mathcal L})^{-1/2}\varphi(x)} \\
&= & \!\int_0^{\pi/2}(-\log\cos\theta)^{-1/2}.\\
&& \quad \int\delta h(y)
    \frac{(\varphi\otimes1)(R_\theta(x,y))
              -(\varphi\otimes1)(R_{-\theta}(x,y))}{2}\cos\theta \mu(dy)d\theta\\
&= &\! \int_W\delta h(y)\int_0^{\pi/2}
     K(\theta)\left({(\varphi\otimes1)(R_\theta(x,y))-(\varphi\otimes1)(R_{-\theta}(x,y))}\right)
     d\theta\mu(dy)\, ,
\end{eqnarray*}

\noindent
where $K(\theta)=\frac{1}{2}\cos\theta(-\log\cos\theta)^{-1/2}$.

\begin{lemma}
\label{Lemma-2}
We have
$$
2K(\theta)-\cot\frac{\theta}{2}\in L^\infty((0,\pi/2]).
$$
\end{lemma}
\proof
The only problem is when $\theta\to0$. To see this let us put
$e^{-t}=\cos\theta$, then
$$
\cot\frac{\theta}{2}=\frac{\sqrt{1+e^{-t}}}{\sqrt{1-e^{-t}}}
\approx\frac{2}{\sqrt{t}}\,
$$ and
$$
K(\theta)=\frac{e^{-t}}{\sqrt{t}}\approx\frac{1}{\sqrt{t}}
$$
hence
$$
2K(\theta)-\cot\frac{\theta}{2}\in L^\infty\left(\left[0,\frac{\pi}{2}\right]\right)\,.
$$
\qed

\noindent
Using Lemma \ref{perturb-lem},  Remark \ref{rem-1} following it and
Corollary \ref{cor-K}, we see that the map
$f\mapsto p.v.\int_0^{\pi/2}(f(x+\theta)-f(x-\theta))K(\theta)d\theta$
is a bounded map from $L^p[0,\pi]$ into itself. Moreover
\begin{lemma}
\label{lemma-3}
Let $F:W\times W\to \R$ be a measurable, bounded function. Define $TF(x,y)$ as
$$
TF(x,y)=p.v.\int_0^{\pi/2}\left[F\circ R_\theta(x,y)-F\circ R_{-\theta}(x,y)\right]
K(\theta)d\theta\,.
$$
Then, for any $p>1$, there exists some $c_p>0$ such that
$$
\|TF\|_{L^p(\mu\times\mu)}\leq c_p\|F\|_{L^p(\mu\times\mu)}\,.
$$
\end{lemma}
\proof
We have
$$
(TF)(R_\beta(x,y))=p.v.\int_0^{\pi/2} (F(R_{\beta+\theta}(x,y))
-F(R_{\beta-\theta}(x,y)))K(\theta)d\theta\,,
$$
this is the Riesz transform for fixed $(x,y)\in W\times W$, hence we have
$$
\int_0^{\pi/2}|TF(R_\beta(x,y))|^pd\beta\leq
c_p\int_0^\pi|F(R_\beta(x,y))|^p d\beta\,,
$$
taking the expectation with respect to $\mu\times\mu$, which is invariant under
$R_\beta\,$, we have
\begin{eqnarray*}
E_{\mu\times\mu}\int_0^\pi|TF(R_\beta(x,y))|^pd\beta
 & = &
   E_{\mu\times\mu}\int_0^\pi |TF(x,y)|^pd\beta \\
& = & \frac{\pi}{2}E[|TF|^p] \\
& \leq & c_p E\int_0^\pi|F(R_\beta(x,y))|^p d\beta\\
& = & \pi c_pE[|F|^p]\,.
\end{eqnarray*}
\qed

\noindent
We have
\begin{theorem}
\label{th-1}
$\nabla\circ (I+{\mathcal L})^{-1/2}:L^p(\mu)\to L^p(\mu,H)$ is a
linear continuous operator  for any $p\in (1,\infty)$.
\end{theorem}
\proof
With the notations of Lemma~\ref{lemma-3}, we have
$$
\nabla_h(I+{\mathcal L})^{-1/2}\varphi= \int_W \delta h(y)\,
T(\varphi\otimes1)(x,y)\mu(dy)\,.
$$
From H\"older  inequality:
\begin{eqnarray*}
|\nabla_h(I+\calL)^{1/2}\phi(x)|&\leq&\|\delta
h\|_q\left(\int_W|T(\phi\otimes 1)(x,y)|^p\mu(dy)\right)^{1/p}\\
&\leq&c_p |h|_H \left(\int_W|T(\phi\otimes 1)(x,y)|^p\mu(dy)\right)^{1/p}\,,
\end{eqnarray*}
where the last inequality follows from the fact that  $y\to \delta
h(y)$ is an $N_1(0,|h|_H^2)$--Gaussian random variable. Hence
$$
|\nabla(I+\calL)^{-1/2}\phi(x)|_H\leq
\left(\int_W|T(\phi\otimes 1)(x,y)|^p\mu(dy)\right)^{1/p}\,
$$
consequently, from Lemma \ref{lemma-3}
\beaa
\|\nabla(I+\calL)^{-1/2}\phi\|_p^p&\leq&\int_{W\times W}|T(\phi\otimes
1)(x,y)|^p\mu(dx)\mu(dy)\\
&\leq&\|\phi\otimes 1\|_{L^p(\mu\times \mu)}^p\\
&=&\|\phi\|_p^p
\eeaa
and this completes the proof.
\qed




\begin{corollary}
\label{cor1}
We have
$$
\|(I+{\mathcal L})^{-1/2}\delta\xi\|_p\leq c_p\|\xi\|_p\,,
$$
for any $\xi\in L^p(\mu;H)$ and  for any  $ p\in (1,\infty)$.
\end{corollary}
\proof
It suffices to  take the adjoint of $\nabla(I+{\mathcal L})^{-1/2}$.
\qed

\begin{corollary}
The following identities are valid for any $\varphi\in \DD$:
\begin{enumerate}
\item
$\|\nabla\varphi\|_p\leq c_p\|(I+{\mathcal L})^{1/2}\varphi\|_p$

\item
$\|(I+{\mathcal L})^{1/2}\varphi\|_p\leq\tilde{c}_p
\left(\|\varphi\|_p+\|\nabla\varphi\|_p\right)$,
\end{enumerate}
where $c_p$ and $\tilde{c}_p$ are two constants independent of $\varphi$.
\end{corollary}
\proof
The first identity follows easily as
\beaa
\|\nabla\varphi\|_p&=&\|\nabla(I+{\mathcal L})^{-1/2}(I+{\mathcal
  L})^{1/2}\varphi\|_p \\
&\leq& c_p\|(I+{\mathcal L})^{1/2}\varphi\|_p\,.
\eeaa
To prove the second we have
\begin{eqnarray*}
\lefteqn{\|(I+{\mathcal L})^{1/2}\varphi\|_p=\|(I+{\mathcal L})^{-1/2}(I+{\mathcal L})\varphi\|_p}\\
& & = \|(I+{\mathcal L})^{-1/2}(I+\delta  \nabla)\varphi\|_p \\
& & \leq \|(I+{\mathcal L})^{-1/2}\varphi\|_p+ \|(I+{\mathcal L})^{-1/2}
     \delta\nabla\varphi\|_p \\
& & \leq \|\varphi\|_p+c_p\|\nabla\varphi\|_p,
\end{eqnarray*}
where the last inequality follows from Corollary  \ref{cor1}.
\qed

\section*{Notes and suggested reading}
{\footnotesize{ The inequalities studied in this chapter are due to
  P. A. Meyer in his seminal paper \cite{PAM}. He discusses at the
  last part of it already about the space of test functions defined by
  the Ornestein-Uhlenbeck operator and proves that this space is an
  algebra. Then the classical  duality  results give birth immediately
  to the space of the distributions on the Wiener space, and this is
  done in \cite{W1}. Later the proof of P. A. Meyer has been
  simplified by several people. Here we have followed an idea of
  D. Feyel, cf. \cite{F}.

}}


\chapter{Hypercontractivity\index[sub]{Hypercontractivity}}
\label{ch.hyper}
\markboth{Hypercontractivity}{}

\section*{Introduction}
We know that the semi-group  of Ornstein-Uhlenbeck is a bounded
operator on $L^p(\mu)$, for any $p\in[1,\infty]$. In fact for
$p\in(1,\infty)$, it is more than bounded. It increases the degree of
integrability, this property is called
{\bf{hypercontractivity}} and it  is used to show the continuity of
linear operators on $L^p(\mu)$-spaces  defined via the Wiener chaos
decomposition or the spectral decomposition of the Ornstein-Uhlenbeck
operator. We shall use it in the next chapter to complete the proof of the
 Meyer inequalities. Hypercontractivity  has been first discovered by
 E. Nelson, here we follow  the proof given by \cite{JN}. We complete
 the chapter by an analytic proof of the
{\bf {logarithmic Sobolev inequality}} of Leonard Gross (cf. \cite{LG},
\cite{D-ST}) for which we shall give another proof in the fifth chapter.

\section{Hypercontractivity via It\^o Calculus}
In the sequel we shall show that this result can be proved using the Ito
formula. Let $(\Omega,{\mathcal A},P)$ be a probability space with
$({\calB}_t;t\in\R_+)$ being a filtration. We take two Brownian
motions $(X_t;t\geq0)$ and $(Y_t;t\geq0)$ which are not necessarily
independent, i.e., $X$ and $Y$ are two continuous, real martingales
such that  $(X_t^2-t)$ and $(Y_t^2-t)$ are again martingales (with
respect  to $({\calB}_t)$) and that $X_t-X_s$ and $Y_t-Y_s$ are
independent  of ${\calB}_s$, for $t>s$. Moreover there exists
($\rho_t;t\in\R_+)$, progressively measurable with values in $[-1,1]$
such that
$$
(X_tY_t-\int_0^t\rho_s ds,t\geq0)
$$
is again a $({\calB}_t)$-martingale. Let us denote by
$$
{\Xi}_t=\sigma(X_s;s\leq t),\quad {\calY}_t=\sigma(Y_s;s\leq t)\,
$$
i.e., the corresponding filtrations of $X$ and $Y$ and by $\Xi$
 and by $\calY$ their respective supremum.

\begin{lemma}
\begin{enumerate}
\item  For any $\varphi\in L^1(\Omega,{\Xi},P)$, $t\geq0$, we have
$$
E[\varphi|{\calB}_t]=E[\varphi|{\Xi}_t]\;\mbox{a.s.}
$$
\item For any $\psi\in L^1(\Omega,{\calY},P)$, $t\geq0$, we have
$$
E[\psi|{\calB}_t]=E[\psi|{\calY}_t]\;\mbox{a.s.}
$$
\end{enumerate}
\end{lemma}
\proof
Since the two claims are similar, we shall prove only the first one.
From Paul  L\'evy's theorem, we have also that $(X_t)$ is an
$({\Xi}_t)$-Brownian motion. Hence
$$
\varphi=E[\varphi]+\int_0^\infty H_s dX_s
$$
where $H$ is $({\Xi}_t)$-adapted process. Hence
$$
E[\varphi|{\calB}_t]=E[\varphi]+ \int_0^t H_sdX_s=
E[\varphi|{\Xi}_t]\,.
$$
\qed

\noindent
Let $T$ be  the operator
$T:L^1(\Omega,{\Xi},P)\to L^1(\Omega,{\calY},P)$ defined as  the
restriction of $E[\,\cdot\,|{\calY}]$ to the space
$L^1(\Omega,{\Xi},P)$. We know that $T:L^p({\Xi})\to L^p({\calY})$ is
a contraction for any $p\geq1$. If we impose supplementary
conditions to $\rho$, then we have more:
\begin{proposition}
If $|\rho_t(w)|\leq r$ $(dt\times \,dP$ a.s.) for some $r\in[0,1]$,
then $T:L^p({\Xi})\to L^q({\calY})$ is a bounded operator, where
$$
p-1\geq r^2(q-1)\,.
$$
\end{proposition}
\proof
$p=1$ is already known. So suppose $p,q\in]1,\infty[\,$. Since $L^\infty({\Xi})$ is dense in $L^p({\Xi})$, it is enough to prove that
$\|TF\|_q\leq\|F\|_p$ for any $F\in L^\infty({\Xi})$. Moreover, since $T$ is
a positive operator, we have $|T(F)|\leq T(|F|)$, hence we can work as well
with $F\in L_+^\infty({\Xi})$.
Due to  the duality between $L^p$-spaces, it suffices to show that
$$
E[T(F)G]\leq\|F\|_p\|G\|_{q'}\,,\qquad
\Big( \frac{1}{q'}+\frac{1}{q}=1\Big),
$$
\noindent
for any $F\in L_+^\infty({\Xi})$, $G\in
L_+^\infty({\calY})$. Since bounded and positive  random variables are
dense  in all $L_+^p$  for any $p>1$, we can suppose without loss of
generality that $F,G\in[a,b]$ almost surely for some  $0<a<b<\infty$.
Let
\beaa
M_t&=&E[F^p|{\Xi}_t]\\
N_t&=&E[G^{q'}|{\calY}_t]\,.
\eeaa
Then, from the Ito representation theorem we have
$$M_t=M_0+\int_0^t\phi_s dX_s$$
$$N_t=N_0+\int_0^t\psi_s dY_s$$

\noindent
where $\phi$ is $\Xi$-adapted, $\psi$ is $\calY$-adapted, $M_0=E[F^p]$,
$N_0=E[G^{q'}]$. From the Ito formula, we have
\begin{eqnarray*}
M_t^\alpha N_t^\beta=M_0^\alpha N_0^\beta & + & \int_0^t \alpha
   M_s^{\alpha-1}N_s^\beta dM_s +
             \beta\int_0^t M_s^\alpha N_s^{\beta-1}dN_s+\\
& + & \frac{1}{2}\int_0^t M_s^\alpha N_s^\beta A_s ds
\end{eqnarray*}
where
$$A_t=\alpha(\alpha-1) \Big(\frac{\phi_t}{M_t}\Big)^2
+2\alpha\beta\frac{\phi_t}{M_t}\frac{\psi_t}{N_t}\rho_t+\beta(\beta-1)
\Big(\frac{\psi_t}{N_t}\Big)^2$$
and $\alpha=\frac{1}{p}\,$, $\beta=\frac{1}{q'}\,$.
To see this it suffices to use  the Ito formula as
\begin{eqnarray*}
M_t^\alpha & = & M_0^\alpha+\alpha\int_0^t M_s^{\alpha-1}\phi_s dX_s+
         \frac{\alpha(\alpha-1)}{2}\int_0^t M_s^{\alpha-2}\phi_s^2ds \\
N_t^\beta & = & \cdots
\end{eqnarray*}
and then as
\begin{eqnarray*}
\lefteqn{M_t^\alpha N_t^\beta - M_0^\alpha N_0^\beta}\\
& = & \int_0^t M_s^\alpha
dN_s^\beta + \int_0^t N_s^\beta dM_s^\alpha + \alpha\beta
\int_0^t M_s^{\alpha-1}N_s^{\beta-1}\phi_s\psi_s\rho_s ds \\
& = &\int_0^t M_s^\alpha\left(\beta N_s^{\beta-1}\psi_s
         dY_s+\frac{\beta(\beta-1)}{2}N_s^{\beta-2}\psi_s^2 ds\right) \\
& & +\int_0^t N_s^\beta\left(\alpha M_s^{\alpha-1}\phi_s dX_s +
        \frac{\alpha(\alpha-1)}{2} M_s^{\alpha-2}\phi_s^2 ds\right) \\
& & + \alpha\beta\int_0^t M_s^{\alpha-1}
         N_s^{\beta-1}\phi_s\psi_s\rho_s ds
\end{eqnarray*}
and finally to pick up  together all the  integrands integrated with
respect to the Lebesgue measure $ds$.
\newline
As everything is square integrable, it comes
\begin{eqnarray*}
E[M_ \infty^\alpha N_ \infty^\beta] & = & E\Big[ E[F^p|{\Xi}_
\infty]^\alpha
      \cdot E[G^{q'}|{\calY}_ \infty]^\beta\Big] \\
& = & E[F\cdot G] \\
& = & \frac{1}{2}\int_0^\infty E[N_t^\beta M_t^\alpha A_t]dt +E
    M_0^\alpha N_0^\beta \\
& = & E[F^p]^\alpha E[G^{q'}]^\beta+\frac{1}{2} \int_0^\infty
    E[M_t^\alpha N_t^\beta A_t]dt\,.
\end{eqnarray*}
Consequently
$$
E[FG]-\|F\|_p\|G\|_{q'}=\frac{1}{2}\int_0^\infty
E\left[M_t^\alpha N_t^\beta A_t\right]dt\,.
$$
Look now  at $A_t$ as a quadratic form  of  with respect  to
$x=\frac{\phi}{M}\,, \,y=\frac{\psi}{N}$:
$$
A_t=\alpha(\alpha-1)x^2+2\alpha\beta \rho_txy+\beta(\beta-1)y^2.
$$
Clearly $(x,y)=(0,0)$ is a stationary point of this quadratic form,
moreover it can not be a minimum, hence it is either a maximum or a
saddle point. For it to be a maximum, the second derivative of $A_t$,
say $D^2A_t(x,y)$, 
should be a negative definite matrix;  in particular the eigenvalues of
this second derivative should be of the same sign and this happens
provided that the determinant of $D^2A_t(x,y)$ is positive. This
latter is implied by the hypothesis $p-1\geq r^2(q-1)$ and hence
$A_t\leq 0$ almost surely and  we obtain
$$
E[FG]=E[T(F)G]\leq\|F\|_p\|G\|_{q'}
$$
which achieves the proof.
\qed

\begin{lemma}
Let $(w,z)=W\times W$ be independent Brownian paths. For $\rho\in[0,1]$, define
$x=\rho w+\sqrt{1-\rho^2}\,z$, ${\Xi}_ \infty$ the $\sigma$-algebra
associated to the paths $x$. Then we have
$$E[F(w)|{\Xi}_\infty]=\int_W F\left(\rho
  x+\sqrt{1-\rho^2}\,z\right)\mu(dz).
$$
\end{lemma}
\proof
For any $G\in L^\infty({\Xi}_\infty)$, we have
\begin{eqnarray*}
E[F(w)\cdot G(x)] & = &
  E\left[F(w)G\left(\rho w+\sqrt{1-\rho^2}\,z\right)\right] \\
& = & E\left[F\left(\rho w+\sqrt{1-\rho^2}\,z\right)G(w)\right] \\
& = & \int\!\!\int F\left(\rho\bar{w}+\sqrt{1-\rho^2}\,\bar{z}\right)
     G(\bar{w})\cdot\mu(d\bar{w})\mu(d\bar{z}) \\
&=&E\left[G(x)\int F\left(\rho
    x+\sqrt{1-\rho^2}\,\bar{z}\right)\cdot\mu(d\bar{z})\right]
\end{eqnarray*}

\noindent
where $\bar{w},\bar{z}$ represent the dummy variables  of integration.
\qed

\begin{corollary}
Under the hypothesis of the above lemma, we have
$$
\left\| \int_W F\left(\rho
x+\sqrt{1-\rho^2}\,\bar{z}\right)\mu(d\bar{z})\right\|_{L^q(\mu)}\leq
\|F\|_{L^p(\mu)}
$$
for any
$$
(p-1)\geq\rho^2(q-1)\,.
$$
\end{corollary}
\section{Logarithmic Sobolev Inequality}
\index[sub]{logarithmic Sobolev inequality}
Let $(P_t,t\geq 0)$ be the Ornstein-Uhlenbeck semigroup. The
commutation relation (cf. \ref{com-rel-1}) 
$$
\nabla P_tf=e^{-t}P_t\nabla f
$$
is directly related to the logarithmic Sobolev inequality of L. Gross:
$$
E\left[f^2\log f^2\right]-E[f^2]\log E[f^2]\leq 2E\left[|\nabla
  f|_H^2\right]\,.
$$
In fact, suppose that $f$ is strictly  positive and lower and upper
bounded.  We have
\bea
\label{inek-log}
E[f\log f]-E[f]\log E[f]&=&-\int_0^\infty
 E\left[ \frac{d}{dt}P_tf\log P_tf\right]dt\nonumber\\
&=&\int_0^\infty E\left[\L P_tf\,\log P_tf\right]dt\nonumber\\
&=&\int_0^\infty E\left[\frac{|\nabla P_tf|_H^2}{ P_tf}\right]dt\nonumber\\
&=&\int_0^\infty e^{-2t}E\left[\frac{|P_t\nabla f|_H^2}{P_tf}\right]dt\,.
\eea
Now insert in \ref{inek-log}  the following,
\beaa
|P_t(\nabla f)|_H^2&=&\left|P_t\left(f^{1/2}\frac{\nabla
      f}{f^{1/2}}\right)\right |_H^2\\
&\leq& (P_tf)\,P_t\left(\frac{|\nabla f|_H^2}{f}\right)
\eeaa
which is a consequence of the H\"older  inequality, to obtain
\beaa
E[f\log f]-E[f]\log E[f]&\leq&\int_0^\infty
e^{-2t}E\left[P_t\left(\frac{|\nabla f|_H^2}{f}\right)\right]dt\\
&=&\int_0^\infty e^{-2t}4 E[|\nabla\sqrt{f}|_H^2]dt\\
&=&2E[|\nabla\sqrt{f}|_H^2]\,,
\eeaa
replacing $f$ by $f^2$ completes the proof of the inequality.

\begin{remarkk}{\rm
Here we have used the fact that if $f>0$ almost surely, then $P_tf>0$
also. In fact one can prove, using the Cameron Martin theorem, that, if
$\mu\{g>0\}>0$, then $P_tg>0$ almost surely. $P_t$ is called a
 positivity improving semi-group (cf. Corollary \ref{pos-imp}).
\index[sub]{positivity improving}}
\end{remarkk}

\section*{Notes and suggested reading}
{\footnotesize{The hypercontractivity property of the
    Ornstein-Uhlenbeck semigroup is due to E. Nelson. The proof given
    here follows the lines given by  J. Neveu, cf. \cite{JN}. For the
    logarithmic Sobolev  inequality and its relations to
    hypercontractivity  cf. \cite{LG,LG-1,G1,G2}.

}}


\chapter{$L^p$-Multipliers Theorem, Meyer Inequalities and
  Distributions}
\label{ch.multip}
\markboth{Multipliers and Inequalities}{Distributions}
\section{$L^p$-Multipliers Theorem}
\index[sub]{multiplier}
$L^p$-Multipliers Theorem gives us a tool to perform some sort of symbolic
calculus to study the continuity of the operators defined via the Wiener
chaos decomposition of the Wiener functionals. With the help of this calculus
 we will complete the proof of the Meyer's inequalities.

 Almost all of these results have been discovered by P. A. Meyer
(cf. \cite{PAM}) and they are consequences of the Nelson's
hypercontractivity theorem (\cite{Nel}).

First let us give  first the following simple and important result:

\begin{theorem}
Let $F\in L^p(\mu)$, $p>1$, denote by $I_n(F_n)$ the projection of $F$
on the $n$-th Wiener  chaos, $n\geq 1$.
Then the map $F\to I_n(F_n)$ is continuous on $L^p(\mu)$.
\end{theorem}
\proof
Suppose first $p>2$. Let $t$ be such that  $p=e^{2t}+1$, then we have
$$
\|P_t F\|_p\leq\|F\|_2\,.
$$
Moreover
$$
\|P_t I_n(F_n)\|_p\leq\|I_n(F_n)\|_2\leq\|F\|_2\leq\|F\|_p
$$
but $P_t I_n(F_n)=e^{-nt}I_n(F_n)$, hence
$$
\|I_n(F_n)\|_p\leq e^{nt}\|F\|_p\,.
$$
For $1<p<2$ we use the duality: let $F\to I_n(F_n)=J_n(F)$. Then
\begin{eqnarray*}
\|I_n(F)\|_p & = & \sup_{\|G\|_q\leq1}| \langle G,J_n(F)\rangle|\\
& =& \sup|\langle J_n(G),F\rangle| \\
& = & \sup|\langle J_nG,J_nF\rangle|\\
&\leq& \sup e^{nt}\|G\|_q\|F\|_p \\
& = & e^{nt}\|F\|_p\,.
\end{eqnarray*}

\vspace*{-3ex}
\qed

\begin{proposition}[Meyer's Multipliers theorem]
Let the function $h$ be defined as
$$
h(x)=\sum_{k=0}^\infty a_k x^k
$$
be  analytic around the origin with
$$
\sum_{k=1}^\infty |a_k|\Big(\frac{1}{n^\alpha}\Big)^k<\infty
$$
for $n\geq n_0$, for some
$n_0\in\N$. Let $\phi(x)=h(x^{-\alpha})$ and define $T_\phi$ on $L^p(\mu)$ as
$$
T_\phi F=\sum_{n=0}^\infty \phi(n)I_n(F_n)\,.
$$
Then the operator to $T_\phi$ is bounded on $L^p(\mu)$ for any $p>1$.
\end{proposition}
\proof
Suppose first $\alpha=1$. Let $T_\phi=T_1+T_2$ where
$$
T_1F=\sum_{n=0}^{n_0-1}\phi(n)I_n(F_n),\qquad T_2F=(I-T_1)F\,.
$$
From the hypercontractivity, $F\mapsto T_1F$ is continuous on $L^p(\mu)$. Let
$$
\Delta_{n_0}F=\sum_{n=n_0}^\infty I_n(F_n).
$$
 Since
$$
(I-\Delta_{n_0})(F)=\sum_{n=0}^{n_0-1}I_n(F_n),
$$
$\Delta_{n_0}:L^p\to L^p$ is
continuous, hence $P_t\Delta_{n_0}:L^p\to L^p$ is also continuous. Applying
Riesz-Thorin interpolation theorem, which says that if
$P_t\Delta_{n_0}$ is $L^q\to
L^q$ and $L^2\to L^2$ then it is $L^p\to L^p$ for any $p$ such that
$\frac{1}{p}$ is in the interval $\Big[\frac{1}{q},\frac{1}{2}\Big]$,
we obtain
$$
\|P_t\Delta_{n_0}\|_{p,p}\leq \|P_t\Delta_{n_0}\|_{2,2}^\theta\,
\|P_t\Delta_{n_0}\|_{q,q}^{1-\theta}\leq \|P_t\Delta_{n_0}\|_{2,2}^\theta\,
\|\Delta_{n_0}\|_{q,q}^{1-\theta}
$$
where $\frac{1}{p}=\frac{\theta}{2}+\frac{1-\theta}{q}\,$, $\theta\in(0,1)\,$.
Choose $q$ large enough such that  $\theta\approx1$ (if necessary). Hence we have
$$
\|P_t\Delta_{n_0}\|_{p.p}\leq e^{-n_0t\theta}K,\qquad K=K(n_0,\theta)\,.
$$
A similar argument  holds for $p\in(1,2)\,$ by duality.

We then  have
\begin{eqnarray*}
T_2(F) & = & \sum_{n\geq n_0}\phi(n)I_n(F_n)\\
& = & \sum_{n\geq n_0}\Big(\sum_k a_k\Big(\frac{1}{n}\Big)^k\Big) I_n(F_n) \\
& = & \sum_k a_k\sum_{n\geq n_0}\Big(\frac{1}{n}\Big)^kI_n(F_n) \\
& = & \sum_k a_k\sum_{n\geq n_0}{\mathcal L}^{-k}I_n(F_n) \\
& = & \sum_k a_k{\mathcal L}^{-k}\Delta_{n_0}F\,.
\end{eqnarray*}
We also have
\begin{eqnarray*}
\|{\mathcal L}^{-1}\Delta_{n_0}F\|_p & = & \Big\|\int_0^\infty
    P_t\Delta_{n_0} Fdt\Big\|_p \leq K \int_0^\infty e^{-n_0\theta t}
              \|F\|_p dt\leq K\cdot\frac{\|F\|_p}{n_0\theta} \\
\|{\mathcal L}^{-2}\Delta_{n_0}F\|_p & = & \Big\|\int_0^\infty
\int_0^\infty
     P_{t+s}\Delta_{n_0} Fdt\,ds\Big\|_p \leq K\cdot\frac{\|F\|_p}{(n_0\theta)^2}, \\
\cdots\\
\|{\mathcal L}^{-k}\Delta_{n_0}F\|_p &\leq &K\|F\|_p \frac{1}{(n_0\theta)^k}\,.
\end{eqnarray*}
Therefore
$$
\|T_2(F)\|_p\leq\sum_k K\|F\|_p\frac{1}{n_0^k\theta^k} \cong \sum
K\|F\|_p\frac{1}{n_0^k}
$$

\noindent
by the hypothesis (take $n_0+1$ instead of $n_0$ if necessary).

For the case $\alpha\in(0,1)\,$, let $\theta_t^{(\alpha)}(ds)$ be the measure
on $\R_+$, defined by
$$
\int_{\R_+}  e^{-\lambda
s}\theta_t^{(\alpha)}(ds)=e^{-t\lambda^\alpha}.
$$
Define
$$
Q_t^\alpha F=\sum_n e^{-n^\alpha t}I_n(F_n)=\int_0^\infty P_s
F\theta_t^{(\alpha)}(ds)\,.
$$
Then
\begin{eqnarray*}
\|Q_t^\alpha\Delta_{n_0}F\|_p&\leq&\|F\|_p \int_0^\infty e^{-n_0\theta
s} \theta_t^{(\alpha)}(ds)\\
&=&\|F\|_p e^{-t(n_0\theta)^{\alpha .}},
\end{eqnarray*}
\noindent
the rest of the proof goes as in the case $\alpha=1$.
\qed

\subsubsection{Examples:}
1) \ Let
\begin{eqnarray*}
\phi(n) & = & \left( \frac{1+\sqrt{n}}{\sqrt{1+n}}\right)^s
     \qquad s\in(-\infty,\infty) \\
& = & h\left( \sqrt{\frac{1}{n}}\right),
     \qquad h(x)=\left(\frac{1+x}{\sqrt{1+x^2}}\right)^s\,.
\end{eqnarray*}
Then $T_\phi:L^p\to L^p$ is bounded. Moreover
$\phi^{-1}(n)=\frac{1}{\phi(n)}=h^{-1}\Big(\sqrt{\frac{1}{n}}\Big)$,
$h^{-1}(x)=\frac{1}{h(x)}$ is also analytic near the origin, hence
$ T_{\phi^{-1}}:L^p\to L^p$ is also a  bounded operator.
\bigskip

\noindent
2) \  Let $\phi(n)=\frac{\sqrt{1+n}}{\sqrt{2+n}}$ then
$h(x)=\sqrt{\frac{x+1}{2x+1}}\,$ satisfies also the above hypothesis.
\bigskip

\noindent
3) \  As an application of (2), look at
\begin{eqnarray*}
\|(I+{\mathcal L})^{1/2}\nabla\varphi\|_p & = &
    \|\nabla(2I+{\mathcal L})^{1/2}\varphi\|_p \\
&\leq&
    \|(I+{\mathcal L})^{1/2}(2I+{\mathcal L})^{1/2}\varphi\|_p \\
& = & \|(2I+{\mathcal L})^{1/2}(I+{\mathcal L})^{1/2}\varphi\|_p  \\
& = & \|T_\phi(I+{\mathcal L})^{1/2}(I+{\mathcal L})^{1/2}\varphi\|_p\\
&\leq& c_p\|(I+{\mathcal L})\varphi\|_p\,.
\end{eqnarray*}

\noindent
Continuing this way we can show that
\begin{eqnarray*}
\|\nabla^k\varphi\|_{L^p(\mu,H^{\otimes k})} & \leq & c_{p,k} \|\varphi\|_{p,k}
           (=\|(I+{\mathcal L})^{k/2}\varphi\|_p) \\
& \leq & \tilde{c}_{p,k}(\|\varphi\|_p+
       \|\nabla^k\varphi\|_{L^p(\mu,H^{\otimes k})})
\end{eqnarray*}

\noindent
and this completes the proof of the Meyer inequalities for the scalar-valued
Wiener functionals.  If $\Xi$ is a separable Hilbert space, we denote
 with  $\DD_{p,k}({\Xi})$ the completion of the $\Xi$-valued
polynomials  with respect to the norm
$$
\|\alpha\|_{\DD_{p,k}({\Xi})}=
\|(I+{\mathcal L})^{k/2}\|_{L^p(\mu,{\Xi})}\,.
$$
We define as in the case ${\Xi}=\R$, the Sobolev derivative
$\nabla$, the divergence $\delta$, etc. All we have said for the real case
extend trivially to the vector case, including the Meyer inequalities.
In fact, in the proof of these inequalities the main step
is the Riesz inequality for the Hilbert
transform. However this inequality is also true for any Hilbert space (in fact
it holds also for a class of Banach spaces which contains Hilbert spaces,
called UMD spaces). The rest is almost the
transcription of the real case combined with  the Khintchine inequalities.
We leave hence this passage to the reader.
\qed

\begin{corollary}
For every $p>1$, $k\in\R$, $\nabla$ has a continuous extension as a map
$\DD_{p,k}\to \DD_{p,k-1}(H)$.
\end{corollary}
\proof
We have
\begin{eqnarray*}
\|\nabla\varphi\|_{p,k}& = & \|(I+{\mathcal L})^{k/2}\nabla\varphi\|_p\\
&=&\|\nabla(2I+{\mathcal L})^{k/2}\varphi\|_p \\
& \leq &c_p \|(1+{\mathcal L})^{1/2}(2I+{\mathcal L})^{k/2}\varphi\|_p \\
&\approx &
   \|(I+{\mathcal L})^{(k+1)/2}\varphi\|_p\\
&=&\|\varphi\|_{p,k+1}\,.
\end{eqnarray*}

\vspace*{-2ex}
\qed

\begin{corollary}
\label{extension-cor}
$\delta=\nabla^\ast:\DD_{p,k}(H)\to \DD_{p,k-1}$ is continuous
for all $ p>1$ and $k\in\R$.
\end{corollary}
\proof
The proof follows from the duality.
\qed

\noindent
In particular we have :
\index[not]{d@$\DD$}
\index[not]{da@$\DD'$}
\begin{corollary}
The Sobolev derivative and its adjoint extend to distribution spaces
as explained below:
\begin{itemize}
\item Sobolev  derivative  operates as a continuous operator on
$$
\nabla:\DD=\bigcap_{p,k}\DD_{p,k}\to D(H)=\bigcap_{p,k} \DD_{p,k}(H)
$$
and it  extends continuously as a map
$$
\nabla:\DD'=\bigcup_{p,k} \DD_{p,k}\to \DD'(H)=\bigcup_{p,k} \DD_{p,k}(H).
$$
The elements of the  space $\DD'$ are called Meyer-Watanabe
distributions. \index[sub]{Meyer-Watanabe distributions}
\item Consequently  its adjoint  has similar properties:
$$
\delta:\bigcap_{p,k} \DD_{p,k}(H)=\DD(H)\to \DD
$$
is continuous and this map has a continuous extension to
$$
\delta:\DD'(H)\to \DD'
$$
\end{itemize}
\end{corollary}
\proof
Everything follows from the dualities
$$
(\DD)'=\DD',\,(\DD(H))'=\DD'(H).
$$
\qed

\begin{definition}
For $n\geq1$, we define $\delta^n$ as $(\nabla^n)^\ast$ with respect to $\mu$.
\end{definition}
Here is the generalization of Corollary \ref{I-W-decomp} promised in
Chapter 2:
\begin{proposition}
\label{taylor}
For $\varphi\in L^2(\mu)$, we have
\begin{equation}
\label{I-W-eqn}
\varphi=E[\varphi]+\sum_{n\geq1}\frac{1}{n!}\delta^n(E[\nabla^n\varphi])\,.
\end{equation}
\end{proposition}
\proof If $f$ is a symmetric element of $H^{\otimes n}$, we e shall
denote by $\tilde{I}_n(f)$ the $n$-th order multiple Ito-Wiener
integral of the density of $f$ with respect to the Lebesgue measure on
$[0,1]^n$ (cf. also Corollary \ref{I-W-decomp}). With this notational
convention,  suppose that
$h\mapsto\varphi(w+h)$ is analytic for almost all  $w$. Then we have
$$
\varphi(w+h)=\varphi(w)+
\sum_{n\geq1}\frac{(\nabla^n\varphi(w),h^{\otimes n})_{H^{\otimes n}}}{n!}\,.
$$
Take the expectations:
\begin{eqnarray*}
E[\varphi(w+h)] & = & \!\!E[\varphi\,\,{\mathcal E}(\delta h)] \\
& = & \!\! E[\varphi]+\sum_{n\cdot} \frac{(E[\nabla^n\varphi],h^{\otimes n})}{n!} \\
& = &\!\! E[\varphi]+ \sum_{n\geq1} E\left[
\frac{{\tilde{I}}_n(E[\nabla^n\varphi])}{n!} {\mathcal E}(\delta h)\right].
\end{eqnarray*}

\noindent
Since the finite linear combinations of the elements of the set
$\{{\mathcal E}(\delta h);h\in H\}$ is dense in any $L^p(\mu)$, we obtain
the identity
$$
\varphi(w)=E[\varphi]+\sum_{n\geq1} \frac{{\tilde{I}}_n(E[\nabla^n\varphi])}{n!}\,.
$$

Let $\psi\in \DD$, then we have (with $E[\psi]=0$),
\begin{eqnarray*}
\langle\varphi,\psi\rangle & = & \sum_{n\geq 1}E[{\tilde{I}}_n(\varphi_n){\tilde{I}}_n(\psi_n)] \\
& = & \sum_{n}E\left[ \frac{{\tilde{I}}_n(E[\nabla^n\varphi])}{n!}\cdot
     {\tilde{I}}_n(\psi_n)\right]  \\
& = & \sum_n(E[\nabla^n\varphi],\psi_n) \\
& = & \sum_n \frac{1}{n!} (E[\nabla^n\varphi],E[\nabla^n\psi]) \\
& = & \sum_n\frac{1}{n!}E[(E[\nabla^n\varphi],\nabla^n\psi)] \\
& = & \sum_{n}\frac{1}{n!}E[\delta^n(E[\nabla^n\varphi])\cdot\psi]
\end{eqnarray*}

\noindent
hence we obtain that
$$
\varphi=\sum_n\frac{1}{n!}\delta^n E[\nabla^n\varphi].
$$

\noindent
In particular it holds true that
$$
\delta^n\left\{E[\nabla^n\varphi]\right\}={\tilde{I}}_n(E[\nabla^n\varphi])\,.
$$
Evidently this identity is valid not only for the expectation of the
n-th derivative of a Wiener functional but for any symmetric element
of $H^{\otimes n}$.
\qed
\begin{remarkk}
\label{taylor-rem}
{\rm
Although in the litterature Proposition \ref{taylor} is announced for
the elements of $\DD$, the proof given here shows its validity for the
elements of $L^2(\mu)$. In fact, although $\nabla^n\phi$ is a
distribution, its expectation is an ordinary symmetric tensor of order
$n$, hence the corresponding multiple Wiener integrals are
well-defined. With a small extra work, we can show that in
fact the formula (\ref{I-W-eqn}) holds for any $\phi\in \cup_{k\in
\Z}\DD_{2,k}$. }
\end{remarkk}
Let us give another result important for the applications:
\begin{proposition}
\label{sob-prop}
Let $F$ be in some $L^p(\mu)$ with $p>1$ and suppose that the
distributional derivative $\nabla F$ of $F$,  is  in some
$L^r(\mu,H)$, ($1<r$). Then $F$ belongs to  $\DD_{r\wedge p,1}$.
\end{proposition}
\proof
Without loss of generality, we can assume that $r\leq p$. Let $(e_i;i\in \N)$ be a complete, orthonormal  basis of the Cameron-Martin space $H$. Denote by
$V_n$ the sigma-field generated by $\delta e_1,\ldots, \delta e_n$, and by $\pi_n$ the orthogonal projection of $H$ onto the subspace spanned by $e_1,\ldots, e_n$, $n\in \N$. Let us define $F_n$ by
$$
F_n=E[P_{1/n}F|V_n],
$$
where $P_{1/n}$ is the Ornstein-Uhlenbeck semi-group at $t=1/n$. Then $F_n$ belongs to $\DD_{r,k}$ for any $k\in \N$ and converges to $F$ in $L^r(\mu)$. Moreover, from Doob's lemma, $F_n$ is of the form
$$
F_n(w)=\alpha(\delta e_1,\ldots, \delta e_n),
$$
with $\alpha$ being a Borel function on $\R^n$, which is in the intersection of the Sobolev spaces $\cap_k W_{r,k}(\R^n,\mu_n)$ defined with the Ornstein-Uhlenbeck operator  $L_n=-\Delta+x\cdot \nabla$ on $\R^n$. Since $L_n$ is elliptic, the Weyl lemma \index[sub]{Weyl's lemma} implies that $\alpha$ can be chosen as a $C^\infty$-function. Consequently, $\nabla F_n$ is again $V_n$-measurable and we find , using the very definition of conditional expectation and the Mehler formula, that
$$
\nabla F_n=E[e^{-1/n}\pi_n P_{1/n}\nabla F|V_n].
$$
Consequently, from the martingale convergence theorem and from the
fact that $\pi_n \rightarrow I_H$ in the weak operators  topology, it
follows that
$$
\nabla F_n\rightarrow \nabla F,
$$
in $L^r(\mu,H)$, consequently $F$ belongs to $\DD_{r,1}$.
\qed

\section*{Appendix: Passing from the classical Wiener space to the
  Abstract Wiener Space (or vice-versa):}

Let $(W,H,\mu)$ be an abstract Wiener space. Since, \` a priori, there is no
 notion of time, it seems that we can not define the notion of anticipation,
non-anticipation, etc.
\noindent
{\bf This difficulty can be overcome in the following way:}
\bigskip

Let $(p_\lambda;\lambda\in\Sigma)$, $\Sigma\subset\R$, be a resolution of
identity\index[sub]{resolution of identity} on the separable Hilbert
space $H$, i.e., each $p_\lambda$ is an
\index[not]{p@$p_\lambda$}orthogonal projection, increasing to $I_H$,
in the sense that
$\lambda\mapsto(p_\lambda h,h)$ is an increasing function. Let us denote by
$H_\lambda=\overline{p_\lambda(H)}$, where $\overline{p_\lambda(H)}$
denotes the closure of $p_\la(H)$ in $H$.

\begin{definition}
We will denote by ${\mathcal F}_\lambda$ the $\sigma$-algebra
generated by the real polynomials $\varphi$ on $W$ such that
$\nabla\varphi\in H_\lambda$ $\mu$-almost surely.
\end{definition}
\index[not]{s@$\Sigma$}
\begin{lemma}
We have
$$
\bigvee_{\lambda\in\Sigma}{\mathcal F}_\lambda={\calB}(W)
$$
up  to $\mu$-negligeable sets.
\end{lemma}
\proof
We have already $\bigvee{\mathcal F}_\lambda\subset{\calB}(W)$.
Conversely, if $h\in
H$, then $\nabla\delta h=h$. Since $\bigcup_{\lambda\in\Sigma}H_\lambda$
is dense in $H$, there exists $(h_n)\subset\bigcup_\lambda H_\lambda$
such that  $h_n\to h$ in $H$. Hence $\delta h_n\to\delta h$ in
$L^p(\mu)$, for all
$p\geq1$. Since each $\delta h_n$ is $\bigvee{\mathcal F}_\lambda$-measurable, so
does $\delta h$. Since ${\calB}(W)$ is generated by $\{\delta h;h\in H\}$ the
proof is completed.
\qed

\begin{definition}
A random variable $\xi:W\to H$ is called a simple, adapted vector field if it
can be written as a finite sum:
$$
\xi=\sum_{i<\infty}F_i(p_{\lambda_{i+1}}h_i-p_{\lambda_i}h_i)
$$
where $h_i\in H$, $F_i$ are ${\mathcal F}_{\lambda_i}$-measurable (and
smooth for the time being) random variables.
\end{definition}
\begin{proposition}
For each adapted simple vector field we have
\begin{itemize}
\item[i)]
$\delta\xi=\sum_{i<\infty}
F_i\delta(p_{\lambda_{i+1}}h_i-p_{\lambda_i}h_i)$

\item[ii)] with Ito's isometry:
$$
E\left[|\delta\xi|^2\right]=E\left[|\xi|^2_H\right]\,.
$$
\end{itemize}
\end{proposition}
\proof
The first part follows from the usual identity
$$
\delta[F_i(p_{\lambda_{i+1}}-p_{\lambda_i})h_i]=F_i\delta[(p_{\lambda_{i+1}}
-p_{\lambda_i})h_i]-\left(\nabla
 F_i,(p_{\lambda_{i+1}}-p_{\lambda_i})h_i\right)_H\,
$$
and from  the fact that the second term is null since $\nabla F_i\in
H_\lambda$ almost surely. The verification of the  second relation is  left
to the reader.
\qed

\begin{remarkk}
{\rm
If we denote $\Sigma F_i\,1_{]\lambda_i,\lambda_{i+1}]}(\lambda)\, h_i$ by
$\dot{\xi}(\lambda)$, we have the following relations:
$$
\delta\xi=\delta\int_\Sigma\dot{\xi}(\lambda)dp_\lambda\quad
{\mbox{with}}
\;\; \|\delta\xi\|_2^2=E\int_\Sigma
d(\dot{\xi}_\lambda,p_\lambda\dot{\xi}_\lambda)=\|\xi\|_{L^2(\mu,H)}^2\;,
$$
which are significantly analogous to the relations  that we have seen
before.}
\end{remarkk}

The Ito representation theorem can be stated  in this setting as
follows: suppose that  $(p_\lambda;\lambda\in\Sigma)$ is weakly
continuous. We mean by this that the function
$$
\la\to (p_\la h,k)_H
$$
is continuous for any $h,k\in H$. Then
\begin{theorem}
Let us denote with $\DD_{2,0}^a(H)$ the completion of adapted simple vector
fields with respect to the $L^2(\mu,H)$-norm. Then we have
$$
L_2(\mu)=\R+\{\delta\xi:\xi\in \DD_{2,0}^a(H)\}\,,
$$
i.e., any $\varphi\in L_2(\mu)$ can be written as
v$$\varphi=E[\varphi]+\delta\,\xi
$$
for some $\xi\in \DD_{2,0}^a(H)$. Moreover such $\xi$ is unique up to
$L^2(\mu,H)$-equivalence classes.
\end{theorem}
The following result explains the reason of the existence of the Brownian
motion (cf. also \cite{U5}):

\begin{theorem}
Suppose that there exists some $\Omega_0\in H$ such that  the set
$\{p_\lambda\Omega_0;$ $\lambda\in\Sigma\}$ has a dense span in $H$ (i.e.\ the
linear combinations from it is a dense set). Then the real-valued $({\mathcal
F}_\lambda)$-martingale defined by
$$
b_\lambda=\delta p_\lambda\Omega_0
$$
\index[not]{o@$\Omega_0$}
\noindent
is a Brownian motion with  a deterministic time change and $({\mathcal
F}_\lambda;\lambda\in\Sigma)$ is its canonical filtration completed with the
negligeable sets.
\end{theorem}
\paragraph{Example:}
Let $H=H_1([0,1])$, define $A$ as the operator defined by
$Ah(t)=\int_0^t s\dot{h}(s)ds$. Then $A$ is a self-adjoint operator on
$H$ with a continuous spectrum which is equal to $[0,1]$. Moreover we have
$$
(p_\lambda h)(t)=\int_0^t 1_{[0,\lambda]}(s)\dot{h}(s)ds
$$

\noindent
and $\Omega_0(t)=\int_0^t 1_{[0,1]}(s)ds$ satisfies the hypothesis of
the above theorem. $\Omega_0$ is called the vacuum vector \index[sub]{vacuum vector} (in physics).

This is the main example, since all the (separable) Hilbert spaces are
isomorphic, we can carry this time structure to any abstract Hilbert-Wiener
space as long as we do not need any particular structure of time.
\section{Exercises}
{\footnotesize{
\begin{enumerate}
\item Give a detailed proof of Corollary \ref{extension-cor}, in
particular explain the why of the existence of  continuous extensions
of $\delta$ and $\nabla$.
\item Prove the last claim of Remark \ref{taylor-rem}.
\end{enumerate}
}}

\section*{Notes and suggested reading}
{\footnotesize{To complete the series of the Meyer inequalities, we
    have been obliged to use the hypercontractivity property of the
    Ornstein-Uhlenbeck semigroup as done  in \cite{PAM}. Once this is
    done the extensions of $\nabla$ and $\delta$ to the distributions
    are immediate via the duality techniques. Proposition \ref{taylor}
    is due to Stroock, \cite{DS} with a different proof. The results
    of the appendix are essentially due to the author, cf. \cite{U5}.
    In \cite{UZ7} a stochastic calculus is constructed in more detail.

}}


\chapter{Some Applications}
\label{ch.appli}
\markboth{Applications}{}

\section*{Introduction}

In this chapter we give some applications of the extended versions of the
derivative and the divergence operators. First we give an extension of the
Ito-Clark formula to the space of the scalar distributions. We refer the reader to  \cite{Bi} and \cite{O} for the  developments of this formula in the case
of Sobolev differentiable Wiener functionals. 
 Let us briefly explain the problem:  although, we
know from the Ito representation theorem, that each square integrable Wiener
functional can be represented as the stochastic integral of an adapted
process, without the use of the distributions, we can not calculate
this process, since any square integrable random variable is not
necessarily in $\DD_{2,1}$, hence it is not Sobolev differentiable in
the ordinary sense. As it will be  explained, this problem is
completely solved using the differentiation in the sense of 
distributions. Afterwards we give a straightforward application of
this result to prove a $0-1$ law for the Wiener measure. At the second
section we construct the composition  of the tempered distributions
with non-degenerate Wiener functionals as Meyer-Watanabe
distributions. This construction carries also the information that the
probability density of a non-degenerate random variable is not only
infinitely differentiable but also it is rapidly decreasing. The same
idea is then applied to prove the regularity of the solutions of the
Zakai equation for the filtering of non-linear diffusions.

\section{Extension of the Ito-Clark formula\index[sub]{Ito-Clark formula}}
\markboth{Applications}{Ito-Clark Formula}
Let $F$ be any integrable random variable. Then we the celebrated Ito
Representation Theorem \ref{Ito-rep-thm} tells us that  $F$ can
be represented as
$$
F=E[F]+\int_0^1 H_sdW_s\,,
$$
where $(H_s;s\in[0,1])$ is an adapted process such that, it is unique and
$$
\int_0^1 H_s^2ds<+\infty\;\;\mbox{a.s.}
$$
Moreover, if $F\in L^p$ $(p>1)$, then we also have
$$
E\Big[\Big(\int_0^1|H_s|^2ds\Big)^{p/2}\Big]<\infty\,.
$$

One question is how to calculate the process $H$. In fact, below we will extend
the Ito representation  and answer to the above question for any $F\in \DD'$
(i.e., the Meyer-Watanabe distributions). We begin with
\begin{lemma}
\label{lem-1}
Let $\xi\in \DD(H)$ be represented as $\xi(t)=\int_0^t\dot{\xi}_sds$,
then $\pi\xi$ defined by \index[not]{p@$\pi$} 
$$
\pi\xi(t)=\int_0^t E[\dot{\xi}_s|{\mathcal F}_s]ds
$$
belongs again to $\DD(H)$. In other words  $\pi:\DD(H)\to \DD(H)$ is a linear
continuous operator.
\end{lemma}
\proof Let $(P_t,t\in \R_+)$ be the Ornstein-Uhlenbeck semigroup. Then
it is easy to see that, for any $\tau\in [0,1]$, if $\phi\in L^1(\mu)$
is $\B_\tau$-measurable, then so is also $P_t\phi$ for any $t\in
\R_+$. This implies in particular that 
 ${\mathcal L}\pi\xi=\pi{\mathcal L}\xi$. Therefore 
\begin{eqnarray*}
\|\pi\xi\|_{p,k} & = & E\Big[\Big(\int_0^1 |(I+{\mathcal L})^{k/2}
     E[\dot{\xi}_s|{\mathcal F}_s]|^2ds\Big)^{p/2}\Big]=\\
& = & E\Big[\Big(\int_0^1| E[(I+{\mathcal L})^{k/2}
    \dot{\xi}_s|{\mathcal F}_s]|^2ds\Big)^{p/2}\Big] \\
& \leq & c_pE\Big[\Big( \int_0^1 |(I+{\mathcal L}^{k/2}
  \dot{\xi}_s|^2ds\Big)^{p/2}\Big]\qquad (c_p\cong p)
\end{eqnarray*}
where the last inequality follows from the convexity inequalities of the dual
predictable projections (c.f. \cite{D-M}).
\qed

\begin{lemma}
$\pi:\DD(H)\to \DD(H)$ extends as a continuous mapping to $\DD'(H)\to \DD'(H)$.
\end{lemma}
\proof
Let $\xi\in \DD(H)$, then we have, for $k>0$,
\begin{eqnarray*}
\|\pi\xi\|_{p,-k} & = & \|(I+{\mathcal L})^{-k/2}\pi\xi\|_p\\
     & = &\|\pi(I+{\mathcal L})^{-k/2}\xi\|_p\leq c_p\|(I+{\mathcal L})^{-k/2}\xi\|_p \\
& \leq & c_p\|\xi\|_{p,-k}\,,
\end{eqnarray*}
then the proof follows since $\DD(H)$ is dense in $\DD'(H)$.
\qed

Before going further let us give a notation: if $F$ is in some
$\DD_{p,1}$ then its Gross-Sobolev derivative $\nabla F$ is  an
$H$-valued random variable. Hence $t\mapsto \nabla F(t)$ is absolutely
continuous with respect to the Lebesgue measure on $[0,1]$.  We shall denote
by $D_sF$ its Radon-Nikodym derivative with respect to the Lebesgue
measure. Note that $D_sF$  is $ds\times d\mu$-almost everywhere
well-defined.\index[not]{de@$D_s$}

\begin{lemma}
\label{ilk-lemma}
Let $\varphi\in \DD$, then we have
\begin{eqnarray*}
\varphi & = & E[\varphi]+\int_0^1 E[D_s\varphi|{\mathcal F}_s]dW_s \\
& = & E[\varphi]+\delta\pi\nabla\varphi\,.
\end{eqnarray*}
Moreover $\pi\nabla\varphi\in \DD(H)$.
\end{lemma}
\proof
Let $u$ be an element of $L^2(\mu,H)$ such that
$u(t)=\int_0^t\dot{u}_sds$ with $(\dot{u}_t;t\in[0,1])$ being
an adapted  and bounded process. Then we have, from the Girsanov theorem,
$$
E\left[\varphi(w+\lambda u(w)).\exp\left\{-\lambda\int_0^1 \dot{u}_sdW_s-
\frac{\lambda^2}{2}\int_0^1\dot{u}_sds\right\}\right]=E[\varphi].
$$
Differentiating both sides at $\lambda=0$, we obtain:
$$
E[(\nabla\varphi(w),u)-\varphi\,\int_0^1\dot{u}_sdW_s]=0\,,
$$
i.e.,
$$
E[(\nabla\varphi,u)]=E[\varphi \,\int_0^1 \dot{u}_sdW_s].
$$
Furthermore
\begin{eqnarray*}
E\Big[\int_0^1 D_s\varphi\dot{u}_sds\Big] & = &
E\Big[\int_0^1 E[D_s\varphi|{\mathcal F}_s]\dot{u}_sds\Big] \\
& = & E[(\pi\nabla\varphi,u)_H] \\
& = & E\Big[\Big(\int_0^1 E[D_s\varphi|{\mathcal F}_s]dW_s\Big)
       \Big(\int_0^1\dot{u}_sdW_s\Big)\Big].
\end{eqnarray*}
\noindent
Since the set of the stochastic integrals  $\int_0^1\dot{u}_sdW_s$ of
the processes  $\dot{u}$ as above is dense
in $L^2_0(\mu)=\{F\in L^2(\mu):\,E[F]=0\}$, we see that
$$
\varphi-E[\varphi]=\int_0^1
     E[D_s\varphi|{\mathcal F}_s]dW_s=\delta\pi\nabla\varphi\,.
$$
The rest is obvious from the Lemma \ref{lem-1} .
\qed

Lemma \ref{ilk-lemma} extends to $\DD'$ as:
\begin{theorem}
\label{ust-thm}
For any $T\in \DD'$, we have
$$
T=\langle T,1\rangle+\delta\pi\nabla T\,.
$$
\end{theorem}
\proof
Let $(\varphi_n)\subset \DD$ such that $\varphi_n\to T$ in $D'$. Then  we have
\begin{eqnarray*}
 T&=& \lim_n\varphi_n\\
 & = & \lim_n\left\{E[\varphi_n]+\delta\pi\nabla\varphi_n\right\} \\
& = & \lim_n E[\varphi_n]+\lim_n\delta\pi\nabla\varphi_n \\
& = & \lim_n\langle1,\varphi_n\rangle+\lim_n\delta\pi\nabla\varphi_n \\
& = & \langle1,T\rangle+\delta\pi\nabla T
\end{eqnarray*}
since $\nabla:\DD'\to \DD'(H)$, $\pi:\DD'(H)\to \DD'(H)$ and
$\delta:\DD'(H)\to\DD'$ are all linear, continuous mappings.
\qed

\noindent
Here is a nontrivial application of the Theorem \ref{ust-thm}:
\begin{theorem}{\bf(0--1 law)}\\
\label{0-1-law}
Let $A\in{\calB}(W)$ such that $\nabla_h\won_A=0$, $h\in H$,  where the
derivative is in the sense of the distributions. Then $\mu(A)=0$ or $1$.
\end{theorem}
\remark
In particular, the above  hypothesis is satisfied when $A+H\subset A$.
\proof
Let $T_A=\won_A$, then Theorem \ref{ust-thm} implies that
$$
T_A=\langle T_A,1\rangle=\mu(A)\,,
$$
hence $\mu(A)^2=\mu(A)$. Another proof can be given as follows: let
$T_t$ be defined as $P_t\won_A$, where $(P_t,t\geq 0)$ is the
Ornstein-Uhlenbeck semigroup. Then, from the hypothesis, $\nabla
T_t=e^{-t}P_t\nabla \won_A=0$, consequently $T_t$ is almost surely a
constant for any $t>0$, this implies that $\lim_{t\to 0}T_t=\won_A$ is also a
constant.
\qed

\begin{remarkk}{\rm
From Doob-Burkholder inequalities, it follows via a duality technique, that
\begin{eqnarray*}
\lefteqn{E\left[\left(\int_0^1|(I+{\mathcal{L}})^{k/2}E[D_s\phi|{\mathcal{F}}_s
      ]|^2ds \right)^{p/2}\right]}\\
   &&\leq c_p\|(I+{\mathcal{L}})^{k/2}\phi\|_{L^p(\mu)}^p\\
  &&\leq c_{p,k}\|\phi\|^p_{p,k}\,,
\end{eqnarray*}
for any $p>1$ and $k\in \R$. Consequently, for any $\phi\in \DD_{p,k}$,
$\pi\nabla \phi\in \DD_{p,k}(H)$ and the Ito integral is an
isomorphisme from the adapted elements of $\DD_{p,k}(H)$ onto
$\DD_{p,k}^o=\{\phi\in \DD_{p,k}:\,\langle \phi,1\rangle =0\}$ (cf.  \cite{U0} for
further details).}
\end{remarkk}
\begin{corollary}{\bf{(Positivity improving)}}
\label{pos-imp}
Let $F\in L^p(\mu),\,p>1$ be a non-negative Wiener functional such that
$\mu\{F>0\}>0$, denote by $(P_\tau, \tau\in \reals_+)$ the
Ornstein-Uhlenbeck semi-group. Then, for any $t>0$,  the set
$A_t=\{w:P_tF(w)>0\}$ has full $\mu$-measure,  in fact we have
$$
A_t+H\subset A_t\,.
$$
\end{corollary}
\index[sub]{positivity improving}
\proof
From the Mehler and Cameron-Martin  formulae, we have
\beaa
P_tF(w+h)&=&\int_W F(e^{-t}(w+h)+\sqrt{1-e^{-2t}}y)\mu(dy)\\
&=&\int_W F(e^{-t}w+\sqrt{1-e^{-2t}}y)\rho(\alpha_t\delta h(y))\mu(dy)
\eeaa
where
$$
\alpha_t=\frac{e^{-t}}{\sqrt{1-e^{-2t}}}
$$
and
$$
\rho(\delta h)=\exp\left\{\delta h-1/2|h|_H^2\right\}\,.
$$
This proves the claim about the $H$-invariance of $A_t$ and the proof
follows from Theorem \ref{0-1-law}.
\qed

\section{Lifting of ${\calS}'(\R^d)$ with  random variables}
\markboth{Applications}{Lifting of the Functionals}
\index[not]{se@${\calS}'(\R^d)$}
\index[sub]{lifting}
\markboth{Lifting of ${\calS}'({\R}^d)$}{}
Let $f:\R\to\R$ be a $C_b^1$-function, $F\in \DD$. Then we know that
$$
\nabla(f(F))=f'(F)\nabla F\,.
$$
Now suppose that $|\nabla F|_H^{-2}\in\bigcap L^p(\mu)$, then
$$
f'(F)=\frac{(\nabla(f(F)),\nabla F)_H}{|\nabla F|_H^2}
$$

\noindent
Even if $f$ is not $C^1$, the right hand side  of this equality has a
sense if we look at $\nabla(f(F))$ as an element of $\DD'$. In the
following we will develop this idea:
\begin{definition}
\label{non-deg-def}
Let $F:W\to\R^d$ be a random variable  such that $F_i\in \DD$, for all
$i=1,\ldots,d$, and that
$$
\left[\det(\nabla F_i,\nabla F_j)\right]^{-1}\in\bigcap_{p>1}L^p(\mu).
$$
 Then we say that $F$ is a {\it\bf{
     non-degenerate}\/}\index[sub]{non-degenerate}  random variable.
\end{definition}

\begin{lemma}
Let  us denote by $\sigma_{ij}=(\nabla F_i,\nabla F_j)_H$ and by
$\gamma=\sigma^{-1}$ (as a matrix). Then $\gamma\in
\DD(\R^d\otimes\R^d)$, in particular $\det \gamma \in \DD$.
\end{lemma}
\proof
Formally, we have, using the relation $\sigma\cdot\gamma=Id$,
$$\nabla\gamma_{ij}=\sum_{k,l}\gamma_{ik}\gamma_{jl}\nabla\sigma_{kl}\,.
$$
To justify this we define first
$\sigma_{ij}^\varepsilon=\sigma_{ij}+\varepsilon\delta_{ij}$,
$\varepsilon>0$. Then we can write
$\gamma_{ij}^\varepsilon=f_{ij}(\sigma^\varepsilon)$, where
$f:\R^d\otimes\R^d\to\R^d\otimes\R^d$ is a smooth function of polynomial
growth. Hence $\gamma_{ij}^\varepsilon\in \DD$. Then from the
dominated convergence theorem we have
$\gamma_{ij}^\varepsilon\to\gamma_{ij}$ in $L^p$ and
$\nabla^k\gamma_{ij}^\varepsilon{\longrightarrow}
\nabla^k\gamma_{ij}$ in $L^p(\mu,H^{\otimes k})$ (the latter  follows
again from $\gamma^\varepsilon\cdot\sigma^\varepsilon=Id$).
\qed

\begin{lemma}
\label{calcul-lemma}
Let $G\in \DD$. Then,  for all $f\in{\calS}(\R^d)$, the following
identities are true:
\begin{enumerate}
\item
$$
E[\partial_i f(F).G]=E[f(F)\,\,l_i(G)]\,,
$$
where $G\mapsto l_i(G)$ is linear and for any $1<r<q<\infty$,
$$
\sup_{\|G\|_{q,1}\leq1}\|l_i(G)\|_r<+\infty\,.
$$
\item  Similarly
$$
E[\partial_{i_1\ldots i_k}f\circ F.G]=E[f(F)\cdot l_{i_1\ldots
  i_k}(G)]
$$
and
$$
\sup_{\|G\|_{q,1}\leq1} \|l_{i_1\ldots i_k}(G)\|_r<\infty\,.
$$
\end{enumerate}
\end{lemma}
\proof
We have
$$
\nabla(f\circ F)=\sum_{i=1}^d\partial_i f(F)\nabla F_i
$$
hence
$$
(\nabla(f\circ F),\nabla F_j)_H=\sum_{j=1}^d\sigma_{ij}\partial_i f(F)\,.
$$
Since $\sigma$ is invertible, we obtain:
$$
\partial_i f(F)=\sum_j\gamma_{ij}(\nabla(f\circ F),\nabla F_j)_H\,.
$$
Then
\begin{eqnarray*}
E[\partial_i f(F).G] & = & \sum_j E\left[\gamma_{ij}
         (\nabla(f\circ F),\nabla F_j)_H G\right]  \\
& = & \sum_j E\left[f\circ F\,\,\delta\{\nabla F_j\gamma_{ij}G\}\right]\,,
\end{eqnarray*}

\noindent
hence we see that $l_i(G)=\sum_j\delta\{\nabla
F_j\gamma_{ij}G\}$. Developing this expression gives
\begin{eqnarray*}
l_i(G) & = & -\sum_j\left[(\nabla(\gamma_{ij}G),
           \nabla F_j)_H-\gamma_{ij}G{\mathcal L}F_j\right] \\
& = & -\sum_j\left[\gamma_{ij}(\nabla G,\nabla
  F_j)_H-\sum_{k,l}\gamma_{ik}\gamma_{jl}
(\nabla\sigma_{kl},\nabla F_j)_HG-\gamma_{ij}G{\mathcal L}F_j\right]\,.
\end{eqnarray*}
Hence
\begin{eqnarray*}
\lefteqn{|l_{i}(G)|\leq\sum_j\left[\sum_{kl}|\gamma_{ik}\gamma_{jl}|\,
|\nabla\sigma_{kl}|\,|\nabla F_j|\,|G|\right.}\\
 &&+     |\gamma_{ij}|\,|\nabla F_j|\,|\nabla G| + \left. \vphantom{\sigma_{kl}}|\gamma_{ij}|\,|G|\,|{\mathcal
    L}F_j|\right]\,.
\end{eqnarray*}
Choose $p$ such that  $\frac{1}{r}=\frac{1}{p}+\frac{1}{q}$ and apply
H\"older's inequality:
\begin{eqnarray*}
\|l_i(G)\|_r & \leq & \sum_{j=1}^d\Big[ \sum_{k,l}\|G\|_q
\|\gamma_{ik}\gamma_{jl}|\nabla\sigma_{kl}|_H|\nabla F_j|_H\|_p+ \\
& & \quad
+\,\|\gamma_{ij}|\nabla F_j|\|_p\||\nabla G|\|_q+\|\gamma_{ij}{\mathcal L}F_j\|_p
         \|G\|_q\Big] \\
& \leq & \|G\|_{q,1}\Big[\sum_{j=1}^d \|\gamma_{ik}\gamma_{jl}|\nabla F_{kl}|
          |\nabla F_j|\|_p+ \\
& & \quad
+\,\|\gamma_{ij}|\nabla F_j|\|_p+ \|\gamma_{ij}{\mathcal L}F_j\|_p\Big]\,.
\end{eqnarray*}

\noindent
To prove the last part  we iterate this procedure for $i>1$.
\qed
\bigskip

\noindent
{\bf Remember now  that} ${\calS}(\R^d)$ can be written as the
intersection (i.e., projective limit) of the  Banach
spaces $S_{2k}$ which are defined as below:\\
Let $A=I-\Delta+|x|^2$ and define  $\|f\|_{2k}=\|A^kf\|_\infty$ (the
uniform norm). Then let $S_{2k}$ be the  completion of ${\calS}(\R^d)$
with respect to the norm  $\|\cdot\|_{2k}\,$.

\begin{theorem}
\label{compo-thm}
Let $F\in \DD(\R^d)$ be a  non-degenerate random variable. Then we
have for  $f\in{\calS}(\R^d)$:
$$
\|f\circ F\|_{p,-2k}\leq c_{p,k}\|f\|_{-2k}\,.
$$
\end{theorem}
\proof
Let $\psi=A^{-k}f\in{\calS}(\R^d)$. For $G\in \DD$,from  Lemma
\ref{calcul-lemma},  we know that there
exists some $\eta_{2k}(G)\in \DD$ with   $G\mapsto\eta_{2k}(G)$ being
linear,  such that
$$
E[A^k\psi\circ F\,\,G]=E[\psi\circ F\,\,\eta_{2k}(G)]\,,
$$
i.e.,
$$
E[f\circ F\,\,G]=E[(A^{-k}f)\circ F\,\,\eta_{2k}(G)]\,.
$$
Hence
$$
|E[f\circ F\,\,G]|\leq\|A^{-k}f\|_\infty\|\eta_{2k}(G)\|_{L^1}
$$
and
\begin{eqnarray*}
\sup_{\|G\|_{q,2k}\leq1}|E[f\circ F.G]| & \leq &
   \|A^{-k}f\|_\infty\sup_{\|G\|_{q,2k}\leq1}\|\eta_{2k}(G)\|_1 \\
& = & K\|f\|_{-2k}\,.
\end{eqnarray*}

\noindent
Consequently
$$
\|f\circ F\|_{p,-2k}\leq K\|f\|_{-2k}\,.
$$
\qed

\begin{corollary}
\label{ext-cor}
The linear  map $f\mapsto f\circ F$ from
${\calS}(\R^d)$ into $\DD$
extends continuously to a map from  ${\calS}'(\R^d)$ into $ \DD'$
whenever $F\in \DD(\R^d)$ is non-degenerate.
\end{corollary}

As we have seen in Theorem \ref{compo-thm} and Corollary
\ref{ext-cor}, if $F:W\to\R^d$ is a non-degenerate random variable, then
the map $f\mapsto f\circ F$ from ${\calS}(\R^d)\to \DD$ has a continuous
extension to ${\calS}'(\R^d)\to \DD'$ which we shall  denote by
$T\mapsto T\circ F$.

For $f\in{\calS}(\R^d)$, let us look at the following Pettis integral:
\index[not]{e@${\mathcal E}_x$}
$$
\int_{\R^d}f(x){\mathcal E}_x dx\,,
$$
where ${\mathcal E}_x$ denotes the Dirac measure at $x\in \R^d$. We have,
for any $g\in{\calS}(\R^d)$,
\begin{eqnarray*}
\Big\langle \int f(x){\mathcal E}_x dx,g\Big\rangle & = & \int\langle
f(x){\mathcal E}_x,g\rangle dx\\
&=&\int f(x)\langle{\mathcal E}_x,g\rangle dx \\
& = & \int f(x)g(x)dx=\langle f,g\rangle\,.
\end{eqnarray*}
Hence we have proven:

\begin{lemma}
\label{rep-lem}
The following representation holds in ${\calS}(\R ^d)$:
$$
f=\int_{\R ^d}f(x){\mathcal E}_xdx\, .
$$
\end{lemma}
From Lemma \ref{rep-lem},  we have

\begin{lemma}
We have
$$
\int\langle{\mathcal E}_y(F),\varphi\rangle f(y)dy=E[f(F)\,\varphi],
$$
 for any $\varphi\in \DD$, where $\langle\cdot,\cdot\rangle$ denotes
 the bilinear form of duality between $\DD'$ and $\DD$.
\end{lemma}
\proof
Let $\rho_\epsilon$ be a mollifier. Then ${\mathcal
E}_y\ast\rho_\epsilon\to{\mathcal E}_y$ in ${\calS}'$ on the other hand
\begin{eqnarray*}
\int_{\R^d}({\mathcal E}_y\ast\rho_\epsilon)(F)f(y)dy & = &
  \int_{\R^d}\rho_\epsilon(F+y)\,\,f(y)dy= \\
& = & \int_{\R^d}\rho_\epsilon(y)f(y+F)dy\,
\longrightarrow_{_{\hspace*{-1.5em}\scriptscriptstyle\epsilon\to 0}}\,f(F)\,.
\end{eqnarray*}
On the other hand, for $\varphi\in \DD$,
\begin{eqnarray*}
\lim_{\varepsilon\to 0}\int_{\R^d}
<({\mathcal E}_y\ast\rho_\epsilon)(F),\varphi>f(y)dy & = &
       \int_{\R^d}\lim_{\epsilon\to 0}<({\mathcal
         E}_y\ast\rho_\epsilon)(F),\varphi> f(y)dy \\
& = & \int_{\R^d}<{\mathcal E}_y(F),\varphi>f(y)dy \\
& = &< f(F),\varphi> \\
&=&E[f(F)\varphi]\, .
\end{eqnarray*}

\vspace*{-3ex}
\qed

\begin{corollary}
\label{rapid-cor}
We have
$$
\langle{\mathcal E}_x(F),1\rangle= \frac{d(F^\ast\mu)}{dx}(x)=p_F(x),
$$
moreover  $p_F\in {\calS}(\R^d)$ (i.e., the probability density of $F$
 is not only $C^\infty$ but it is also a rapidly decreasing function).
\end{corollary}
\proof
We know that, for any $\varphi \in \DD$,  the map
$T\to \langle T(F),\varphi\rangle$ is continuous
on ${\calS}'(\R^d)$ hence there exists  some
$p_{F,\varphi}\in{\calS}(\R^d)$ such that
$$
E[T(F)\,.\varphi]= {}_{\mathcal{S}}\langle
p_{F,\varphi},T\rangle_{\mathcal{S}'}\,.
$$
Let $p_{F,1}=p_F$, then it follows from the Lemma \ref{rep-lem}  that
\begin{eqnarray*}
E[f(F)]&=&\int\langle{\mathcal E}_y(F),1\rangle f(y)dy\\
&=&\int p_F(y) f(y)dy\,.
\end{eqnarray*}
\qed

\begin{remarkk}{\rm
From the disintegration of measures, we have
\begin{eqnarray*}
E[f(F)\,\,\varphi]&=&\int_{\R^d} p_F(x) E[\varphi|F=x]f(x)dx  \\
& = & \int_{\R^d} f(x)\langle{\mathcal E}_x(F),\varphi\rangle dx
\end{eqnarray*}
\noindent
hence
$$
E[\varphi|F=x]=\frac{\langle{\mathcal E}_x(F),\varphi\rangle}{p_F(x)}
$$
$dx$-almost surely on the support of the law of $F$.  In fact the
right hand side is an everywhere defined version of this conditional
probability.}
\end{remarkk}
\section{Supports of the laws  of  nondegenerate  Wiener functions }
Recall the local characters of the Sobolev derivative and the
divergence operators, namely: for any $\phi\in \DD_{p,1}$, $p>1$, we
have, for any $c\in\R$,  $\nabla\phi=0$ on the set $\{w:\,\phi(w)=c\}$
$\mu$-almost surely. There is also a similar result for the divergence
operator. A simple consequence of this observation is 
\begin{proposition}
\label{0-1-prop}
For a measurable subset $A$ of $W$ to be in some Sobolev space
$\DD_{p,1}$ $(p> 1)$, it is necessary and sufficient that
$\mu(A)\in\{0,1\}$.
\end{proposition}
\proof
Since the indicator function is idempotent, supposing the sufficiency,
we get $\nabla 1_A=2\,1_A\nabla 1_A$, which implies, from the locality
property of $\nabla$ explained above, that $\nabla
1_A=0$ a.s., hence we get from the Clark formula 
$$
1_A=E[1_A]=\mu(A)\,,
$$ 
which implies  $\mu(A)\in\{0,1\}$.
\qed

\begin{proposition}
\label{connex-prop}
Suppose that $F\in \DD_{p,1}(\R^m),\,p>1$. Then the support of the
measure induced by $F$, denoted as $F(\mu)$, is a connected subset of $\R^m$.
\end{proposition}
\proof
Suppose the contrary, then there exists two closed, disjoint sets
$A,B$ of positive $F(\mu)$ measure. For $k\geq 2$, choose $\psi_k\in
C^\infty(\R^m)$, with values in $[0,1]$, such that $\psi_k(x)=0$ for
$|x|\geq k$ and $\psi_k(x)=1$ for $|x|\leq k-1$ with the property that
$\sup_k\|D\psi_K\|_\infty<\infty$. Let $D_k$ be the closed ball of
radius $k$ in $\R^m$ and define $A_k=D_k\cap A,B_k=D_k\cap B$. For any
$k$, we choose an $f_K\in C^1(\R^m)$ which is equal to one on $A_k$
and to zero on $B_k$. It follows that 
\begin{equation}
\label{lim-pro}
(f_k\psi_k)\circ F\to 1_{\{F\in A\}}
\end{equation}
a.s. as $k\to\infty$. Moreover
\beaa
\nabla((f_k\psi_k)\circ F)&=&\sum_{i=1}^m [(\partial f_k\psi_k)\circ
F+(f_k\partial \psi_k)\circ F]\nabla F^i\\
&=&\sum_{i=1}^m (f_k\partial \psi_k)\circ F \nabla F^i
\eeaa
which implies that the sequence $((f_k\psi_k)\circ F,\,k\geq 2)$ is
bounded in $\DD_{p,1}$. This observation, combined with the limiting
property (\ref{lim-pro})  implies that $1_{\{F\in A\}}\in \DD_{p,1}$,
hence from Proposition \ref{0-1-prop}, we have $\mu\{F\in A\}=1$ which
is a contradiction to the assumption that $F(\mu)(B)>0$.
\qed

\remark Note that if $F$ is real-valued, then its support is an interval.

\begin{proposition}
\label{+-prop}
Assume that $F\in \DD_{p,1},\,p>2$, and that $F(\mu)$ possesses a
locally Lipschitz density $p_F(x)$. If $a\in \R$ is an element of the
interior of the support of $F(\mu)$, then $p_F(a)>0$.
\end{proposition}
\proof
Assume the contrary, i.e., that $p_F(a)=0$, let $r=2p/(2+p)>2$, by the
hypothesis, we should have $\mu\{F>a\}\in (0,1)$. Fix $\eps>0$ and
define 
$$
G_\eps(x)=\frac{1}{2\eps}\int_{-\infty}^x 1_{[a-\eps,a+\eps]}(y) dy\,.
$$
Since $G_\eps$ is a Lipschitz map, $G_\eps\circ F$ belongs to
$\DD_{n,1}$ for any $n\geq 0$, moreover, using the H\"older
inequality, we obtain 
$$
E[|\nabla (G_\eps\circ F)|_H^r]\leq E[|\nabla
F|_H^p]^{\frac{2}{2+p}}\left(\frac{1}{2\eps^2}\int_{a-\eps}^{a+\eps}p_F(y)dy\right)^{p/p+2}\,.
$$
Since we have supposed that $p_F(a)=0$, the integrand of the last
integral can be upper bounded  as $p_F(y)=|p_F(y)-p_F(a)|\leq K|y-a|$,
where $K$ is a Lipschitz constant of $p_F$, corresponding to its
Lipschitz character on  any fixed  compact interval containing
$[a-\eps,a+\eps]$. This upper estimate implies the boundedness of
$\DD_{p,1}$-norm of $G_\eps\circ F$ uniformly w.r. to $\eps$, which
implies that its weak limit, which is equal to the indicator function
of the set $\{F>a\}$, belongs to the Sobolev space $\DD_{p,1}$, which
is clearly a contradiction.
\qed

\subsection{Extension of the Ito Formula}
\markboth{Applications}{Ito Formula}
Let $(x_t)$ be the solution of the following stochastic differential equation:
\beaa
dx_t(w)&=&b_i(x_t(w))dt+\sigma_i(x_t(w))dw_t^i\\
x_0&=&x \,{\mbox{ given}},
\eeaa
where $b:\R^d\to\R^d$ and $\sigma_i:\R^d\to\R^d$ are smooth vector fields
with bounded derivatives. Let us denote by
$$
X_0=\sum_{i=1}^d\tilde{b}_0^i\frac{\partial}{\partial x_i}\,,\qquad
X_j=\sum \sigma_i^j\frac{\partial}{\partial x_j}
$$
\noindent
where
$$
\tilde{b}^i(x)=b^i(x)-\frac{1}{2}
\sum_{k,\alpha}\partial_k\sigma_\alpha^i(x)\sigma_\alpha^k(x).
$$
 If the Lie algebra of vector fields generated by $\{X_0,X_1,\ldots,X_d\}$ has
dimension equal to $d$ at any $x\in\R^d$, then $x_t(w)$ is non-degenerate
 cf. \cite{W1}. In fact it is also uniformly non-degenerate in the following
sense:
$$
E\int_s^t|{\mbox{\rm det}}(\nabla x_r^i,\nabla x_r^j)|^{-p} dr<\infty\,,
$$
forall $ 0<s<t$ and $ p>1$.

As a corollary of this result, combined with the lifting of ${\calS}'$
to $\DD'$, we can show the following:
\begin{theorem}
For any $T\in{\calS}'(\R^d)$, one has the following:
$$
T(x_t)-T(x_s)=\int_s^t AT(x_s)ds+
    \int_s^t\sigma_{ij}(x_s)\cdot\partial_j T(x_s)dW_s^i,
$$
where the Lebesgue integral is a Bochner integral, the stochastic integral
is as defined at the first section of this chapter and we have used
the following notation:
$$
A=\sum b^i\partial_i+\frac{1}{2}\sum a_{ij}(x)
   \frac{\partial^2}{\partial x_i\partial x_j}\,,\qquad
a(x)=(\sigma\sigma^\ast)_{ij},\;\sigma=[\sigma_1,\ldots,\sigma_d]\,.
$$
\end{theorem}
\subsection{Applications to the filtering\index[sub]{filtering} of
  the diffusions}
\markboth{Applications}{Filtering}
Suppose that we are given, for any $t\geq 0$,
$$
y_t=\int_0^t h(x_s)ds+B_t
$$
where $h\in
C_b^\infty(\R^d)\otimes\R^d$, $B$ is another Brownian motion independent of $w$
above. The process $(y_t;t\in[0,1])$ is called an (noisy) observation  of
$(x_t,t\in \reals_+)$. Let ${\mathcal Y}_t=\sigma\{y_s;s\in[0,t]\}$ be
the observed
data till $t$. The filtering problem consists of calculating the
random measure $f\mapsto E[f(x_t)|{\mathcal Y}_t]$. Let $P^0$ be the
probability defined by
$$
dP^0=Z_1^{-1}dP
$$
where
$$
Z_t=\exp\left\{\int_0^t (h(x_s),dy_s)-\frac{1}{2}
  \int_0^t|h(x_s)|^2ds\right\}\,.
$$
Then for any bounded, ${\mathcal Y}_t$-measurable random variable
$Y_t$, we have:
\begin{eqnarray*}
E[f(x_t).Y_t] & = & E\left[\frac{Z_t}{Z_t}f(x_t).Y_t\right]\\
&=&E^0\left[Z_tf(x_t)Y_t\right] \\
& = & E^0\left[E^0[Z_tf(x_t)|{\calY}_t]\cdot Y_t\right] \\
& = & E\left[ \frac{1}{E^0[Z_t|{\calY}_t]}E^0[Z_tf(x_t)|{\calY}_t]\cdot Y_t\right],
\end{eqnarray*}
hence
$$
E[f(x_t)|{\calY}_t]= \frac{E^0[Z_tf(x_t)|{\mathcal
Y}_t]}{E^0\Big[Z_t|{\calY}_t\Big]}\,.
$$

\noindent
If we want to study the smoothness of the measure $f\mapsto E[f(x_t)|{\mathcal
Y}_t]$, then from the above formula, we see that it is sufficient to study the
smoothness of $f\mapsto E^0[Z_tf(x_t)|{\calY}_t]$. The reason for the
use of $P^0$ \index[not]{p@$P^0$} is that
$w$ and $(y_t;t\in[0,1])$ are two independent Brownian motions
\footnote{ This claim follows directly from Paul L\'evy's theorem of
  the characterization of the Brownian motion.} under $P^0$
\begin{remarkk}
\index[sub]{Zakai equation}
{\rm Let us note that the random distribution $f\to \nu_t(f)$ defined
  by
$$
\nu_t(f)= E^0[Z_tf(x_t)|{\calY}_t]
$$
 satisfies the Zakai equation:
$$
\nu_t(f)=\nu_0(f)+\int_0^t\nu_s(Af)ds+\int_0^t\sum_{i}\nu_s(h_if)dy^i_s\,,
$$
where $A$ denotes the infinitesimal generator of the diffusion process
$(x_t,\,t\in \reals_+)$}
\end{remarkk}

After this preliminaries, we can prove the following
\begin{theorem}
Suppose that the map $f\mapsto f(x_t)$ from ${\mathcal{S}}(\R^d)$ into $\DD$
has a continuous extension as a map from ${\calS}'(\R^d)$ into $\DD'$.
Then the measure $f\mapsto E[f(x_t)|{\calY}_t]$ has a
density in ${\calS}(\R^d)$.
\end{theorem}
\proof
As explained above, it is sufficient to prove that the (random) measure
$f\mapsto E^0[Z_tf(x_t)|{\calY}_t]$ has a density in ${\calS}(\R^d)$.
Let ${\mathcal L}_y$ be the Ornstein-Uhlenbeck operator on the space of the
Brownian motion $(y_t;t\in[0,1])$. Then we have
$$
{\mathcal L}_yZ_t=Z_t\left( -\int_0^t h(x_s)dy_s+\frac{1}{2}
\int_0^t|h(x_s)|^2ds\right)\,\in\bigcap_p L^p\,.
$$
It is also easy to see that
$$
{\mathcal L}_w^k Z_t\in\bigcap_p L^p\,.
$$
From these observations we draw the following conclusions:
\begin{itemize}
\item
Hence $Z_t(w,y)\in \DD(w,y)$, where $\DD(w,y)$ denotes the space of
test functions defined on the product Wiener space with respect to the
laws of w and y.
\item
The second point is that the operator $E^0[\,\cdot\,|{\calY}_t]$ is a
continuous
mapping from $\DD_{p,k}(w,y)$ into $\DD_{p,k}^0(y)$, for any $p\geq
1,k\in\Z)$,  since ${\mathcal L}_y$ commutes with $E^0[\,\cdot\,|{\calY}_t]$ .
\item
Hence the map
$$
T\mapsto E^0[T(x_t)Z_t|{\calY}_t]
$$
 is continuous from
${\calS}'(\R^d)\to \DD'(y)$. In particular, for fixed $T\in{\calS}'$,
there exist  $ p>1$ and $k\in\N$ such that  $T(x_t)\in \DD_{p,-k}(w)$. Since
$Z_t\in \DD(w,y)$,
$$
Z_tT(x_t)\in \DD_{p,-k}(w,y)
$$
and
$$
T(x_t).(I+{\mathcal L}_y)^{k/2}Z_t\in \DD_{p,-k}(w,y).
$$
\item Consequently
$$
E^0[T(x_t)\cdot(I+{\mathcal
L}_y)^{k/2}Z_t|{\calY}_t]\in \DD_{p,-k}(y).
$$
\item
Finally it follows from the latter that
$$
(I+{\mathcal L})^{-k/2}E^0[T(x_t)(I+{\mathcal L}_y)^{k/2}Z_t|{\calY}_t]=
E^0[T(x_t)Z_t|{\calY}_t]
$$
belongs to $L^p(y)$. Therefore we see that:
$$
T\mapsto E^0[T(x_t)Z_t|{\calY}_t]
$$
defines a linear, continuous (use the closed graph theorem for instance) map
from
${\calS}'(\R^d)$ into $L^p(y)$.
\end{itemize}
Since ${\calS}'(\R^d)$ is a nuclear
space, the map
$$
T\stackrel{\Theta}{\mapsto}E^0[T(x_t)Z_t|{\calY}_t]
$$
is a nuclear  operator. This implies that $\Theta$  can be represented
as
$$
\Theta=\sum_{i=1}^\infty\lambda_i f_i\otimes\alpha_i
$$
where $(\lambda_i)\in l^1$, $(f_i)\subset{\calS}(\R^d)$ and
$(\alpha_i)\subset L^p(y)$ are bounded sequences. Define
$$
k_t(x,y)=\sum_{i=1}^\infty\lambda_i f_i(x)\alpha_i(y)\in{\calS}(\R^d)
\widetilde{\otimes}_1 L^p(y)\,
$$
where $\widetilde{\otimes}_1$ denotes the projective tensor product topology.
It is easy now to see that, for $g\in{\calS}(\R^d)$
$$
\int_{\R^d}g(x)k_t(x,y)dx=E^0[g(x_t)\cdot Z_t|{\calY}_t]\,
$$
and this completes the proof.
\qed

\section{Some applications of the Clark formula}
\markboth{Applications of Ito Clark Formula}{}
\subsection{Case of non-differentiable functionals}
In this example we use the Clark representation theorem for the
elements of $\DD'$ and the composition of the tempered distributions
with the non-degenerate Wiener functionals:
Let $w\mapsto \kappa(w)$ be the sign of the random variable
$w\mapsto W_1(w)$ where $W_1$ denotes the value of the Wiener path
$(W_t,t\in [0,1])$ at
time $t=1$. We have, using Theorem \ref{compo-thm}
$$
E[D_t\kappa|{\mathcal{F}}_t]=
2\exp\left\{-\frac{W_t^2}{2(1-t)}\right\}\frac{1}{\sqrt{2\pi(1-t)}}\,,
$$
$dt\times d\mu$-almost surely. Hence
$$
\kappa=2\int_0^1\exp\left\{-\frac{W_t^2}{2(1-t)}\right\}\frac{1}{\sqrt{2\pi(1-t)}}dW_t\,,
$$
$\mu$-almost surely. Note that, although $\kappa$ is not strongly Sobolev
differentiable, the integrand of the stochastic integral is an ordinary square
integrable process. This phenomena can be explained by the fact that the
conditional expectation tames the distribution, in such a way that the result
becomes an ordinary random variable.

Here is another application of the Clark formula:
\begin{proposition}
\label{set-independence}
Assume that $A$ is a measurable subset of $W$, then from Theorem
\ref{ust-thm}, there exists an $e_A\in L^2(\mu,H)$ which can be
represented as $e_A(t)=\int_0^t\dot{e}_A(\tau)d\tau,\,t\in [0,1]$,
such that $\dot{e}_A$ is adapted and
$$
\won_A=\mu(A)+\delta e_A\,.
$$
If $B$ is another measurable set, then $A$ and $B$ are independent if and
only if
$$
E\left[\left(e_A,e_B\right)_H\right]=0\,.
$$
\end{proposition}
\proof
It suffices to observe that
\begin{equation}
\label{scal-prod}
\mu(A\cap B)=\mu(A)\mu(B)+E[(e_A,e_B)_H]\,,
\end{equation}
hence $A$ and $B$ is independent if and only if the last term in
(\ref{scal-prod}) is null.
\qed

\subsection{Logarithmic Sobolev Inequality}
\markboth{Applications}{Log-Sobolev Inequality}
\index[sub]{logarithmic Sobolev inequality}
As another application  of the Clark representation theorem, we shall
give a quick proof of the logarithmic Sobolev inequality of L. Gross
{\footnote {The  proof which is given here  is similar  to that of 
    B. Maurey.}} (cf. \cite{LG}).
\begin{theorem}[log-Sobolev inequality]
For any $\phi\in \DD_{2,1}$, we have
$$
E[\phi^2\log \phi^2]\leq E[\phi^2]\log E[\phi^2]+2E[|\nabla
\phi|_H^2]\,.
$$
\end{theorem}
\proof
Clearly it suffices to prove the following inequality 
$$
E[f\log f]\leq \frac{1}{2}E\left[\frac{1}{f}|\nabla f|_H^2\right]\,,
$$
for any $f\in \DD_{2,1}$ which is strictly positive, lower bounded
with some $\eps>0$ and  with $E[f]=1$. Using
the It\^o-Clark representation theorem, we can write 
$$
f=\exp\left(\int_0^1\frac{E[D_sf|\calF_s]}{f_s}dW_s-\frac{1}{2}\int_0^1\left(\frac{E[D_sf|\calF_s]}{f_s}\right)^2ds\right)\,,
$$
where $f_s=E[f|\calF_s]$. It follows from the It\^o formula that 
$$
E[f\log
f]=\frac{1}{2}E\left[f\int_0^1\left(\frac{E[D_sf|\calF_s]}{f_s}\right)^2
ds\right]\,.
$$
Let $\nu$ be the probability defined by $d\nu=f\,d\mu$. Then we have 
\beaa
E[f\log f]&=&\frac{1}{2}E\left[f\int_0^1\left(\frac{E[f\,D_s\log
      f|\calF_s]}{f_s}\right)^2 ds\right]\\
&=&\frac{1}{2}E_\nu\left[\int_0^1 \left(E_\nu[D_s\log
      f|\calF_s]\right)^2 ds\right]\\
&\leq&\frac{1}{2}E_\nu\int_0^1(D_s\log f)^2ds\\
&=&\frac{1}{2}E[f\,|\nabla\log f|_H^2]\\
&=&\frac{1}{2}E\left[\frac{|\nabla f|_H^2}{f}\right]\,,
\eeaa


\qed

\begin{remarkk}{\rm  We have given the proof in the frame of the classical
Wiener space. However this result extends immediately to any abstract
Wiener space by the use of the techniques explained in the Appendix of
the fourth chapter.}
\end{remarkk}
\begin{remarkk}{\rm
A straightforward implication of the Clark representation, as we
have seen in the sequel of the  proof, is the {\bf Poincar\'e}
inequality which says that,  for any $F\in \DD_{2,1}$, one has
$$
E[|F-E[F]|^2]\leq E[|\nabla F|_H^2]\,.
$$
This inequality is the first step  towards  the logarithmic Sobolev
inequality.}
\end{remarkk}
\section*{Exercises}
{\footnotesize{
\begin{enumerate}
\item Assume that $F:W\to X$ is a measurable Wiener function,
where $X$ is a separable Hilbert space. Assume further that
$$
\|F(w+h)-F(w+k)\|_X\leq K|h-k|_H
$$
$\mu$-almost surely, for any $h,\,k\in H$. Prove that there exists
$F'=F$ almost surely such that
$$
\|F'(w+h)-F'(w+k)\|_X\leq K|h-k|_H
$$
for any $w\in W$ and $h,\,k\in H$.
\item Deduce from this result that if $A$ is a measurable subset
of $W$, such that $A+h\subset A$ almost surely, then $A$ has a
modification, say $A'$ such that $A'+H\subset A'$.
\end{enumerate}
}}

\section*{Notes and suggested reading}
{\footnotesize{Ito-Clark formula has been discovered first by Clark in
    the case of Fr\'echet differentiable Wiener functionals. Later its
    connections with the Girsanov theorem has been remarked by
    J.-M. Bismut, \cite{Bi}. D. Ocone has extended it to the Wiener
    functionals in $\DD_{2,1}$, cf. \cite{O}. Its extension to the
    distributions is due to the author, cf. \cite{U0}. Later D. Ocone
    and I. Karatzas have also extended it to the functionals of
    $\DD_{1,1}$.

 Composition of the non-degenerate Wiener functionals with the
 elements of ${\mathcal{S}}'(\R^d)$ is due to Kuo \cite{KUO1}.
 Watanabe has generalized it to more general Wiener functionals,
 \cite{W1,W}. Later it has been observed by the author  that this implies
 automatically the fact that the density of the law of a
 non-degenerate Wiener functional is a rapidly decreasing $C^\infty$
 function. This last result remains true for the conditional density
 of the non-linear filtering as it has been first proven in \cite{U4}.

}}

\chapter{Positive distributions and  applications}
\label{ch.posit}
\markboth{Positive Distributions}{}

\section{Positive Meyer-Watanabe distributions}

If $\theta$ is a positive distribution on $\R^d$, then a well-known theorem
says that $\theta$ is a positive measure, finite on the compact sets. We will
prove an analogous result for the  Meyer-Watanabe distributions  in this
section, show that they are absolutely continuous with respect to the
capacities defined with respect to the scale of the Sobolev spaces on
the Wiener space and  give an application to the construction of the
local time of the Wiener process. We end the chapter by making some
remarks about the Sobolev spaces constructed by the second
quantization of an elliptic operator on the Cameron-Martin space.

We will work on the classical Wiener space $C_0([0,1])=W$. First we have the
following:

\begin{proposition}
\label{basic-prop}
Suppose $(T_n)\subset \DD'$ and each $T_n$ is also a probability on $W$. If
$T_n\to T$ in $ \DD'$, then $T$ is also a probability and $T_n\to T$
in the weak-star  topology of measures on $W$.
\end{proposition}

\noindent
For the proof of this proposition, we shall need the following result
whose proof can be found in \cite{S-V}
\begin{lemma}[Garsia-Rademich-Ramsey lemma]
\label{lemma-GRR}
Let $p,\psi$ be two continuous, stritly increasing functions on $\R_+$
such that $\psi(0)=p(0)=0$ and that
$\lim_{t\to\infty}\psi(t)=\infty$. Let $T>0$ and $f\in
C([0,T],\R^d)$. If
$$
\int_{[0,T]^2}\psi\left(\frac{|f(t)-f(s)|}{p(|t-s|)}\right)ds \,dt\leq
B\,,
$$
then for any $0\leq s\leq t\leq T$, we have 
$$
|f(t)-f(s)|\leq
8\int_0^{t-s}\psi^{-1}\left(\frac{4B}{u^2}\right)p(du)\,.
$$
\end{lemma}
\proof[{\bf{of the Proposition}}]
It is sufficient to prove that the sequence  of  probability measures
$(\nu_n, n\geq 1)$
associated to $(T_n,n\geq 1)$, is tight. In fact, let $S= \DD \cap C_b(W)$,
if the tightness holds, then we would  have, for $\nu=w-\lim \nu_n$
(taking a subsequence if necessary),
where $w-\lim$ denotes the limit in the weak-star topology of measures, 
$$
\nu(\varphi)=T(\varphi)\quad {\mbox{on}}\; S\,.
$$
Since the mapping  $w\to e^{i\langle w,w^\ast\rangle}$ ($w^\ast\in
W^\ast$) belongs to $S$, $S$ separates the probability measures on
$(W,\B(W))$  and the proof would follow.

In order to realize this program, let $G:W\to \R$ be defined as
$$
G(w)=\int_0^1\!\!\int_0^1
\frac{|w(t)-w(s)|^8}{|t-s|^3}ds\,dt.
$$
Then  $G\in \DD$ and
$A_{\lambda}=\{G(w)\leq\lambda\}$ is a compact
subset of $W$. In fact, from  Lemma \ref{lemma-GRR} , the
inequality $G(w)\leq \la$ implies the existence of a constant
$K_\la$ such that  
$$
|w(s)-w(t)|\leq K_\la |t-s|^{\frac{15}{4}}\,,
$$
for $0\leq s<t\leq 1$, hence $A_\la$ is equicontinuous, then the
Arzela-Ascoli Theorem implies that the set $\{w:\,G(w)\leq \la\}$ is
relatively compact in $W$, moreover it is a  closed set  since $G$ 
is a lower semi-continuous function by the Fatou Lemma. In particular,
it is measurable with respect to the non-completed Borel sigma algebra
of $W$.
Moreover, we have $\bigcup_{\lambda\geq0}A_\lambda=W$ almost surely. Let
$\varphi\in C^\infty(\R)$ such that  $0\leq\varphi\leq1$; $\varphi(x)=1$ for
$x\geq0$, $\varphi(x)=0$ for $x\leq-1$. Let
$\varphi_\lambda(x)=\varphi(x-\lambda)$.
We have
$$
\nu_n(A_\lambda^c)\leq\int_W\varphi_\lambda(G(w))\,\nu_n(dw)\,.
$$

\noindent
We claim that
$$
\int_W\varphi_\lambda(G)d\nu_n=\langle\varphi_\lambda(G),T_n\rangle.
$$

To see this, for $\varepsilon>0$, write
$$
G_\varepsilon(w)=\int_{[0,1]^2}
\frac{|w(t)-w(s)|^8}{(\varepsilon+|t-s|)^3}ds\,dt\,.
$$
Then $\varphi_\lambda(G_\varepsilon)\in S$ (but not $\varphi_\lambda(G)$, since
$G$ is not continuous on $W$)
and we have
$$
\int\varphi_\lambda(G_\varepsilon)d\nu_n=
\langle\varphi_\lambda(G_\varepsilon),T_n\rangle.
$$

\noindent
Moreover  $\varphi_\lambda(G_\varepsilon)\to\varphi_\lambda(G)$ in $\DD$, hence
$$
\lim_{\varepsilon\to  0}\langle\varphi_\lambda(G_\varepsilon),
T_n\rangle=\langle\varphi_\lambda(G),T_n\rangle\,.
$$
From the dominated convergence theorem, we have also
$$
\lim_{\varepsilon\to 0}\int\varphi_\lambda(G_\varepsilon)d\nu_n=
\int\varphi_\lambda(G)d\nu_n.
$$
Since $T_n\to T$ in
$\DD'$, there  exist  some $ k>0$ and  $p>1$ such that  $T_n\to T$ in
$\DD_{p,-k}$.
Therefore
\begin{eqnarray*}
\left\langle\varphi_\lambda(G),T_n\right\rangle & = &
    \langle(I+{\mathcal L})^{k/2}\varphi_\lambda(G),(I+{\mathcal L})^{-k/2}T_n\rangle \\
& \leq &\left \|(I+{\mathcal L})^{k/2}\varphi_\lambda(G)\right\|_q\,\,
\sup_n\left\|(I+{\mathcal L})^{-k/2}T_n\right\|_p\,.
\end{eqnarray*}
From the Meyer inequalities, we see that
$$
\lim_{\lambda\to\infty}\left\|(I+{\mathcal L})^{k/2}\varphi_\lambda(G)\right\|_q=0\,,
$$
in fact, it is sufficient to see that $\nabla^i(\varphi_\lambda(G))\to0$ in
$L^p$ for all $i\leq [k]+1$, but this is obvious from the choice of
$\varphi_\lambda$. We  have proven that
\beaa
\lefteqn{\lim_{\lambda\to\infty}\sup_n \nu_n(A_\lambda^c)}\\
&\leq& \sup_n\left\|(I+{\mathcal
L})^{-k/2} T_n\right\|_p\lim_{\lambda\to\infty}\left\|(I+{\mathcal
L})^{k/2}\varphi_\lambda (G)\right\|_p=0\,,
\eeaa
which implies the tightness and the proof is completed.
\qed

\begin{corollary}
\label{pos-dist}
Let $T\in \DD'$ such that $\langle T,\varphi\rangle\geq0$, for all positive
$\varphi\in \DD$. Then $T$ is a Radon measure on $W$.
\end{corollary}
\proof
Let $(h_i)\subset H$ be a complete, orthonormal basis  of $H$. Let
$V_n=\sigma\{\delta h_1,\ldots,\delta h_n\}$.
Define $T_n$ as $T_n=E[P_{1/n}T|V_n]$ where $P_{1/n}$
is the Ornstein-Uhlenbeck  semi-group on $W$. Then $T_n\geq0$ and it is
a random variable in some $L^p(\mu)$. Therefore it defines a measure on
$W$ (it is even  absolutely continuous with respect to $\mu$). Moreover
$T_n\to T$ in $\DD'$, hence the proof follows from Proposition
\ref{basic-prop}.
\qed

Another application  is the following:
\begin{proposition}
\label{++-prop}
Let $F\in\DD(\R^m)$ be a nondegenerate random vector and denote by
$p_F$ the density of its law, which is described as in Corollary
\ref{rapid-cor}. If $p_F(a)=0$ for some $a\in \R^m$, then any derivative
$\partial^\alpha p_F$ is again zero at this point $a\in \R^m$.
\end{proposition}
\proof
Assume then $p_F(a)=0$, since $p_F(a)=<\E_a(F),1>$ and we know
that $\E_a\circ F$ is a positive measure on the Wiener space, and
$p_F(a)=0$ implies that this measure is equal to zero. We can also
write 
$$
\nabla \E_a(F)=\sum_i\partial_i\E_a(F)\nabla F_i\,,
$$
hence 
$$
\partial_i\E_a(F)=\sum_{j\leq m}(\nabla \E_a(F),\nabla
F^j)_H(\gamma_F^{-1})^{ji}=0\,,
$$
where $\gamma_F$ is the Gramm-Malliavin matrix associated to
$F$. Since 
$$
\partial^\alpha p_F(x)=<(\partial \E_x)(F),1>\,,
$$
for any multi-index $\alpha$, the proof follows by induction.
\qed

\section{Capacities\index[sub]{capacity} and positive Wiener
  func\-tio\-nals}
\markboth{Capacities}{Positive Distributions}
We begin with the following definitions:
\begin{definition}
Let  $p\in[1,\infty)$ and $k>0$. If  $O\subset W$ is an open set, we define the
$(p,k)$-capacity of $O$ as \index[not]{ca@$C_{p,k}$}
$$
C_{p,k}(O)=\inf\{\|\varphi\|_{p,k}^p:\varphi\in \DD_{p,k},
\varphi\geq 1 \mu-{\mbox{a.e.  on }}O\}\,.
$$
If $A\subset W$ is any set, define its $(p,k)$-capacity as
$$
C_{p,k}(A)=\inf\{C_{p,k}(O);O\;{\mbox{is open}}\;O\supset A\}\,.
$$
\end{definition}
\begin{itemize}
\item We say that some property takes place
\index[sub]{quasi-everywhere} $(p,k)$-quasi everywhere if the set on
which it does not hold has $(p,k)$-capacity zero.
\item  We say $N$ is a slim set\index[sub]{slim set} if
$C_{p,k}(N)=0$, for all $ p>1$, $k>0$.
\item A function is called $(p,k)${\bf{-quasi continuous\/}} if
for any  $\varepsilon>0$ \index[sub]{quasi-continuous}, there exists
an  open set $O_\varepsilon$  such that
$C_{p,k}(O_\varepsilon)<\varepsilon$ and the function is continuous on
$O_\varepsilon^c $.
\item  A function is  called $\infty$-quasi continuous if it is $(p,k)$-quasi
continuous for any  $p>1,\,k\in \NN$.
\end{itemize}

\noindent
The  results contained in the next lemma  are proved by Fukushima \&
Kaneko (cf. \cite{Fu-K}):
\begin{lemma}
\label{F-K-lemma}
\begin{enumerate}
\item
If $F\in \DD_{p,k}$, then there exists a $(p,k)$-quasi continuous function
$\tilde{F}$ such that $F=\tilde{F}$ $\mu$-a.e. and $\tilde{F}$ is
$(p,k)$-quasi everywhere defined, i.e.\ if $\tilde{G}$ is another such
function, then $C_{p,k}(\{\tilde{F}\not=\tilde{G}\}))=0$.
\item
If $A\subset W$ is arbitrary, then
$$
C_{p,k}(A)=\inf\{\|\varphi\|_{p,k}:\varphi\in \DD_{p,k}\,,\quad
\tilde{\varphi}\geq 1\;(p,r)-q.e.\;{\mbox{on}}\;A\}
$$
\item
There exists a unique element $U_A\in \DD_{p,k}$ such that  $\tilde{U}_A\geq1$
$(p,k)$-quasi everywhere  on $A$ with  $C_{p,k}(A)=\|U_A\|_{p,k}$, and
 $\tilde{U}_A\geq0$ $ (p,k)$-quasi everywhere.  $U_A$ is called the $(p,k)$-equilibrium potential of $A$. \index[sub]{equilibrium potential}
\end{enumerate}
\end{lemma}
\begin{theorem}
Let $T\in \DD'$ be a positive distribution and suppose that $T\in \DD_{q,-k}$ for
some $q>1$, $k\geq0$. Then, if we denote by $\nu_T$ the measure associated to
$T$, we have
$$
\bar{\nu}_T(A)\leq\|T\|_{q,-k}(C_{p,k}(A))^{1/p},
$$
for any set $A\subset W$, where $\bar{\nu}_T$ denotes the outer measure
with respect to
$\nu_T$. In particular $\nu_T$ does not charge the slim sets.
\end{theorem}
\proof
Let $V$ be an open set in $W$ and let $U_V$ be its equilibrium potential of
order $(p,k)$. We have
\begin{eqnarray*}
\langle P_{1/n}T,U_V\rangle & = &
\int P_{1/n}T\, U_V d\mu\\
&\geq&\int_V P_{1/n}T\, U_V d\mu \\
& \geq& \int_V P_{1/n}T d\mu\\
&=&\nu_{P_{1/n}T}(V)\,.
\end{eqnarray*}

\noindent
Since $V$ is open, we have, from the fact that $\nu_{P_{1/n}T}\to\nu_T$ weakly,
$$
\liminf_{n\to\infty}\nu_{P_{1/n}T}(V)\geq\nu_T(V)\,.
$$
On the other hand
\begin{eqnarray*}
\lim_{n\to\infty}\langle P_{1/n}T,U_V\rangle&=&\langle T,U_V\rangle \\
& \leq & \|T\|_{q,-k}\|U_V\|_{p,k} \\
& = & \|T\|_{q,-k}C_{p,k}(V)^{1/p}.
\end{eqnarray*}
\qed
\section{Some  Applications}
Below we  use  the characterization of the  positive distributions to
give a different interpretation of the local times. Afterwards the
$0-1$ law is revisited via the capacities.
\subsection{Applications to Ito formula and local times}
\markboth{Ito Formula}{Local Time}
 \ Let $f:\R^d\to\R$ be a function from ${\calS}'(\R^d)$ and suppose that
$(X_t,t\geq 0)$ is a hypoelliptic  diffusion on $\R^d$
which  is constructed as the solution of the following
stochastic differential equation with smooth coefficients:
\begin{eqnarray}
\label{SDE}
dX_t&=&\sigma(X_t)dW_t+b(X_t)dt\\
X_0&=&x\in \R^d\,. \nonumber
\end{eqnarray}
We denote by  $L$  the  infinitesimal generator of the diffusion
process $(X_t,t\geq 0)$. For any $t>0$, $X_t$  is a non-degenerate
random variable in the sense of Definition
\ref{non-deg-def}. Consequently  we have the extension of the  Ito
formula 
$$
f(X_t)-f(X_u)=\int_u^t Lf(X_s)ds+
\int_u^t\sigma_{ij}(X_s)\partial_i f(X_s)dW_s^j\,,
$$
for $0<u\leq t\leq 1$. Note that, since we did not make any
differentiability hypothesis about $f$, the above integrals are to be
regarded  as the elements of $\DD'$. Suppose that $Lf$ is a bounded
measure on $\R^d$, from our result about the positive distributions,
we see that $\int_u^t Lf(X_s)ds$ is a measure on $W$ which does not
charge the slim sets. By difference, so does the term
$\int_u^t\sigma_{ij}(X_s)\partial_if(X_s)dW_s^j$.

\noindent
As a particular case, we can take  $d=1$, $L=\frac{1}{2}\Delta$
(i.e. $\sigma=1$), $f(x)=|x|$ and this gives
$$
|W_t|-|W_u|=\frac{1}{2}\int_u^t\Delta|x|(W_s)ds+
\int_u^t\frac{d}{dx}|x|(W_s)dW_s\,.
$$
As $\frac{d}{dx}|x|={\mbox{sign}}(x)$, we have
$$
\int_u^t\frac{d}{dx}|x|(W_s)dW_s= \int_u^t{\mbox{sign}}(W_s)dW_s=M_t^u
$$
is a measure absolutely continuous with respect to $\mu$.
Since $\lim_{u\to0}M_t^u=N_t$ exists in all $L^p$, so does
$$
\lim_{u\to0}\int_u^t\Delta|x|(W_s)ds
$$
in $L^p$ for any $p\geq1$. Consequently $\int_0^t\Delta|x|(W_s)ds$ is
absolutely continuous with respect to $\mu$, i.e., it is a random variable.
It is easy to see that
$$
\Delta|x|(W_s)=2{\mathcal E}_0(W_s)\,,
$$
where ${\mathcal E}_0$ denotes the Dirac measure at zero, hence  we obtain
\begin{eqnarray*}
\int_0^t2{\mathcal E}_0(W_s)ds&=&\int_0^t\Delta|x|(W_s)ds\\
&=&2l_t^0
\end{eqnarray*}
which is the local time\index[sub]{local time} of Tanaka. Note that, although
${\mathcal E}_0(W_s)$ is singular
with respect to $\mu$, its Pettis integral is absolutely continuous with
respect to  $\mu$.

\begin{remarkk}{\rm  If $F:W\to\R^d$ is a non-degenerate random
    variable, then for any 
$S\in{\mathcal{S}}'(\R^d)$ with $S\geq0$ on ${\calS}_+(\R^d)$,
$S(F)\in \DD'$ is a
positive distribution, hence it is a positive  Radon measure on $W$. In
particular ${\mathcal E}_x(F)$ is a positive  Radon measure.
}
\end{remarkk}
\subsection{Applications to $0-1$ law and to  the gauge functionals of
  sets}
\markboth{Applications}{$0-1$-Law}
In Theorem \ref{0-1-law} we have seen that an $H$-nvariant subset of
$W$ has measure which is equal either to zero or to one. In this
section we shall refine this result using the capacities. Let us first
begin by defining the gauge function of a measurable subset of $W$: if
$A\in \calB(W)$, define
\begin{equation}
\label{lip-con}
q_A(w)=\inf\bigl[|h|_H:\,h\in (A-w)\cap H\bigr]\,,
\end{equation}
where the infimum is defined as to be infinity on the empty set. We
have
\begin{lemma}
\label{gauge-lem}
For any $A\in \calB(W)$, the map $q_A$ is measurable with respect to
the $\mu$-completion of $\calB(W)$. Moreover
\begin{equation}
\label{lip-property}
|q_A(w+h)-q_A(w)|\leq |h|_H
\end{equation}
almost surely, for any $h\in H$ and $\mu\{q_A<\infty\}=0$ or $1$.
\end{lemma}
\proof
Without loss of generality, we may assume that $A$ is a compact subset
of $W$ with $\mu(A)>0$. Then the set $K(w)=(A-w)\cap H\neq \emptyset$
almost surely. Therefore $w\to K(w)$ is a multivalued map with values
in the non-empty subsets of $H$ for almost all $w\in W$. Let us denote
by $G(K)$ its graph, i.e.,
$$
G(K)=\{(h,w):h\in K(w)\}\,.
$$
Since $(h,w)\mapsto h+w$ is measurable from $H\times W$ to $W$ when
the first space is equipped with the product sigma algebra, due to the
continuity of the map $(h,w)\to w+h$, it follows that $G(K)$ is a
measurable subset of $H\times W$. From a theorem about the measurable
multi-valued maps, it follows that $w\to K(w)$ is measurable with
respect to the $\mu$- completed sigma field $\calB(W)$
(cf. \cite{C-V}).
Hence  there
is a countable sequence of $H$-valued measurable  selectors $(u_i,i\in
\NN)$ of $K$ (i.e., $u_i:W\to H$ such that $u_i(w)\in K(w)$ almost
surely) such that $(u_i(w),i\in \NN)$ is dense in $K(w)$ almost
surely. To see the measurability, it suffices to remark that
$$
q_A(w)=\inf(|u_i(w)|_H:i\in \NN)\,.
$$
The relation \ref{lip-property} is evident from the definition of $q_A$. To
complete the proof it suffices to remark that the set
$Z=\{w:\,q_A(w)<\infty\}$ is $H$-invariant, hence from Theorem
\ref{0-1-law}, $\mu(Z)=0$ or $1$. Since $Z$ contains $A$ and $\mu(A)>0$,
 $\mu(Z)=1$.
\qed

The following result refines the $0-1$--law (cf. \cite{F-P},  \cite{Kus1}):
\begin{theorem}
\label{0-1-cap}
Assume that $A\subset W$ is an $H$-invariant set of zero Wiener
measure. Then
$$
C_{r,1}(A)=0
$$
for any $r>1$.
\end{theorem}
\proof
Choose a compact  $K\subset A^c$ with $\mu(K)>0$. Denote by $B_n$ the
ball of radius $n$ of $H$ and define $K_n=K+B_n$. It is easy to see
that $\cup_nK_n$ is an $H$-invariant set. Moreoever
$$
\left(\cup_n K_n\right)\cap A=\emptyset\,,
$$
otherwise, due to the $H$-invariance of of $A$, we would have $A\cap
K\neq \emptyset$. We also have $\mu(K_n)\to 1$. Let
$$
p_n(w)=q_{K_n}(w)\wedge 1\,.
$$
From Proposition \ref{sob-prop}, we see that $p_n\in
\cap_p\DD_{p,1}$. Moreover $p_n(w)=1$ on $K_{n+1}^c$ (hence on $A$) by
construction. Since $p_n=0$ on $K_n$, from  Lemma \ref{der-loc}
$\nabla p_n=0$ almost surely  on  $K_n$. Consequently
\beaa
C_{r,1}(A)&\leq&\int (|p_n|^r+|\nabla p_n|_H^r)d\mu\\
          &=&\int_{K_n^c} (|p_n|^r+|\nabla p_n|_H^r)d\mu\\
          &\leq &2\mu(K_n^c)\to 0
\eeaa
as $n\to \infty$.
\qed

\section{Local Sobolev spaces}
\markboth{Applications}{Local Sobolev Spaces}
In Chapter II we have observed the local character of the Sobolev
derivative and the divergence operator. This permits us to define the
local Sobolev spaces as follows:
\begin{definition}
We say that a Wiener functional $F$ with values in some separable
Hilbert space $X$ belongs to $\Dll_{p,1}(X)$, $p>1$,  if there exists a
sequence $(\Om_n,n\geq 1)$ of  measurable subsets of $W$ whose union
is equal to $W$ almost surely and
$$
F=F_n \,\,{\mbox{ a.s. on }}\,\Om_n\,,
$$
where $F_n\in \DD_{p,1}(X)$ for any $n\geq 1$. We call
$((\Om_n,F_n),n\geq 1)$ a localizing sequence for $F$.
\index[sub]{localizing sequence}
\index[not]{d@$\Nl$}
\index[not]{div@$\dl$}
\end{definition}
Lemma \ref{der-loc} and Lemma \ref{div-loc} of Section
\ref{locality-section} permit us to define the local Sobolev
derivative and local divergence of the Wiener functionals. In fact, if
$F\in \Dll_{p,1}(X)$, then we define $\Nl F$ as
$$
\Nl F=\nabla F_n\,\,{\mbox{on }}\,\Om_n\,.
$$
Similarly, if $\xi\in \Dll_{p,1}(X\otimes H)$, then we define
$$
\dl \xi=\delta \xi_n\,\,{\mbox{ on}}\,\Om_n\,.
$$
From the lemmas quoted  above $\Nl F$ and $\dl \xi$ are independent of
the choice of their localizing sequences.

\remark Note that we can define also the spaces $\Dll_{p,k}(X)$
similarly.

\noindent
The most essential property of the Sobolev derivative and the
divergence operator is the fact that the latter is the adjoint of the
former under the Wiener measure. In other words they satisfy the
integration by parts formula:
$$
E[(\nabla \phi,\xi)_H]=E[\phi\,\delta \xi]\,.
$$
In general this important formula is no longer valid when we replace
$\nabla$ and  $\delta$ with $\Nl$ and $\dl$ respectively. The theorem
given below gives the exact  condition when the local derivative or
divergence of a vector field is in fact equal to the global one.
\begin{theorem}
\label{derloc-thm}
Assume that $\phi\in \Dll_{p,1}(X)$, and let $((\phi_n,\Om_n),n\in
\NN)$ be a localizing sequence of $\phi$. A neccessary and sufficient
condition for $\phi\in \DD_{p,1}(X)$ and for  $\nabla \phi=\Nl\phi$ almost
surely, is
\begin{equation}
\label{cap-cond}
\lim_{n\to \infty}C_{p,1}(\Om^c_n)=0\,.
\end{equation}
\end{theorem}
\proof
The neccessity is trivial since, from Lemma \ref{F-K-lemma}. To prove
the sufficiency we can assume without loss of generality that $\phi$
is bounded. In
fact, if the theorem is proved for the bounded functions, then to
prove the general case,  we can replace $\phi$ by
$$
\phi_k=\left(1+\frac{1}{k}\|\phi\|_X\right)^{-1}\phi\,,
$$
which converges in $\DD_{p,1}(X)$ as $k\to \infty$ due to the
closedness of the Sobolev derivative. Hence we shall assume that
$\phi$ is bounded. Let $\eps>0$ be arbitrary, since
$C_{p,1}(\Om_n^c)\to 0$, by Lemma \ref{F-K-lemma}, there exists some
$F_n\in \DD_{p,1}$ such that $F_n\geq 1$ on $\Om_n^c$ quasi-everywhere
and $\|F_n\|_{p,1}\leq C_{p,1}(\Om_n^c)+\eps 2^{-n}$, for any $n\in
\NN$. Evidently, the sequence $(F_n,n\in \NN)$ converges to zero in
$\DD_{p,1}$. Let $f:\R\to [0,1]$ be a smooth function such that
$f(t)=0$ for $|t|\geq 3/4$ and $f(t)=1$ for $|t|\leq 1/2$. Define
$A_n=f\circ F_n$, then $A_n=0$ on $\Om_n^c$ quasi-everywhere and the
sequence $(A_n,n\in \NN)$ converges to the constant $1$ in
$\DD_{p,1}$. As a consequence  of this observation
$\phi\,A_n=\phi_n\,A_n$ almost surely and by the dominated convergence
theorem, $(\phi_nA_n,n\in \NN)$ converges to $\phi$ in
$L^p(\mu,X)$. Moreover
\beaa
\nabla(\phi\,A_n)&=&\nabla(\phi_nA_n)\\
                &=&A_n\Nl\phi+\phi\nabla A_n\to \Nl\phi
\eeaa
in $L^p(\mu,X\otimes H)$ since $(A_n,n\in \NN)$ and $\phi$ are
bounded. Consequently $\nabla(\phi_nA_n)\to \Nl \phi$ in
$L^p(\mu,X\otimes H)$, since $\nabla$ is a closed operator on
$L^p(\mu,X)$ the convergence takes place also in $\DD_{p,1}(X)$ and
the proof is completed.
\qed

We have also a similar result for the divergence operator:
\begin{theorem}
\label{divloc-thm}
Let  $\xi$ be in $\Dll_{p,1}(H)$ with a localizing sequence
$((\xi_n,\Om_n),n\in \NN)$ such  that $\xi\in L^p(\mu,H)$ and
$\dl\xi\in L^p(\mu)$. Assume moreover
\begin{equation}
\label{cap-q-con}
\lim_{n\to \infty}C_{q,1}(\Om_n^c)=0\,,
\end{equation}
where $q=p/(p-1)$. Then $\xi\in \Dom_p(\delta)$ and $\dl\xi=\delta\xi$
almost surely.
\end{theorem}
\proof
Due to  the hypothesis (\ref{cap-q-con}), we can construct a sequence
$(A_n,n\in \NN)$ as in the proof of Theorem \ref{derloc-thm}, which is
bounded in $L^\infty(\mu)$, converging to the constant function $1$ in
$\DD_{q,1}$ such that $A_n=0$ on $\Om_n^c$. Let $\gamma\in \DD$ be
bounded, with a bounded Sobolev derivative. We have
\beaa
E\left[A_n(\dl \xi)\gamma\right]&=&E\left[A_n\delta\xi_n\,\gamma\right]\\
  &=&E\left[A_n(\xi_n,\nabla\gamma)_H\right]+E\left[(\nabla
    A_n,\xi_n)_H\gamma\right]\\
&=&E\left[A_n(\xi,\nabla\gamma)_H\right]+E\left[(\nabla
  A_n,\xi)_H\gamma\right]\\
&&\to E\left[(\xi,\nabla\gamma)_H\right]\,.
\eeaa
Moreover, from the dominated convergence theorem we have
$$
\lim_nE[A_n(\dl\xi)\gamma]=E[(\dl\xi)\gamma]\,,
$$
hence
$$
E\left[(\dl\xi)\,\gamma\right]=E\left[(\xi,\nabla\gamma)_H\right]\,.
$$
Since the set of functionals $\gamma$ with the above prescribed
properties is dense in $\DD_{q,1}$, the proof is completed.
\qed

\section{Distributions associated to
  $\Gamma(A)$\index[not]{g@$\Gamma(A)$}}
\markboth{Applications}{$\Gamma(A)$-Distributions}
It is sometimes useful to have a scale of distribution spaces which
are defined with a more ``elliptic'' operator than the
Ornstein-Uhlenbeck semigroup. In this way objects which are more
singular than Meyer distributions can be interpreted as the elements
of the dual space. This is important essentially for the constructive
Quantum field theory, cf. \cite{Si}.
We begin with an abtract Wiener space  $(W,H,\mu)$. Let $A$ be a
self-adjoint operator on $H$, we suppose that
its spectrum lies in $(1,\infty)$, hence $A^{-1}$ is bounded and
$\|A^{-1}\|<1$. Let
$$
H_\infty=\bigcap_n {\mbox{\rm Dom}}(A^n)\,,
$$

\noindent
hence $H_\infty$ is dense in $H$ and $\alpha\mapsto(A^\alpha h,h)_H$ is
increasing. Denote by $H_\alpha$ the completion of $H_\infty$ with respect to
the norm $|h|_\alpha^2=(A^\alpha h,h)$; $\alpha\in\R$. Evidently
$H_\alpha'\cong H_{-\alpha}$ (isomorphism). If $\varphi:W\to\R$ is a nice
Wiener functional with $\varphi=\sum_{n=0}^\infty I_n(\varphi_n)$,
define the second quantization of $A$\index[sub]{second quantization}
$$
\Gamma(A)\varphi=E[\varphi]+
    \sum_{n=1}^ \infty I_n(A^{\otimes n}\varphi_n)\,.
$$

\begin{definition}
For $p>1$, $k\in\Z$, $\alpha\in\R$, we define $\DD_{p,k}^\alpha$ as the
completion of polynomials based on $H_\infty$, with respect to the norm:
$$
\|\varphi\|_{p,k;\alpha}=\|(I+{\mathcal L})^{k/2}
       \Gamma(A^{\alpha/2})\varphi\|_{L^p(\mu)}\,,
$$
where $\varphi(w)=p(\delta h_1,\ldots,\delta h_n)$, $p$ is a
polynomial on $\reals^n$ and $h_i\in H_\infty\,$.
If $\Xi$ is a separable Hilbert space,
$\DD_{p,k}^\alpha({\Xi})$ is defined
similarly except that $\varphi$ is taken as an $\Xi$-valued polynomial.
\end{definition}

\begin{remarkk}{\rm

 If $\varphi=\exp(\delta h-\frac{1}{2}|h|^2)$ then we have
$$
\Gamma(A)\varphi=\exp\left\{\delta(Ah)-
{\textstyle\frac{1}{2}}|Ah|^2\right\}\,.
$$}
\end{remarkk}
\begin{remarkk}{\rm
 $\DD_{p,k}^\alpha$ is decreasing with respect to $\alpha,p$ and
  $k$.}
\end{remarkk}

\begin{theorem}
\label{H1-th}
Let $(W^\alpha,H_\alpha,\mu_\alpha)$ be the abstract Wiener space
corresponding to the
Cameron-Martin space $H_\alpha$. Let us denote by $\DD_{p,k}^{(\alpha)}$ the
Sobolev space on $W^\alpha$ defined by
$$
\|\varphi\|_{\DD_{p,k}^{(\alpha)}}=
      \|(I+{\mathcal L})^{k/2}\varphi\|_{L^p(\mu_\alpha,W^\alpha)}
$$

\noindent
Then $\DD_{p,k}^{(\alpha)}$ and $\DD_{p,k}^\alpha$ are isomorphic.
\end{theorem}
\remark
This isomorphism is not algebraic, i.e., it does not commute with the
point-wise multiplication fion.

\proof
We have
$$
E[e^{i\delta(A^{\alpha/2}h)}]= \exp{\textstyle\frac{1}{2}}
|A^{\alpha/2}h|^2=\exp\frac{|h|_\alpha^2}{2}
$$
which is the characteristic function of $\mu_\alpha$ on $W^\alpha$.
\qed

\begin{theorem}
\label{H2-th}
\begin{enumerate}
\item   For $p>2$, $\alpha\in\R$, $k\in\Z$, there exists some
$\beta>\frac{\alpha}{2}$ such that
$$
\|\varphi\|_{\DD_{p,k}^\alpha}\leq\|\varphi\|_{\DD_{2,k}^\beta}
$$
consequently
$$
\bigcap_{\alpha,k}\DD_{2,k}^\alpha=
\bigcap_{\alpha,p,k} \DD_{p,k}^\alpha\,.
$$
\item  Moreover, for some $\beta>\alpha$ we have
$$
\|\varphi\|_{\DD_{p,k}^\alpha}\leq\|\varphi\|_{\DD_{2,0}^\beta}\,,
$$
\end{enumerate}
hence we have also
$$
\bigcap_{\alpha}\DD_{2,0}^\alpha=
\bigcap_{\alpha,p,k} \DD_{p,k}^\alpha\,.
$$

\end{theorem}
\proof
1)  We have
\begin{eqnarray*}
\|\varphi\|_{\DD_{p,k}^\alpha} & = &
    \|\sum_n(1+n)^{k/2}I_n((A^{\alpha/2})^{\otimes n}\varphi_n)\|_{L^p}\\
& = & \Big\| \sum(1+n)^{k/2}e^{nt}e^{-nt}I_n
  ((A^{\alpha/2})^{\otimes n}\varphi_n)\Big\|_{L^p}\,.
\end{eqnarray*}
From the hypercontractivity of $P_t$, we can choose $t$ such that
$p=e^{2t}+1$ then
$$
\Big\|\sum(1+n)^{k/2}e^{nt}e^{-nt} I_n(\ldots)\Big\|_p
   \leq \Big\|\sum(1+n)^{k/2}e^{nt}I_n(\ldots)\Big\|_2\,.
$$
Choose $\beta>0$ such that  $\|A^{-\beta}\|\leq e^{-t}$, hence
\begin{eqnarray*}
\lefteqn{\Big\|\sum(1+n)^{k/2}e^{nt}I_n(\ldots)\Big\|_2}\\
 & \leq&
   \Big\|\sum(1+n)^{k/2}\Gamma(A^\beta)\Gamma(A^{-\beta})
        e^{nt}I_n((A^{\alpha/2})^{\otimes n}\varphi_n)\Big\|_2 \\
& \leq & \sum(1+n)^{k/2}\|I_n((A^{\beta+\alpha/2})^{\otimes n}\varphi_n))\|_2 \\
& = & \|\varphi\|_{\DD_{2,k}^{2\beta+\alpha}}\,.
\end{eqnarray*}

\noindent
2) \  If we choose $\|A^{-\beta}\|<e^{-t}$ then the difference suffices to
absorb the action of the multiplicator $(1+n)^{k/2}$ which is of polynomial
growth and the former gives an exponential decrease.
\qed

\begin{corollary}
We have similar relations for  any separable Hilbert space valued
functionals.
\end{corollary}
\proof
This statement follows easily from  the Khintchine inequality.
\qed

As another corollary we have

\begin{corollary}
Let us denote by $\Phi(H_\infty)$ the space
$\bigcap_\alpha\Phi(H_\alpha)$. Then
\begin{enumerate}
\item   $\nabla:\Phi\to\Phi(H_\infty)$ and
$\delta:\Phi(H_\infty)\to\Phi$ are linear  continuous
operators. Consequently $\nabla$ and $\delta$ have continuous
extensions as linear operators $\nabla:\Phi'\to\Phi'(H_{-\infty})$ and
$\delta:\Phi'(H_{-\infty})\to\Phi'$.
\item  $\Phi$ is an algebra.
\item  For any $T\in\Phi'$, there exists some $\zeta\in\Phi'(H_{-\infty})$
 such that
$$
T=\langle T,1\rangle +\delta\zeta\,.
$$
\end{enumerate}
\end{corollary}
\proof
The first claim follows from Theorems  \ref{H1-th} and \ref{H2-th}.
To prove the second one it is sufficient to show that
$\varphi^2\in\Phi$ if $\varphi\in\Phi$. This follows from the
multiplication formula of the multiple Wiener integrals.
 (cf. Lemma \ref{mult-lemma}).
To prove the last one let us observe that  if $T\in\Phi'$, then there
exists some $\alpha>0$ such that
$T\in \DD_{2,0}^{-\alpha}$, i.e., $T$ under the isomorphism of Theorem
\ref{H1-th} is in $L^2(\mu_\alpha,W^\alpha)$ on which we have Ito
representation (cf.\ Appendix to the Chapter IV).
\qed

\begin{proposition}
Suppose that $A^{-1}$ is $p$-nuclear, i.e., there exists some $p\geq1$
such that  $A^{-p}$ is nuclear. Then $\Phi$ is a nuclear Fr\'echet
space.
\end{proposition}
\proof
This goes as in the classical  white noise case, except that the
eigenvectors of $\Gamma(A^{-1})$ are of the form
$H_{\vec{\alpha}}(\delta h_{\alpha_1},\ldots,\delta h_n)$ with
$h_{\alpha_i}$ are the eigenvectors of $A$.
\qed

\section{Applications to positive distributions}

Let $T\in\Phi'$ be a positive distribution. Then, from the
construction of the distribution spaces,  there exists some
$\DD_{p,-k}^{-\alpha}$ such that $T\in \DD_{p,-k}^{-\alpha}$ and $\langle
T,\varphi\rangle\geq0$ for any $\varphi\in \DD_{q,k}^\alpha$, $\varphi\geq0$.
Hence $i_\alpha(T)$ is a positive functional on $\DD_{1,k}^{(\alpha)}$
which is the Sobolev space on $W^\alpha$. Therefore $ i_\alpha(T)$ is
a Radon measure on $W^{-\alpha}$ and  we find  in fact that  the
support of $T$ is $W^{-\alpha}$ which is much smaller than
$H_{-\infty}$.
Let us give an example of such a positive distribution:

\begin{proposition}
\label{d-girsanov}
Assume that $u\in L^2(\mu,H)$ such that
\begin{equation}
\label{fin-sum}
\sum_{n=0}^\infty
\frac{1}{\sqrt{n!}}E\left[|u|_H^n\right]^{1/2}<\infty\,.
\end{equation}
Then the  mapping defined by
$$
\phi\to E[\phi(w+u(w))]=<L_u,\phi>
$$
is a positive distribution and it can be expressed as
$$
\sum_{n=0}^\infty\frac{1}{n!} \delta^n u^{\otimes n}\,.
$$
Moreover this sum is weakly uniformly convergent in
$\DD_{2,0}^{-\alpha}$, for any $\alpha>0$ such that
$\|A^{-1}\|^{2\alpha}<\frac{1}{2}$.
\end{proposition}
\proof
It follows trivially from the Taylor formula and from the definition
of $\delta^n$ as the adjoint of $\nabla^n$ with respect to $\mu$,  that
$$
<L_u,\phi>=\sum_{n=0}^\infty \frac{1}{n!}E\left[\phi\,\delta^n
  u^{\otimes n}\right]
$$
for any cylindrical, analytic function $\phi$. To complete the proof
it suffices to show that $\phi\to <L_u,\phi>$ extends continuously to
$\Phi$. If $\phi $ has the chaos decomposition
$$
\phi=\sum_{k=0}^\infty I_k(\phi_k)\,,
$$
with $\phi_k\in H_\infty^{\circ k}$, then
\beaa
E\left[\|\nabla^n\phi\|_{H^{\circ k}}^2\right]&=&\sum_{k\geq
  n}\frac{\left(k!\right)^2}{(k-n)!} \left\|\phi_k\right\|_{H^{\circ k}}^2\\
&\leq& \sum_{k\geq n}c^{-\alpha k}\frac{(k!)^2}{(k-n)!}
\|\phi_k\|_{H_\alpha^{\circ  k}}^2\,,
\eeaa
where $c^{-\alpha}$ is an upper bound for the norm of
$A^{-\alpha/2}$. Hence we the following a priori bound:
\beaa
|<L_u,\phi>|&\leq&\sum_{n=0}^\infty\frac{1}{n!}|<\nabla^n\phi,u^{\otimes
  n}>|\\
&\leq&\sum_{n=0}^\infty
\frac{1}{\sqrt{n!}}E\left[|u|_H^n\right]^{1/2}\left(\sum_{k\geq n}c^{-\alpha k}\frac{(k!)^2}{(k-n)!}
\|\phi_k\|_{H_\alpha^{\circ  k}}^2\right)^{1/2}\\
&\leq&\sum_{n=0}^\infty
\frac{1}{\sqrt{n!}}E\left[|u|_H^n\right]^{1/2}\left(\sum_{k\geq n}
(2c^{-\alpha})^k
\|\phi_k\|_{H_\alpha^{\circ  k}}^2\right)^{1/2}\,.
\eeaa
Choose now $\alpha$ such that $c^{-\alpha}<1/2$, then the sum inside
the square root is dominated by
$$
\sum_{k=0}^\infty k!\|\phi_k\|_{H_\alpha^{\circ
    k}}^2=\|\phi\|_{\DD^\alpha_{2,0}}^2\,.
$$
Hence the sum is absolutely convergent provided that $u$ satisfies the
condition (\ref{fin-sum}).
\qed
\section{Exercises}
{\footnotesize{
\begin{enumerate}
\item Let $K$ be a closed vector subspace of $H$ and denote by $P$ the
  orthogonal projection associated to it. Denote by $\calF_K$ the
  sigma algebra generated by $\{\delta k,\,k\in K\}$. Prove that
$$
\Gamma(P)f=E[f|\calF_K]\,,
$$
for any $f\in L^2(\mu)$.
\item Assume that $M$ and $N$ are two closed  vector subspaces of the 
Cameron-Martin space $H$, denote by $P$ and $Q$ respectively the 
corresponding orthogonal projections. For any $f,g\in L^2(\mu)$ prove the 
following inequality:
$$
\left|E\left[(f-E[f])(g-E[g])\right]\right|\leq 
\|PQ\|\|f\|_{L^2(\mu)}\,\|g\|_{L^2(\mu)}\,,
$$
where $\|PQ\|$ is the operator norm of $PQ$.

\item Prove that $\Gamma(e^{-t}I_H)=P_t$, $t\geq 0$, where $P_t$ denotes
the Ornstein-Uhlenbeck semi-group.
\item Let $B$ is a bounded operator on $H$, define $d\Gamma(B)$ as
$$
d\Gamma(B)f=\frac{d}{dt}\Gamma(e^{tB})f|_{t=0}\,.
$$
Prove that
$$
d\Gamma(B)f=\delta\{ B\nabla f\}
$$
and that
$$
d\Gamma(B)(fg)=f\,d\Gamma(B) g+g\,d\Gamma(B) f\,,
$$
(i.e., $d\Gamma(B)$ is a derivation)
for any $f,\,g\in \DD$ whenever $B$ is {\bf{skew-symmetric}}.
\end{enumerate}
}}


\section*{Notes and suggested reading}
{\footnotesize{
The fact that a positive Meyer distribution defines a Radon measure on
the Wiener space has been indicated for the first time in
\cite{Mal}. The notion of the capacity in an abstract frame has been
studied by several people, cf. in particular \cite{B-H},
\cite{Mall-1} and the references there. Application  to the local times
is original, the capacity version of $0-1$--law is taken from
\cite{Kus1}. Proposition \ref{d-girsanov} is taken from \cite{KU}, for
the more general distribution spaces we refer the
reader to \cite{KU, K-T, M-Y} and to the references there.
}}


\chapter{Characterization of independence of some Wiener functionals}
\label{ch.indep}
\markboth{Independence}{}

\section*{Introduction \index[sub]{independence}}

In probability theory, probably the   most important concept is the
independence since it is the basic property which differentiates the
probability theory from the abstract measure theory or from the
functional analysis. Besides it is almost always difficult to verify
 the independence of random variables. In fact, even in the
elementary probability, the tests required to verify the independence
of three or more random
variables get very quickly quite cumbersome. Hence it is very tempting
to try to  characterize the independence of random variables via the local
 operators as $\nabla$ or $\delta$ that we have studied  in the
 preceding chapters.

Let us begin with two random variables: let $F,G\in \DD_{p,1}$ for some $p>1$.
They are independent if and only if
$$
E[e^{i\alpha F}e^{i\beta G}]=E[e^{i\alpha F}]E[e^{i\beta G}]
$$
for any $\alpha,\beta\in\R$, which is equivalent to
$$
E[a(F)b(G)]=E[a(F)]E[b(G)]
$$
for any $a,b\in C_b(\R)$.

Let us denote by $\tilde{a}(F)=a(F)-E[a(F)]$, then we have:\\
 $F$ and $G$ are independent if and only if
$$
E[\tilde{a}(F)\cdot b(G)]=0\,,\qquad \forall a,b\in C_b(\R)\,.
$$

\noindent
Since $e^{i\alpha x}$ can be approximated point-wise with smooth functions, we
can suppose as well that $a,b\in C_b^1(\R)$ (or $C_0^\infty(\R)$). Since
$\mathcal L$ is invertible on the centered random variables, we have
\begin{eqnarray*}
E[\tilde{a}(F)b(G)] & = & E[{\mathcal {L L}}^{-1}\tilde{a}(F)\cdot b(G)] \\
& = & E[\delta\nabla{\mathcal L}^{-1}\tilde{a}(F)\cdot b(G)] \\
& = & E[(\nabla{\mathcal L}^{-1}\tilde{a}(F),\nabla(b(G)))_H] \\
& = & E[((I+{\mathcal L})^{-1}\nabla a(F),\nabla(b(G)))] \\
& = & E[((I+{\mathcal L})^{-1}(a'(F)\nabla F),b'(G)\nabla G)_H] \\
& = & E[b'(G)\cdot((I+{\mathcal L})^{-1}(a'(F)\nabla F),\nabla G)_H] \\
& = & E[b'(G)\cdot
          E[((I+{\mathcal L})^{-1}(a'(F)\nabla F,\nabla G)_H|\sigma(G)]]\,.
\end{eqnarray*}

\noindent
In particular choosing  $a=e^{i\alpha x}$, we find that

\begin{proposition}
\label{ind-char}
$F$ and $G$ (in $\DD_{p,1})$ are independent if and only if
$$
E\left[((I+{\mathcal L})^{-1}(e^{i\alpha F}
      \nabla F),\nabla G)_H\Big|\sigma(G)\right]=0\;\;{\mbox{a.s.}}
$$
\end{proposition}
\section{The case of multiple Wiener integrals}
Proposition \ref{ind-char}  is not very useful,  because of the non-localness
property of the operator ${\mathcal L}^{-1}$. Let us however look at the case
of multiple Wiener integrals:

First recall the following multiplication formula
\index[sub]{multiplication formula } of the multiple Wiener
integrals:

\begin{lemma}
\label{mult-lemma}
Let $f\in\hat{L}^2([0,1]^p)$, $g\in\hat{L}^2([0,1]^q)$. Then we have
\beaa
I_p(f)\,\, I_q(g)&=&\sum_{m=0}^{p\wedge q} \frac{p!\,q!}{m!(p-m)!(q-m)!}
I_{p+q-2m}(f\otimes_m g)\\
&=&\sum_{m=0}^{p\wedge q} \frac{p!\,q!}{m!(p-m)!(q-m)!}
I_{p+q-2m}(f\hat{\otimes}_m g)\,,
\eeaa
where $f\otimes_m g$\index[not]{f@$f\otimes_m g$} denotes the
contraction \index[sub]{contraction} of  order $m$ of the tensor 
$f\otimes g$, i.e., the partial scalar product of $f$ and $g$ in
$L^2([0,1]^m)$ and $f\hat{\otimes}_m g$ is its symmetrization.
\end{lemma}

\noindent
To prove above lemma we need to prove the Leibniz formula whose proof
follows from its finite dimensional version:
\begin{lemma}
\label{Leibniz}
Assume that $F, \,G$ are in $\DD$, then, for any $n\in \NN$,
we have 
$$
\nabla^n(F\,G)=\sum_{i=0}^n \frac{n!}{i!\,(n-i)!}\nabla^i
F\hat{\otimes}\nabla^{n-i}G\,.\,
$$
almost surely.
\end{lemma}

\noindent
{\bf{Proof of Lemma \ref{mult-lemma}:}}
Suppose that $p>q$ and let $\phi\in \DD$, using the identity
$\delta^pf=I_p(f)$ and the fact that $\delta^p$ is the adjoint of the
operator $\nabla^p$, we get, from Lemma \ref{Leibniz}
\beaa
E[I_p(f)I_q(g)\phi]&=&E[(f,\nabla^p(I_q(g)\phi)) _{H^{\hat{\otimes}p}}]\\
&=&E\left[\sum_{i=0}^p C_{p,i}(f,\nabla^i
I_q(g)\hat{\otimes}\nabla^{p-i}\phi) _{H^{\hat{\otimes}p}}\right]\\
&=&E\left[\sum_{i=0}^p C_{p,i} \frac{q!}{(q-i)!}                
 (f,I_{q-i}(g)\otimes\nabla^{p-i}\phi) _{H^{\otimes p}}\right]\\
&=&\sum_{i=0}^p
C_{p,i}\frac{q!}{(q-i)!} \,E\left[(f,\,I_{q-i}(g)\otimes\nabla^{p-i}\phi)_{H^{\otimes p}}\right]\\
&=&\sum_{i=0}^p
C_{p,i}\frac{q!}{(q-i)!} \,E\left[(I_{q-i}(g)\otimes_i\,f,\nabla^{p-i}\phi)_{H^{\otimes(p-i)}}\right]\\
&=&\sum_{i=0}^p
C_{p,i}\frac{q!}{(q-i)!} \,E[(g\otimes_i\,f,\nabla^{q-i}\nabla^{p-i}\phi)_{H^{\otimes(p+q-2i)}}]\\
&=&\sum_{i=0}^p
C_{p,i}\frac{q!}{(q-i)!} \,E[(g\otimes_i\,f,\nabla^{p+q-2i}\phi)_{H^{\otimes(p+q-2i)}}]\\
&=&\sum_{i=0}^p
C_{p,i}\frac{q!}{(q-i)!} \,E[I_{p+q-2i}(g\hat{\otimes}_i\,f),\phi]\,,
\eeaa
where $C_{p,i}=p!/i!(p-i)!$ and the proof of the lemma follows.
\qed

By the help of this lemma we will prove:

\begin{theorem}
\label{ind-thm}
$I_p(f)$ and $I_q(g)$ are independent if and only if
$$
f\otimes_1 g=0\quad {\mbox{a.s.~on}}\; [0,1]^{p+q-2}\,.
$$
\end{theorem}
\proof
$(\Rightarrow):$ By independence, we have
$$
E[I_p^2I_q^2]=p!\|f\|^2q!\|g||^2=p!q!||f\otimes g\|^2\,.
$$
On the other hand
$$
I_p(f)I_q(g)=\sum_0^{p\wedge q}m!C_p^mC_q^m\;I_{p+q-2m}(f\otimes_m g)\,,
$$
hence
\begin{eqnarray*}
\lefteqn{E[(I_p(f)I_q(g))^2]} \\
& =&    \sum_0^{p\wedge q}(m!C_p^mC_q^m)^2(p+q-2m)!\|f\hat{\otimes}_m g\|^2 \\
& \geq & (p+q)!\|f\hat{\otimes}g||^2\quad
\mbox{(dropping the terms with $m\geq 1 $)}\,.
\end{eqnarray*}
We have, by definition:
\begin{eqnarray*}
\|f\hat{\otimes}g\|^2 & = & \Big\| \frac{1}{(p+q)!} \sum_{\sigma\in S_{p+q}}
      f(t_{\sigma(1)},\ldots,t_{\sigma(p)})
      g(t_{\sigma(p+1)},\ldots,t_{\sigma(p+q)})\Big\|^2 \\
& = & \frac{1}{((p+q)!)^2}\sum_{\sigma,\pi\in S_{p+q}}\lambda_{\sigma,\pi}\,,
\end{eqnarray*}
where $S_{p+q}$ denotes the group of permutations of order $p+q$ and
\begin{eqnarray*}
\lambda_{\sigma,\pi} & = & \int_{[0,1]^{p+q}}
   f(t_{\sigma(1)},\ldots,t_{\sigma(p)})
     g(t_{\sigma(p+1)},\ldots,t_{\sigma(p+q)})\cdot \\
& & f(t_{\pi(1)},\ldots,t_{\pi(p)})
g(t_{\pi(p+1)},\ldots,t_{\pi(p+q)})dt_1\ldots dt_{p+q}\,.
\end{eqnarray*}

Without loss of generality, we may suppose that $p\leq q$. Suppose now that
$(\sigma(1),\ldots,\sigma(p))$ and $(\pi(1),\ldots,\pi(p))$ has $k\geq0$
elements in common. If we use the block notations, then
\begin{eqnarray*}
& & (t_{\sigma(1)},\ldots,t_{\sigma(p)})=(A_k,\tilde{A}) \\
& & (t_{\sigma(p+1)},\ldots,t_{\sigma(p+q)})=B \\
& & (t_{\pi(1)},\ldots,t_{\pi(p)})=(A_k,\tilde{C}) \\
& & (t_{\pi(p+1)},\ldots,t_{\pi(p+q)})=D
\end{eqnarray*}
where $A_k$ is the sub-block containing elements common to
$(t_{\pi(1)},\ldots,t_{\pi(p)})$ and $(t_{\sigma(1)},\ldots,t_{\sigma(p)})$.
Then we have
$$\lambda_{\sigma,\pi}=\int_{[0,1]^{p+q}}f(A_k,\tilde{A})g(B)\cdot
f(A_k,\tilde{C})g(D)dt_1\ldots dt_{p+q}\,.$$

\noindent
Note that $A_k\cup\tilde{A}\cup B=A_k\cup\tilde{C}\cup
D=\{t_1,\ldots,t_{p+q}\}$, $\tilde{A}\cap\tilde{C}=\emptyset$. Hence we have
$\tilde{A}\cup B=\tilde{C}\cup D$. Since $\tilde{A}\cap\tilde{C}=\emptyset$, we
have $\tilde{C}\subset B$ and $\tilde{A}\subset D$. From the fact that
$(\tilde{A},B)$ and $(\tilde{C},D)$ are the partitions of the same set, we have
$D\backslash\tilde{A}=B\backslash\tilde{C}$. Hence we can write, with the
obvious notations:
\begin{eqnarray*}
\lefteqn{\lambda_{\sigma,\pi}  = }\\
& = &\int_{[0,1]^{p+q}}f(A_k,\tilde{A})g(\tilde{C},B\backslash\tilde{C})\cdot
    f(A_k,\tilde{C})g(\tilde{A},D\backslash\tilde{A})dt_1\ldots dt_{p+q} \\
& = & \int_{[0,1]^{p+q}}f(A_k,\tilde{A})g(\tilde{C},B\backslash
   \tilde{C})f(A_k,\tilde{C})g(\tilde{A},B\backslash \tilde{C})dA_k
d\tilde{A} d\tilde{C}d(B\backslash\tilde{C})\\
& = & \int_{[0,1]^{q-p+2k}}(f\otimes_{p-k}g) (A_k,B\backslash\tilde{C})
(f\otimes_{p-k}g)(A_k,B\backslash\tilde{C})\cdot dA_kd(B\backslash\tilde{C}) \\
& = & \|f\otimes_{p-k}g\|^2_{L^2([0,1]^{q-p+2k})}\,
\end{eqnarray*}
where we have used the relation $D\backslash\tilde{A}=B\backslash\tilde{C}$
in the second line of the above equalities.
\noindent
Note that for $k=p$ we have $\lambda_{\sigma,\pi}=\|f\otimes g\|_{L^2}^2$.
Hence we have
\begin{eqnarray*}
\lefteqn{E[I_p^2(f)I_q^2(g)]= p!\|f\|^2\cdot q!\|g\|^2} \\
&& \geq (p+q)! \left[
\frac{1}{((p+q)!)^2}\Big[ \sum_{\sigma,\pi}\lambda_{\sigma,\pi}(k=p)+
\sum_{\sigma,\pi}\lambda_{\sigma,\pi}(k\not= p))\Big]\right].
\end{eqnarray*}
The number  of $\lambda_{\sigma,\pi}$ with $(k=p)$ is exactly
${p+q\choose p}(p!)^2(q!)^2$, hence we have
$$
p!q!\|f\|^2\|g\|^2\geq p!q!\|f\otimes g\|^2+ \sum_{k=0}^{p-1}
c_k\|f\otimes_{p-k}g\|^2_{L^2([0,1]^{q-p+2k})}
$$
with $c_k>0$. For this relation to hold we should have
$$
\|f\otimes_{p-k}g\|=0\,,\qquad k=0,\ldots,p-1
$$
in particular for $k=p-1$, we have
$$
\|f\otimes_1g\|=0\,.
$$
$(\Leftarrow)$: \ From the Proposition \ref{ind-char}, we see that it
is sufficient to prove
$$
((I+{\mathcal L})^{-1}e^{i\alpha F}\nabla F,\nabla I_q(g))=0\;\;{\mbox{a.s.}}
$$
with $F=I_p(f)$, under the hypothesis $f\otimes_1g=0$ a.s.  Let us write
$$
e^{i\alpha I_p(f)}=\sum_{k=0}^\infty I_k(h_k)\,,
$$
then
\begin{eqnarray*}
e^{i\alpha I_p(f)}\nabla I_p(f) & = & p\sum_{k=0}^\infty I_k(h_k)\cdot I_{p-1}(f) \\
& = & p\sum_{k=0}^\infty\sum_{r=0}^{k\wedge(p-1)} \alpha_{p,k,r}
I_{p-1+k-2r}(h_k\otimes_r f)\,.
\end{eqnarray*}
Hence
$$
(I+{\mathcal L})^{-1}e^{i\alpha F}\nabla F=
p\sum_k\!\!\sum_{r=0}^{k\wedge(p-1)}(1+p+k-1-2r)^{-1}
       I_{p-1+k-2r}(h_k\otimes_p f)\,.
$$

\noindent
When we take the scalar product with $\nabla I_q(g)$, we will have terms of the
type:
\begin{eqnarray*}
\lefteqn{\left(I_{p-1+k-2r}(h_k\otimes_r f),I_{q-1}(g)\right)_H=}\\
& & =\sum_{i=1}^\infty I_{p-1+k-2r}(h_k\otimes_r f(e_i))I_{q-1}(g(e_i))\,.
\end{eqnarray*}

\noindent
If we use the multiplication formula to calculate each term, we find the terms
as
\begin{eqnarray*}
\lefteqn{\sum_{i=1}^\infty\int(h_k\otimes_r
f(e_i))(t_1,\ldots,t_{p+k-2r-1})g(e_i)(t_1,\ldots,t_{q-1})dt_1dt_2\ldots} \\
& & =\int\!\!\int_{\theta=0}^1(h_k\otimes_r f(\theta))(t_1,\ldots,
t_{p+k-2r-1})g(\theta,t_1,\ldots,t_{q-1})d\theta\,dt_1\ldots
\end{eqnarray*}
From the hypothesis we have
$$
\int_0^1 f(\theta ,t_1\ldots)g(\theta ,s_1\ldots ,)d\theta =0\quad
{\mbox{a.s.}}\,,
$$
hence the Fubini theorem completes the proof.
\qed

\remark
For a more elementary proof of the sufficiency of Theorem
\ref{ind-thm} cf. \cite{OK}.
\begin{remarkk}
\label{ind-remark}
{\rm 
In the proof of the necessity we have used only the fact that
$I_p(f)^2$ and $I_q(g)^2$ are independent. Hence, as a byproduct we
obtain also the fact that $I_p$ and $I_q$ are independent if and only
if their squares are independent.}
\end{remarkk}

\begin{corollary}
Let $f$ and $g$ be symmetric $L^2$-kernels respectively on $[0,1]^p$ and
$[0,1]^q$. Let
$$
S_f={\mbox{\rm span}}\{f\otimes_{p-1}h:h\in L^2([0,1])^{p-1}\}
$$
and
$$
S_g={\mbox{\rm span}}\{g\otimes_{q-1}k;k\in L^2(]0,1]^{q-1})\}\,.
$$
Then the following are equivalent:
\begin{itemize}
\item[i)]
$I_p(f)$ and $I_q(g)$ are independent,
\item[ii)] $I_p(f)^2$ and $I_q(g)^2$ are independent, 
\item[iii)]
$S_f$ and $ S_g$ are orthogonal  in $H$,
\item[iv)]
the Gaussian-generated $\sigma$-fields $\sigma\{ I_1(k);k\in S_f\}$ and
$\sigma\{I_1(l);l\in S_g\}$ are independent.
\end{itemize}
\end{corollary}
\proof
As it is indicated in  Remark \ref{ind-remark}, the independence of
$I_p$ and $I_q$ is equivalent to the independence of their squares.\\ 
(i$\Rightarrow$iii): The hypothesis  implies that$f\otimes_1g=0$ a.s. If $a\in
S_f$, $b\in S_g$ then they can be written as finite linear combinations of the
vectors  $f\otimes_{p-1}h$ and $g\otimes_{q-1}k$ respectively. Hence,
it suffices to assume, by linearity, that $a=f\otimes_{p-1}h$ and
$b=g\otimes_{q-1}k$. Then it follows from the Fubini theorem
\begin{eqnarray*}
(a,b)=(f\otimes_{p-1}h,g\otimes_{q-1}k) & = &
(f\otimes_1g,h\otimes k)_{(L^2)^{\otimes p+q-2}} \\
& = & 0\,.
\end{eqnarray*}
(iii$\Rightarrow$i) \ If $(f\otimes_1 g,h\otimes k)=0$ for all $ h\in
L^2([0,1]^{p-1}), k\in L^2([0,1]^{q-1})$, then $f\otimes_1g=0$ a.s.\ since
finite combinations of $h\otimes k$ are dense in $L^2([0,1]^{p+q-2})$.
Finally, the equivalence of (iii) and (iv) is obvious.
\qed

\begin{proposition}
Suppose that $I_p(f)$ is independent of $I_q(g)$ and $I_p(f)$ is independent of
$I_r(h)$. Then $I_p(f)$ is independent of $\{I_q(g),I_r(h)\}$.
\end{proposition}
\proof
We have $f\otimes_1g=f\otimes_1h=0$ a.s. This implies the independence of
$I_p(f)$ and $\{I_g(g),I_r(h)\}$ from the  calculations similar to those  of
the proof of sufficiency of the theorem.
\qed

\noindent
In a similar way we have
\begin{proposition}
Let $\{I_{p_\alpha}(f_\alpha);\alpha\in J\}$ and $I_{q_\beta}(g_\beta);\beta\in
K\}$ be two arbitrary families of multiple Wiener integrals. The two families
are independent if and
only if $I_{p_\alpha}(f_\alpha)$ is independent of $I_{q_\beta}(g_\beta)$
for all $(\alpha,\beta)\in J\times K$.
\end{proposition}
\begin{corollary}
If $I_p(f)$ and $I_q(g)$ are independent, so are  also $I_p(f)(w+\tilde{h})$
and
$I_q(g)(w+\tilde{k})$ for any $\tilde{h},\tilde{k}\in H$.
\end{corollary}
\proof
Let us denote, respectively, by $h$ and $k$ the Lebesgue densities of
$\tilde{h}$ and $\tilde{k}$. We have then
$$
I_p(f)(w+\tilde{h})=\sum_{i=0}^p{p\choose i}(I_{p-i}(f),h^{\otimes i})
_{H^{\otimes i}}\,.
$$
Let us define  $f[h^{\otimes i}]\in L^2[0,1]^{p-i}$ by
$$
I_{p-i}(f[h^{\otimes i}])=(I_{p-i}(f),h^{\otimes i}).
$$
If $f\otimes_1g=0$ then it is easy to see that
$$
f[h^{\otimes i}]\otimes_1g[k^{\otimes j}]=0\,,
$$
hence the corollary follows from Theorem \ref{ind-thm}.

\qed

\noindent
From the corollary it follows

\begin{corollary}
 $I_p(f)$ and $I_q(g)$ are independent if and only if  the germ $\sigma$-fields
$$
\sigma\{I_p(f),\nabla I_p(f),\ldots,\nabla^{p-1}I_p(f)\}
$$
and
$$
\sigma\{I_q(g),\ldots,\nabla^{q-1}I_q(g)\}
$$
are independent.
\end{corollary}
\begin{corollary}
Let $X,Y\in L^2(\mu)$, $Y=\sum_0^\infty I_n(g_n)$. If
$$
\nabla X\otimes_1g_n=0\quad{\mbox{a.s.}}\; \forall n\,,
$$
then $X$ and $Y$ are independent.
\end{corollary}
\proof
This follows from Proposition \ref{ind-char}.
\qed

\begin{corollary}
In particular, if $\tilde{h}\in H$, then $\nabla_{\tilde{h}}\varphi=0$
a.s.\ implies that $\varphi$ and $I_1(h)=\delta\tilde{h}$  are independent.
\end{corollary}
\section{Exercises}{\footnotesize{
\begin{enumerate}
\item Let$f\in {\hat{L}}^2([0,1]^p)$ and $h\in L^2([0,1])$. Prove the
  product formula
\begin{equation}
\label{prod-1-form}
I_p(f)I_1(h)=I_{p+1}(f\otimes h)+pI_{p-1}(f\otimes_1h)\,.
\end{equation}
\item Prove  by induction and with the help of (\ref{prod-1-form}),
  the general multiplication  formula
$$
I_p(f)\,\, I_q(g)=\sum_{m=0}^{p\wedge q} \frac{p!\,q!}{m!(p-m)!(q-m)!}
I_{p+q-2m}(f\otimes_m g)\,, 
$$
where $f\in {\hat{L}}^2([0,1]^p)$ and $g\in {\hat{L}}^2([0,1]^q)$.
\end{enumerate}

}}
\section*{Notes and suggested reading}
{\footnotesize{All the results of this chapter are taken from
    \cite{UZ1,UZ2}, cf. also \cite{OK} for some simplification of the
    sufficiency of Theorem \ref{ind-thm}. Note that, in Theorem
    \ref{ind-thm}, we have used only the independence of $I_p(f)^2$
    and $I_q(g)^2$. Hence two multiple Ito-Wiener integrals are
    independent if and only if their squares are independent. The
    important Lemma \ref{mult-lemma} is proven by Shigekawa,
    cf. \cite{Sh} using the induction, the proof that we give here is
    totally original and it is more in harmony with the spirit of
    Malliavin Calculus.

}}


\chapter{Moment inequalities for Wiener functionals}
\label{ch.ineq}
\markboth{Moment Inequalities}{}

\section*{Introduction}
In several applications, as limit theorems, large deviations, degree
theory of Wiener maps, calculation of the Radon-Nikodym densities,
etc., it is
important to control the (exponential) moments of Wiener functionals by those
of their derivatives. In this chapter we will give two  results on this
subject. The first one concerns the tail probabilities of the Wiener
functionals with essentially bounded Gross-Sobolev derivatives. This
result is a straightforward generalization of the celebrated
Fernique's lemma which says that the square of the  supremum of the
Brownian path on any bounded interval has an exponential moment
provided that it is multiplied with a sufficiently small, positive
constant. The second inequality says that for a Wiener functional
$F\in \DD_{p,1}$, we have
\begin{equation}
E_w\times E_z[U(F(w)-F(z))]\leq E_w\times
              E_z\left[U\left(\frac{\pi}{2}I_1(\nabla F(w)\right)(z)\right],
\end{equation}
where $w$ and $z$ represent two independent Wiener paths, $E_w$ and
$E_z$ are the corresponding expectations, and $I_1(\nabla F(w))(z)$ is
the first order Wiener integral with respect to $z$ of $\nabla F(w)$
and $U$ is any lower bounded, convex function on $\R$.
 Then combining  these two inequalities  we will obtain some interesting
majorations.

In the next section we show that the log-Sobolev inequality implies
the exponential integrability of the square of the Wiener functionals
whose derivatives are essentially bounded. In this section we study
with general measures which satisfy a logarithmic Sobolev inequality.

The next  inequality is an interpolation inequality which says that the
Sobolev norm of first order can be upper bounded  by the product  of
the  second order and of the  zero-th order Sobolev norms.

In the last part we study the exponential integrability of the Wiener
functionals in the divergence form, a problem which has gained
considerable importance due to the degree theorem on the Wiener space
as it is explained in more detail in the notes at the end of this chapter.

\section{Exponential tightness \index[sub]{exponential tightness}}
\markboth{Moment Inequalities}{Exponential Tightness}
First we will show the following result which is a consequence of the
Doob inequality:
\begin{theorem}
\label{exp-int-thm}
Let $\varphi\in \DD_{p,1}$ for some $p>1$. Suppose that $\nabla\varphi\in
L^\infty(\mu,H)$. Then we have
$$
\mu\{|\varphi|>c\}\leq2\exp\left\{-
\frac{(c-E[\varphi])^2}{2\|\nabla\varphi\|^2_{L^\infty(\mu,H)}}\right\}
$$
for any $c\geq 0$.
\end{theorem}
\proof
Suppose that $E[\varphi]=0$. Let $(e_i)\subset H$ be a complete, orthonormal
basis of $H$. Define
$V_n=\sigma\{\delta e_1,\ldots,\delta e_n\}$ and let
$\varphi_n=E[P_{1/n}\varphi|\nabla_n]$, where $P_t$ denotes the
Ornstein-Uhlenbeck semi-group on $W$.  Then, from Doob's Lemma,
$$
\varphi_n  =  f_n(\delta e_1,\ldots,\delta e_n).
$$
Note that, since
$f_n\in\bigcap_{p,k}W_{p,k}(\R^n,\mu_n)$, the Sobolev embedding
theorem implies that after a modification on a  set of null Lebesgue
measure, $ f_n$ can be chosen in $C^\infty(\R^n)$. Let
$(B_t;t\in[0,1])$ be an  $\R^n$-valued Brownian motion. Then
\begin{eqnarray*}
\mu\{|\varphi_n|>c\} & = & \P\{|f_n(B_1)|>c\} \\
& \leq & \P\{\sup_{t\in[0,1]}|E[f_n(B_1)|{\calB}_t]|>c\} \\
& = & \P\{\sup_{t\in[0,1]}|Q_{1-t}f_n(B_t)|>c\}\,,
\end{eqnarray*}
where $\P$ is the canonical Wiener measure on $C([0,1],\R^n)$  and  $Q_t$ is
the heat kernel associated to $(B_t)$, i.e.
$$
Q_t(x,A)=\P\{B_t+x\in A\}\,.
$$

\noindent
From the Ito formula, we have
$$
Q_{1-t}f_n(B_t)=Q_1f_n(B_0)+\int_0^t(DQ_{1-s}f_n(B_s),dB_s)\,.
$$
By definition
\begin{eqnarray*}
Q_1f_n(B_0) & = & Q_1f_n(0)=\int f_n(y)\cdot Q_1(0,dy) \\
& = & \int_{\R^n} f_n(y)e^{-\frac{1}{2}|y|^2}\frac{dy}{(2\pi)^{n/2}}\\
& = &E\left[E[P_{1/n}\varphi|V_n]\right]\\
&=&E\left[P_{1/n}\varphi\right]\\
&=&E[\varphi]\\
&=&0\,.
\end{eqnarray*}
Moreover we have $DQ_tf=Q_tDf$, hence
$$
Q_{1-t}f_n(B_t)=\int_0^t(Q_{1-s}Df_n(B_s),dB_s)=M_t^n\,.
$$
The Doob-Meyer process $(\langle M^n,M^n\rangle_t,t\in \reals_+)$ of
the martingale $M^n$ can be controlled as
\begin{eqnarray*}
\langle M^n,M^n\rangle_t & = & \int_0^t|DQ_{1-s}f_n(B_s)|^2ds \\
& \leq & \int_0^t\|Df_n\|_{C_b}^2 ds=t\|\nabla f_n\|_{C_b}^2\\
&=&t\|\nabla f_n\|_{L^\infty(\mu_n)} \\
&\leq& t\|\nabla\varphi\|_{L^\infty(\mu,H)}^2\,.
\end{eqnarray*}
Hence from the exponential Doob inequality, we obtain
$$
\P\left\{\sup_{t\in[0,1]}|Q_{1-t}f_n(B_t)|>c\right\}\leq 2\exp\left[-
\frac{c^2}{2\|\nabla\varphi\|_{L^\infty(\mu,H)}^2}\right]\,.
$$
Consequently
$$
\mu\{|\varphi_n|>c\}\leq2\exp
\left[-\frac{c^2}{2\|\nabla\varphi\|_{L^\infty(\mu,H)}^2}\right]\,.
$$
Since $\varphi_n\to\varphi$ in probability the proof is completed.
\qed

\begin{corollary}
\label{exp-cor}
Under the hypothesis of the theorem, for any
$$
\lambda<\left[2\|\nabla\varphi\|_{L^\infty(\mu,H)}\right]^{-1}\,,
$$
 we have
$$
E\left[\exp\lambda|\varphi|^2\right]<\infty.
$$
\end{corollary}
\proof
The first part follows from the fact that, for $F\geq 0$ a.s.,
$$
E[F]=\int_0^\infty P\{F>t\}dt\,.
$$
\qed

\remark In the next  sections  we will give more  precise estimate for
$E[\exp \lambda F^2]$.\\
\newline
In the applications, we encounter random variables $F$ satisfying
$$
|F(w+h)-F(w)|\leq c |h|_H,
$$
almost surely, for any $h$ in the Cameron-Martin space $H$ and a fixed
constant $c>0$, without any hypothesis of integrability. For example,
$F(w)=\sup_{t\in [0,1]}|w(t)|$, defined on $C_0[0,1]$ is such a
functional. In fact the above
hypothesis contains the integrability and Sobolev differentiability  of $F$.
We begin first by proving that under the integrability hypothesis,
such a functional is in the domain of $\nabla$:
\begin{lemma}
\label{lemma-bir}
Suppose that $F:W\mapsto \R$ is a measurable random variable in
$\cup_{p>1}L^p(\mu)$, satisfying
\begin{equation}
\label{diff-ineq}
|F(w+h)-F(w)|\leq c |h|_H,
\end{equation}
almost surely, for any $h\in H$, where $c>0$ is a fixed constant. Then $F$
belongs to $\DD_{p,1}$ for any $p>1$.
\end{lemma}

\remark If in (\ref{diff-ineq})  the negligeable set on which the
inequality is satisfied is independent of $h\in H$,  then the
functional $F$ is called  H-Lipschitz.
\index[sub]{H-Lipschitz}

\proof
Since, for some $p_0>1$, $F\in L^{p_0}$, the distributional derivative of $F$,
$\nabla F$ exists . We have $\nabla_k F\in \DD'$ for any $k\in
H$. Moreover, for $\phi \in \DD$, from the integration by parts
formula
\begin{eqnarray*}
E[\nabla_kF\,\phi]&=&-E[F\, \nabla_k \phi]+E[F \delta k\, \phi]\\
  &=&-\frac{d}{dt}|_{t=0}E[F\,\phi(w+tk)]+E[F \delta k\, \phi]\\
&=&-\frac{d}{dt}|_{t=0}E\left[F(w-tk)\,\phi\,\,\varepsilon(t\delta k)\right]+
   E\left[F \delta k\, \phi\right]\\
&=&\lim_{t\rightarrow 0}-E\left[\frac{F(w-tk)-F(w)}{t}\phi\right],
\end{eqnarray*}
where $\varepsilon(\delta k)$ denotes the Wick exponential of the
Gaussian random variable $\delta k$, i.e., 
$$
\varepsilon(\delta k)=\exp\left\{\delta k-\frac{1}{2}|k|^2\right\}\,.
$$
Consequently,
\begin{eqnarray*}
|E[\nabla_k F\, \phi]|&\leq&c|k|_H E[|\phi|]\\
                     &\leq& c|k|_H \|\phi\|_q,
\end{eqnarray*}
for any $q>1$, i.e., $\nabla F$ belongs to $L^p(\mu,H)$ for any $p>1$.
Let now $(e_i;i\in \N)$ be a complete, orthonormal basis of $H$, denote by
$V_n$ the sigma-field generated by $\delta e_1,\ldots, \delta e_n$, $n\in \N$
and let $\pi_n$ be the orthogonal projection onto the the subspace of $H$
spanned by $e_1,\ldots, e_n$. Let us define
$$
F_n=E[P_{1/n}F|V_n],
$$
where $P_{1/n}$ is the Ornstein-Uhlenbeck semi-group at the instant $t=1/n$.
Then $F_n\in\cap_{k} \DD_{p_0,k}$  and  it is immediate, from the martingale
convergence theorem and from the fact that $\pi_n$ tends to the identity
operator of $H$ pointwise, that
$$
\nabla F_n=E[e^{-1/n}\pi_n P_{1/n}\nabla F|V_n]\rightarrow \nabla F,
$$
in $L^p(\mu,H)$, for any $p>1$, as $n$ tends to infinity. Since, by
construction, $(F_n;n\in \N)$ converges also to $F$ in  $L^{p_0}(\mu)$,
$F$ belongs to $ \DD_{p_0,1}$. Hence we can apply the Corollary \ref{exp-cor}.

\qed

\begin{lemma}
\label{integ-lem}
Suppose that $F:W\mapsto \R$ is a measurable random variable satisfying
$$
|F(w+h)-F(w)|\leq c |h|_H,
$$
almost surely, for any $h\in H$, where $c>0$ is a fixed constant. Then $F$
belongs to $\DD_{p,1}$ for any $p>1$.
\end{lemma}
\proof
Let $ F_n=|F|\wedge n$, $n\in \N$. A simple calculation shows that
$$
|F_n(w+h)-F_n(w)|\leq c |h|_H,
$$
hence $F_n\in \DD_{p,1}$ for any $p>1$ and $|\nabla F_n|\leq c$ almost
surely from  Lemma \ref{lemma-bir}.
We have from the Ito-Clark formula (cf. Theorem \ref{ust-thm}),
$$
F_n=E[F_n]+\int_0^1 E[D_sF_n|{\mathcal F}_s] dW_s.
$$
From the definition of the stochastic integral, we have
\begin{eqnarray*}
E\left[\left(\int_0^1 E[D_sF_n|{\mathcal F}_s] dW_s\right)^2\,\right]&=&
       E\left[\int_0^1| E[D_sF_n|{\mathcal F}_s]|^2 ds\right]\\
&\leq & E\left[\int_0^1| D_sF_n|^2 ds\right]\\
&=&E[|\nabla F_n|^2]\\
&\leq & c^2\,.
\end{eqnarray*}
Since $F_n$ converges to $|F|$ in probability, and the stochastic
integral is bounded in $L^2(\mu)$, by taking the difference, we see
that $(E[F_n], n\in \N)$ is a sequence of (degenerate) random
variables bounded in the space of random variables under the topology
of convergence in probability, denoted by $L^0(\mu)$. Therefore
$\sup_n\mu\{E[F_n]>c\}\rightarrow 0$ as $c\rightarrow \infty$. Hence
$\lim_nE[F_n]=E[|F|]$ is finite. Now we apply the dominated
convergence theorem to obtain that $F\in L^2(\mu)$. Since the
distributional derivative of $F$
is a square integrable random variable, $F\in \DD_{2,1}$. We can now apply the
Lemma \ref{lemma-bir} which implies that $F\in \DD_{p,1}$ for any $p$.

\qed

\remark Although we have used the classical Wiener space structure in
the proof, the case of the Abstract Wiener space can be reduced to
this case using the method explained in the appendix of Chapter IV.

\begin{corollary}[Fernique's Lemma]
For any
$\lambda<\frac{1}{2}\,$, we have
$$
E[\exp\lambda\|w\|_W^2]<\infty,
$$
where $\|w\|$ is the norm of the Wiener path $w\in W$.
\end{corollary}
\proof
It suffices to remark that
$$
|\|w+h\|-\|w\||\leq |h|_H
$$
for any $h\in H$ and $w\in W$.
\qed

\section{Coupling  inequalities\index[sub]{coupling inequality}}
\markboth{Moment Inequalities}{Coupling Inequalities}
We begin with the following elementary  lemma (cf. \cite{P}):

\begin{lemma}
Let $X$ be a Gaussian random variable with values in  $\R^d$. Then for
any convex function $U$ on $\R$
and $C^1$-function $V:\R^d\to\R$, we have the following inequality:
$$
E[U(V(X)-V(Y))]\leq E\Big[U\Big(\frac{\pi}{2}(V'(X),Y)_{\R^d}\Big)\Big],
$$
where $Y$ is an independent copy of $X$ and $E$ is the expectation with
respect to  the product measure.
\end{lemma}
\proof
Let $X_\theta=X\sin\theta+Y\cos\theta$. Then
\begin{eqnarray*}
V(X)-V(Y) & = & \int_{[0,\pi/2]}\frac{d}{d\theta}V(X_\theta)d\theta \\
& = & \int_{[0,\pi/2]}(V'(X_\theta),X'_\theta)_{\R^d}\,d\theta \\
& = & \frac{\pi}{2}\int_{[0,\pi/2]}
     (V'(X_\theta),X'_\theta)_{\R^d}\, d\tilde{\theta}
\end{eqnarray*}

\noindent
where $d\tilde{\theta}=\frac{d\theta}{\pi/2}\,$. Since $U$ is convex, we have
$$
U(V(X)-V(Y))\leq\int_0^{\pi/2}U\Big(\frac{\pi}{2}(V'(X_\theta),
X'_\theta)\Big)d\tilde{\theta}\,.
$$

\noindent
Moreover $X_\theta$ and $X'_\theta$ are two independent Gaussian random
variables with the same law as the one of $X$. Hence
\begin{eqnarray*}
E[U(V(X)-V(Y))] & \leq & \int_0^{\pi/2} E\Big[U\Big(\frac{\pi}{2}
(V'(X),Y)\Big)\Big]d\tilde{\theta} \\
& = & E\Big[U\Big(\frac{\pi}{2}(V'(X),Y)\Big)\Big].
\end{eqnarray*}

\vspace*{-3ex}
\qed

Now we will extend this result to the Wiener space:
\begin{theorem}
Suppose that $\varphi\in \DD_{p,1}$, for some $p>1$ and $U$ is a lower bounded,
 convex function (hence lower semi-continuous) on $\R$. We have
$$
E[U(\varphi(w)-\varphi(z))]\leq E\Big[ U\Big(
\frac{\pi}{2}I_1(\nabla\varphi(w))(z)\Big)\Big]
$$
where $E$ is taken with respect to $\mu(dw)\times\mu(dz)$ on $W\times W$ and
on the classical Wiener space, we have
$$
I_1(\nabla\varphi(w))(z)=\int_0^1\frac{d}{dt}
\nabla\varphi(w,t)dz_t\,.
$$
\end{theorem}
\label{th2}
\proof
Suppose first that
$$
\varphi=f(\delta h_1(w),\ldots,\delta h_n(w))
$$
with $f$ smooth on $\R^n$, $h_i\in H$, $(h_i,h_j)=\delta_{ij}$. We have
\begin{eqnarray*}
I_1(\nabla\varphi(w))(z) & = & I_1 \Big( \sum_{i=1}^n \partial_if(\delta
h_1(w),\ldots,\delta h_n(w))h_i\Big) \\
& = & \sum_{i=1}^{n}\partial_if(\delta h_1(w),
       \ldots,\delta h_n(w))I_1(h_i)(z) \\
& = & (f'(X),Y)_{\R^n}
\end{eqnarray*}
where $X=(\delta h_1(w),\ldots,\delta h_n(w))$ and $Y=(\delta
h_1(z),\ldots,\delta h_n(z))$. Hence the inequality is trivially true in this
case.

For general $\varphi$, let $(h_i)$ be a complete, orthonormal basis  in $H$,
$$V_n=\sigma\{\delta h_1,\ldots,\delta h_n\}$$ and let
$$
\varphi_n=E[P_{1/n}\varphi|V_n]\,,
$$
where $P_{1/n}$ is the Ornstein-Uhlenbeck  semi-group on $W$.
We have then
$$
E[U(\varphi_n(w)-\varphi_n(z))]\leq E\Big[U\Big(\frac{\pi}{2}
I_1(\nabla\varphi_n(w))(z)\Big)\Big].
$$
Let $\pi_n$ be the orthogonal projection from $H$ onto span
$\{h_1,\ldots,h_n\}$. We have
\begin{eqnarray*}
I_1(\nabla\varphi_n(w))(z) & = &
I_1(\nabla_w E_w[P_{1/n}\varphi|V_n])(z) \\
& = & I_1(E_w[e^{-1/n}P_{1/n}\pi_n\nabla\varphi|V_n])(z) \\
& = & I_1(\pi_n E_w[e^{-1/n}P_{1/n}\nabla\varphi|V_n])(z) \\
& = & E_z[I_1^z(E_w[e^{-1/n}P_{1/n}^w\nabla\varphi|V_n])|\tilde{V}_n]
\end{eqnarray*}
where $\tilde{V}_n$ is the copy of $V_n$ on the second Wiener space. Then
\begin{eqnarray*}
\lefteqn{E\Big[U\Big(\frac{\pi}{2}I_1(\nabla\varphi_n(w))(z)\Big)\Big]} \\
& \leq &
       E\Big[U\Big(\frac{\pi}{2}
        I_1( E_w[e^{-1/n}P_{1/n}\nabla\varphi|V_n])(z)\Big)\Big] \\
& = &
E\Big[U\Big(\frac{\pi}{2}e^{-1/n} E_w[
        I_1(P_{1/n}\nabla\varphi(w))(z)|V_n]\Big)\Big] \\
& \leq & E\Big[U\Big(
     \frac{\pi}{2}e^{-1/n}I_1(P_{1/n}\nabla\varphi(w))(z)\Big)\Big] \\
& = & E\Big[U\Big(\frac{\pi}{2}e^{-1/n}P_{1/n}^w
I_1(\nabla\varphi(w))(z)\Big)\Big] \\
& \leq & E\Big[U\Big(\frac{\pi}{2}e^{-1/n}I_1(\nabla\varphi(w))(z)\Big)\Big]\\
& = & E\Big[U\Big(\frac{\pi}{2}
        P_{1/n}^{(z)}I_1(\nabla\varphi(w))(z)\Big)\Big] \\
& \leq & E\Big[U\Big(\frac{\pi}{2}I_1(\nabla\varphi(w))(z)\Big)\Big].
\end{eqnarray*}
Now Fatou's lemma completes the proof.
\qed

Let us give some consequences of this result:
\begin{theorem}
\label{th3}
The following Poincar\'e  inequalities are valid:
\begin{itemize}
\item[i)]
$E[\exp(\varphi-E[\varphi])]\leq E
\Big[\exp\left(\frac{\pi^2}{8}|\nabla\varphi|_H^2\right)\Big]$,
\item[ii)]
$E[|\varphi-E[\varphi]|]\leq \frac{\pi}{2}E[|\nabla\varphi|_H]$.
\item[iii)]
$E[|\varphi-E[\varphi]|^{2k}]\leq \Big(\frac{\pi}{2}\Big)^{2k}
\frac{(2k)!}{2^kk!}E[|\nabla\varphi|_H^{2k}]$, $k\in\N$.
\end{itemize}
\end{theorem}
\begin{remarkk}{\rm
Let us note that the result of (ii) can not be obtained with the classical
methods, such as the Ito-Clark representation theorem, since the optional
projection is not a continuous map in $L^1$-setting. Moreover, using the
H\"older inequality and the Stirling formula, we deduce the following set
of inequalities:
$$
\|\varphi-E[\varphi]\|_{p}\leq p\,\frac{\pi}{2} \|\nabla\varphi\|_{L^p(\mu,H)},
$$
for any $p\geq 1$ . To compare this result with those already known,
let us recall that using first the Ito-Clark formula, then the
Burkholder-Davis-Gundy inequality combined with the convexity
inequalities for the dual projections and some duality techniques, we
obtain, only for $p>1$ the  inequality 
$$
\|\varphi-E[\varphi]\|_{p}\leq K p^{3/2} \|\nabla\varphi\|_{L^p(\mu,H)},
$$
where $K$ is some positive constant.}
\end{remarkk}
\proof
 Replacing  the function  $U$ of Theorem \ref{th2}  by  the
 exponential function, we have
\begin{eqnarray*}
E[\exp(\varphi-E[\varphi])] & \leq &
 E_w\times E_z[\exp(\varphi(w)-\varphi(z))]\leq \\
& \leq & E_w\left[{
    E_z\left[{[\exp\frac{\pi}{2}I_1(\nabla\varphi(w))(z)}\right]}\right]\\ 
& = & E\left[\exp\frac{\pi^2}{8}|\nabla\varphi|_H^2\right].
\end{eqnarray*}

\noindent
 (ii) and (iii) are  similar provided that we take  $U(x)=|x|^{k}$, $k\in \N$.
\qed

\begin{theorem}
Let $\varphi\in \DD_{p,2}$ for some $p>1$ and that $\nabla|\nabla\varphi|_H\in
L^\infty(\mu,H)$. Then there exists some $\lambda>0$
such that
$$
E[\exp\lambda|\varphi|]<\infty\,.
$$
In particular, this hypothesis  is satisfied if
$\|\nabla^2\varphi\|_{\mbox{\rm{op}}}\in L^\infty(\mu)$, where
$\|\cdot\,\|_{\mbox{\rm{op}}}$ denotes the operator norm.
\end{theorem}
\proof
From Theorem \ref{th3} (i), we know that
$$E[\exp\lambda|\varphi-E[\varphi]|] \leq
2E\left[\exp\frac{\lambda^2\pi^2}{8}|\nabla\varphi|^2\right]\,.
$$
 Hence it is
sufficient to prove that
$$
E\left[\exp\lambda^2|\nabla\varphi|^2\right]<\infty
$$
for some $\lambda>0$. However Theorem \ref{exp-int-thm} applies  since
$\nabla|\nabla\varphi|\in L^\infty(\mu,H)$. The last claim is obvious
since $|\nabla|\nabla\phi|_H|_H\leq \|\nabla^2\varphi\|_{\mbox{op}}$
almost surely.
\qed

\begin{corollary}
\label{cor3}
Let $F\in \DD_{p,1}$ for some $p>1$ such that $|\nabla F|_H \in L^{\infty}(\mu)$.
We then have
\begin{equation}
\label{est-1}
E[\exp\lambda F^2]\leq
  E\left[\frac{1}{\sqrt{1-\frac{\lambda\pi^2}{4}|\nabla F|_H^2}}\exp\left(\frac{\lambda E[F]^2}{1-\frac{\lambda \pi^2}{4}|\nabla F|^2}\right)\right],
\end{equation}
for any $\lambda >0$ such that $\||\nabla F|_H\|_{L^\infty(\mu)}^2 \frac{\lambda\pi^2}{4}<1$.
\end{corollary}
\proof
Ley $Y$ be an auxiliary, real-valued Gaussian random variable, living
on a separate probability space $(\Omega, {\mathcal U},P)$ with
variance one and zero expectation. We have, using Theorem \ref{th3} :
\begin{eqnarray*}
E[\exp \lambda F^2]&=&E\otimes E_P[\exp\sqrt{2\lambda}FY]\\
       &\leq&E\otimes E_P\left[\exp\left\{\sqrt{2\lambda}E[F]Y+
                 |\nabla F|^2 Y^2\frac{\lambda \pi^2}{4}\right\}\right]\\
     &=& E\left[\frac{1}{\sqrt{1-\frac{\lambda\pi^2}{2}|\nabla F|_H^2}}\exp
\left(\frac{\lambda E[F]^2}{1-\frac{\lambda \pi^2}{2}|\nabla F|^2}\right)
\right],
\end{eqnarray*}
where $E_P$ denotes the expectation with respect to the probability $P$.
\qed

\remark
In the next section we shall obtain a better estimate then the one
given by (\ref{est-1}).
\section{Log-Sobolev inequality and exponential integrability}
\markboth{Log-Sobolev}{Exponential integrability}

There is a close relationship between the probability measures
satisfying the log-Sobolev inequality and the exponential
integrability of the random variables having essentially bounded
Sobolev derivatives. We shall explain this in the frame of the Wiener
space: let $\nu$ be a probability measure on $(W,\calB(W))$ such that
the operator $\nabla$ is a closable operator on $L^2(\nu)$.
Assume that we have
$$
E_\nu[\calH_\nu(f^2)]\leq KE_\nu[|\nabla f|_H^2]
$$
for any cylindrical $f:W\to \reals$, where $\calH_\nu(f^2)=f^2(\log
f^2-\log  E_\nu[f^2])$. Since $\nabla$ is a closable operator, of
course this inequality extends immediately
to the extended $L^2$- domain of it.
\begin{lemma}
Assume now that $f$ is  in the
extended $L^2$-domain of $\nabla$ such that  $|\nabla f|_H$ is
$\nu$-essentially bounded by one. Then
\begin{equation}
\label{ineq-0}
E_\nu[e^{tf}]\leq \exp\left\{tE_\nu[f]+\frac{Kt^2}{4}\right\}\,,
\end{equation}
for any $t\in \reals$.
\end{lemma}
\proof
Let $f_n=\min(|f|,n)$, then it is easy to see that $|\nabla f_n|_H\leq
|\nabla f|_H$ $\nu$-almost surely.
Let $t\in \reals$ and define $g_n$ as
to be $e^{\frac{t}{2}f_n}$. Denote by $\theta(t)$ the function
$E[e^{tf_n}]$. Then it follows from the above inequality that
\begin{equation}
\label{ineq-1}
t\theta'(t)-\theta(t)\log\theta(t)\leq \frac{Kt^2}{4}\theta(t)\,.
\end{equation}
If we write $\beta(t)=\frac{1}{t}\log \theta(t)$, then $\lim_{t\to
0}\beta(t)=E[f_n]$, and (\ref{ineq-1}) implies that $\beta'(t)\leq K/4$,
hence we have
$$
\beta(t)\leq E_\nu[f_n]+\frac{Kt}{4}\,,
$$
therefore
\begin{equation}
\label{ineq-2}
\theta(t)\leq \exp\left(tE_\nu[f_n]+\frac{Kt^2}{4}\right)\,.
\end{equation}
It follows from the monotone convergence theorem that
$E[e^{tf}]<\infty$, for any $t\in \reals$. Hence the function
$\theta(t)=E[e^{tf}]$ satisfies also the inequality (\ref{ineq-1})
which implies the inequality (\ref{ineq-0}).
\qed

Using now the inequality (\ref{ineq-0})  and an auxillary Gaussian
random variable as in Corollary \ref{cor3}, we can show easily:
\begin{proposition}
Assume that $f\in L^p(\nu)$  has $\nu$-essentially bounded Sobolev
derivative and that this bound is equal to one. Then we have, for any
$\eps>0$,
$$
E_\nu[e^{\eps f^2}]\leq
\frac{1}{\sqrt{1-\eps K}}\exp\left(\frac{2\eps E_\nu[f]^2}{1-\eps K}\right)\,,
$$
provided $\eps K<1$.
\end{proposition}


\section{An interpolation inequality}
\markboth{Moment Inequalities}{Interpolation Inequality}
Another useful inequality for the Wiener functionals \footnote{This
  result has been proven as an answer to a question posed by
  D. W. Stroock, cf. also \cite{D-H-U}.} is the following
interpolation inequality which helps to control the $L^p$- norm of
$\nabla F$ with the help of the $L^p$-norms of $F$ and
$\nabla^2 F$.
\begin{theorem}
\label{inter-th}
For any $p>1$, there exists a constant $C_p$, such that, for any $F\in \DD_{p,2}$, one has
$$
\|\nabla F\|_p\leq C_p\left[\|F\|_p+\|F\|_p^{1/2}\|\nabla^2F\|_p^{1/2}\right].
$$
\end{theorem}
Theorem \ref{inter-th}  will be proven, thanks to  the Meyer
inequalities, if we can  prove the following
\begin{theorem}
\label{th6}
For any $p>1$, we have
$$
\|(I+{\mathcal L})^{1/2}F\|_p\leq \frac{4}{\Gamma(1/2)}
  \|F\|_p^{1/2}\|(I+{\mathcal L})F\|^{1/2}_p.
$$
\end{theorem}
\proof
Denote by $G$ the functional $(I+{\mathcal L})F$. Then we have
$F=(I+{\mathcal L})^{-1}G$. Therefore it suffices to show that
$$
\|(I+{\mathcal L})^{-1/2}G\|_p\leq  \frac{4}{\Gamma(1/2)}
\|G\|^{1/2}_p\|(I+{\mathcal L})^{-1}G\|^{1/2}_p.
$$
We have
$$
(I+{\mathcal L})^{-1/2}G=\frac{\sqrt{2}}{\Gamma(1/2)}\int_0^\infty
t^{-1/2}e^{-t}P_t G dt,
$$
where $P_t$ denotes the semi-group of Ornstein-Uhlenbeck. For any
$a>0$, we can write
$$
(I+{\mathcal L})^{-1/2}G=\frac{\sqrt{2}}{\Gamma(1/2)}\left[\int_0^a
t^{-1/2}e^{-t}P_t G dt+\int_a^\infty t^{-1/2}e^{-t}P_t G dt\right].
$$
Let us denote the two terms at the right hand side of the above
equality, respectively, by $I_a$ and $II_a$. We have
$$
\|(I+{\mathcal L})^{-1/2}G\|_p\leq \frac{\sqrt{2}}{\Gamma(1/2)}
[\|I_a\|_p+\|II_a\|_p].
$$
The first term at the right hand side can be upper bounded  as
\begin{eqnarray*}
\|I_a\|_p&\leq & \int_0^a t^{-1/2}\|G\|_p dt\\
&=& 2\sqrt{a}\|G\|_p.
\end{eqnarray*}
Let $g=(I+{\mathcal L})^{-1}G$. Then
\begin{eqnarray*}
\int_a^\infty t^{-1/2}e^{-t}P_t G dt&=&\int_a^\infty t^{-1/2}e^{-t}P_t
(I+{\mathcal L}) (I+{\mathcal L})^{-1}G dt\\
&=&\int_a^\infty t^{-1/2}e^{-t}P_t (I+{\mathcal L})g  dt\\
&=&\int_a^\infty t^{-1/2}\frac{d}{dt}(e^{-t}P_t) dt\\
&=&-a^{-1/2}e^{-a}P_a g+\frac{1}{2}\int_a^\infty t^{-3/2} e^{-t}P_tg dt,
\end{eqnarray*}
where the third equality follows from the integration by parts formula.
Therefore
\begin{eqnarray*}
\|II_a\|_p &\leq & a^{-1/2}\|e^{-a}P_ag\|_p+\frac{1}{2}\int_a^\infty t^{-3/2}
\|e^{-t}P_tg\|_pdt\\
 &\leq & a^{-1/2}\|g\|_p+\frac{1}{2}\int_a^\infty t^{-3/2}
\|g\|_pdt\\
&=&2a^{-1/2}\|g\|_p\\
&=&2a^{-1/2}\|(I+{\mathcal L})^{-1}G\|_p.
\end{eqnarray*}
Finally we have
$$
\|(I+{\mathcal L})^{-1/2}G\|_p\leq \frac{2}{\Gamma(1/2)}\left[\,a^{1/2}\|G\|_p+
a^{-1/2}\|(I+{\mathcal L})^{-1}G\|_p\right].
$$
This expression attains its minimum when we take
$$
a=\frac{\|(I+{\mathcal L})^{-1}G\|_p}{\|G\|_p}.
$$
\qed

Combining Theorem \ref{inter-th}  with Meyer inequalities, we have
\begin{corollary}
Suppose that $(F_n, n\in \N)$ converges to zero in $\DD_{p,k}$,
$p>1,k\in\Z$, and that it is bounded in $\DD_{p,k+2}$. Then the
convergence takes place also in  $\DD_{p,k+1}$.
\end{corollary}
\section{Exponential integrability of the divergence}
\markboth{Moment Inequlaties}{Divergence}
We begin with two lemmas which are of some interest:
\begin{lemma}
\label{ident-1}
Let $\phi \in L^p(\mu)$, $p>1$, then, for any $h\in H,\ t>0$, we have
$$
\nabla_h P_t \phi(x)=\frac{e^{-t}}{\sqrt{1-e^{-2t}}}\int_W
\phi(e^{-t}x+\sqrt{1-e^{-2t}}y)\delta h(y)\mu(dy)
$$
almost surely, where $\nabla_h P_t\phi$ represents $(\nabla P_t\phi,h)_H$.
\end{lemma}
\proof
From the  Mehler  formula (cf. \ref{mehler-formula}), we have
\begin{eqnarray*}
\nabla_h P_t\phi(x)&=&\frac{d}{d\lambda}\Big|_{\lambda=0}
       \int_W \phi\left(e^{-t}(x+\lambda h)+\sqrt{1-e^{-2t}}y\right)\mu(dy)\\
 &=&\frac{d}{d\lambda}\Big|_{\lambda=0}
     \int_W \phi\left(e^{-t}x+\sqrt{1-e^{-2t}}\left(y+\frac{\lambda e^{-t}}
    {\sqrt{1-e^{-2t}}}h\right)\right)\mu(dy)\\
&=&\frac{d}{d\lambda}\Big|_{\lambda=0}\int_W
     \phi\left(e^{-t}x+\sqrt{1-e^{-2t}}y\right)
    \varepsilon\left(\frac{\lambda e^{-t}}{\sqrt{1-e^{-2t}}}\delta
      h\right)(y)\mu(dy)\\
 &=&\int_W
 \phi\left(e^{-t}x+\sqrt{1-e^{-2t}}y\right)\frac{e^{-t}}{\sqrt{1-e^{-2t}}}
        \delta h(y) \mu(dy),
\end{eqnarray*}
 where $\varepsilon(\delta h)$ denotes $\exp(\delta h-1/2|h|_{H}^2)$.
\qed

\begin{lemma}
\label{iden-2}
Let $\xi\in L^p(\mu,H),\ p>1$ and for $(x,y)\in W\times W,\ t\geq 0$, define
$$
R_t(x,y)=e^{-t}x+(1-e^{-2t})^{1/2}y
$$
and
$$
S_t(x,y)=(1-e^{-2t})^{1/2}x-e^{-t}y\,.
$$
Then  $S_t(x,y)$ and $R_t(x,y)$ are independant,
identically distributed Gaussian random variables on
 $(W\times W,\mu(dx)\times \mu(dy))$. Moreover  the following
 identity holds true:
$$
P_t\delta\xi(x)=\frac{e^{-t}}{\sqrt{1-e^{-2t}}}
     \int_W I_1(\xi(R_t(x,y)))(S_t(x,y))\mu(dy),
$$
where
$$
I_1(\xi(R_t(x,y)))(S_t(x,y))
$$
denotes the first order Wiener
integral of $\xi(R_t(x,y))$ with respect to the independent path
$S_t(x,y)$ under the product measure $\mu(dx)\times \mu(dy)$.
\end{lemma}
\proof
The first part of the lemma is a well-known property  of the Gaussian
random variables and left to the reader. In the proof of  the second part,
for the typographical facility, we shall denote in the sequel by $e(t)$ the
function $(\exp-t)/(1-\exp-2t)^{1/2}$. Let now $\phi$ be an element of
$\DD$, we have, via duality and using  Lemma \ref{ident-1}
\begin{eqnarray*}
<P_t\delta\xi,\phi>&=&<\xi,\nabla P_t\phi>\\
                   &=&\sum_{i=1}^{\infty}<\xi_i,\nabla_{h_i}P_t\phi>\\
              &=&\sum_i e(t)E\left[\int_W
                \xi_i(x)\phi\left(R_t(x,y)\right)\delta h_i(y)
                          \mu(dy)\right]
\end{eqnarray*}
where $(h_i;i\in \N)\subset W^{\star}$ is a complete orthonormal basis
of  $H$, $\xi_i$ is the component of $\xi$ in the direction of $e_i$
and $<.,.>$  represents the duality bracket corresponding to the dual
pairs $(\DD,\DD')$
or $(\DD(H),\DD'(H))$. Let us make the following change of variables,
which preserves $\mu\times \mu$ :
\begin{eqnarray*}
 x &\mapsto & e^{-t}x+\sqrt{1-e^{-2t}}y\\
 y & \mapsto & \sqrt{1-e^{-2t}}x-e^{-t}y.
\end{eqnarray*}
We then obtain
$$
<P_t(\delta\xi),\phi>=e(t)\int_W
\phi(x)I_1\left(\xi(R_t(x,y))\right)\left(S_t(x,y)\right)\mu(dx)
\mu(dy),
$$
for any $\phi\in \DD$ and the lemma follows from the density of $\DD$ in all
$L^p$-spaces.
\qed

We are now ready to prove the following
\begin{theorem}
Let $\beta >1/2$ and suppose that $\eta\in \DD_{2,2\beta}(H)$. Then we have
$$
E\left[\exp\delta\eta\right]\leq E\left[\exp\left(\alpha
  |(2I+\L)^{\beta}\eta|_{H}^2\right)\right],
$$
for any $\alpha$ satisfying
$$
\alpha\geq \frac{1}{2}\left[\frac{1}{\Gamma(\beta)}
\int_{\R_+}\frac{t^{\beta-1} e^{-2t}}{\sqrt{1-e^{-2t}}}dt\right]^{-2},
$$
where $\L$ denotes the Ornstein-Uhlenbeck or the number operator on $W$.
\end{theorem}
\proof
Let $\xi=(2I+\L)^{\beta}\eta$, then the above inequality is equivalent to
$$
E\left[\exp\left((I+\L)^{-\beta}\delta\xi\right)\right]\leq
E\left[\exp\alpha|\xi|_{H}^2\right]\,,
$$
where we have used the identity
$$
(I+\L)^{-\beta}\delta\xi=\delta\left((2I+\L)^{-\beta}\xi\right).
$$
We have from the resolvent identity  and from the Lemma \ref{iden-2},
\begin{eqnarray*}
\lefteqn{(I+\L)^{-\beta}\delta\xi=\frac{1}{\Gamma(\beta)}
                            \int_{\R_+}t^{\beta-1}e^{-t}P_t\delta\xi dt}\\
 &=&\int_{\R_+\times W}\frac{e^{-t}}{\Gamma(\beta)\sqrt{1-e^{-2t}}}
              t^{\beta-1}e^{-t}I_1(\xi(R_t(x,y)))(S_t(x,y))\mu(dy) dt.
\end{eqnarray*}
Let
$$
\lambda_0=\frac{1}{\Gamma(\beta)}
\int_{\R_+}\frac{t^{\beta-1} e^{-2t}}{\sqrt{1-e^{-2t}}}dt
$$
and
$$
\nu(dt)=\won_{\R_+}(t)\frac{1}{\lambda_0 \Gamma(\beta)}
\frac{t^{\beta-1} e^{-2t}}{\sqrt{1-e^{-2t}}}dt.
$$
Then, from the H\"older inequality
\begin{eqnarray*}
\lefteqn{E\left[\exp\left\{(I+\L)^{-\beta}\delta\xi\right\}\right]}\\
&=&E\left[\exp
   \left\{ \lambda_0
  \int_{\R_+}\int_W I_1(\xi(R_t(x,y)))(S_t(x,y))\mu(dy)\nu(dt)\right\}\right]\\
&\leq &\int_{\R_+}\int_W \int_W
        \exp\left\{\lambda_0
          I_1(\xi(R_t(x,y)))(S_t(x,y))\right\}\mu(dx)\mu(dy)\nu(dt)\\
&=&E\left[\exp\left\{ \frac{\lambda_0^2}{2}|\xi|_{H}^2\right\} \right]\,,
\end{eqnarray*}
which completes the proof.
\qed

\noindent
In the applications, we need also to control the moments like
$E[\exp\|\nabla\eta\|_2^2]$ (cf. \cite{UZ7}), where $\eta$ is an
$H$-valued random variable and  $\|\,.\,\|_2$ denotes
the Hilbert-Schmidt norm. The following result gives an answer to
this question:
\begin{proposition}
Suppose that $\beta>1/2$ and that $\eta\in \DD_{2,2\beta}(H)$. Then we have
$$
E[\exp\|\nabla\eta\|_2^2]\leq E[\exp c|(I+\L)^{\beta}\eta|_{H}^2],
$$
for any
$$
c\geq c_0= \left[\frac{1}{\Gamma(\beta)}\int_{\R_+}t^{\beta-1}e^{-2t}
             (1-e^{-2t})^{-1/2}dt\right]^2\,.
$$
In particular, for $\beta=1$ we have $c\geq 1/4$.
\end{proposition}
\proof
Setting $\xi=(I+\L)^{\beta}\eta$, it is sufficient to show that
$$
E\left[\exp\|\nabla (I+\L)^{-\beta}\xi\|_2^2\right]\leq E\left[\exp
  c|\xi|_{H}^2\right]\,.
$$
Let $(E_i, i\in \N)$ be a complete, orthonormal basis of $H\otimes H$
which is  the completion of the tensor product of
$H$ with itself under the Hilbert-Schmidt topology.  Then
$$
\|\nabla \eta\|_2^2=\sum_i K_i(\nabla\eta,E_i)_2,
$$
where $(.,.)_2$ is the scalar product in $H\otimes H$ and
 $K_i=(\nabla \eta,E_i)_2$. Let   $\theta(t)$ be the function
$$
 \frac{1}{\Gamma(\beta)}t^{\beta-1}e^{-2t}(1-e^{-2t})^{-1/2}
$$
and let $\gamma_0=\int_0^\infty \theta(t)dt$. From  Lemmas
\ref{ident-1} and \ref{iden-2}, we have
\begin{eqnarray*}
\|\nabla\eta(x)\|_2^2&=&\|\nabla(I+\L)^{-\beta}\xi(x)\|_2^2\\
    &=&\sum_i K_i(x)\int_{\R_+}\theta(t)
           \int_W\left(I_1(E_i)(y),\xi(R_t(x,y))\right)_H \mu(dy)dt\\
    &=&\int_{\R_+\times W}\theta(t)
           \left(I_1(\nabla\eta(x))(y),\xi(R_t(x,y))\right)_H \mu(dy)dt\\
  &\leq &\int_{\R_+}\theta(t)
       \left(\int_W|I_1(\nabla\eta(x))(y)|_H^2\mu(dy)\right)^{1/2}\\
  && \hspace{1.5cm} \left(\int_W|\xi(R_t(x,y))|_H^2\mu(dy)\right)^{1/2}dt\\
&=&\int_{\R_+}\theta(t)
           \|\nabla\eta(x)\|_2(P_t(|\xi|_H^2))^{1/2}dt,
\end{eqnarray*}
where $I_1(\nabla \eta(x))(y)$ denotes the first order Wiener integral of
$\nabla\eta(x)$ with respect to the independent path (or variable) $y$.
Consequently we have the following inequality:
$$
\|\nabla\eta\|_2\leq  \int_{\R_+}\theta(t)\left(P_t(|\xi|_H^2)\right)^{1/2}dt.
$$
Therefore 
\begin{eqnarray*}
E[\exp\|\nabla \eta\|_2^2]&\leq & E\left[\exp\left\{\int_{\R_+}\gamma_0^2
  P_t(|\xi|_H^2)
     \frac{\theta(t)}{\gamma_0}dt\right\}\right]\\
  &\leq & E\int_{\R_+}\frac{\theta(t)}{\gamma_0}
     \exp\left\{\gamma_0^2 P_t(|\xi|_H^2)\right\}dt\\
 &\leq &  E\int_{\R_+}\frac{\theta(t)}{\gamma_0}
     \exp\left\{\gamma_0^2 |\xi|_H^2\right\} dt\\
&=&E\left[\exp \gamma_0^2|\xi|_H^2\right].
\end{eqnarray*}
\qed

\noindent
As an example of application of these results let us give the
following theorem of  the degree theory of the Wiener maps (cf. \cite{UZ7}):
\begin{corollary}
\label{deg-appli}
Suppose that $\eta\in \DD_{2,2\beta}(H)$, $\beta>1/2$,  satisfies
$$
E\left[\exp a\left|(2I+\L)^\beta\eta\right|_H^2\right]<\infty,
$$
for some $a>0$. Then for any $\lambda\leq \sqrt{\frac{a}{4c_0}}$ and 
$h\in H$,  we have
$$
E\left[e^{i\left(\delta
      h+\lambda(h,\eta\right)_H)}\Lambda\right]=\exp-\frac{1}{2}|h|_H^2, 
$$
where $\Lambda$ is defined by
$$
\Lambda=\dett(I_H+\lambda\nabla\eta)\exp\left\{-\lambda\delta
h-\frac{\lambda^2}{2}|\eta|_H^2\right\}.
$$
In particular, if we deal with the classical Wiener space, the path defined by
$$
T_\lambda(w)=w+\lambda \eta(w),
$$
is a Brownian motion under the new probability measure
$E[\Lambda|\sigma(T_\lambda)]d\mu$, where $\sigma(T_\lambda)$ denotes
the sigma field generated by the mapping $T_\lambda$.
\end{corollary}
\proof
This result follows from the degree theorem for Wiener maps
(cf. \cite{UZ2}). In fact from the Theorem 3.2 of \cite{UZ2} (cf. also
\cite{UZ7}), it follows that $E[\Lambda]=1$. On the other hand, from
the Theorem 3.1 of the same reference, we have
$$
E[F\circ T_\lambda \Lambda]=E[F]E[\Lambda].
$$
Hence the proof follows.
\qed

\section*{Notes and suggested reading}
{\footnotesize{The results about the exponential tightness go back
    till to the celebrated Lemma of X. Fernique about the
    exponential integrability of the square of semi-norms
    (cf. \cite{KUO}). It is also proven
    by B. Maurey in the finite dimensional case for the Lipschitz
    continuous maps with the same  method that we have used here
    (cf. \cite{P}). A similar
    result in the abstract Wiener space case has been given by
    S. Kusuoka under the hypothesis of $H$-continuity, i.e.,  $h\to
    \phi(w+h)$ is continuous for any $w\in W$. We have proven the actual result
    without this latter hypothesis. However, it has been proven later
    that the essential boundedness of the Sobolev derivative implies
    the existence of a version which is $H$-continuous by Enchev and
    Stroock (cf. \cite{E-S}). Later it has been discovered that the
    exponential integrability is implied by the logarithmic Sobolev
    inequality (cf. \cite{A-M-S}). The derivation of the inequality
    (\ref{ineq-2}) is attributed to Herbst (cf. \cite{Led}).

    In any case the exponential
    integrability of the square of the Wiener functionals has found
    one of its most important applications in the analysis of
    non-linear Gaussian functionals. In fact in the proof of the Ramer
    theorem and its extensions this property plays an important role
    (cf. Chapter X, \cite{UZ5}, \cite{UZ6}  and  \cite{UZ7}).
    Corollary \ref{deg-appli} uses some results about the degree
    theory of the Wiener maps which are explained below:
\begin{theorem}
\label{degree-thm}
   Assume that $\gamma$ and $r$ be fixed strictly positive numbers
   such that $r>(1+\gamma)\gamma{-1}$. Let $u\in\DD_{r,2}(H)$ and assume
   that
\begin{enumerate}
\item $\La_u\in L^{1+\gamma}(\mu)$,
\item $\La_u(I_H+\nabla u)^{-1}h\in L^{1+\gamma}(\mu,H)$ for any $h\in H$,
\end{enumerate}
where
$$
\La_u=\dett(I_H+\nabla u)\exp\left\{-\delta u-\frac{1}{2}|u|_H^2\right\}\,.
$$
Then, for any $F\in C_b(W)$, we have
$$
E\left[F(w+u(w))\La_u\right]=E[\La_u]E[F]\,.
$$
\end{theorem}
In particular, using a homotopy argument, one can show that, if
$$
\exp\left\{-\delta u+\frac{1+\eps}{2}\|\nabla u\|_2^2\right\}\in
L^{1+\alpha}(\mu)\,,
$$
for some $\alpha>0,\,\eps>0$, then $E[\La_u]=1$.
We refer the reader to \cite{UZ7} for further information about this topic.

}}


\chapter{Introduction to the Theorem of Ramer}
\label{ch.ramer}
\index[sub]{Ramer's Theorem}
\markboth{Theorem of Ramer}{}

\section*{Introduction}
The Girsanov theorem tells us that if $u:W\mapsto H$ is a Wiener
functional such that $\frac{du}{dt}=\dot{u}(t)$ is an adapted process
such that
$$
E\left[\exp\left\{ -\int_0^1 \dot{u}(s)
  dW_s-\frac{1}{2}\int_0^1|\dot{u}(s)|^2ds\right\}\right]=1,
$$
then under the new probability $Ld\mu$, where
$$
L=\exp\left\{ -\int_0^1 \dot{u}(s)
  dW_s-\frac{1}{2}\int_0^1|\dot{u}(s)|^2ds\right\},
$$
$w\to w+u(w)$ is a Brownian motion. The theorem of Ramer studies the
same problem
without hypothesis of adaptedness of the process $\dot{u}$. This problem has
been initiated by Cameron and Martin. Their work has been  extended by Gross
and others. It was  Ramer \cite{R} who gave a main impulse to the problem by
realizing that the ordinary determinant can be replaced by the modified
Carleman-Fredholm determinant via defining a Gaussian divergence instead of
the ordinary Lebesgue
divergence. The problem has been further studied by Kusuoka \cite{Kus} and the
final solution in the case of (locally) differentiable shifts in the
Cameron-Martin space direction  has been given by \"Ust\"unel and
Zakai \cite{UZ5}.  In this chapter we will give a partial
( however indispensable for the proof of the  general ) result.

To understand the problem, let us consider first the finite dimensional case:
let $W=\R^n$ and let $\mu_n$ be the standard Gauss measure on
$\R^n$. If $u:\R^n\mapsto \R^n$ is a differentiable mapping such that
$I+u$ is a diffeomorphism of $\R^n$, then the theorem of Jacobi tells
us that, for any smooth function $F$ on $\R^n$, we have 
$$
\int_{\R^n}F(x+u(x))|\det(I+\partial
u(x))|\exp\left\{-<u(x),x>-\frac{1}{2}|u|^2\right\}\mu_n(dx) 
$$
$$
 = \int_{\R^n}F(x)\mu_n(dx),
$$
where $\partial u$ \index[not]{d@$\partial u$} denotes the derivative of $u$.
The natural idea now is to pass to the infinite dimension. For this,
note that, if we define $\dett(I+\partial\, u)$ by
\begin{eqnarray*}
\dett(I+\partial u(x))&=&{\mbox{\rm det}}(I+\partial\, u(x))
\,\,e^{-\t\, [\partial\, u(x)]}\\
  &=&\prod_i (1+\lambda_i)\exp-\lambda_i,
\end{eqnarray*}
where $(\lambda_i)$ are the eigenvalues of $\partial\,u(x)$ counted
with respect to their multiplicity,
then the density of the left hand side can be written as
\index[not]{c@$\dett(I+\nabla u)$}
$$
\Lambda=|\dett(I+\partial\, u(x))|\exp\left\{-<u(x),x>+\t \, \partial\,
u(x)-\frac{1}{2}|u|^2\right\}
$$
and let us remark that
$$
<u(x),x>-\t \, \partial\, u(x)=\delta u(x),
$$
where $\delta$ is the adjoint of the $\partial$ with respect to the
Gaussian measure $\mu_n$.
Hence, we can express the density $\Lambda$  as
$$
\Lambda=|\dett(I+\partial u(x))|\exp\left\{-\delta
  u(x)-\frac{|u(x)|^2}{2}\right\}\,. 
$$
 As remarked first by Ramer, cf. \cite{R}, this expression has two
advantages: first $\dett(I+\partial u)$, called Carleman-Fredholm determinant,
\index[sub]{Carleman-Fredholm determinant} can be defined for the
mappings $u$ such that $\partial u(x)$ is with values in the space of
Hilbert-Schmidt operators rather than nuclear operators (the latter is
a smaller class than the former), secondly, as we have already seen,
$\delta u$ is well-defined for a large class of mappings meanwhile
$<u(x),x>$ is a highly singular object in the Wiener space.

\section{Ramer's Theorem}
After  these preliminaries, we can announce,  using our standard
notations,  the main result of this
chapter:
\begin{theorem}
\label{ramer-th}
Suppose that $u:W\mapsto H$ is a measurable map belonging to
$\DD_{p,1}(H)$ for some $p>1$. Assume that there are constants $c$ and
$d$ with {\bf$c<1$} such that for almost all $w\in W$, 
$$
\|\nabla u\|\leq c<1
$$
and
$$
\|\nabla u\|_2\leq d<\infty,
$$
where $\| \cdot \| $ denotes the operator norm and
$\| \cdot \|_2 $ denotes the Hilbert-Schmidt norm for the linear
operators on $H$. 
Then:
\begin{itemize}
\item
Almost surely $w\mapsto T(w) = w + u(w)$ is bijective. The inverse of
$T$, denoted by $S$ is of the form $S(w)=w+v(w)$, where $v$ belongs to
$\DD_{p,1}(H)$ for any $p>1$, moreover 
$$
\|\nabla v\|\leq \frac{c}{1-c}\,\,{\mbox{and}}\,\, \|\nabla
v\|_2\leq\frac{d}{1-c}, 
$$
$\mu$-almost surely.

\item For all bounded and measurable $F$, we have
$$
E[F(w)] = E[F(T(w)) \cdot | \Lambda_u(w) |]
$$
 and in particular
$$
E|\Lambda_u| =1,
$$
where
$$
\Lambda_u=|\dett(I+\nabla u)|\exp-\delta u-\frac{1}{2}|u|_H^2,
$$
and $\det_2(I+\nabla u)$ denotes the Carleman-Fredholm determinant of
$I+\nabla u$.
\item
The measures
$\mu$,  $ T^{\star}\mu$ and $S^*\mu$ are mutually absolutely continuous, where
$T^{\star}\mu$ (respectively $S^*\mu$) denotes the image of $\mu$
under $T$ (respectively $S$). We have 
\begin{eqnarray*}
\frac{dS^*\mu}{d\mu}&=&|\Lambda_u|\,,\\
\frac{dT^*\mu}{d\mu}&=&|\Lambda_v|,
\end{eqnarray*}
where $\Lambda_v$ is defined similarly.
\end{itemize}
\end{theorem}
\begin{remarkk}{\rm
If $\|\nabla u\|\leq 1$ instead of
$\|\nabla u\|\leq c< 1$, then taking $u_\epsilon=(1-\epsilon)u$ we see that
 the hypothesis of the theorem are satisfied for $u_\epsilon$. Hence using the
   Fatou lemma, we obtain
$$
E[F\circ T \,|\Lambda_u|]\leq E[F]
$$
for any positive $F\in C_b(W)$. Consequently, if $\Lambda_u\neq 0$
almost surely, then $T^*\mu$ is absolutely continuous with respect to
$\mu$.} 
\end{remarkk}

The proof of  Theorem \ref{ramer-th}  will be done in several steps. As we have
indicated above, the main idea is to pass to the limit from finite to
infinite dimensions. The key point in this procedure will be the use
of the Theorem 1 of the preceding chapter which will imply the
uniform integrability of the finite dimensional densities.
We shall first prove  the same theorem in the  cylindrical case:

\begin{lemma}
\label{lemma1}
Let $\xi:W\mapsto H$ be a shift of the following form:
$$
\xi(w)=\sum_{i=1}^{n}\alpha_i(\delta h_1,\ldots,\delta h_n)h_i,
$$
with $\alpha_i\in C^\infty(\R^n)$ with bounded first derivative,
$h_i\in W^*$ are orthonormal\footnote{In fact
$h_i \in W^*$ should be distinguished from its image in $H$, denoted
by $j(h)$. For notational simplicity, we denote both by $h_i$, as long
as there is no ambiguity.} in $H$. 
Suppose furthermore that $\|\nabla \xi\|\leq c<1$ and that
$\|\nabla\xi \|_2\leq d$ as above. Then we have
\begin{itemize}
\item
Almost surely $w\mapsto U(w) = w + \xi(w)$ is bijective.
\item
The measures
$\mu$ and $ U^{\star}\mu$ are mutually absolutely continuous.
\item For all bounded and measurable $F$, we have
$$
E[F(w)] = E[F(U(w)) \cdot | \Lambda_\xi(w) |]
$$
for all bounded and measurable $F$ and in particular
$$
E[|\Lambda_\xi|] =1,
$$
where
$$
\Lambda_\xi=|\dett(I+\nabla \xi)|\exp-\delta \xi-\frac{1}{2}|\xi|_H^2.
$$
\item The inverse of $U$, denoted by $V$ is of the form $V(w)=w+\eta(w)$,
where $$
\eta(w)=\sum_{i=1}^{n}\beta_i(\delta h_1,\ldots,\delta h_n)h_i,
$$
such that  $\|\nabla \eta\|\leq \frac{c}{1-c}$ and
$\|\nabla\eta \|_2\leq\frac{d}{1-c}$.
\end{itemize}
\end{lemma}
\proof
Note first that due to the Corollary \ref{exp-cor} of the Chapter
VIII, $E[\exp \lambda|\xi|^2]<\infty$ for any $\lambda
<\frac{1}{2c}$. We shall construct the inverse of $U$ by imitating the
fixed point techniques: let
\begin{eqnarray*}
\eta_0(w)&=&0\\
\eta_{n+1}(w)&=&-\xi(w+\eta_n(w)).
\end{eqnarray*}
We have
\begin{eqnarray*}
|\eta_{n+1}(w)-\eta_n(w)|_H&\leq & c|\eta_{n}(w)-\eta_{n-1}(w)|_H\\
                            &\leq & c^n|\xi(w)|_H.
\end{eqnarray*}
Therefore $\eta(w)=\lim_{n\rightarrow \infty}\eta_n(w)$ exists and it
is bounded  by $\frac{1}{1-c}|\xi(w)|_H$. By the triangle inequality
\begin{eqnarray*}
\left|\eta_{n+1}(w+h)-\eta_{n+1}(w)\right|_H&\leq &\left|\xi(w+h+\eta_n(w+h))-
                             \xi(w+\eta_n(w))\right|_H\\
                 &\leq &c|h|_H+c|\eta_n(w+h)-\eta_n(w)|_H.
\end{eqnarray*}
Hence passing to the limit, we find
$$
|\eta(w+h)-\eta(w)|_H\leq \frac{c}{1-c}|h|_H.
$$
We also have
\begin{eqnarray*}
U(w+\eta(w))&=&w+\eta(w)+\xi(w+\eta(w))\\
            &=&w+\eta(w)-\eta(w)\\
           &=&w,
\end{eqnarray*}
hence $U\circ (I_W+\eta)=I_W$, i.e., $U$ is an onto map. If $U(w)=U(w')$, then
\begin{eqnarray*}
|\xi(w)-\xi(w')|_H&=&|\xi(w'+\xi(w')-\xi(w))-\xi(w')|_H\\
                &\leq &c|\xi(w)-\xi(w')|_H,
\end{eqnarray*}
which implies that $U$ is also injective.
 To show the Girsanov identity, let us complete the sequence
 $(h_i, i\leq n)$ to a complete orthonormal basis whose elements are
 chosen from $W^*$. From a theorem of Ito-Nisio  \cite{I-N}, we can
 express the Wiener path $w$ as
$$
w=\sum_{i=1}^\infty \delta h_i(w) h_i,
$$
where the sum converges almost surely in the norm topology of $W$. Let
$F$ be a nice function on $W$, denote by $\mu_n$ the image of the
Wiener measure $\mu$ under the map $w\mapsto \sum_{i\leq n} \delta
h_i(w) h_i$ and by $\nu$ the image of $\mu$ under $w\mapsto \sum_{i>
  n} \delta h_i(w) h_i$. Evidently 
 $\mu=\mu_n\times \nu$. Therefore
\begin{eqnarray*}
E[F\circ U\, |\Lambda_\xi|]&=&
\!\int_{\R^n}E_{\nu}\left[F\left(w+\sum_{i\leq
      n}(x_i+\alpha_i(x_1\ldots,x_n))h_i\right)|\Lambda_\xi|\right]
\mu_{\R^n}(dx)\\ 
  &=&E[F],
\end{eqnarray*}
where $\mu_{\R^n}(dx)$ denotes the standard Gaussian measure on $\R^n$
and the equality follows from the Fubini theorem. In fact by changing
the order of integrals, we reduce the problem to a finite dimensional
one and then the result is immediate from the theorem of Jacobi as
explained above. From the construction of $V$, it is trivial to see
that 
$$
\eta(w)=\sum_{i\leq n}\beta_i(\delta h_1,\ldots,\delta h_n)h_i,
$$
for some  vector field  $(\beta_1,\ldots,\beta_n)$ which is a
$C^\infty$ mapping from $\R^n$ into itself due to the finite
dimensional inverse mapping theorem. Now it is routine to verify that 
$$
\nabla \eta=-(I+\nabla \eta )^\star  \nabla\xi \circ V,
$$
hence
\begin{eqnarray*}
\|\nabla \eta\|_2&\leq &\|I+\nabla \eta\|\|\nabla \xi\circ V\|_2\\
                 &\leq & (1+\|\nabla \eta\|)\|\nabla \xi\circ V\|_2\\
                &\leq &d\,\left(1+\frac{c}{1-c}\right)\\
                &=&\frac{d}{1-c}\,\,.
\end{eqnarray*}
\qed

\begin{lemma}
\label{lemma2}
With the notations and hypothesis of Lemma \ref{lemma1}, we have
$$
\delta \xi\circ V=-\delta \eta-|\eta|_H^2+{\mbox{\rm trace}}\left[(\nabla
\xi\circ V)\cdot \nabla \eta\right]\,,
$$
almost surely.
\end{lemma}
\proof
We have
$$
\delta \xi=\sum_{i=1}^\infty \left\{(\xi,e_i)_H\delta
  e_i-\nabla_{e_i}(\xi,e_i)_H\right\},
$$
where the sum converges in  $L^2$ and the result  is independent
of the choice of the orthonormal basis $(e_i;\,i\in \N)$. Therefore we
can choose as basis $h_1,\ldots, h_n$ that we have already used in
Lemma \ref{lemma1}, completed  with the elements of $W^*$ to form an
orthonormal basis of $H$, denoted by $(h_i;\, i\in \N)$. Hence
$$
\delta \xi=\sum_{i=1}^n \left\{(\xi,h_i)_H\delta
  h_i-\nabla_{h_i}(\xi,h_i)_H\right\}\,.
$$
From the Lemma \ref{lemma1}, we have $\xi\circ V=-\eta$ and since, $h_i$ are
originating from $W^*$, it is immediate to see that $\delta h_i\circ
V=\delta h_i+(h_i,\eta)_H$. Moreover, from the preceding lemma we know
that $\nabla (\xi\circ V)=(I+\nabla\eta)^*\nabla \xi\circ V$.
 Consequently, applying all this, we obtain
\begin{eqnarray*}
\delta \xi\circ V&=&\sum_1^n(\xi\circ V,h_i)_H(\delta h_i+(h_i,\eta)_H)-
 (\nabla_{h_i}(\xi,h_i)_H)\circ V\\
&=&(\xi\circ V,\eta)_H+\delta(\xi\circ V)+
 \sum_1^n\nabla_{h_i}(\xi\circ V,h_i)_H-\nabla_{h_i}(\xi,h_i)_H\circ V\\
&=&-|\eta|_H^2-\delta \eta+
\sum_1^n\left(\nabla\xi\circ V\, [h_i],\nabla \eta\,[h_i]\right)_H\\
&=&-|\eta|_H^2-\delta \eta+{\mbox{\rm trace}}(\nabla \xi\circ V\cdot
\nabla \eta),
\end{eqnarray*}
where $\nabla\xi\,[h]$ denotes the Hilbert-Schmidt operator
$\nabla\xi$ applied to the vector $h\in H$.
\qed
\begin{remarkk}
{\sl
Since $\xi$ and $\eta$ are symmetric, we have $\eta\circ U=-\xi$ and
consequently
$$
\delta \eta\circ U=-\delta \xi-|\xi|_H^2+{\mbox{\rm trace}}\left[(\nabla
\eta\circ U)\cdot \nabla \xi\right]\,.
$$ }
\end{remarkk}
\begin{corollary}
\label{identification-cor}
For any cylindrical function $F$ on $W$, we have
\begin{eqnarray*}
E[F\circ V]&=&E\left[F\, |\Lambda_{\xi}|\right].\\
E[F\circ U]&=&E\left[F\, |\Lambda_{\eta}|\right].
\end{eqnarray*}
\end{corollary}
\proof
The first part follows from the identity
\begin{eqnarray*}
E\left[F\, |\Lambda_{\xi}|\right]&=&E\left[F\circ V\circ U\,
  |\Lambda_\xi|\right]\\ 
&=&E[F\circ V].
\end{eqnarray*}
To see the second part, we have
\begin{eqnarray*}
E[F\circ U]=&=&E\left[F\circ U \frac{1}{|\Lambda_\xi|\circ V}\circ U\,
  |\Lambda_\xi|\right]\\ 
&=&E\left[F\, \frac{1}{|\Lambda_\xi|\circ V}\right].
\end{eqnarray*}
From  Lemma \ref{lemma2}, it follows that
\begin{eqnarray*}
\frac{1}{|\Lambda_\xi|\circ V}&=&\frac{1}{|{\mbox{det}}_2(I+\nabla\xi)\circ V|}
\exp\left\{\delta\xi+1/2|\xi|_H^2\right\}\circ V\\
&=&\frac{1}{|{\mbox{det}}_2(I+\nabla\xi)\circ
  V|}\\
&&\hspace{1cm}\exp\left\{-\delta\eta-1/2|\eta|_H^2+
{\mbox{\rm trace}}((\nabla\xi\circ V)\cdot \nabla \eta)\right\}\\
&=&|\Lambda_\eta|,
\end{eqnarray*}
since, for general Hilbert-Schmidt maps $A$ and $B$, we have
\begin{equation}
\label{det-iden}
\dett(I+A)\cdot
\dett(I+B)=\exp\left\{{\mbox{\rm trace}}(AB)\right\}
\dett((I+A)(I+B))
\end{equation}
and in our case we have
$$
(I+\nabla\xi\circ V)\cdot (I+\nabla \eta)=I\,.
$$
\qed

\remark
In fact the  equality  (\ref{det-iden})  follows from the
multiplicative property of
the ordinary determinants and from the formula (cf.
\cite{D-S}, page 1106, Lemma 22):
$$
\dett(I+A)=\prod_{i=1}^\infty (1+\la_i)e^{-\la_i}\,,
$$
where $(\la_i,i\in \NN)$ are the eigenvalues of $A$ counted with
respect to their multiplicity.

\paragraph{Proof of Theorem \ref{ramer-th}:} Let $(h_i,i\in \N)\subset
W^*$ be a complete
orthonormal basis of $H$. For $n\in \N$, let $V_n$ be the sigma algebra on $W$
generated by $\{\delta h_1,\ldots, \delta h_n\}$, $\pi_n$ be the
orthogonal projection of $H$ onto the subspace spanned by
$\{h_1,\ldots,h_n\}$. Define
$$
\xi_n=E\left[\pi_n\,P_{1/n}u|V_n\right]\,,
$$
where $P_{1/n}$ is the Ornstein-Uhlenbeck semi-group on $W$ with $t=1/n$. Then
$\xi_n\rightarrow \xi$ in $\DD_{p,1}(H)$ for any $p>1$ (cf., Lemma
\ref{integ-lem} of Chapter IX). Moreover $\xi_n$ has the following
form:
$$
\xi_n=\sum_{i=1}^n \alpha_i^n(\delta h_1,\ldots,\delta h_n)h_i,
$$
where $\alpha_i^n$ are $C^\infty$-functions due to the finite
dimensional  Sobolev embedding theorem. We have
$$
\nabla \xi_n=E\left[\pi_n\otimes \pi_n e^{-1/n}P_{1/n}\nabla u|V_n\right]\,,
$$
hence
$$
\|\nabla \xi_n\|\leq e^{-1/n}E\left[P_{1/n}\|\nabla u\||V_n\right],
$$
and the same inequality holds also with the Hilbert-Schmidt
norm. Consequently, we have
$$
\|\nabla \xi_n\|\leq c\,\,,\,\|\nabla \xi_n\|_2\leq d\,,
$$
$\mu$-almost surely. Hence, each $\xi_n$ satisfies the hypothesis of Lemma
\ref{lemma1}. Let us denote by $\eta_n$ the shift corresponding to the
inverse of  $U_n=I+\xi_n$ and let $V_n=I+\eta_n$. Denote by
$\Lambda_n$ and $L_n$ the densities corresponding,  respectively, to
$\xi_n$ and $\eta_n$, i.e., with the old notations
$$
\Lambda_n=\Lambda_{\xi_n}\,\,\,{\mbox{and}}\,\,\, L_n=\Lambda_{\eta_n}.
$$
We will prove that the sequences of densities
$$
\{\Lambda_n\,:n\in \N\}\,\,\,{\mbox{and}}\,\,\,\{L_n\,:n\in\N\}
$$
are uniformly integrable. In fact we will do this only for the first
sequence since the proof for the second is very similar to the proof
of the first case.
To prove the uniform integrability, from the lemma of de la
Vall\'e-Poussin, it suffices to show
$$
\sup_n E\left[|\Lambda_n||\log\Lambda_n|\right]<\infty,
$$
which amounts to show, from the Corollary \ref{identification-cor}, that
$$
\sup_n E\left[|\log\Lambda_n\circ V_n|\right]<\infty\,.
$$
Hence we have to control
$$
E\left[|\log{\mbox{det}}_2(I+\nabla\xi_n\circ V_n)|+
|\delta\xi_n\circ V_n|+1/2|\xi_n\circ V_n|^2\right]\,.
$$
From the Lemma \ref{lemma2}, we have
$$
\delta \xi_n\circ V_n=-\delta \eta_n-|\eta_n|_H^2+{\mbox{\rm
    trace}}(\nabla \xi_n\circ V_n)\cdot \nabla \eta_n,
$$
hence
\begin{eqnarray*}
E[|\delta \xi_n\circ V_n|]&\leq&\|\delta \eta_n\|_{L^2(\mu)}+E[|\eta_n|^2]+
E[\|\nabla\xi_n\circ V_n\|_2\,\|\nabla \eta_n\|_2]\\
&\leq& \|\eta_n\|_{L^2(\mu,H)}+\|\eta_n\|_{L^2(\mu,H)}^2+
 \|\nabla \eta_n\|_{L^2(\mu,H\otimes H)}+
\frac{d^2}{1-c}\\
&\leq& \|\eta_n\|_{L^2(\mu,H)}+\|\eta_n\|_{L^2(\mu,H)}^2+ \frac{d(1+d)}{1-c}\,,
\end{eqnarray*}
where the second inequality follows from
$$
\|\delta \gamma\|_{L^2(\mu)}\leq \|\nabla\gamma\|_{L^2(\mu,H\otimes H)}+
               \|\gamma\|_{L^2(\mu,H)}.
$$
From the Corollary \ref{exp-cor} of Chapter IX, we have
$$
\sup_n E\left[\exp \alpha |\eta_n|_H^2\right]<\infty,
$$
for any $\alpha<\frac{(1-c)^2}{2d^2}$, hence
$$
\sup_n E[|\eta_n|^2]<\infty\,\,.
$$
We have a well-known inequality (cf. \cite{UZ7}, Appendix), which says that
$$
|\dett(I+A)|\leq \exp \frac{1}{2}\|A\|^2_2
$$
for any Hilbert-Schmidt operator $A$ on $H$. Applying this inequality
to our case, we obtain
$$
\sup_n\left|\log \dett(I+\nabla \xi_n\circ V_n)\right|\leq \frac{d^2}{2}
$$
and this proves the uniform integrability of $(\La_n,n\in \NN)$.
Therefore the sequence  $(\Lambda_n,\,n\in \N)$
converges to $\Lambda_u$ in $L^1(\mu)$ and  we have
$$
E[F\circ T\,\,|\Lambda_u|]=E[F],
$$
for any $F\in C_b(W)$, where $T(w)=w+u(w)$.

To prove  the existence  of the inverse transformation  we begin with
\begin{eqnarray*}
\left|\eta_n-\eta_m\right|_H&\leq &\left|\xi_n\circ V_n-\xi_m\circ
  V_n\right|_H+ \left|\xi_m\circ V_n-\xi_m\circ V_m\right|_H\\
&\leq &\left|\xi_n\circ V_n-\xi_m\circ
  V_n\right|_H+c\left|\eta_n-\eta_m\right|_H,
\end{eqnarray*}
since $c <1$, we obtain:
$$
(1-c)|\eta_n-\eta_m|_H\leq |\xi_n\circ V_n-\xi_m\circ V_n|_H\,.
$$
Consequently, for any $K>0$,
\begin{eqnarray*}
\mu\left\{|\eta_n-\eta_m|_H>K\right\}&\leq & \mu\left\{|\xi_n\circ V_n-\xi_m\circ V_n|_H>(1-c)K\right\}\\
         &=& E\left[|\Lambda_n|\won_{\{|\xi_n-\xi_m|>(1-c)K\}}\right]
     \rightarrow 0,
\end{eqnarray*}
as $n$ and $m$ go to infinity, by the uniform integrability of
$(\Lambda_n;n\in \N)$ and by the convergence in probability of
$(\xi_n; n\in \N)$. As the sequence $(\eta_n;\,n\in\N)$ is bounded in all
$L^p$ spaces, this result implies the existence of an $H$-valued
random variable, say $v$ which is the limit of $(\eta_n;\,n\in\N)$ in
probability.
By uniform integrability, the convergence takes place in $L^p(\mu,H)$ for any
$p>1$ and since the sequence $(\nabla \eta_n;\,n\in \N)$ is bounded in
$L^\infty(\mu,H\otimes H)$, also the convergence takes place in $\DD_{p,1}(H)$
for any $p>1$. Consequently, we have
$$
E[F(w+v(w))\,\,|\Lambda_v|]=E[F],
$$
and
$$
E[F(w+v(w))]=E[F\,\,|\Lambda_u|]\,,
$$
for any $F\in C_b(W)$.

Let us show that $S:W\rightarrow W$, defined by $S(w)=w+v(w)$ is the inverse
of $T$ : let $a>0$ be any number, then
\begin{eqnarray*}
\mu\left\{\|T\circ S(w)-w\|_W>a\right\}&=&\mu\left\{\|T\circ
  S-U_n\circ S\|_W>a/2\right\}\\
         & &  +\mu\left\{\|U_n\circ S-U_n\circ V_n\|_W>a/2\right\}\\
        &=&E\left[|\Lambda_u|\won_{\{\|T-U_n\|_W>a/2\}}\right]\\
        &
        &+\mu\left\{|\xi_n(w+v(w))-\xi_n(w+\eta_n(w))|_H>\frac{a}{2}\right\}\\
       &\leq&E\left[|\Lambda_u|\won_{\{|u-\xi_n|_H>a/2\}}\right]\\
       & &+\mu\left\{|v-\eta_n|_H>\frac{a}{2c}\right\}\rightarrow 0,
\end{eqnarray*}
as $n$ tends to infinity, hence  $\mu$-almost surely $T\circ S(w)=w$.
Moreover
\begin{eqnarray*}
\mu\left\{\|S\circ T(w)-w\|_W>a\right\}&=&\mu\left\{\|S\circ T-S\circ
  U_n\|_W>a/2\right\}\\
      & &  +\mu\left\{\|S\circ U_n-V_n\circ U_n\|_W>a/2\right\}\\
&\leq&\mu\left\{|u-\xi_n|_H>\frac{a(1-c)}{2c}\right\}\\
& &+E\left[|\Lambda_{\eta_n}|\won_{\{|v-\eta_n|_H>a/2\}}\right]\rightarrow 0,
\end{eqnarray*}
by the uniform integrability of $(\Lambda_{\eta_n};\,n\in \N)$, therefore
$\mu$-almost surely, we have $S\circ T(w)=w$.
\qed

\section{Applications}
\markboth{Theorem of Ramer}{Applications}
In the sequel we shall give two applications. The first one consists
of  a very simple case  of the Ramer
formula which is  used in Physics litterature (cf. \cite{D-U} for more
details). The second one concerns the logarithmic Sobolev inequality
for the measures $T^*\mu$ for the shifts $T$ studied in this chapter.
\subsection{Van-Vleck formula}
\markboth{Theorem of Ramer}{Van-Vleck Formula}
\begin{lemma}
\label{lem_inverse}
  Let $K\in L^2(H)$ be a symmetric   Hilbert--Schmidt  operator on $H$
 such that $-1$ does not belong to its spectrum.
         Set
        $T_K(w)=w+\delta  K(w)$, then  $T_K:W\rightarrow W $ is almost surely
        invertible and
$$
T_K^{-1}(w)=w-\delta  [(I+K)^{-1} K](w)\,,
$$
almost surely.
\end{lemma}
\proof
By the properties of the divergence operator (cf. Lemma \ref{lemma2})
\begin{eqnarray*}
\lefteqn{ T_K(w-\delta ((I+K)^{-1} K)(w))}&&\\
 & = & w-\delta ((I+K)^{-1}K)(w)+\delta K(w)
        -\langle \delta ((I+K)^{-1}K)(w), \ K \rangle _H\\
&=&w-\delta ((I+K)^{-1} K)(w)+\delta K(w)-K(\delta ((I+K)^{-1}K)(w))\\
&=&w+\delta K(w)-(I+K)\delta ((I+K)^{-1}K)(w)\\
&=&w,
\end{eqnarray*}
and this proves the lemma.

\qed
\begin{lemma}
\label{normedeltaK}
Let $K$  be a  symmetric  Hilbert--Schmidt operator on $H$. We have
\begin{displaymath}
 \|\delta K\|_H^2=\delta^{(2)}K^2+\trace K^2,
 \end{displaymath}
\end{lemma}
where $\delta^{(2)}$ denotes the second order divergence, i.e.,
$\delta^{(2)}=(\nabla^2)^\star$ with respect to $\mu$.

\proof
 Let  $\{e_i, \ i\ge 0\}$ be the complete, orthonormal basis of $H$
 corresponding to the eigenfunctions of $K$ and denote by
   $\{\alpha_i, \ i\ge 0\}$ its eigenvalues. We can represent $K$ as
$$
 K=\sum_{i=0}^{\infty} \alpha_i e_i\otimes e_i
$$
 and
$$
 K^2=\sum_{i=0}^{\infty} \alpha_i^2 e_i\otimes e_i\ .
$$
Since $\delta K=\sum_i \alpha_i \delta e_i\, e_i$, we have
        \begin{eqnarray*}
                \|\delta K\|_H^2&=&\sum_{i=0}^{\infty} \alpha_i^2 \delta e_i^2 \\
               & =&\sum_{i=0}^{\infty} \alpha_i^2 (\delta e_i^2-1 )
                +\sum_{i=0}^{\infty} \alpha_i^2 \\
               & =&\sum_{i=0}^{\infty} \alpha_i^2 \delta(\delta e_i .e_i)
                +\trace K^2\\
               & =&\delta^{(2)}K^2+\trace K^2 .
        \end{eqnarray*}
\qed
\begin{theorem}
\label{vanvleck}
        Let $K\in L^2(H)$ be a  symmetric   Hilbert--Schmidt  operator such
        that $(I+K)$ is invertible and let  $h_1,\dots, h_n$
        be $n$ linearly independent elements
        of $H$. Denote by $\delta\vec{h}$ the random vector $(\delta
        h_1,\ldots,\delta h_n)$.  Then we have, for any
        $x=(x_1,\ldots,x_n)\in \reals^n$
        \begin{eqnarray*}
\lefteqn{E\left[\exp\left(-\delta^{(2)}\left\{K+\frac{1}{2}K^2\right\}\right)
    \Big| \delta \vec{h}=x\right]}\\
&=&\exp\left(\frac{1}{2}\trace K^2 \right)
           \Big|\dett(I+K)\Big|^{-1}\frac{q_K(x)}{q_0(x)},
        \end{eqnarray*}
where $q_0(x)$ and  $q_K(x)$ denote respectively
the  densities  of the laws  of the Gaussian vectors
$(\delta h_1,\dots, \delta  h_n)$ and
 $$
\left(\delta(I+K)^{-1} h_1,\dots, \delta (I+K)^{-1}h_n\right)\,.
$$
\end{theorem}
\proof
 By the Ramer formula (cf. Theorem \ref{ramer-th}), for any nice
 function   $f$ on $\reals^n$, we have
 \begin{eqnarray*}
 \lefteqn{
   E\left[f(\delta\vec{h})|\dett(I+K)|\exp\left(-\delta^{(2)}K-\frac{1}{2}\|
\delta K\|_H^2\right)\right]}\\
&=&E\left[f(\delta \vec{h})\circ T_K^{-1}(w)\right].
\end{eqnarray*}
Hence
\begin{eqnarray*}
\lefteqn{\int_{\reals^n} E\left[\exp\left(-\delta
    ^{(2)}\left(K+\frac{1}{2}K^2\right)\right)| \delta \vec{h}=x\right]
 f(x)q_0(x)\ dx}\\
&=&\exp\left(\frac{1}{2}\trace K^2 \right)
                \Big|\dett(I+K)\Big|^{-1}
\int_{\reals^n} f(x) q_K(x) dx.
\end{eqnarray*}
\qed
\begin{corollary}
       Suppose that  $A$ is  a symmetric  Hilbert--Schmidt  operator whose
      spectrum is included  in $(-1/2,1/2)$.
   Let $h_1,\dots, h_n$ be $n$ linearly independent  elements
        of $H$ and define   the symmetric, Hilbert-Schmidt  operator $K$  as
       $K=(I+2A)^{1/2}-I$. Then the following identity holds:
\begin{equation}
\label{vanvleckavecA}
E\left[\exp(-\delta^{(2)} A)\ |\ \delta\vec{h}=x\right]
=\frac{1}{\sqrt{\det_2(I+2A)}}\,\frac{q_K(x)}{q_0(x)}\,,
\end{equation}
for any $x=(x_1,\ldots,x_n)\in \reals^n$.
\end{corollary}
\proof
Since the spectrum of $A$ is included in $(-1/2,1/2)$, the operator
$I+2A$ is symmetric and definite.
 It is easy to see that  the operator $K$ is Hilbert-Schmidt. We have
$K+K^2/2=A$, hence  the result follows  by Theorem \ref{vanvleck}.
\qed
\subsection{Logarithmic Sobolev inequality}
\markboth{Theorem of Ramer}{Log-Sobolev Inequality}
\index[sub]{logarithmic Sobolev inequality}
Recall that the logarithmic Sobolev inequality for the Wiener measure  says
\be
\label{log-sob1}
E\left[f^2\log\frac{f^2}{E[f^2]}\right]\leq 2E[|\nabla f|_H^2]\,,
\ee
for any $f\in \DD_{2,1}$.
We can extend this inequality easily to the measures $\nu=T^*\mu$,
where $T=I_W+u$ satisfies the hypothesis of Theorem \ref{ramer-th}
\begin{theorem}
Assume that $\nu$ is a measure given by $\nu=T^*\mu$, where $T=I_W+u$
satisfies the hypothesis of Theorem \ref{ramer-th}, in particular
$\|\nabla u\|\leq c$ almost surely for some $c\in (0,1)$. Then, we
have
\be
\label{log-sob2}
E_\nu\left[f^2\log\frac{f^2}{E[f^2]}\right]\leq
2\left(\frac{c}{1-c}\right)^2E_\nu[|\nabla f|_H^2]
\ee
for any cylindrical Wiener functional $f$, where $E_\nu[\cdot\,]$
represents the expectation with respect to $\nu$.
\end{theorem}
\proof
Let us denote by $S=I_W+v$ the inverse of $T$ whose existence has been
proven in Theorem \ref{ramer-th}. Apply now  the inequality
\ref{log-sob1}  to $f\circ T$:
\beaa
E\left[(f\circ T)^2\log\frac{(f\circ T)^2}{E[(f\circ T)^2]}\right]
&\leq& 2E\left[|\nabla (f\circ T)|_H^2\right]\\
&\leq&2E\left[|\nabla f\circ T|_H^2\|I_H+\nabla u\|^2\right]\\
&=&2E\left[|\nabla f\circ T|_H^2\|I_H+\nabla u\circ S\circ T\|^2\right]\\
&=&2E_\nu\left[|\nabla f|_H^2\|I_H+\nabla u\circ S\|^2\right]\\
&=&2E_\nu\left[|\nabla f|_H^2\|(I_H+\nabla v)^{-1}\|^2\right]\\
&\leq&2\left(\frac{c}{1-c}\right)^2E_\nu\left[|\nabla f|_H^2\right]\,
\eeaa
and this completes the proof.
\qed

We have also the following:
\index[sub]{closable}
\begin{theorem}
The operator $\nabla$ is closable in $L^p(\nu)$ for any $p>1$.
\end{theorem}
\proof
Assume that $(f_n,n\in \NN)$ is a sequence of cylindrical Wiener
functionals, converging in $L^p(\nu)$ to zero, and assume also that
$(\nabla f_n,n\in \NN)$ is Cauchy in $L^p(\nu,H)$, denote its limit by
$\xi$. Then, by
definition, $(f_n\circ T,n\in \NN)$ converges to zero in $L^p(\mu)$,
hence $(\nabla(f_n\circ T),n\in \NN)$ converges to zero in
$\DD_{p,-1}(H)$. Moreover
$$
\nabla (f_n\circ T)=(I_H+\nabla u)^*\nabla f_n\circ T\,,
$$
hence for any cylindrical  $\eta\in \DD(H)$, we have
\beaa
\lim_nE[(\nabla(f_n\circ T),\eta)_H]
&=&\lim_nE[(\nabla f_n\circ T,(I_H+\nabla u)\eta)_H]\\
&=&E[(\xi\circ T,(I_H+\nabla u)\eta)_H]\\
&=&0\,.
\eeaa
Since $T$ is invertible, the sigma algebra generated by $T$ is equal
to the Borel sigma algebra of $W$ upto the negligeable
sets. Consequently, we have
$$
(I_H+\nabla u)^*\xi\circ T=0
$$
$\mu$-almost surely. Since $I_H+\nabla u$ is almost surely invertible,
$\mu$-almost surely we have $\xi\circ T=0$  and this amounts up to saying
$\xi=0$ $\nu$-almost surely.
\qed

\section*{Notes and suggested reading}
{\footnotesize{The Ramer theorem has been proved, with some stronger
    hypothesis (Fr\'echet regularity of $u$) in \cite{R}, later some
    of its hypothesis  have been relaxed in \cite{Kus}. The version
    given here has been proved in \cite{UZ5}. We refer the reader to
    \cite{UZ7} for its further extensions and applications to the
    degree theory of Wiener maps (cf. \cite{UZ6} also). The Van-Vleck
    formula is well-known in Physics, however the general approach
    that we have used here as well as the logarithmic Sobolev
    inequalities with these new measures are  original.}}

%

\chapter{Convexity on Wiener space}
\label{ch.convex}
\markboth{Convexity}{Convexity}
\section*{Introduction}
On an infinite dimensional vector space $W$  the notion of convex or concave
function is well-known. Assume now that this space is equipped with a
probability measure. Suppose that there are two measurable functions  on
this vector space, say $F$ and $G$ such that $F=G$ almost surely. If
$F$ is a  convex function, then from the probabilistic point of view,
we would like to say  that $G$ is also convex. However this is false;  since
in general the underlying probability measure is not (quasi)
invariant under the translations by the elements of the vector
space. If $W$ contains a dense  subspace $H$ such that $w\to w+h$
($h\in H$) induces a measure which is equivalent to the initial
measure or absolutely continuous with respect to it, then   we can
define a  notion of ``$H$--convexity'' or ``$H$--concavity  in the
direction  of $H$ for the equivalence classes of real random
variables.  Hence these notions will be  particularly useful
for  the probabilistic calculations.

The notion of $H$-convexity  has
been used  in \cite{UZ7} to study the
absolute continuity of the image of the Wiener measure under the
monotone shifts. In this chapter we study further  properties of
such functions  and some additional ones  in the frame of an abstract
Wiener space, namely  $H$-convex, $H$-concave, log $H$-concave and log
$H$-convex Wiener functions,
where $H$ denotes the associated  Cameron-Martin space. In particular
we extend some finite dimensional results of  \cite{PRE71} and
\cite{BLI}  to this setting and prove that some finite dimensional
convexity-concavity inequalities have their counterparts in infinite
dimensions.
\section{Preliminaries}
\markboth{Convexity}{Preliminaries}
In the sequel $(W,H,\mu)$ denotes an abstract Wiener space, i.e., $H$
is a separable Hilbert space, called the Cameron-Martin space.
It is identified with its continuous dual.  $W$
is a Banach or a Fr\'echet space into which $H$ is injected
continuously and densely. $\mu$ is the standard  cylindrical  Gaussian
measure on $H$ which is concentrated in $W$ as a Radon probability measure.

In the sequel we shall use the notion of second quantization of
bounded operators on $H$; although this is a well-known subject, we
give a brief outline below for the reader's convenience
(cf. \cite{BAD}, \cite{P-F}, \cite{Si}).
Assume that $A:H\to H$ is a bounded, linear operator, then it has a
unique, $\mu$-measurable (i.e., measurable with respect to the
$\mu$-completion of $\calB(W)$) extension, denoted by $\tilde{A}$,  as
a linear map on $W$ (cf. \cite{BAD,P-F}). Assume in particular
that $\|A\|\leq 1$ and define $S=(I_H-A^*A)^{1/2}$,
$T=(I_H-AA^*)^{1/2}$ and $U:H\times H\to H\times H$ as
$U(h,k)=(Ah+Tk,-Sh+A^*k)$. $U$ is then  a  unitary operator on $H\times H$,
hence its $\mu\times\mu$-measurable linear extension to $W\times W$
preserves the Wiener measure $\mu\times\mu$ (this is called the
rotation associated
to $U$, cf. \cite{UZ7}, Chapter VIII). Using this observation, one can define
the second quantization of $A$ via the generalized Mehler formula as
$$
\Gamma(A)f(w)=\int_Wf(\tilde{A^*}w+\tilde{S}y)\mu(dy)\,,
$$
which happens to be a Markovian contraction on $L^p(\mu)$ for any
$p\geq 1$. $\Gamma (A)$ can be calculated explicitly for the Wick
exponentials as
$$
\Gamma(A)\exp\left\{\delta h-1/2|h|_H^2\right\}=\exp\left\{\delta A
  h-1/2|Ah|_H^2\right\}\,\,(h\in H)\,.
$$
This identity implies that $\Gamma (AB)=\Gamma(A)\Gamma(B)$ and that
for any sequence $(A_n,n\in \NN)$ of operators whose norms are bounded
by one, $\Gamma(A_n)$ converges strongly to $\Gamma(A)$ if
$\lim_nA_n=A$ in the strong operator topology.
 A particular case of interest is  when we take $A=e^{-t}I_H$, then
 $\Gamma(e^{-t}I_H)$ equals to  the Ornstein-Uhlenbeck semigroup
 $P_t$. Also  if
 $\pi$ is the  orthogonal projection of $H$ onto a closed vector subspace
 $K$, then $\Gamma(\pi)$ is the conditional expectation with respect
 to the sigma field generated  by $\{\delta k,\,k\in K\}$.



\section{$H$-convexity and its properties}
\markboth{Convexity}{$H$-convexity}
Let us give  the notion of $H$-convexity on the  Wiener space $W$:
\begin{definition}
Let  $F:W\to \reals\cup\{\infty\}$ be a measurable function.
It is called $H$-convex if for any
$h,k\in H$, $\al\in [0,1]$
\begin{equation}
\label{con-def}
F(w+\al h+(1-\al)k)\leq \al F(w+h)+(1-\al)F(w+k)
\end{equation}
almost surely.
\end{definition}
\remarks
\begin{itemize}
\item This definition is more general than the one given in
  \cite{UZ0,UZ7} since $F$ may be   infinite on a set of  positive
  measure.
\item Note that the negligeable  set on which the relation
  (\ref{con-def}) fails may depend on the choice of $h,k$ and of $\al$.
\item If $G:W\to\reals\cup\{\infty\}$ is a measurable  convex
function, then it is  necessarily $H$-convex.
\item To conclude the $H$-convexity, it suffices to verify the
  relation (\ref{con-def}) for $k=-h$ and $\alpha=1/2$.
\end{itemize}

\noindent
The following properties of $H$-convex Wiener functionals have been
proved in  \cite{UZ0,USZ98,UZ7}:
\begin{theorem}
\label{char-thm}
\begin{enumerate}
\item If $(F_n,n \in \NN)$ is a sequence  of $H$-convex
functionals converging in probability, then the limit is also
  $H$-convex.
\item If $F\in L^p(\mu)$ ($p>1$) is $H$-convex if and only if
  $\nabla^2 F$ is positive and symmetric  Hilbert-Schmidt operator
  valued distribution on $W$.
\item If $F\in L^1(\mu)$ is $H$-convex, then $P_tF$ is also $H$-convex
  for any $t\geq 0$, where $P_t$ is the Ornstein-Uhlenbeck semi-group
  on $W$.
\end{enumerate}
\end{theorem}
\noindent
The following result is immediate from Theorem \ref{char-thm} :
\begin{corollary}
\label{char-cor}
$F\in \cup_{p>1}L^p(\mu)$ is $H$-convex if and only if
$$
E\left[\varphi\,\left(\nabla^2F(w),h\otimes h\right)_2\right]\geq 0
$$
for any $h\in H$ and $\varphi\in \DD_+$, where $(\cdot\,,\,\cdot)_2$
denotes the scalar
product for the Hilbert-Schmidt operators on $H$ .
\end{corollary}
We have also
\begin{corollary}
If $F\in L^p(\mu)$, $p>1$, is $H$-convex and if $E[\nabla^2F]=0$, then
$F$ is of the form
$$
F=E[F]+\delta\left(E[\nabla F]\right)\,.
$$
\end{corollary}
\proof
Let $(P_t,t\geq 0)$ denote the Ornstein-Uhlenbeck semigroup, then $P_tF$ is
again $H$-convex and Sobolev differentiable.  Moreover $\nabla^2
P_tF=e^{-2t}P_t\nabla^2F$. Hence $E[\nabla^2P_tF]=0$, and  the
positivity of $\nabla^2P_tF$ implies that $\nabla^2 P_tF=0$ almost
surely, hence $\nabla^2 F=0$. This implies that $F$ is in the first
two Wiener chaos.
\qed

\remark It may be worth-while  to note that the random variable which
represents  the share price of the Black and
Scholes model in financial mathematics  is $H$-convex.
\newline

We shall need also the  concept of $\calC$-convex functionals:
\begin{definition}
\label{C-defn}
 Let $(e_i,i\in \NN)\subset W^*$ be any complete, orthonormal basis of
 $H$. For $w\in W$, define  $w_n=\sum_{i=1}^n\delta e_i(w)e_i$ and
 $w_n^\perp=w-w_n$, then a  Wiener functional $f:W\to \R$ is called
 $\calC$-convex if, for any such basis $(e_i,i\in \NN)$, for almost
 all $w_n^\perp$,  the partial map
$$
w_n\to f(w_n^\perp+w_n)
$$
has a modification which is convex on the space
${\rm span}\{e_1,\ldots,e_n\}\simeq  \R^n$.
\end{definition}
\remark It follows from Corollary \ref{char-cor}  that, if $f$ is
$H$-convex and   in
some $L^p(\mu)$ $(p>1)$, then it is $\calC$-convex. We shall prove
that this is also true without any integrability hypothesis.

We begin with the following lemma whose proof is obvious:

\begin{lemma}
\label{L1}
If $f$ is $\calC$-convex then it is $H$-convex.
\end{lemma}
\noindent
In order to  prove the validity of the  converse of Lemma
\ref{L1}  we need some technical results from  the
harmonic analysis on finite dimensional Euclidean spaces that we shall
state as separate lemmas:
\begin{lemma}
\label{L2}
Let $B\in {\calB}(\reals^n)$ be a set of positive Lebesgue
measure. Then $B+B$ contains a non-empty open set.
\end{lemma}
\proof
Let $\phi(x)=1_B\star 1_B(x)$, where ``$\star$''  denotes the
convolution of functions  with respect to the Lebesgue measure. Then
$\phi$ is a non-negative,  continuous function, hence the set $O=\{x\in
\reals^n:\,\phi(x)>0\}$ is an open set. Since $B$ has positive
measure, $\phi$ can not be identically zero, hence $O$ is
non-empty. Besides, if $x\in O$, then the set of $y\in \reals^n$ such
that $y\in B$ and $x-y\in B$ has positive Lebesgue measure, otherwise
$\phi(x)$ would have been null. Consequently $O\subset B+B$.
\qed

The following lemma gives a more precise statement than Lemma \ref{L2}:
\begin{lemma}
\label{L3}
Let $B\in {\calB}(\reals^n)$ be a set of positive Lebesgue
measure and  assume  that  $A\subset \reals^n\times \reals^n$ with
$B\times B=A$ almost surely with respect to the Lebesgue measure of
$\reals^n\times \reals^n$. Then the set $\{x+y:\,(x,y)\in A\}$
contains almost surely an open subset of $\reals^n$.
\end{lemma}
\proof
It follows from an obvious change of variables that
$$
1_A(y,x-y)=1_B(y)1_B(x-y)
$$
almost surely, hence
$$
\int_{\reals^n}1_A(y,x-y)dy=\phi(x)
$$
almost surely, where $\phi(x)=1_B\star 1_B(x)$. Consequently, for
almost all $x\in \reals^n$ such that $\phi(x)>0$, one has $(y,x-y)\in
A$, this means that
$$
\{x\in \reals^n:\,\phi(x)>0\}\subset
\{u+v:\,(u,v)\in A\}
$$
almost surely.
\qed

The following lemma is particularly important for the sequel:
\begin{lemma}
\label{L4}
Let $f:\reals^n\to \reals_+\cup\{\infty\}$ be a Borel function which
is finite on a set of positive Lebesgue measure. Assume that, for any
$u\in \reals^n$,
\begin{equation}
\label{C1}
f(x)\leq \frac{1}{2}[f(x+u)+f(x-u)]
\end{equation}
$dx$-almost surely (the negligeable set on which the inequality
(\ref{C1}) fails may depend on $u$). Then there exists a non-empty,
open convex subset $U$ of $\reals^n$ such that $f$ is locally essentially
bounded on $U$. Moreover let $D$ be the set consisting of
 $x\in \reals^n$ such  that  any neighbourhood
 of $x\in D$ contains a Borel  set of positive Lebesgue measure on
 which $f$ is finite, then $D\subset \overline{U}$, in particular
 $f=\infty$ almost surely  on the complement of $\overline{U}$.
\end{lemma}
\proof
From the theorem of Fubini, the inequality (\ref{C1}) implies that
\begin{equation}
\label{C2}
2f\left(\frac{x+y}{2}\right)\leq f(x)+f(y)
\end{equation}
$dx\times dy$-almost surely. Let $B\in {\calB}(\reals^n)$ be a set of
positive Lebesgue measure on which $f$ is  bounded by some constant
$M>0$. Then from Lemma
\ref{L2}, $B+B$ contains an open set $O$. Let $A$ be the set
consisting of the elements of $B\times B$ for which the inequality
(\ref{C2}) holds. Then $A=B\times B$ almost surely, hence from Lemma
\ref{L3}, the set $\Gamma=\{x+y:\,(x,y)\in A\}$ contains almost surely
the open set $O$. Hence for almost all $z\in \frac{1}{2}O$, $2z$
belongs to the set $\Gamma$, consequently $z=\frac{1}{2}(x+y)$, with
$(x,y)\in A$. This implies, from (\ref{C2}),  that $f(z)\leq
M$. Consequently $f$ is essentially bounded on the open set
$\frac{1}{2}O$.

Let now $U$ be set of points  which have neighbourhoods on which $f$ is
essentially bounded. Clearly $U$ is open and non-empty by what we have
shown above. Let $S$ and $T$ be two balls of radius $\rho$, on which
$f$ is bounded by some $M>0$. Assume that they are centered at the
points $a$ and $b$ respectively. Let $u=\frac{1}{2}(b-a)$, then for
almost all $x\in \frac{1}{2}(S+T)$, $x+u\in T$ and $x-u\in S$, hence,
from the inequality (\ref{C1})
$f(x)\leq M$, which shows that $f$ is essentially bounded on the set
$\frac{1}{2}(S+T)$ and this proves the convexity of $U$.

To prove the
last  claim, let $x$ be any element of $D$ and let $V$ be any
neighbourhood of $x$; without loss of generality, we may assume that
$V$ is convex. Then there exists a Borel set $B\subset V$ of positive
measure on which $f$ is bounded, hence from the first part of the
proof, there exists an open  neighbourhood $O\subset B+B$ such that $f$ is
essentially bounded on $\frac{1}{2}O\subset \frac{1}{2}(V+V)\subset
V$, hence $\frac{1}{2}O\subset U$. Consequently $V\cap U\neq
\emptyset$, and this implies that $x$ is in the closure of $U$,
i.e. $D\subset \overline{U}$. The fact that $f=\infty$ almost surely
on the complement of $\overline{U}$ is obvious from the definition of $D$.
\qed

\begin{theorem}
\label{con-char}
Let $g:\reals^n\to \reals\cup\{\infty\}$ be a measurable mapping such
that, for almost all $u\in \reals^n$,
\begin{equation}
\label{CC}
g(u+\alpha x+\beta y)\leq \alpha g(u+x)+\beta g(u+y)
\end{equation}
for any $\alpha,\beta \in [0,1]$ with $\alpha+\beta =1$ and for any
$x,y\in \R^n$, where the
negligeable set on which the relation (\ref{CC}) fails may depend on
the choice of $x,y$ and of $\alpha$ . Then $g$ has
a modification $g'$ which is a convex function.
\end{theorem}
\proof
Assume first that $g$ is positive, then with the notations of Lemma
\ref{L4}, define $g'=g$ on the open, convex set $U$ and as $g'=\infty$
on $U^c$. From the relation (\ref{CC}), $g'$ is a distribution on $U$
whose second derivative is positive, hence it is convex on $U$, hence
it is convex on the whole space $\reals^n$. Moreover we have $\{g'\neq
g\}\subset \partial U$ and $\partial U$ has zero Lebesgue measure,
consequently $g=g'$ almost surely. For general $g$, define
$f_\eps=e^{\eps g}$ ($\eps>0$), then, from what is proven above,
  $f_\eps$ has a modification $f'_\eps$ which is convex (with the same
  fixed open and convex set $U$), hence $\lim\sup_{\eps\to
    0}\frac{f'_\eps-1}{\eps}=g'$ is also convex and $g=g'$ almost
  surely.
\qed

\begin{theorem}
\label{H->C}
A Wiener functional $F:W\to \reals\cup\{\infty\}$ is $H$-convex if and
only if it is $\calC$-convex.
\end{theorem}
\proof
We have already proven the sufficiency. To prove the necessity, with
the notations of Definition \ref{C-defn}, $H$-convexity implies that
$h\to F(w_n^\perp+w_n+h)$ satisfies the hypothesis of Theorem
\ref{con-char} when $h$ runs in any $n$-dimensional Euclidean subspace
of $H$, hence the partial mapping $w_n\to F(w_n^\perp+w_n)$ has a
modification which is convex on the vector space spanned by
$\{e_1,\ldots,e_n\}$.
\qed

\section{Log $H$-concave and $\calC$-$\log$ concave  Wiener functionals}
\markboth{Convexity}{log $H$-concave}

\begin{definition}
Let $F$ be a  measurable mapping  from $W$ into $\reals_+$ with
$\mu\{F>0\}>0$.
\begin{enumerate}
\item $F$ is  called log $H$-concave, if for any $h,\,k\in H$,
  $\alpha\in [0,1]$, one has
\begin{equation}
\label{log-cv-def}
F\bigl(w+\alpha h+(1-\alpha)k\bigr)\geq F(w+ h)^\alpha\,  F(w+k)^{1-\alpha}
\end{equation}
almost surely, where the negligeable set on which the relation
(\ref{log-cv-def})
fails may depend on $h,\,k$ and on $\alpha$.
\item We shall say that $F$ is $\calC$-log concave, if for any
  complete, orthonormal basis  $(e_i,i\in \NN)\subset W^*$ of $H$, the
  partial map $w_n\to F(w_n^\perp+w_n)$ is log-concave (cf. Definition
  \ref{C-defn} for the notation),  up to a
  modification,  on ${\rm span}\{e_1,\ldots,e_n\}\simeq \R^n$.
\end{enumerate}
\end{definition}

Let us  remark immediately that if $F=G$ almost surely then $G$ is also log
$H$-concave. Moreover, any limit in probability of log $H$-concave random
variables is again log $H$-concave. We shall prove below some less immediate
properties. Let us begin with the following observation which is a
direct consequence of Theorem \ref{H->C}:

\remark
\label{H-concave=C-concave}
$F$ is log $H$-concave if and only if $-\log F$ is $H$-convex
 (which may be infinity with a positive probability), hence if and
 only if $F$ is $\calC$-log concave.

\begin{theorem}
\label{prekopa1-the}
Suppose that $(W_i,H_i,\mu_i)$, $i=1,2$, are two abstract Wiener
spaces. Consider $(W_1\times W_2, H_1\times H_1,\mu_1\times \mu_2)$ as an
abstract Wiener space. Assume that $F:W_1\times W_2\to \reals_+$ is log
$H_1\times H_2$-concave. Then the map
$$
w_2\to \int_{W_1}F(w_1,w_2)\,d\mu_1(w_1)
$$
is log $H_2$-concave.
\end{theorem}
\proof
If $F$ is log $H\times H$-concave, so is also $F\wedge c$ ($c\in
\reals_+$), hence  we may suppose without
loss of generality that $F$ is bounded.
Let $(e_i,i\in \NN)$ be a complete, orthonormal basis in $H_2$. It suffices to
prove that
$$
E_1[F](w_2+\alpha h+\beta l)\geq\left(E_1[F](w_2+h)\right)^\alpha
  \left(E_1[F](w_2+l)\right)^\beta
$$
almost surely, for any $h,\,l\in {\mbox{span}}\{e_1,\ldots,e_k\}$,
$\alpha,\beta\in [0,1]$ with $\alpha+\beta=1$, where $E_1$ denotes the
expectation with respect to $\mu_1$. Let $(P_n,n\in \NN)$ be a sequence of
orthogonal projections of finite rank  on $H_1$ increasing to the
identity map of it. Denote by
$\mu_1^n$ the image of $\mu_1$ under the map $w_1\to {\tilde{P}}_n w_1$ and by
$\mu_1^{n\perp}$ the image of $\mu_1$ under $w_1\to w_1-{\tilde{P}}_n w_1$. We
have, from the martingale convergence theorem,
$$
\int_{W_1}F(w_1,w_2)\,d\mu_1(w_1)=\lim_n\int
F(w_1^{n\perp}+w_1^n,w_2)\,d\mu_1^n(w^n_1)
$$
almost surely. Let $(Q_n,n\in \NN)$ be a sequence of orthogonal
projections of finite rank  on
$H_2$ increasing to the identity, corresponding to the basis $(e_n,n\in
\NN)$. Let $w_2^k={\tilde{Q}}_kw_2$ and $w_2^{k\perp}=w_2-w_2^k$. Write
\beaa
F(w_1,w_2)&=&F(w_1^{n\perp}+w_1^n,w_2^k+w_2^{k\perp})\\
       &=&F_{w_1^{n\perp},w_2^{k\perp}}(w_1^n,w_2^k)\,.
\eeaa
From the hypothesis
$$
(w_1^n,w_2^k)\to F_{w_1^{n\perp},w_2^{k\perp}}(w_1^n,w_2^k)
$$
has a  log concave modification  on the $(n+k)$-dimensional Euclidean
space. From the theorem of Pr\'ekopa (cf. \cite{PRE71}), it follows that
$$
w_2^k\to \int F_{w_1^{n\perp},w_2^{k\perp}}(w_1^n,w_2^k)\,d\mu_1^n(w_1^n)
$$
is log concave on $\R^k$  for any $k\in \NN$ (upto a modification), hence
$$
w_2\to \int F(w_1^{n\perp}+w_1^n,w_2)\, d\mu(w_1^n)
$$
is log $H_2$-concave for any $n\in \NN$, then the proof follows by passing to
the limit with respect to $n$.
\qed

\begin{theorem}
\label{Gamma-thm}
Let $A:H\to H$ be a linear operator with $\|A\|\leq 1$, denote by
$\Gamma(A)$ its second  quantization as explained in the
preliminaries. If $F:W\to \reals_+$ is a log $H$-concave Wiener
functional, then  $\Gamma(A)F$  is also log $H$-concave.
\end{theorem}
\proof
Replacing $F$ by $F\wedge c=\min(F,c),\,c>0$, we may suppose that $F$
is bounded. It is
easy  to see that the mapping
$$
(w,y)\to F(\tilde{A^*}w+\tilde{S}y)
$$
is log $H\times H$-concave on $W\times W$. In fact, for any
$\alpha+\beta=1$, $h,k,u,v\in H$, one has
\bea
\label{formul}
\lefteqn{F(\tilde{A^*}w+\tilde{S}y+\alpha(A^*h+Sk)+\beta(A^*u+Sv))}\\
&&\geq F(\tilde{A^*}w+\tilde{S}y+A^*h+Sk)^\alpha\,
F(\tilde{A^*}w+\tilde{S}y+A^*u+Sv)^\beta\,,\nonumber
\eea
$d\mu\times d\mu$-almost surely. Let us recall that, since the image
of $\mu\times \mu$ under the map $(w,y)\to \tilde{A^*}w+\tilde{S}y$
is $\mu$, the terms in the inequality (\ref{formul}) are defined without
ambiguity. Hence
$$
\Gamma(A)F(w)=\int_W F(\tilde{A^*}w+\tilde{S}y)\mu(dy)
$$
is log $H$-concave on $W$ from Theorem \ref{prekopa1-the}.
\qed

\begin{corollary}
\label{con-cor}
Let $F:W\to \reals_+$ be a log $H$-concave functional. Assume that $K$
is any closed vector subspace of $H$ and denote by $V(K)$ the sigma algebra
generated by $\{\delta k,\,k\in K\}$. Then the  conditional
expectation of $F$ with respect to $V(K)$, i.e., $E[F|V(K)]$  is again
log $H$-concave.
\end{corollary}
\proof
The proof follows from Theorem \ref{Gamma-thm} as soon as we remark
that $\Gamma(\pi_K)F=E[F|V(K)]$, where $\pi_K$ denotes the orthogonal
projection associated to $K$.
\qed

\begin{corollary}
\label{semi-con}
Let $F$  be log $H$-concave. If $P_t$ denotes the
Ornstein-Uhlenbeck semigroup on $W$, then $w\to P_tF(w)$ is log $H$-concave.
\end{corollary}
\proof
Since $P_t=\Gamma(e^{-t}I_H)$, the proof follows from Theorem
\ref{Gamma-thm}.
\qed

Here is an important application of these results:
\begin{theorem}
\label{version-thm}
Assume that $F:W\to \reals\cup \{\infty\}$ is an $H$-convex Wiener
functional, then $F$ has a modification $F'$ which is a Borel
measurable  convex
function on $W$.  Any log $H$-concave
functional $G$ has a modification $G'$ which is Borel measurable and
log-concave on $W$.
\end{theorem}
\proof
Assume first that $F$ is positive, let $G=\exp-F$, then $G$ is a
positive, bounded $\calC$-log concave function. Define $G_n$ as
$$
G_n=E[P_{1/n}G|V_n]\,,
$$
where $V_n$ is the sigma algebra generated by $\{\delta
e_1,\ldots,\delta e_n\}$, and $(e_i,i\in \NN)\subset W^*$ is a
complete orthonormal basis of $H$. Since
$P_{1/n}E[G|V_n]=E[P_{1/n}G|V_n]$, the positivity improving property
of the Ornstein-Uhlenbeck semigroup implies that $G_n$ is almost
surely strictly positive  (even quasi-surely). As we have attained the
finite dimensional case, $G_n$ has a modification
$G'_n$ which is continuous on $W$ and, from Corollary \ref{con-cor}
and Corollary \ref{semi-con}, it satisfies
\begin{equation}
\label{f1}
G'_n(w+ah+bk)\geq G'_n(w+h)^a G'_n(w+k)^b
\end{equation}
almost surely, for any $h,k\in H$ and $a+b=1$. The continuity of
$G'_n$ implies that  the relation (\ref{f1}) holds for any $h,k\in H$,
$w\in W$ and $a\in [0,1]$. Hence $G'_n$ is log-concave on $W$ and
this implies that $-\log G'_n$ is convex on $W$. Define
$F'=\lim\sup_n(-\log G'_n)$, then $F'$ is convex and Borel measurable
on $W$ and $F=F'$ almost surely.

For general $F$, define $f_\eps=e^{\eps F}$, then from above, there
exists a modification of $f_\eps$, say  $f'_\eps$ which is convex and
Borel measurable  on
$W$. To complete the proof it suffices to define $F'$ as
$$
F'=\lim\sup_{\eps\to 0}\frac{f'_\eps-1}{\eps}\,.
$$
The rest is   now obvious.
\qed

Under the light of Theorem \ref{version-thm}, the following definition
is natural:
\begin{definition}
A Wiener functional $F:W\to \reals\cup\{\infty\}$ will be called
almost surely convex if it has a modification $F'$ which is convex and
Borel measurable on $W$. Similarly, a non-negative functional $G$ will
be called almost surely log-concave if it has a modification $G'$
which is log-concave on $W$.
\end{definition}

The following proposition summarizes the main results of this section:
\begin{theorem}
\label{equiv-prop}
Assume that $F:W\to \reals\cup\{\infty\} $ is a Wiener functional such that
$$
\mu\{F<\infty\}>0\,.
$$
Then the following are equivalent:
\begin{enumerate}
\item $F$ is $H$-convex,
\item $F$ is $\calC$-convex,
\item $F$ is almost surely convex.
\end{enumerate}
Similarly, for $G:W\to \reals_+$, with $\mu\{G>0\}>0$, the following
properties are equivalent:
\begin{enumerate}
\item $G$ is log $H$-concave,
\item $G$ is log $\calC$-concave,
\item $G$ is almost surely log-concave.
\end{enumerate}
\end{theorem}

The notion of a convex set can be extended as
\begin{definition}
Any measurable subset $A$ of $W$ will be called $H$-convex  if its indicator
function $1_A$ is log $H$-concave.
\end{definition}
\remark
Evidently any measurable  convex subset of $W$ is
$H$-convex. Moreover, if $A=A'$ almost surely and if $A$ is
$H$-convex, then $A'$ is also $H$-convex.

\remark
If $\phi$ is an $H$-convex Wiener functional, then
the set
$$
\{w\in W: \phi(w)\leq t\}
$$
is $H$-convex for any $t\in
\reals$.

We have the following result about the characterization of the
$H$-convex sets:
\begin{theorem}
Assume that $A$ is an $H$-convex set, then there exists a convex set
$A'$, which is Borel measurable such that $A=A'$ almost surely.
\end{theorem}
\proof
Since, by definition, $1_A$ is a log $H$-concave Wiener functional,
from Theorem \ref{version-thm}, there exists a log-concave Wiener
functional $f_A$ such that $f_A=1_A$ almost surely. It suffices to
define $A'$ as the set
$$
A'=\{w\in W:\,f_A(w)\geq 1\}\,.
$$
\qed
\paragraph{Example:}
Assume that $A$ is an $H$-convex subset of $W$ of positive
measure. Define $p_A$ as
$$
p_A(w)=\inf\left(|h|_H:h\in (A-w)\cap H\right)\,.
$$
Then $p_A$ is $H$-convex, hence almost surely convex (and $H$-Lipschitz
c.f. \cite{UZ7}). Moreover, the  $\{w:\,p_A(w)\leq \alpha\}$ is an
$H$-convex set for any $\alpha\in \reals_+$.


\section{Extensions and some applications}
\markboth{Convexity}{Extensions}
\begin{definition}
Let $(e_i,i\in \NN)$ be any  complete orthonormal  basis of  $H$. We shall
denote, as before, by $w_n=\sum_{i=1}^n\delta e_i(w)\,e_i$ and
$w_n^\perp=w-w_n$.
Assume now that $F:W\to \reals\cup\{\infty\}$ is a measurable mapping with
$\mu\{F<\infty\}>0$.
\begin{enumerate}
\item We say that it is $a$-convex ($a\in \reals$), if the partial map
$$
w_n\to \frac{a}{2}|w_n|^2+F(w_n^\perp+w_n)
$$
is almost surely convex for any $n\geq 1$, where $|w_n|$ is the
Euclidean norm of $w_n$.
\item We call $G$ $a$-log-concave if
$$
w_n\to \exp\left\{-\frac{a}{2}|w_n|^2\right\}G(w_n^\perp+w_n)
$$
\end{enumerate}
is almost surely log-concave for any $n\in \NN$.
\end{definition}
\remark
$G$ is $a$-log-concave if and only if $-\log G$ is $a$-convex.
\newline

\noindent
The following theorem gives  a practical method to verify
$a$-convexity or log-concavity:
\begin{theorem}
\label{a-con-thm}
Let $F:W\to \reals\cup\{\infty\}$ be a measurable map such that
$\mu\{F<\infty\}>0$. Define the map $F_a$ on $H\times W$ as
$$
F_a(h,w+h)=\frac{a}{2}|h|_H^2+F(w+h)\,.
$$
Then $F$ is $a$-convex if and only if, for any $h,k\,\in H$ and
$\alpha,\beta\in[0,1]$ with $\alpha+\beta=1$, one has
\begin{equation}
\label{H-cond}
F_a(\alpha h+\beta k,w+\alpha h+\beta k)\leq \alpha\, F_a(h,w+h)+\beta\,
F_a(k,w+k)
\end{equation}
$\mu$-almost surely, where the negligeable set on which the inequality
(\ref{H-cond}) fails may depend on the choice of $h,k$ and of
$\alpha$.

Similarly a measurable mapping $G:W\to \reals_+$ is
$a$-log-concave if and only if the map defined by
$$
G_a(h,w+h)=\exp\left\{-\frac{a}{2}|h|_H^2\right\}G(w+h)
$$
satisfies the inequality
\begin{equation}
\label{a-cond}
G_a(\alpha h+\beta k,w+\alpha h+\beta k)\geq G_a( h,w+h)^\alpha
G_a( k,w+k)^\beta\,,
\end{equation}
$\mu$-almost surely, where the negligeable set on which the inequality
(\ref{a-cond}) fails may depend on the choice of $h,k$ and of
$\alpha$.
\end{theorem}
\proof
Let us denote by $h_n$ its projection on
the vector space spanned by $\{e_1,\ldots,e_n\}$, i.e.  $h_n=\sum_{i\leq
  n}(h,e_i)_H e_i$. Then, from Theorem \ref{equiv-prop}, $F$ is
$a$-convex if and only if the map
$$
h_n\rightarrow
\frac{a}{2}\left[|w_n|^2+2(w_n,h_n)+|h_n|^2\right]+F(w+h_n)
$$
satisfies a convexity inequality like (\ref{H-cond}). Besides the
term $|w_n|^2$ being kept constant in this operation,   it can be
removed from the both sides of the inequality. Similarly, since
$h_n\to (w_n,h_n)$ is being affine, it also cancels from the both
sides of this inequality. Hence $a$-convexity is equivalent to
$$
F_a(\alpha h_n+\beta k_n,w+\alpha h_n+\beta k_n)\leq
 \alpha\, F_a(h_n,w+h_n)+\beta\,F_a(k_n,w+k_n)
$$
where $k_n$ is defined as $h_n$ from a $k\in H$.

The second part of the theorem is obvious since $G$ is $a$-log-concave
if and only if $-\log G$ is $a$-convex.
\qed
\begin{corollary}
\label{car-cor}
\begin{enumerate}
\item Let $\hat{L}^0(\mu)$ be the space of the $\mu$-equivalence
  classes of   $\R\cup\{\infty\}$-valued
  random variables regarded as a topological semi-group under addition
  and convergence in probability. Then $F\in \hat{L}^0(\mu)$ is
  $\beta$-convex if and only if the mapping
$$
h\to \frac{\beta}{2}|h|_H^2+F(w+h)
$$
is a convex and continuous mapping from $H$ into $\hat{L}^0(\mu)$.
\item $F\in L^p(\mu),\,p>1$ is $\beta$-convex if and only if
$$
E\left[\left((\beta I_H+\nabla^2F)h,h\right)_H\,\phi\right]\geq 0
$$
for any $\phi\in \DD$ positive and $h\in H$, where $\nabla^2F$ is to
be understood in the sense of the distributions $\DD'$.
\end{enumerate}
\end{corollary}

\paragraph{Example:}
Note for  instance that $\sin \delta h$ with $|h|_H=1$, is a 1-convex
random variale and that $\exp(\sin \delta h)$ is $1$-log-concave.

The following result is a direct consequence of Prekopa's theorem:
\begin{proposition}
\label{c-exp-prop}
Let $G$ be an $a$-log concave Wiener functional, $a\in [0,1]$,  and
assume that $V$ is
any sigma algebra generated by the elements of the first Wiener
chaos. Then $E[G|V]$ is again $a$-log-concave.
\end{proposition}
\proof
From Corollary \ref{car-cor}, it suffices to prove the case $V$ is
generated by $\{\delta e_1,\ldots,\delta e_k\}$, where $(e_n,n\in
\NN)$ is an orthonormal basis of $H$. Let
\beaa
w_k&=&\sum_{i\leq k}\delta e_i(w)e_i\\
z_k&=&w-w_k\\
z_{k,n}&=&\sum_{i=k+1}^{k+n}\delta e_i(w)e_i
\eeaa
 and let  $z_{k,n}^\perp=z_k-z_{k,n}$. Then we have
\beaa
E[G|V]&=&\int G(z_k+w_k)d\mu(z_k)\\
&=&\lim_n\frac{1}{(2\pi)^{n/2}}\int_{\R^n}G(z_{k,n}^\perp+z_{k,n}+w_k)
e^{-\frac{|z_{k,n}|^2}{2}}dz_{k,n}\,.
\eeaa
Since
$$
(z_{k,n},w_k)\to
\exp\left\{-\frac{1}{2}(a|w_k|^2+|z_{n,k}|^2)\right\}
G(z_{k,n}^\perp+z_{k,n}+w_k)
$$
is almost surely log-concave, the proof follows from Prekopa's
theorem (cf. \cite{PRE71}).
\qed

\noindent
The following theorem extends Theorem  \ref{Gamma-thm} :
\begin{theorem}
\label{general-O-U}
Let $G$ be an $a$-log-concave Wiener functional, where $a\in
[0,1)$. Then   $\Gamma(A)G$ is
$a$-log-concave, where $A\in L(H,H)$ (i.e. the space of bounded linear
operators on $H$) with $\|A\|\leq 1$. In particular $P_tG$ is
$a$-log-concave for any $t\geq 0$, where $(P_t,t\geq 0)$ denotes the
Ornstein-Uhlenbeck semi-group on $W$.
\end{theorem}
\proof
Let $(e_i,i\in \NN)$ be a complete, orthonormal basis of $H$, denote
by $\pi_n$ the orthogonal projection from $H$ onto the linear  space spanned
by $\{e_1,\ldots,e_n\}$ and by $V_n$ the sigma algebra generated by
$\{\delta e_1,\ldots,\delta e_n\}$. From Proposition \ref{c-exp-prop}
and from the fact that $\Gamma(\pi_nA\pi_n)\to \Gamma(A)$ in the
strong operator topology as $n$ tends to infinity, it suffices to
prove the theorem when $W=\reals^n$. We may then assume that $G$ is
bounded and of compact support. Define $F$ as
\beaa
G(x)&=&F(x)e^{\frac{a}{2}|x|^2}\\
    &=&F(x)\int_{\reals^n}e^{\sqrt{a}(x,\xi)}d\mu(\xi)\,.
\eeaa
From the hypothesis, $F$ is almost surely log-concave. Then, using the
notations explained in Section 2:
\beaa
\lefteqn{e^{-a\frac{|x|^2}{2}}\Gamma(A)G(x)}\\
&=&\int\int
F(A^*x+Sy)\exp\left\{-a\frac{|x|^2}{2}+\sqrt{a}(A^*x+Sy,\xi)\right\}
d\mu(y)d\mu(\xi)\\
&=&(2\pi)^{-n}\int\int F(A^*x+Sy)\exp-\frac{\Theta(x,y,\xi)}{2}\, dyd\xi\,,
\eeaa
where
\beaa
\Theta(x,y,\xi)&=&a|x|^2-2\sqrt{a}(A^*x+Sy,\xi)+|y|^2+|\xi|^2\\
  &=&|\sqrt{a}x-A\xi|^2+|\sqrt{a}y-S\xi|^2+(1-a)|y|^2\,,
\eeaa
which is a convex function of $(x,y,\xi)$.  Hence the proof follows
from Pr\'ekopa's theorem (cf. \cite{PRE71}).
\qed

The following proposition  extends a well-known  finite dimensional
inequality (cf. \cite{Hu}):
\begin{proposition}
Assume that  $f$ and $g$ are $H$-convex Wiener functionals such that
$f\in L^p(\mu)$ and $g\in L^q(\mu)$ with  $p>1,\,p^{-1}=1-q^{-1}$. Then
\begin{equation}
\label{cov-ineq}
E[f\,g]\geq E[f] E[g]+\left(E[\nabla f],E[\nabla g]\,\right)_H\,.
\end{equation}
\end{proposition}
\proof
Define the smooth and convex functions $f_n$ and $g_n$ on $W$ by
\beaa
P_{1/n}f&=&f_n\\
P_{1/n}g&=&g_n\,.
\eeaa
Using the fact that $P_t=e^{-t\L}$, where $\L$ is the number operator
$\L=\delta\circ \nabla$ and the commutation relation $\nabla
P_t=e^{-t}P_t\nabla$, for any $0\leq t\leq T$,  we have
\begin{eqnarray}
\label{yosida}
E\left[P_{T-t}f_n\,g_n\right]&=&E[P_Tf_n\,g_n]+
\int_0^tE\left[\L P_{T-s}f_n\,g_n\right]ds\nonumber\\
  &=&E[P_Tf_n\,g_n]+\int_0^te^{-(T-s)}E\left[\left(P_{T-s}\nabla f_n,\nabla
    g_n\right)_H\right]ds\nonumber\\
&=&E[P_Tf_n\,g_n]+\int_0^te^{-(T-s)}E\left[\left(P_T\nabla f_n,\nabla
g_n\right)_H\right]ds\nonumber\\
&&+e^{-2T}\int_0^t\int_0^se^{s+\tau}
E\left[\left(P_{T-\tau}\nabla^2f_n,\nabla^2g_n\right)_2\right]
d\tau ds\nonumber\\
&\geq &E[P_Tf_n\,g_n]\nonumber\\
&&+E\left[\left(P_T\nabla f_n,\nabla g_n\right)_H\right]\,e^{-T}(e^t-1)
\end{eqnarray}
where $(\cdot,\cdot)_2$ denotes the Hilbert-Schmidt scalar product
and the inequality  (\ref{yosida}) follows from the convexity of
$f_n$ and $g_n$. In fact their convexity implies that $P_t
\nabla^2 f_n$ and $\nabla^2 g_n$ are positive operators, hence
their Hilbert-Schmidt tensor product is positive. Letting $T=t$ in
the above inequality we have
\begin{equation}
\label{yosida-2}
E[f_n\,g_n]\geq
E\left[P_Tf_n\,g_n\right]+(1-e^{-T})E\left[\left(P_T\nabla f_n,\nabla
    g_n\right)_H\right]\,.
\end{equation}
Letting  $T\rightarrow \infty$ in (\ref{yosida-2}), we obtain, by the
ergodicity of $(P_t,t\geq 0)$, the
claimed inequality for $f_n$ and $g_n$. It suffices then to  take the
limit of this inequality  as $n$ tends to infinity.
\qed


\begin{proposition}
\label{con-con}
Let $G$ be a (positive) $\gamma$-log-concave Wiener functional with $\gamma\in
[0,1]$. Then the map $h\to E[G(w+h)]$ is a log-concave mapping on $H$.
In particular, if $G$ is symmetric, i.e., if
$G(w)=G(-w)$, then
$$
E[G(w+h)]\leq E[G]\,.
$$
\end{proposition}
\proof
Without loss of generality, we may suppose that $G$ is bounded. Using
the usual notations, we have, for any $h$ in any finite dimensional
subspace $L$ of $H$,
$$
E[G(w+h)]=\lim_n\frac{1}{(2\pi)^{n/2}}\int_{W_n}G(w_n^\perp+w_n+h)
\exp\left\{-\frac{|w_n|^2}{2}\right\}dw_n\,,
$$
from the hypothesis, the integrand is almost surely log-concave on
$W_n\times L$, from Prekopa's theorem, the integral is log-concave on
$L$, hence the limit is  also log-concave. Since $L$ is arbitrary,
the first part of the  proof follows. To prove the second part, let
$g(h)=E[G(w+h)]$, then, from the log-concavity of $g$ and symmetry of
$G$, we have
\beaa
E[G]&=&g(0)\\
&=&g\left(1/2(h)+1/2(-h)\right)\\
          &\geq&g(h)^{1/2}g(-h)^{1/2}\\
          &=&g(h)\\
   &=&E[G(w+h)]\,.
\eeaa
\qed

\remark In fact, with a little bit more attention, we can see that
the map $h\to \exp\{\frac{1}{2}(1-\gamma)|h|^2_H\}E[G(w+h)]$ is
log-concave on $H$.
\newline

\noindent
We have the following immediate corollary:
\begin{corollary}
Assume that $A\subset W$ is an $H$-convex  and  symmetric set. Then we
have
$$
\mu(A+h)\leq \mu(A)\,,
$$
for any $h\in H$.
\end{corollary}
\proof
Since $1_A$ is log $H$-concave, the proof  follows from Proposition
\ref{con-con}.
\qed

\begin{proposition}
Let $F\in L^p(\mu)$ be a  positive  log $H$-convex function. Then for
any $u\in \DD_{q,2}(H)$, we have
$$
E_F\left[\left(\delta u-E_F[\delta u]\right)^2\right]\geq
 E_F\left[|u|_H^2+2\delta(\nabla_u u)+\trace (\nabla u\cdot \nabla u)\right]\,,
$$
where $E_F$ denotes the mathematical expectation with respect to the
probability defined as
$$
\frac{F}{E[F]}d\mu\,.
$$
\end{proposition}
\proof Let $F_\tau$ be $P_\tau F$, where $(P_\tau,\tau\in \reals_+)$
denotes the Ornstein-Uhlenbeck semi-group. $F_\tau$ has a
modification, denoted again by the same letter, such that the mapping
$h\mapsto F_\tau(w+h)$ is real-analytic on $H$ for all $w\in W$
(cf. \cite{UZ7}). Suppose first also that $\|\nabla u\|_{_2}\in
L^\infty (\mu, H\otimes H)$ where $\|\cdot\|_{_2}$ denotes the
Hilbert-Schmidt norm. Then, for any $r>1$, there exists some $t_r>0$
such that, for any $0\leq t<t_r$,  the image of the Wiener measure
under $w\mapsto w+tu(w)$ is equivalent to $\mu$ with the Radon-Nikodym
density $L_t\in L^r(\mu)$.
Hence $w\mapsto F_\tau(w+tu(w))$ is a well-defined
mapping on $W$ and it is in some $L^r(\mu)$ for small $t>0$
(cf. \cite{UZ7}, Chapter~3 and  Lemma B.8.8). Besides
$t\mapsto F(w+tu(w))$ is log convex on $\reals$
since $F_\tau $ is log $H$-convex. Consequently $t\mapsto
E[F_\tau(w+tu(w))]$ is log convex and strictly positive. Then the
second derivative of its logarithm  at $t=0$ should be positive. This
implies immediately the claimed inequality for $\nabla u$ bounded. We
then pass to the limit with respect to $u$ in $\DD_{q,2}(H)$ and then
let $\tau \to 0$ to complete the proof.
\qed
\markboth{Poincar\'e inequality}{log-Sobolev inequality}
\section{Poincar\'e and logarithmic Sobolev inequalities}
The following theorem  extends the Poincar\'e- Brascamp-Lieb
inequality:
\begin{theorem}
\label{B-L-thm}
Assume that $F$ is a Wiener functional in
$\cup_{p>1}\DD_{p,2}$ with $e^{-F}\in L^1(\mu)$ and assume also  that  there
exists a constant $\eps>0$ such that
\begin{equation}
\label{growth-cond}
\left((I_H+\nabla^2F)h,h\right)_H\geq \eps |h|_H^2
\end{equation}
almost surely, for any $h\in H$, i.e. $F$ is $(1-\eps)$-convex.
Let us denote by $\nu_F$ the
probability measure on $(W,\calB(W))$ defined by
$$
d\nu_F=\exp\left\{-F-\log E\left[e^{-F}\right]\right\}d\mu\,.
$$
Then for any smooth cylindrical Wiener
functional $\phi$, we have
\begin{equation}
\label{poincare-ineq}
\int_W|\phi-E_{\nu_F}[\phi]|^2d\nu_F\leq \int_W\left((I_H+\nabla^2F)^{-1}\nabla
\phi,\nabla \phi\right)_Hd\nu_F\,.
\end{equation}
In particular, if $F$ is an $H$-convex Wiener functional, then the
condition (\ref{growth-cond}) is satisfied with $\eps=1$.
\end{theorem}
\proof
Assume first that $W=\R^n$ and that $F$ is a smooth  function
on $\R^n$ satisfying the inequality (\ref{growth-cond}) in this
setting. Assume also for the typographical facility that
$E[e^{-F}]=1$. For any smooth function function $\phi$ on $\R^n$, we have
\begin{equation}
\label{fin-dim}
\int_{\R^n}\left|\phi-E_{\nu_F}[\phi]\right|^2d\nu_F=\frac{1}{(2\pi)^{n/2}}
\int_{\R^n} e^{-F(x)-|x|^2/2}\left|\phi(x)-E_F[\phi]\right|^2dx\,.
\end{equation}
The function $G(x)=F(x)+\frac{1}{2}|x|^2$ is a strictly convex smooth
function. Hence Brascamp-Lieb inequality (cf. \cite{BLI})
implies that:
\beaa
\int_{\R^n}\left|\phi-E_{\nu_F}[\phi]\right|^2d\nu_F&\leq&\int_{\R^n}\left
(\left({\mbox{\rm
        Hess}}\,G(x)\right)^{-1}\nabla\phi(x),\nabla\phi(x)\right)_{\R^n}d\nu_F(x)\\
&=&\int_{\R^n}\left((I_{\R^n}+\nabla^2F)^{-1}\nabla\phi,\nabla\phi\right)_
{\R^n}d\nu_F\,.
\eeaa
To prove the general case we proceed by approximation as before:
indeed let $(e_i,i\in \NN)$ be a complete,
orthonormal basis of $H$, denote by $V_n$ the sigma algebra generated
by $\{\delta e_1,\ldots,\delta e_n\}$. Define $F_n$ as to be
$E[P_{1/n}F|V_n]$, where $P_{1/n}$ is the Ornstein-Uhlenbeck semigroup
at $t=1/n$. Then from the martingale convergence theorem and
the fact that $V_n$ is a smooth sigma algebra, the sequence $(F_n,n\in
\NN)$ converges to $F$ in some $\DD_{p,2}$. Moreover $F_n$
satisfies the hypothesis (with  a better constant in the inequality
 (\ref{growth-cond})) since $\nabla^2F_n=e^{-2/n}E[Q_n^{\otimes
   2}\nabla^2F|V_n]$, where $Q_n$ denotes the orthogonal projection
 onto the vector space spanned by $\{e_1,\ldots,e_n\}$. Besides $F_n$
 can be represented as $F_n=\theta(\delta e_1,\ldots,\delta e_n)$,
 where $\theta$ is a smooth  function on $\R^n$ satisfying
$$
((I_{\reals^n}+\nabla^2\theta(x))y,y)_{\reals^n}\geq \eps
|y|_{\reals^n}^2\,,
$$
for any $x,y\in \reals^n$. Let $w_n={\tilde{Q}}_n(w)=\sum_{i\leq
  n}(\delta e_i)e_i$, $W_n={\tilde{P}}_n(W)$ and
$W_n^\perp=(I_W-{\tilde{Q}}_n)(W)$ as before.
Let us denote by $\nu_n$ the probability measure corresponding to
$F_n$. Let us also denote by $V_n^\perp$ the sigma algebra generated
by $\{\delta e_k,k>n\}$.  Using the finite dimensional result that we
have derived,  the Fubini theorem and the inequality $2|ab|\leq \kappa
a^2+\frac{1}{\kappa} b^2$, for any $\kappa>0$,  we obtain
\bea
\label{uzun-eqn}
\lefteqn{E_{\nu_{n}}\left[\left|\phi-E_{\nu_{n}}[\phi]\right|^2\right]}\nonumber\\
&=&\int_{W_n\times
  W_n^\perp}e^{-F'_n(w_n)}|\phi(w_n+w_n^\perp)-E_{\nu_n}[\phi]|^2d\mu_n(w_n)d\mu_n^\perp(w_n^\perp)\nonumber\\
&\leq&
(1+\kappa)\int_We^{-F'_n}|\phi-E[e^{-F'_n}\phi|V_n^\perp]|^2d\mu\nonumber\\
&&+\left(1+\frac{1}{\kappa}\right)\int_We^{-F'_n}|E[e^{-F'_n}\phi|V_n^\perp]-E_{\nu_n}[\phi]|^2d\mu\nonumber\\
&\leq&
(1+\kappa)E_{\nu_n}\left[\left((I_H+\nabla^2F_n)^{-1}\nabla\phi,\nabla\phi\right)_H\right]\nonumber\\
&&+\left(1+\frac{1}{\kappa}\right)\int_We^{-F'_n}|E[e^{-F'_n}\phi|V_n^\perp]-E_{\nu_n}[\phi]|^2d\mu\,,
\eea
where $F_n'$ denotes $F_n-\log E[e^{-F_n}]$. Since $V_n$ and
$V_n^\perp$ are independent sigma algebras, we have
\beaa
|E[e^{-F'_n}\phi|V_n^\perp]|&=&\frac{1}{E[e^{-F_n}]}
|E[e^{-F'_n}\phi|V_n^\perp]|\\
&\leq&\frac{1}{E[e^{-F_n}]}E[e^{-F_n}|V_n^\perp]\|\phi\|_\infty\\
&=&\|\phi\|_\infty\,,
\eeaa
hence, using the triangle inequality  and the dominated convergence
theorem, we realize that the last term in (\ref{uzun-eqn}) converges
to zero as $n$ tends to infinity. Since the sequence of  operator valued random
variables $((I_H+\nabla^2F_n)^{-1},n\in \NN)$ is  essentially bounded
in the strong operator norm, we can pass to the limit on both sides
and this gives the claimed inequality with a factor $1+\kappa$, since
$\kappa>0$ is arbitrary, the  proof is completed.
\qed


\remark
Let $T:W\to W$ be a shift  defined as $T(w)=w+u(w)$, where $u:W\to H$
is a measurable map satisfying $(u(w+h)-u(w),h)_H\geq -\eps|h|^2$. In
\cite{UZ0} and in \cite{UZ7}, Chapter 6, we have studied such
transformations, called $\eps$-monotone shifts. Here  the hypothesis of Theorem
\ref{B-L-thm} says that the shift $T=I_W+\nabla F$ is $\eps$-monotone.

The Sobolev regularity hypothesis can be omitted if we are after a
Poincar\'e inequality with another  constant:
\begin{theorem}
Assume that $F\in \cup_{p>1}L^p(\mu)$ with
$E\left[e^{-F}\right]$ is finite and that, for some constant $\eps>0$,
$$
E\left[\left((I_H+\nabla^2F)h,h\right)_H\,\psi\right]\geq
\eps\,|h|_H^2E[\psi]\,,
$$
for any $h\in H$ and positive test function $\psi\in\DD$, where
$\nabla^2F$ denotes the second order derivative in the sense of the
distributions. Then we have
\begin{equation}
\label{poincare-2-ineq}
E_{\nu_F}\left[|\phi-E_F[\phi]|^2\right]\leq
\frac{1}{\eps}E_{\nu_F}[|\nabla\phi|_H^2]
\end{equation}
for any cylindrical Wiener functional $\phi$. In particular, if $F$ is
$H$-convex, then we can take $\eps=1$.
\end{theorem}
\proof
Let $F_t$ be defined as $P_tF$, where $P_t$ denotes the
Ornstein-Uhlenbeck semigroup. Then $F_t$ satisfies the hypothesis of
Theorem \ref{B-L-thm}, hence we have
$$
E_{\nu_{F_t}}\left[\left|\phi-E_{F_t}[\phi]\right|^2\right]\leq
\frac{1}{\eps}E_{\nu_{F_t}}\left[|\nabla\phi|_H^2\right]
$$
for any $t>0$. The claim follows when we take the limits of both sides
as $t\to 0$.
\qed
\paragraph{\bf Example:}
Let $F(w)=\|w\|+\frac{1}{2}\sin(\delta h)$ with $|h|_H\leq 1$, where
$\|\cdot\|$ denotes the norm of the Banach
space $W$. Then in general $F$ is not in $\cup_{p>1}\DD_{p,2}$,
however the Poincar\'e  inequality (\ref{poincare-2-ineq}) holds with
$\eps=1/2$.
\begin{theorem}
\label{log-sob-thm}
Assume that $F$ is a Wiener functional in
$\cup_{p>1}\DD_{p,2}$ with $E[\exp-F]<\infty$. Assume  that  there
exists a constant $\eps>0$ such that
\begin{equation}
\label{growth-cond-1}
\left((I_H+\nabla^2F)h,h\right)_H\geq \eps |h|_H^2
\end{equation}
almost surely, for any $h\in H$.
Let us denote by $\nu_F$ the
probability measure on $(W,\calB(W))$ defined by
$$
d\nu_F=\exp\left\{-F-\log E\left[e^{-F}\right]\right\}d\mu\,.
$$
Then for any smooth cylindrical Wiener
functional $\phi$, we have
\begin{equation}
\label{log-sob-ineq}
E_{\nu_F}\left[\phi^2\left\{\log \phi^2-\log \|\phi\|_{L^2(\nu_F)}^2\right\}\right]\leq
\frac{2}{\eps} E_{\nu_F}\left[|\nabla \phi|_H^2\right]\,.
\end{equation}
In particular, if $F$ is an $H$-convex Wiener functional, then the
condition (\ref{growth-cond-1}) is satisfied with $\eps=1$.
\end{theorem}
\proof
We shall proceed as in the proof of Theorem \ref{B-L-thm}. Assume then
that $W=\R^n$ and that $F$ is a smooth function satisfying the
inequality (\ref{growth-cond-1})  in this frame. In this case it is
immediate to see that  function $G(x)=\frac{1}{2}|x|^2+F(x)$ satisfies
the Bakry-Emery condition (cf. \cite{B-E}, \cite{D-S}), which is known
as a sufficient condition for the inequality (\ref{log-sob-ineq}). For
the infinite dimensional case we define as in the proof of Theorem
\ref{B-L-thm}, $F_n,\nu_n,V_n,V_n^\perp$. Then, denoting by $E_n$ the
expectation with respect to the probability $\exp\{-F'_n\}d\mu$, where
$F_n'=F_n-\log E[e^{-F_n}]$, we have
\bea
\label{l-s-1}
\lefteqn{E_n\left[\phi^2\left\{\log \phi^2-\log
    \|\phi\|_{L^2(\nu_F)}^2\right\}\right]}\nonumber\\
&=&E_n\left[\phi^2\left\{\log \phi^2-\log
    E[e^{-F_n'}\phi^2|V_n^\perp]\right\}\right]\nonumber\\
&&+E_n\left[\phi^2\left\{\log E[e^{-F_n'}\phi^2|V_n^\perp]-\log
    E_n[\phi^2]\right\}\right]\nonumber\\ 
&\leq&\frac{2}{\eps}E_n\left[|\nabla
  \phi|_H^2\right]\nonumber\\
  &&+E_n\left[\phi^2\left\{\log
    E[e^{-F_n'}\phi^2|V_n^\perp]-\log E_n[\phi^2]\right\}\right]\,,
\eea
where we have used, as in the proof of Theorem \ref{B-L-thm}, the
finite dimensional log-Sobolev inequality to obtain the inequality
(\ref{l-s-1}).
Since in the above inequalities everything is squared, we can assume
that $\phi$ is positive, and adding a constant $\kappa>0$, we can also
replace $\phi $ with  $\phi_\kappa=\phi+ \kappa$. Again by the
independance of $V_n$ and $V_n^\perp$, we  can pass to the limit with
respect to $n$ in
the inequality (\ref{l-s-1}) for $\phi=\phi_\kappa$ to obtain
$$
E_{\nu_F}\left[\phi_\kappa^2\left\{\log \phi_\kappa^2-\log
    \|\phi_\kappa\|_{L^2(\nu_F)}^2\right\}\right]\leq 
\frac{2}{\eps} E_{\nu_F}\left[|\nabla \phi_\kappa|_H^2\right]\,.
$$
To complete the proof it suffices to pass to the limit as $\kappa\to
0$.
\qed

The following theorem fully   extends Theorem \ref{log-sob-thm} and it
is useful for the applications:
\begin{theorem}
\label{l-s-3}
Assume that $G$ is a (positive) $\gamma$-log-concave Wiener functional for
some $\gamma\in [0,1)$ with $E[G]<\infty$. Let us denote by $E_G[\cdot\,]$
the expectation with respect to the probability measure defined by
$$
d\nu_G=\frac{G}{E[G]}d\mu\,.
$$
Then we have
\begin{equation}
\label{l-s-2}
E_G\left[\phi^2\left\{\log \phi^2-\log E_G[\phi^2]\right\}\right]\leq
\frac{2}{1-\gamma}E_G[|\nabla \phi|_H^2]\,,
\end{equation}
for any cylindrical Wiener functional $\phi$.
\end{theorem}
\proof
Since $G\wedge c$, $c>0$,  is again $\gamma$-log-concave, we may
suppose without loss of generality that $G$ is bounded. Let now
$(e_i,i\in \NN)$ be a
complete, orthonormal basis for $H$, denote by $V_n$ the sigma algebra
generated by $\{\delta e_1,\ldots,\delta e_n\}$. Define $G_n$ as to be
$E[P_{1/n}G|V_n]$. From Proposition \ref{c-exp-prop} and Theorem
\ref{general-O-U}, $G_n$ is again a  $\gamma$-log-concave, strictly positive
Wiener functional. It can be represented  as
$$
G_n(w)=g_n(\delta e_1,\ldots,\delta e_n)
$$
and due to the Sobolev embedding  theorem, after a modification on a
set of zero Lebesgue measure, we can assume that  $g_n$ is  a smooth
function on $\reals^n$. Since it is
strictly positive, it is of the form $e^{-f_n}$, where $f_n$ is a
smooth, $\gamma$-convex function. It follows then from Theorem
\ref{log-sob-thm} that the inequality (\ref{l-s-2}) holds when we
replace $G$ by $G_n$, then the proof follows by taking the limits  of
both  sides as $n\to \infty$.
\qed
\paragraph{Example:}
Assume that $A$ is a measurable subset of $W$ and let $H$ be a
measurable Wiener functional with values in
$\reals\cup\{\infty\}$. If $G$ defined by  $G=1_A\,H$ is
$\gamma$-log-concave with $\gamma\in [0,1)$, then the hypothesis of Theorem
\ref{l-s-3} are satisfied.

\begin{definition}
Let $T\in \DD'$ be a positive distribution. We say that it is
$a$-log-concave if $P_tT$ is an $a$-log-concave Wiener functional. If
$a=0$, then we call $T$ simply log-concave.
\end{definition}
\remark
From Corollary \ref{pos-dist},  to any positive
distribution on $W$, it
corresponds a positive Radon measure $\nu_T$ such that
$$
<T,\phi>=\int_W {\tilde{\phi}}(w)d\nu_T(w)
$$
 for any $\phi\in \DD$, where $\tilde{\phi}$ represents a
 quasi-continuous version of $\phi$.
\paragraph{Example:}
Let $(w_t,t\in [0,1])$ be the one-dimensional Wiener process and
denote by  $p_\tau$ the heat kernel on $\reals$.  Then the
distribution defined as $\veps_0(w_1)=\lim_{\tau\to 0}p_\tau(w_1)$ is
log-concave, where $\veps_0$ denotes the Dirac measure at zero.

The following result is a Corollary of Theorem \ref{l-s-3}:

\begin{theorem}
\label{dist-case}
Assume that $T\in \DD'$ is a positive, $\beta$-log-concave distribution
with $\beta\in [0,1)$. Let $\gamma$ be the probability  Radon  measure
defined by
$$
\gamma=\frac{\nu_T}{<T,1>}\,.
$$
Then we have
\begin{equation}
\label{l-s-4}
E_\gamma\left[\phi^2\left\{\log \phi^2-\log E_\gamma[\phi^2]\right\}\right]\leq
\frac{2}{1-\beta}E_\gamma[|\nabla \phi|_H^2]\,,
\end{equation}
for any smooth cylindrical function $\phi:W\to \reals$.
\end{theorem}

Here is  an application of this result:
\begin{proposition}
\label{non-deg}
Let $F$ be a Wiener functional in $\DD_{r,2}$ for some $r>1$. Suppose
that it is $p$-non-degenerate in the sense that
\begin{equation}
\label{non-deg-con}
\delta\left\{\frac{\nabla F}{|F|^2}\phi\right\}\in L^p(\mu)
\end{equation}
for any $\phi\in \DD$, for some $p>1$. Assume furthermore that, for
some $x_0\in \reals$,
\begin{equation}
\label{x0-cond}
(F-x_0)\nabla^2F+\nabla F\otimes \nabla F\geq 0
\end{equation}
almost surely. Then we have
$$
E\left[\phi^2\left\{\log \phi^2-\log
    E\left[\phi^2|F=x_0\right]\right\}|F=x_0\right]\leq
 2\,E\left[|\nabla \phi|_H^2|F=x_0\right]
$$
for any smooth cylindrical $\phi$.
\end{proposition}
\proof
Note that the non-degeneracy  hypothesis (\ref{non-deg-con}) implies
the existence of a continuous
density of the law of $F$ with respect to the Lebesgue
measure (cf. \cite{Mall-1} and the references there). Moreover it
implies also the fact that
$$
\lim_{\tau\to 0}p_\tau(F-x_0)=\veps_{x_0}(F)\,,
$$
in $\DD'$, where $\veps_{x_0}$ denotes the Dirac measure at $x_0$ and
$p_\tau$ is  the heat kernel on $\reals$. The
inequality (\ref{x0-cond}) implies that the distribution
defined by
$$
\phi\to E[\phi|F=x_0]=\frac{<\veps_{x_0}(F),\phi>}{<\veps_{x_0}(F),1>}
$$
is log-concave, hence the conclusion follows from Theorem \ref{dist-case}.
\qed
\section{Change of variables formula and log-Sobolev inequality}
\markboth{Change of Variables}{log-Sobolev}
In this section we shall derive a different kind of logarithmic
Sobolev inequality using the change of variables formula for the
monotone shifts studied in \cite{UZ0} and in more detail in
\cite{UZ7}. An analogous approach to derive log-Sobolev-type
inequalities  using the Girsanov theorem has
been employed  in \cite{UST00}.

\begin{theorem}
\label{log-sob-tordu}
Suppose that $F\in L^p(\mu)$, for some $p>1$, is an $a$-convex Wiener
functional, $a\in [0,1)$ with $E[F]=0$. Assume that
\begin{equation}
\label{growth2-cond}
E\left[\exp\left\{ c\,\|\nabla^2 \L^{-1}F\|_{_2}^2\right\}\right]<\infty\,,
\end{equation}
for some
$$
c>\frac{2+(1-a)}{2(1-a)}\,,
$$
where $\|\cdot\|_{_2}$ denotes the Hilbert-Schmidt norm on $H\otimes
H$ and $\L^{-1}F=\int_{\R_+}P_tF\,dt$. Denote by $\nu$ the
probability measure defined by
$$
d\nu=\La\,d\mu\,,
$$
where
$$
\La=\dett(I_H+\nabla^2\L^{-1}F)\exp\left\{-F-\frac{1}{2}|\nabla
  \L^{-1}F|_H^2\right\}\,
$$
and $\dett(I_H+\nabla^2\L^{-1}F)$ denotes the modified Carleman-Fredholm
determinant. Then we have
\begin{equation}
\label{eqn-tordu}
E_\nu\left[f^2\log\left(\frac{f^2}{\|f\|^2_{L^2(\nu)}}\right)\right]\leq
2 E_\nu\left[|(I_H+\nabla^2 \L^{-1}F)^{-1}\nabla f|_H^2\right]
\end{equation}
and
\begin{equation}
\label{poincare-tordu}
E_\nu[|f-E_\nu[f]|^2]\leq E_\nu\left[|(I_H+\nabla^2
  \L^{-1}F)^{-1}\nabla f|_H^2\right]
\end{equation}
for any smooth, cylindrical $f$.
\end{theorem}
\proof
Let $F_n=E[P_{1/n}F|V_n]$, where $V_n$ is the sigma algebra generated
by $\{\delta e_1,\ldots,\delta e_n\}$ and let  $(e_n,n\in \NN)$ be  a complete,
orthonormal basis of $H$. Define $\xi_n$ by $\nabla \L^{-1}F_n$, then
$\xi_n$ is $(1-a)$-strongly monotone (cf. \cite{UZ0} or \cite{UZ7}) and
smooth. Consequently, the shift $T_n:W\to W$, defined by
$T_n(w)=w+\xi_n(w)$ is a bijection of $W$ (cf. \cite{UZ7} Corollary
6.4.1), whose inverse is of the form $S_n=I_W+\eta_n$, where
$\eta_n(w)=g_n(\delta e_1,\ldots,\delta e_n)$ such that
$g_n:\reals^n\to \reals^n$ is a smooth function. Moreover the images of
$\mu$ under $T_n$ and $S_n$, denoted by $T_n^*\mu$ and $S_n^*\mu$
respectively,  are equivalent to $\mu$ and we have
\beaa
\frac{dS_n^*\mu}{d\mu}&=&\La_n\\
\frac{dT_n^*\mu}{d\mu}&=&L_n
\eeaa
where
\beaa
\La_n&=&\dett(I_H+\nabla \xi_n)\exp\left\{-\delta
  \xi_n-\frac{1}{2}|\xi_n|_H^2\right\}\\
L_n&=&\dett(I_H+\nabla \eta_n)\exp\left\{-\delta
  \eta_n-\frac{1}{2}|\eta_n|_H^2\right\}\,.
\eeaa
The hypothesis (\ref{growth2-cond}) implies the uniform integrability
of the densities $(\La_n,n\geq 1)$ and $(L_n,n\geq 1)$
(cf. \cite{USZ98,UZ7}).
For any probability $P$ on $(W,\calB(W))$  and any positive, measurable
function $f$,  define  $\calH_P(f)$ as
\begin{equation}
\label{entropy}
\calH_P(f)=f(\log f-\log E_P[f]).
\end{equation}
Using    the logarithmic Sobolev inequality of L. Gross  for $\mu$
(cf. \cite{LG}) and the relation
$$
(I_H+\nabla \eta_n)\circ T_n=(I_H+\nabla \xi_n)^{-1}\,,
$$
we have
\begin{eqnarray}
\label{neweq}
E[\La_n \calH_{\La_nd\mu}(f^2)]&=&E[\calH_\mu(f^2\circ S_n)]\nonumber\\
&\leq& 2E[|\nabla (f\circ S_n)|_H^2]\nonumber\\
&=&2E[|(I_H+\nabla \eta_n)\nabla f\circ S_n|_H^2]\nonumber\\
&=&2E[\La_n|(I_H+\nabla \xi_n)^{-1}\nabla f|_H^2]\,.
\end{eqnarray}
 It follows by
the $a$-convexity  of $F$
that
$$
\|(I_H+\nabla \xi_n)^{-1}\|\leq \frac{1}{1-a}
$$
almost surely for any $n\geq 1$, where $\|\cdot\|$ denotes the operator norm.
Since the sequence $(\La_n,n\in \NN)$ is uniformly integrable,
 the limit of (\ref{neweq}) exists  in $L^1(\mu)$ and the proof of
 (\ref{eqn-tordu}) follows. The proof of the inequality
 (\ref{poincare-tordu}) is now trivial.
\qed

\begin{corollary}
Assume that $F$ satisfies the hypothesis of Theorem
\ref{log-sob-tordu}. Let $Z$ be the functional defined by
$$
Z=\dett(I_H+\nabla^2 \L^{-1}F)\exp\frac{1}{2}|\nabla
\L^{-1}F|_H^2
$$
and assume that $Z,\,Z^{-1}\in L^\infty(\mu)$.
Then we have
\begin{equation}
\label{last-ineq}
E\left[e^{-F}f^2\log\left\{\frac{f^2}{E[e^{-F}f^2]}\right\}\right]\leq
2K
E\left[e^{-F}\left|(I_H+\nabla^2 \L^{-1}F)^{-1}\nabla f\right|_H^2\right]
\end{equation}
and
\begin{equation}
\label{last-poincare}
E\left[e^{-F}\left|f-E[e^{-F}f]\right|^2\right]\leq K E\left[e^{-F}
\left|(I_H+\nabla^2 \L^{-1}F)^{-1}\nabla f\right|_H^2\right]
\end{equation}
for any smooth, cylindrical $f$, where
$K=\|Z\|_{L^\infty(\mu)}\|Z^{-1}\|_{L^\infty(\mu)}$.
\end{corollary}
\proof
Using the identity remarked by  Holley and
Stroock (cf. \cite{H-S}, p.1183)
$$
E_P\left[\calH_P(f^2)\right]=\inf_{x>0}
E_P\left[f^2\log\left(\frac{f^2}{x}\right)-(f^2-x)\right]\,,
$$
where $P$ is an arbitrary probability measure, and $\calH$ is defined
by the relation  (\ref{entropy}),
we see that the inequality (\ref{last-ineq}) follows from Theorem
\ref{log-sob-tordu} and  the inequality (\ref{last-poincare}) is trivial.
\qed

\section*{Exercises}
{\footnotesize{
\begin{enumerate}
\item Assume that $A_1,\ldots,A_n$ are almost surely convex and symmetric
sets. Prove the following inequality:
\begin{equation}
\label{mcon-ineq}
 \mu\left(\bigcap_{i=1}^n(A_i+h_i)\right)\leq
\mu\left(\bigcap_{i=1}^nA_i\right)\,,
\end{equation}
for any $h_1,\ldots,h_n\in H$.
\item Assume that $F$ is a positive,  symmetric, almost surely log-concave
Wiener functional such that $\mu\{F>0\}>0$. Denote by $\mu_F$ the
probability defined by
$$
d\mu_F=\frac{F}{E[F]}\,d\mu\,.
$$
Prove the inequality (\ref{mcon-ineq}) when $\mu$ is replaced by
$\mu_F$.
\item Let $A$ and $B$ be two almost surely convex sets. For
$\alpha\in [0,1]$, define the map $(\alpha,w)\to f(\alpha,w)$ as
$$
f(\alpha,w)=\won_{C_\alpha}(w)\,,
$$
where $C_\alpha=\alpha A+(1-\alpha)B$. Prove that $(\alpha,w)\to
f(\alpha,w)$ is almost surely log-concave. Deduce from that and
from Pr\'ekopa's theorem the inequality:
$$
\mu(C_\alpha)\geq \mu(A)^\alpha\,\mu(B)^{1-\alpha}\,.
$$
\item Let $F$ and $G$ be two almost surely convex, symmetric
Wiener functionals from $\DD_{2,2}$. Prove that
$$
E[(\nabla F,\nabla G)_H]\geq 0\,.
$$
\item Let $W$ be the classical Wiener space $C_0([0,1],\R)$ and let
  $f$ and $g$ be two $H$-convex functions in $L^2(\mu)$. With the help
  of the Clark's formula, prove that 
$$
E[E[D_tf|\F_t]\,E[D_tg|\F_t]]\geq E[D_tf]E[D_tg]\,,
$$
$dt$-almost surely.
\end{enumerate}

}}

\section*{Notes and references}
{\footnotesize{ The notion of convexity for the equivalence
classes of Wiener random variables is a new subject. It has been
studied for the first time in \cite{F-U}. Even in the finite
dimensional case it is not evident to find a result about the
$H$-convexity.

The log-Sobolev inequalities given here are well-known in the
finite dimensional case except the content of the last section.
The fact that log-concavity is preserved under the action of
certain semi-groups and especially its implications
concerning log-concave distributions  seem to be  novel.

}}

\chapter{Monge-Kantorovitch Mass  Transportation}
\markboth{Monge-Kantorovitch}{Mass Transportation}
\section{Introduction}
\label{intro}

In 1781, Gaspard Monge has published his celebrated memoire about the
most economical way of earth-moving \cite{Monge}. The configurations
of excavated earth and remblai were  modelized as two measures of
equal mass, say $\rho$ and $\nu$, that Monge had supposed absolutely
continuous with respect to the volume measure. Later Amp\`ere has
studied an analogous question about the electricity current in a media
with varying conductivity. In modern language of measure theory we can
express the problem in the following terms: let $W$ be a Polish space
on which are given two positive measures $\rho$ and $\nu$,  of finite,
equal mass. Let $c(x,y)$ be a cost function on $W\times W$, which is,
usually, assumed  positive. Does there  exist a map $T:W\to W$ such that
$T\rho=\nu$ and $T$ minimizes the integral 
$$
\int_W c(x,T(x))d\rho(x)
$$
between all such maps? The problem has been further studied by Appell
 \cite{App-1,App-2} and by Kantorovitch \cite{Kan}. Kantarovitch has
 succeeded to transform this highly nonlinear problem of Monge into a
 linear problem by replacing the search for $T$ with the search of a
 measure $\ga$ on $W\times W$ with marginals $\rho$ and $\nu$  such that the
 integral 
$$
\int_{W\times W}c(x,y)d\ga(x,y)
$$
is the minimum  of all the integrals 
$$
\int_{W\times W}c(x,y)d\beta(x,y)
$$
where $\beta$ runs in the set of measures on $W\times W$ whose
marginals are $\rho$ and $\nu$. Since then the problem adressed above
is called the  Monge problem and the quest of the optimal
measure is called the  Monge-Kantorovitch problem.

In this chapter we study the Monge-Kantorovitch and the 
Monge   problem in the frame of an abstract  Wiener
space with a singular cost. In other words, let  $W$ be a separable
Fr\'echet space with its Borel sigma algebra 
$\calB(W)$ and assume that there is a separable  Hilbert space $H$ which is
injected densely and continuously into $W$, hence in general the
topology of $H$ is stronger than the topology induced by $W$. The cost
function $c:W\times W\to \reals_+\cup\{\infty\}$ is defined  as 
$$
c(x,y)=|x-y|_H^2\,,
$$
we suppose that $c(x,y)=\infty$ if $x-y$ does not belong to
$H$. Clearly, this
choice of the function $c$ is not arbitrary, in fact it is closely
related to Ito Calculus, hence also to the problems originating from 
Physics, quantum chemistry, large deviations, etc. Since
for all the interesting measures on $W$, the Cameron-Martin space is a
negligeable set, the cost function will be infinity very frequently. 
Let  $\Sigma(\rho,\nu)$ denote  the set  of probability measures 
on $W\times W$ with given marginals $\rho$ and $\nu$. It is a convex,  compact
 set under the weak topology $\sigma(\Sigma,C_b(W\times W))$. As
 explained above, the problem of Monge consists of finding a measurable map 
$T:W\to W$, called the optimal transport of $\rho$ to $\nu$, i.e.,
$T\rho=\nu$\footnote{We denote the
  push-forward of $\rho$ by $T$, i.e., the image of $\rho$ under $T$,  by $T\rho$.} which minimizes
the cost
$$
U\to  \int_W|x-U(x)|_H^2d\rho(x)\,,
$$
between all the maps $U:W\to W$ such that $U\rho=\nu$. The 
Monge-Kantorovitch  problem will   consist of
finding a  measure on $W\times W$,  which minimizes the function
$\theta\to J(\theta)$, defined by
\begin{equation}
\label{J-defn}
J(\theta)=\int_{W\times W}\big|x-y\big|_H^2d\theta(x,y)\,,
\end{equation}
where $\theta$ runs in $\Sigma(\rho,\nu)$. Note that 
$\inf\{J(\theta):\,\theta\in \Sigma(\rho,\nu)\}$ is the square of  Wasserstein
metric $d_H(\rho,\nu)$ with respect to the Cameron-Martin space $H$. 

Any solution $\ga$ of the
 Monge-Kantorovitch problem will give a solution to  the
Monge problem  provided that its support is included in
the graph of a map. Hence our work consists of realizing this program.
Although in  the finite dimensional case this problem is well-studied
in the path-breaking papers of Brenier \cite{BRE} and McCann \cite{Mc1,Mc2}
the things do not come up easily in our setting  and the  difficulty
is due to the 
fact that the cost function is not continuous with respect to the
Fr\'echet topology of $W$, for instance  the weak convergence of the
probability measures does not imply the convergence of the integrals  of
the cost function.  In other words the
function $|x-y|_H^2$ takes the value plus  infinity ``very
often''. On the other hand the results we obtain  seem to
have important applications to several problems of stochastic
analysis that we shall explain while enumerating the contents of this 
chapter. 

Section \ref{inequalities} is
devoted to the derivation of some inequalities which control the
Wasserstein distance. In particular, with the help of  the Girsanov
theorem, we give a very simple proof of an
inequality, initially  discovered by Talagrand (\cite{Tal}); this
facility gives already an idea about the efficiency of 
the infinite dimensional techniques for the Monge-Kantorovitch
problem{\footnote{In Section \ref{equation} we shall see another
    illustration of this phenomena.}}.  We indicate some simple
consequences of this inequality to 
control the measures of subsets of the Wiener space with respect to  second
moments of their  gauge functionals defined with  the
Cameron-Martin distance. These inequalities  are quite useful in the
theory of large deviations. Using a different representation of the
target measure, namely by constructing a flow of diffeomorphisms of
the Wiener space (cf. Chapter V of \cite{UZ7}) which maps the Wiener
measure to the target measure, we obtain also  a new control of the
Kantorovitch-Rubinstein 
metric of order one. The method we employ for this inequality
generalizes directly to a more general class of measures, namely those
for which one can define a reasonable divergence operator.

In Section \ref{gaussian}, we solve directly  the original problem of
Monge when the first measure is the Wiener measure and the second one is
given with a density, in such a way that the Wasserstein distance
between these two measures is finite. We prove the existence and the
uniqueness  of a
transformation of $W$ of the form $T=I_W+\nabla\phi$, where $\phi$ is a
$1$-convex function in the Gaussian Sobolev space $\DD_{2,1}$ such
that the measure $\ga=(I_W\times T)\mu$ is the unique solution of the
 problem of Monge-Kantorovitch.  This result gives
a new insight to the question  of representing an  integrable, positive
random variable whose expectation is unity,  as the Radon-Nikodym
derivative  of the image of the Wiener measure under a map which is a
perturbation of identity, a problem which has been studied by X. Fernique
and by one of us with M. Zakai (cf., \cite{Fer1,Fer2,UZ7}).  In \cite{UZ7},
Chapter II,  it is shown that such random 
variables are dense in $L^1_{1,+}(\mu)$ (the lower index $1$ means
that the expectations are equal to one), here we prove that this set
of random variables contains  the random variables who are at finite
Wasserstein distance from the Wiener measure. In fact even if this
distance is infinite, we show that there is a solution to this problem
if we  enlarge  $W$ slightly by taking $\NN\times W$. 

Section \ref{factorization} is devoted to the immediate
implications of the existence and the  uniqueness  of the solutions of 
Monge-Kantorovitch and Monge  problems constructed in Section \ref{gaussian}. 
Indeed the uniqueness  implies at once that the absolutely continuous
transformations of the Wiener space, at finite (Wasserstein) distance,
have a unique decomposition in the sense that they can be written as
the composition of a measure preserving map in the form of the
perturbation of identity with another one which is the perturbation of
identity with the Sobolev derivative of a $1$-convex function. This
means in particular that the class  of  $1$-convex functions is  as
basic  as the class  of adapted processes in the setting  of Wiener space.

In Section \ref{general-case} we prove the existence and the
uniqueness of solutions of the   Monge-Kantorovitch and Monge 
problems  for the measures which are at finite Wasserstein distance  from
each other. The fundamental  hypothesis  we use  is that the regular
conditional probabilities which are obtained by the disintegration of
one of the measures along the orthogonals of a sequence of regular, finite
dimensional projections vanish on the sets of co-dimension one. In
particular, this hypothesis  is satisfied if the measure under
question is absolutely continuous with respect to the Wiener measure.
The method  we use in this section is totally  different from the one
of Section \ref{gaussian}; it is  based on the notion of 
cyclic monotonicity of the supports of the regular conditional
probabilities obtained through some specific  disintegrations of the
optimal measures. The importance of cyclic monotonicity has first been
remarked by McCann and used abundently in \cite{Mc1} and in
\cite{G-Mc} for the finite dimensional case. Here the things are much
more complicated due to the singularity of the cost function, in
particular, contrary to the finite dimensional case,  the cyclic
monotonicity  is not compatible with the weak convergence of
probability measures.  A curious reader may ask why we did not treat
first the general case and then attack the subject of Section
\ref{gaussian}. The answer is twofold: even if we had done so, we
would have needed similar calculations as in Section \ref{gaussian} in
order to show the Sobolev regularity of the transport map, hence
concerning the volume, the order that we have chosen does not change
anything. Secondly, the construction used in  Section \ref{gaussian}
has  an interest by itself since it explains interesting  relations
between the transport map and its inverse and the optimal measure in a
more detectable situation, in this sense this construction  is rather
complementary to the material of Section \ref{general-case}.

Section \ref{equation} studies the Monge-Amp\`ere equation for the
measures which are absolutely continuous with respect to the Wiener
measure. First we briefly indicate the notion of  second order
Alexandroff derivative   and the Alexandroff version of the
Ornstein-Uhlenbeck operator applied to a $1$-convex function   in the
finite dimensional case. With the help of these observations, we write the
corresponding Jacobian using the modified Carleman-Fredholm
determinant which is natural in the infinite dimensional case (cf.,
\cite{UZ7}). Afterwards we attack the infinite dimensional case by
proving that the absolutely continuous part of the Ornstein-Uhlenbeck
operator applied to the finite rank conditional expectations of the
transport function is a submartingale which converges almost
surely. Hence the only difficulty lies in the calculation of the limit
of the Carleman-Fredholm determinants. Here we have a major difficulty
which originates from the pathology  of the
Radon-Nikodym derivatives of the vector measures with respect to a
scalar measure as explained in \cite{THOMAS}: in fact even if the
second order Sobolev  derivative of a Wiener function is a vector
measure with values in the space of Hilbert-Schmidt operators, its
absolutely continuous part has no reason to be Hilbert-Schmidt. Hence
the Carleman-Fredholm determinant may not exist, however 
due  to the $1$-convexity, the detereminants of the approximating
sequence are all with values in the interval $[0,1]$. Consequently  we
can construct the subsolutions with the help of the  Fatou lemma. 

Last but not the least, in section \ref{sub-equation}, we prove that
all these difficulties can be overcome thanks to the natural
renormalization of the Ito stochastic calculus. In fact using the Ito
representation theorem and the Wiener space analysis extended to the
distributions,  we can give the explicit solution of the
Monge-Amp\`ere equation. This is a remarkable result in the sense that
such techniques do not exist in the finite dimensional case.


\section{Preliminaries and notations}
\label{preliminaries}
Let $W$ be a separable Fr\'echet space equipped with  a Gaussian
measure $\mu$ of zero mean whose support is the whole space. The
corresponding Cameron-Martin space is denoted by $H$. Recall that the
injection $H\hookrightarrow W$ is compact and its adjoint is the
natural injection $W^\star\hookrightarrow H^\star\subset
L^2(\mu)$. The triple $(W,\mu,H)$ is called 
an abstract Wiener space. Recall that $W=H$ if and only if $W$ is
finite dimensional. A subspace $F$ of $H$ is called regular if the
corresponding orthogonal projection 
has a continuous extension to $W$, denoted again  by the same letter.
It is well-known that there exists an increasing sequence of regular
subspaces $(F_n,n\geq 1)$, called total,  such that $\cup_nF_n$ is
dense in $H$ and in $W$. Let $\sigma(\pi_{F_n})${\footnote{For the notational
  simplicity, in the sequel we shall denote  it by  $\pi_{F_n}$.}}  be the
$\sigma$-algebra generated by $\pi_{F_n}$, then  for any  $f\in
L^p(\mu)$, the martingale  sequence 
$(E[f|\sigma(\pi_{F_n})],n\geq 1)$
converges to $f$ (strongly if 
$p<\infty$) in $L^p(\mu)$. Observe that the function
$f_n=E[f|\sigma(\pi_{F_n})]$ can be identified with a function on the
finite dimensional abstract Wiener space $(F_n,\mu_n,F_n)$, where
$\mu_n=\pi_n\mu$.

Let us recall some facts from the convex analysis. Let $K$ be a 
Hilbert space, a subset $S$ of $K\times K$ is called cyclically
monotone if  any finite subset 
$\{(x_1,y_1),\ldots,(x_N,y_N)\}$ of
$S$ satisfies the following algebraic condition:
$$
\langle y_1,x_2-x_1\rangle+\langle y_2,x_3-x_2\rangle+\cdots+\langle
y_{N-1},x_N-x_{N-1}\rangle+\langle y_N,x_1-x_N\rangle\leq 0\,,
$$
where $\langle\cdot,\cdot\rangle$ denotes the inner product of
$K$. It turns out  that $S$ is
cyclically monotone if and only if 
$$
\sum_{i=1}^N(y_i,x_{\sigma(i)}-x_i)\leq 0\,,
$$
for any permutation $\sigma$ of $\{1,\ldots,N\}$ and for any finite
subset $\{(x_i,y_i):\,i=1,\ldots,N\}$ of $S$.
Note that  $S$ is
cyclically monotone if and only if any translate of it is cyclically
monotone.  By a theorem of Rockafellar,  any cyclically monotone set is
contained in the graph of the subdifferential  of a convex function in the
sense of convex analysis (\cite{ROC}) and even if the function may not
be unique its subdifferential  is unique. 

\noindent
Let now  $(W,\mu,H)$ be an abstract Wiener space;  a measurable  function
$f:W\to \reals\cup\{\infty\}$  is called $1$-convex if the map 
$$
h\to f(x+h)+\frac{1}{2}|h|_H^2=F(x,h)
$$
is convex on the Cameron-Martin space $H$ with values in
$L^0(\mu)$. Note that this notion is compatible with the
$\mu$-equivalence classes of random variables thanks to the
Cameron-Martin theorem. It is proven in Chapter  \ref{ch.convex} that 
this definition  is equivalent  the following condition:
Let $(\pi_n,n\geq 1)$ be a sequence of regular, finite dimensional,
 orthogonal projections of 
  $H$,  increasing to the identity map
  $I_H$. Denote also  by $\pi_n$ its  continuous extension  to $W$ and
  define $\pi_n^\bot=I_W-\pi_n$. For $x\in W$, let $x_n=\pi_nx$ and
  $x_n^\bot=\pi_n^\bot x$.   Then $f$ is $1$-convex if and only if 
$$
x_n\to \frac{1}{2}|x_n|_H^2+f(x_n+x_n^\bot)
$$ 
is  $\pi_n^\bot\mu$-almost surely convex.


\section{Some Inequalities}
\label{inequalities}
\begin{definition}
Let $\xi$ and $\eta$ be two probabilities on $(W,\calB(W))$. We say
that a probability $\ga$ on $(W\times W,\calB(W\times W))$ is a
solution of the Monge-Kantorovitch problem associated to the
couple $(\xi,\eta)$ if the first marginal of $\ga$ is $\xi$, the
second one is $\eta$ and if 
$$
J(\ga)=\int_{W\times W}|x-y|_H^2d\ga(x,y)=\inf\left\{\int_{W\times
  W}|x-y|_H^2d\beta(x,y):\,\beta\in \Sigma(\xi,\eta)\right\}\,,
$$
where $\Sigma(\xi,\eta)$ denotes the set of all the probability
measures on $W\times W$ whose first and second marginals are
respectively $\xi$ and $\eta$. We shall denote the Wasserstein
distance between $\xi$ and $\eta$, which is  the positive
square-root of  this infimum, with $d_H(\xi,\eta)$. 
\end{definition}
\remark
Since the set of probability measures on $W\times W$ is weakly compact
and since the integrand in the definition is lower semi-continuous and
strictly convex,  the infimum in the definition is always  attained
even if the functional $J$ is identically infinity.

\noindent
The following result is an extension of an inequality due to Talagrand
\cite{Tal} and it gives a sufficient condition for the Wasserstein
distance to be finite:
\begin{theorem}
\label{ineq-thm}
Let $L\in \LL\log\LL(\mu)$ be a positive random variable with
$E[L]=1$ and let  $\nu$ be  the measure $d\nu=Ld\mu$.  We then  have  
\begin{equation}
\label{tal-ineq}
d_H^2(\nu,\mu)\leq 2E[L\log L]\,.
\end{equation}
\end{theorem}
\proof
Without loss of generality, we may suppose  that $W$ is equipped with
a filtration of sigma algebras in 
such a way that it becomes a classical Wiener space as
$W=C_0(\reals_+,\R^d)$. Assume first that  $L$ is  a strictly positive
and bounded  random variable. We can represent
it as 
$$
L=\exp\left[-\int_0^\infty(\dot{u}_s,dW_s)-\frac{1}{2}|u|_H^2\right]\,,
$$
where $u=\int_0^\cdot \dot{u}_sds$ is an $H$-valued, adapted random
variable. Define $\tau_n$ as
$$
\tau_n(x)=\inf\left\{t\in \reals_+:\,\int_0^t|\dot{u}_s(x)|^2ds>n\right\}\,.
$$
$\tau_n$ is a stopping time with respect to the canonical filtration
$(\calF_t,t\in \reals_+)$ of the Wiener process $(W_t,t\in \reals_+)$
and $\lim_n\tau_n=\infty$ almost surely. Define $u^n$ as
$$
u^n(t,x)=\int_0^t\won_{[0,\tau_n(x)]}(s)\dot{u}_s(x)ds\,.
$$
Let  $U_n:W\to W$ be  the map
$U_n(x)=x+u^n(x)$, then the Girsanov 
theorem says that $(t,x)\to U_n(x)(t)=x(t)+\int_0^t\dot{u}^n_sds$ is a
Wiener process under the measure $L_nd\mu$, where
$L_n=E[L|\calF_{\tau_n}]$. Therefore 
\beaa
E[L_n\log L_n]&=&E\left[L_n\,\left\{-\int_0^\infty
      (\dot{u}^n_s,dW_s)-\frac{1}{2}|u^n|_H^2\right\}\right]\\
&=&\frac{1}{2}E[L_n|u^n|_H^2]\\
&=&\frac{1}{2}E[L|u^n|_H^2]\,.
\eeaa
Define now the measure $\beta_n$ on $W\times W$ as 
$$
\int_{W\times W} f(x,y)d\beta_n(x,y)=\int_W f(U_n(x),x)L_n(x)d\mu(x)\,.
$$
Then the first marginal of $\beta_n$ is $\mu$ and the second one is
$L_n.\mu$. Consequently
\beaa
\lefteqn{\inf\left\{\int_{W\times
  W}|x-y|_H^2d\theta:\pi_1\theta=\mu,\,\pi_2\theta=L_n.\mu\right\}}\\
&\leq& \int_W|U_n(x)-x|_H^2L_nd\mu\\
&=&2E[L_n\log L_n]\,.
\eeaa
Hence we obtain 
$$
d_H^2(L_n.\mu,\mu)=J(\ga_n)\leq 2E[L_n\log L_n]\,,
$$
where $\ga_n$ is a solution of the Monge-Kantorovitch  problem in
$\Sigma(L_n.\mu,\mu)$. Let now $\ga$ be any cluster point of the
sequence $(\ga_n,n\geq 1)$, since $\ga\to J(\ga)$ is lower
semi-continuous with respect to the weak topology of probability
measures, we have 
\beaa
J(\ga)&\leq& \lim\inf_nJ(\ga_n)\\
&\leq& \sup_n2E[L_n\log L_n] \\
&\leq&2E[L\log L]\,,
\eeaa
since $\ga\in \Sigma(L.\mu,\mu)$, it follows that 
$$
d_H^2(L.\mu,\mu)\leq 2E[L\log L]\,.
$$
For the general case we stop the martingale $E[L|\calF_t]$
appropriately to obtain a bounded density $L_n$, then  replace it  by
$P_{1/n} L_n$ to improve the positivity, where $(P_t,t\geq 0)$ denotes
the Ornstein-Uhlenbeck semigroup. Then, from the Jensen inequality,
$$
E[P_{1/n}L_n\log P_{1/n}L_n]\leq E[L\log L]\,,
$$
therefore, using the same reasoning as above 
\beaa
d_H^2(L.\mu,\mu)&\leq&\lim\inf_nd_H^2(P_{1/n}L_n.\mu,\mu)\\
&\leq& 2E[L\log L]\,,
\eeaa
and this completes the proof.
\qed

\begin{corollary}
\label{tri-cor}
Assume that $\nu_i\,(i=1,2)$ have Radon-Nikodym  densities
$L_i\,(i=1,2)$ with respect to the Wiener measure $\mu$ which are in
$\LL\log\LL$. Then
$$
d_H(\nu_1,\nu_2)<\infty\,.
$$
\end{corollary}
\proof
This is a simple consequence of the triangle  inequality (cf. \cite{B-F}):
$$
d_H(\nu_1,\nu_2)\leq d_H(\nu_1,\mu)+d_H(\nu_2,\mu)\,.
$$
\qed

Let us give a simple application of the above result in the lines of
\cite{Mar}: 
\begin{corollary}
\label{iso-cor}
Assume that $A\in \calB(W)$ is any set of positive Wiener
measure. Define the $H$-gauge function of $A$ as 
$$
q_A(x)=\inf(|h|_H:\,h\in (A-x)\cap H)\,.
$$
Then we have 
$$
E[q_A^2]\leq 2\log\frac{1}{\mu(A)}\,,
$$
in other words
$$
\mu(A)\leq \exp\left\{-\frac{E[q_A^2]}{2}\right\}\,.
$$
Similarly if $A$ and $B$ are $H$-separated, i.e., if $A_\eps\cap
B=\emptyset$,  for some $\eps>0$, where  $A_\eps=\{x\in
W:\,q_A(x)\leq\eps\}$, then  
$$
\mu(A_\eps^c)\leq \frac{1}{\mu(A)}e^{-\eps^2/4}
$$
and consequently 
$$
\mu(A)\,\mu(B)\leq \exp\left(-\frac{\eps^2}{4}\right)\,.
$$
\end{corollary}
\remark We already know that, from the $0-1$--law,  $q_A$ is almost
surely finite, besides it satisfies $|q_A(x+h)-q_A(x)|\leq |h|_H$, hence
$E[\exp\la q_A^2]<\infty$ for any $\la<1/2$ (cf. \cite{UZ7}). In fact
all these assertions can also be proved with the technique used below.

\proof
Let $\nu_A$ be the measure defined by 
$$
d\nu_A=\frac{1}{\mu(A)}1_A d\mu\,.
$$
Let $\ga_A$ be the solution of the Monge-Kantorovitch problem, it is
easy to see that the support of $\ga_A$ is included in $W\times A$,
hence 
$$
|x-y|_H\geq \inf\{|x-z|_H:\,z\in A\}=q_A(x)\,,
$$
$\ga_A$-almost surely. This implies in particular that $q_A$ is almost
surely finite. It follows now from the inequality
(\ref{tal-ineq})
$$
E[q_A^2]\leq -2\log\mu(A)\,,
$$
hence the proof of the first inequality  follows. For the second let
$B=A_\eps^c$ and let $\ga_{AB}$ be the solution of the
Monge-Kantorovitch problem corresponding to $\nu_A,\nu_B$. Then we
have from the Corollary \ref{tri-cor}, 
$$
d^2_H(\nu_A,\nu_B)\leq -4\log\mu(A)\mu(B)\,.
$$
Besides the support of the measure $\ga_{AB}$ is in $A\times B$, hence
$\ga_{AB}$-almost surely $|x-y|_H\geq \eps$ and the proof follows.
\qed

\noindent
For the distance defined by 
$$
d_1(\nu,\mu)=\inf\left\{\int_{W\times
  W}|x-y|_H d\theta:\pi_1\theta=\mu,\,\pi_2\theta=\nu\right\}
$$
we have the following  control:

\begin{theorem}
\label{our-thm}
Let $L\in \LL_+^1(\mu)$ with $E[L]=1$. Then we have
\begin{equation}
\label{our-ineq}
d_1(L.\mu,\mu)\leq E\left[\left|(I+\calL)^{-1}\nabla L\right|_H\right]\,.
\end{equation}
\end{theorem}
\proof To prove the theorem we shall use a technique developed in
\cite{Da-M}. 
Using the conditioning with respect to the sigma algebra
$V_n=\sigma\{\delta e_1,\ldots,\delta e_n\}$, where $(e_i,i\geq 1)$ is
a complete, orthonormal basis of $H$, 
we reduce the problem to the finite dimensional case. Moreover, we can
assume that $L$ is a smooth, strictly positive  function on
$\reals^n$. Define now $\sigma=(I+\calL)^{-1}\nabla L$ and 
$$
\sigma_t(x)=\frac{\sigma(x)}{t+(1-t)L}\,,
$$
for $t\in [0,1]$. Let $(\phi_{s,t}(x),s\leq t\in[0,1])$ be the flow of
diffeomorphisms defined by the following differential equation:
$$
\phi_{s,t}(x)=x-\int_s^t\sigma_\tau(\phi_{s,\tau}(x))d\tau\,.
$$
From the standart results (cf. \cite{UZ7}, Chapter V), it follows that
$x\to\phi_{s,t}(x)$ is Gaussian under the probability $\La_{s,t}.\mu$,
where 
$$
\La_{s,t}=\exp\int_s^t(\delta\sigma_\tau)(\phi_{s,\tau}(x))d\tau
$$
is the Radon-Nikodym density of $\phi_{s,t}^{-1}\mu$ with respect to $\mu$.
Define 
$$
H_s(t,x)=\La_{s,t}(x)\left\{t+(1-t)L\circ\phi_{s,t}(x)\right\}\,.
$$
It is easy to see that 
$$
\frac{d}{dt}H_s(t,x)=0
$$
for  $t\in(s,1)$. Hence the map $t\to H_s(t,x)$ is a constant, this
implies that 
$$
\La_{s,1}(x)=s+(1-s)L(x)\,.
$$
We have, as in the proof of Theorem \ref{ineq-thm}, 
\beaa
d_1(L.\mu,\mu)&\leq&E[|\phi_{0,1}(x)-x|_H\La_{0,1}]\\
&\leq&E\left[\La_{0,1}\int_0^1|\sigma_t(\phi_{0,t}(x))|_Hdt\right]\\
&=&E\left[\int_0^1\left|\sigma_t(\phi_{0,t}\circ\phi_{0,1}^{-1})(x)\right|_Hdt\right]\\
&=&E\left[\int_0^1\left|\sigma_t(\phi_{t,1}^{-1}(x))\right|_Hdt\right]\\
&=&E\left[\int_0^1|\sigma_t(x)|_H\La_{t,1}dt\right]\\
&=&E[|\sigma|_H]\,,
\eeaa
and the general case follows via the usual approximation procedure.
\qed


\section{Construction of the transport map}
\label{gaussian}
In this section we give  the construction of the transport map in  the
Gaussian case. We begin with the following lemma:
\begin{lemma}
\label{integ-lemma}
Let $(W,\mu,H)$ be an abstract Wiener space, assume that $f:W\to
\reals$ is a measurable function such that it is G\^ateaux  differentiable
in the direction of the Cameron-Martin space $H$, i.e., there exists
some $\nabla f:W\to H$ such that 
$$
f(x+h)=f(x)+\int_0^1(\nabla f(x+\tau h),h)_Hd\tau\,,
$$
$\mu$-almost surely, for any $h\in H$.
If $|\nabla f|_H\in L^2(\mu)$, then $f$ belongs to the Sobolev space
$\DD_{2,1}$. 
\end{lemma}
\proof 
Since $|\nabla |f||_H\leq |\nabla f|_H$, we can assume that $f$ is
positive. Moreover, for any $n\in \NN$, the function $f_n=\min(f,n)$
has also a G\^ateaux  derivative such that $|\nabla f_n|_H\leq |\nabla
f|_H$ $\mu$-almost surely. It follows from the Poincar\'e inequality
that the sequence $(f_n-E[f_n],n\geq 1)$ is bounded in $L^2(\mu)$,
hence it is also bounded in $L^0(\mu)$. Since $f$ is almost surely
finite, the sequence $(f_n, n\geq 1)$ is bounded in $L^0(\mu)$,
consequently the deterministic sequence $(E[f_n],n\geq 1)$ is also
bounded in $L^0(\mu)$. This means that $\sup_nE[f_n]<\infty$, 
hence  the monotone convergence theorem implies that  $E[f]<\infty$
and the proof is completed.
\qed

\begin{theorem}
\label{gaussian-case}
Let $\nu$ be the measure $d\nu=Ld\mu$, where $L$ is a positive random variable,
with $E[L]=1$. Assume that $d_H(\mu,\nu)<\infty$ (for instance 
$L\in \LL\log \LL$). Then there exists a  $1$-convex function $\phi\in
\DD_{2,1}$, unique upto a constant,  such that  the map
$T=I_W+\nabla \phi$ is the unique solution of the original problem of
Monge. Moreover, its graph supports  the  unique
solution of the 
Monge-Kantorovitch problem $\ga$. Consequently    
$$
(I_W\times T)\mu=\ga
$$
In particular  $T$ maps $\mu$ to $\nu$ and  $T$ is almost surely
invertible, i.e., there exists some $T^{-1}$ such that $T^{-1}\nu=\mu$
and that 
\beaa
1&=&\mu\left\{x:\,T^{-1}\circ T(x)=x\right\}\\
&=&\nu\left\{y\in W:\,T\circ T^{-1}(y)=y\right\}\,.
\eeaa
\end{theorem}
\proof Let $(\pi_n,n\geq 1)$ be a sequence of regular, finite
dimensional orthogonal projections of $H$ increasing to $I_H$. Denote
their continuous extensions to $W$ by the same letters. For $x\in W$,
we define $\pi_n^\bot x=:x_n^\bot=x-\pi_nx$. Let  $\nu_n$ be  the measure
$\pi_n\nu$. Since $\nu$ is absolutely continuous with respect to
$\mu$, $\nu_n$ is absolutely continuous with respect to
$\mu_n:=\pi_n\mu$ and 
$$
\frac{d\nu_n}{d\mu_n}\circ\pi_n=E[L|V_n]=:L_n\,,
$$
where $V_n$ is the sigma algebra  $\sigma(\pi_n)$ and the conditional
expectation is taken with 
respect to $\mu$. On the space $H_n$, the  Monge-Kantorovitch
problem, which consists of finding the probability measure which
realizes the following infimum 
$$
d_H^2(\mu_n,\nu_n)=\inf\left\{J(\beta):\,\beta\in M_1(H_n\times
H_n)\,,p_1\beta=\mu_n,p_2\beta=\nu_n\right\} 
$$
where
$$
J(\beta)=\int_{H_n\times H_n}|x-y|^2d\beta(x,y)\,,
$$
has a unique solution $\ga_n$, where $p_i,\,i=1,2$ denote  the
projections  $(x_1,x_2)\to x_i,\,i=1,2$ from $H_n\times H_n$ to $H_n$
and $M_1(H_n\times H_n)$ denotes the set of probability measures on
$H_n\times H_n$.  The measure $\ga_n$ may be regarded
as a measure on $W\times W$, by taking its image under the injection
$H_n\times H_n\hookrightarrow W\times W$ which  we shall denote again by
$\ga_n$. It results from the finite dimensional results of Brenier and
of McCann(\cite{BRE}, \cite{Mc1}) that 
there are two convex continuous  functions (hence almost everywhere
differentiable) $\Phi_n$ and $\Psi_n$ on 
$H_n$ such that 
$$
\Phi_n(x)+\Psi_n(y)\geq (x,y)_H
$$
for all $x,y\in H_n$ and that 
$$
\Phi_n(x)+\Psi_n(y)= (x,y)_H
$$
$\ga_n$-almost everywhere. Hence the 
support of $\ga_n$ is included in the graph of the derivative
$\nabla\Phi_n$  of $\Phi_n$, hence $\nabla\Phi_n\mu_n=\nu_n$ and the inverse of
$\nabla\Phi_n$ is equal to $\nabla\Psi_n$. Let 
\beaa
\phi_n(x)&=&\Phi_n(x)-\frac{1}{2}|x|_H^2\\
\psi_n(y)&=&\Psi_n(y)-\frac{1}{2}|y|_H^2\,.
\eeaa
Then $\phi_n$ and $\psi_n$ are $1$-convex functions and they satisfy
the following relations:
\begin{equation}
\label{+-ineq}
\phi_n(x)+\psi_n(y)+\frac{1}{2}|x-y|_H^2\geq 0\,,
\end{equation}
for all $x,y\in H_n$ and 
\begin{equation}
\label{=-eq}
\phi_n(x)+\psi_n(y)+\frac{1}{2}|x-y|^2_H=0\,,
\end{equation}
$\ga_n$-almost everywhere. From what we have said above, it follows
that $\ga_n$-almost surely  $y=x+\nabla\phi_n(x)$, consequently
\begin{equation}
\label{energy-iden}
J(\ga_n)=E[|\nabla\phi_n|_H^2]\,.
\end{equation}
Let $q_n:W\times W\to H_n\times H_n$ be defined as
$q_n(x,y)=(\pi_nx,\pi_ny)$. If $\ga$ is any solution of the 
Monge-Kantorovitch problem, then $q_n\ga\in \Sigma(\mu_n,\nu_n)$,
hence
\begin{equation}
\label{1st-estimate}
J(\ga_n)\leq J(q_n\ga)\leq J(\ga)=d_H^2(\mu,\nu)\,.
\end{equation}
Combining the relation (\ref{energy-iden}) with the inequality
(\ref{1st-estimate}), we obtain the following bound
\bea
\label{control-ineq}
\sup_nJ(\ga_n)&=&\sup_nd_H^2(\mu_n,\nu_n)\nonumber\\
&=&\sup_nE[|\nabla\phi_n|_H^2]\nonumber\\
&\leq&d_H^2(\mu,\nu)=J(\ga)\,.
\eea
For $m\leq n$, $q_m\ga_n\in\Sigma(\mu_m,\nu_m)$, hence  we should have 
\beaa
J(\ga_m)&=&\int_{W\times W}|\pi_mx-\pi_my|_H^2d\ga_m(x,y)\\
&\leq&\int_{W\times W}|\pi_mx-\pi_my|_H^2d\ga_n(x,y)\\
&\leq&\int_{W\times W}|\pi_nx-\pi_ny|_H^2d\ga_n(x,y)\\
&=&\int_{W\times W}|x-y|_H^2d\ga_n(x,y)\\
&=&J(\ga_n)\,,
\eeaa
where the third equality follows from the fact that we have denoted
the $\ga_n$ on $H_n\times H_n$ and its image in $W\times W$ by the
same letter.
Let now $\ga$ be a weak cluster point of the sequence of measures
$(\ga_n,n\geq 1)$, where the word `` weak''\footnote{To prevent the
  reader against the trivial errors let us emphasize  that $\ga_n$ is
  not the projection of $\ga$ on $W_n\times W_n$.} refers to the weak
convergence of 
measures on $W\times W$. Since $(x,y)\to |x-y|_H$ is lower
semi-continuous, we have 
\beaa
J(\ga)&=&\int_{W\times W}|x-y|_H^2d\ga(x,y)\\
&\leq&\lim\inf_n\int_{W\times W}|x-y|_H^2d\ga_n(x,y)\\
&=&\lim\inf_nJ(\ga_n)\\
&\leq&\sup_nJ(\ga_n)\\
&\leq& J(\ga)=d_H^2(\mu,\nu)\,,
\eeaa
from the relation (\ref{control-ineq}). Consequently 
\be
\label{lim-eq}
J(\ga)=\lim_nJ(\ga_n)\,.
\ee
Again from (\ref{control-ineq}), if we replace $\phi_n$ with $\phi_n-E[\phi_n]$
and $\psi_n$ with $\psi_n+E[\phi_n]$ we obtain a bounded sequence
$(\phi_n,n\geq 1)$ in $\DD_{2,1}$, in particular it is bounded in the
space $L^2(\ga)$ if we inject it into  latter  by $\phi_n(x)\to
\phi_n(x)\otimes 1(y)$. Consider now the sequence of the 
positive, lower semi-continuous functions $(F_n,n\geq 1)$ defined on
$W\times W$ as   
$$
F_n(x,y)=\phi_n(x)+\psi_n(y)+\frac{1}{2}|x-y|_H^2\,.
$$
We have, from the relation (\ref{=-eq})
\beaa
\int_{W\times
  W}F_n(x,y)d\ga(x,y)&=&\int_W\phi_nd\mu+\int_W\psi_n(y)d\nu+\frac{1}{2}J(\ga)
\\ 
&=&\frac{1}{2}\left(J(\ga)-J(\ga_n)\right)\to 0\,.
\eeaa
Consequently the sequence $(F_n,n\geq 1)$ converges to zero in
$L^1(\ga)$, therefore it is uniformly integrable. Since $(\phi_n,n\geq
1)$ is uniformly integrable as explained above and since $|x-y|^2$ has
a finite expectation with respect to $\ga$, it follows that
$(\psi_n,n\geq 1)$ is also uniformly integrable in $L^1(\ga)$ hence
also in $L^1(\nu)$. Let $\phi'$ be a weak cluster point of
$(\phi_n,n\geq 1)$, then there exists a sequence $(\phi_n',n\geq 1)$
whose elements are the convex combinations of some  elements of
$(\phi_k,k\geq n)$ such that $(\phi_n',n\geq 1)$ converges in the norm
topology of $\DD_{2,1}$ and $\mu$-almost everywhere. Therefore the
 sequence $(\psi_n',n\geq 1)$, constructed from $(\psi_k,k\geq n)$,
 converges in $L^1(\nu)$ and $\nu$-almost surely. Define $\phi$ and
 $\psi$ as  
\beaa
\phi(x)&=&\lim\sup_n\phi_n'(x)\\
\psi(y)&=&\lim\sup_n\psi_n'(y)\,,
\eeaa
hence we have 
$$
G(x,y)=\phi(x)+\psi(y)+\frac{1}{2}|x-y|_H^2\geq 0
$$
for all $(x,y)\in W\times W$, also the  equality holds $\ga$-almost
everywhere. Let now $h$ be any element of $H$, since $x-y$ is in $H$
for $\ga$-almost all $(x,y)\in W\times W$, we have 
$$
|x+h-y|_H^2=|x-y|_H^2+|h|_H^2+2(h,x-y)_H
$$
$\ga$-almost surely. Consequently 
$$
\phi(x+h)-\phi(x)\geq -(h,x-y)_H-\frac{1}{2}|h|_H^2
$$
$\ga$-almost surely and  this implies that 
$$
y=x+\nabla \phi(x)
$$
$\ga$-almost everywhere. Define now the map $T:W\to W$ as
$T(x)=x+\nabla\phi(x)$, then 
\beaa
\int_{W\times W}f(x,y)d\ga(x,y)&=&\int_{W\times W}f(x,T(x))d\ga(x,y)\\
&=&\int_{W}f(x,T(x))d\mu(x)\,,
\eeaa
for any $f\in C_b(W\times W)$, consequently $(I_W\times T)\mu=\ga$, in
particular $T\mu=\nu$.

Let us notice
that any weak cluster point of $(\phi_n,n\geq 1)$, say $\tilde{\phi}$,
satisfies 
$$
\nabla\tilde{\phi}(x)=y-x
$$
$\ga$-almost surely, hence $\mu$-almost surely we have
$\tilde{\phi}=\phi$. This implies that $(\phi_n,n\geq 1)$ has a unique
cluster point $\phi$, consequently the sequence $(\phi_n,n\geq 1)$
converges weakly in $\DD_{2,1}$ to $\phi$. Besides we have
\beaa
\lim_n\int_W|\nabla\phi_n|_H^2d\mu&=&\lim_nJ(\ga_n)\\
&=&J(\ga)\\
&=&\int_{W\times W}|x-y|_H^2d\ga(x,y)\\
&=&\int_W|\nabla\phi|_H^2d\mu\,,
\eeaa
hence $(\phi_n,n\geq 1)$ converges to $\phi$ in the norm topology of
$\DD_{2,1}$.  Let us recapitulate what we
have done till here: we have taken an arbitrary optimal $\ga\in
\Sigma(\mu,\nu)$ and an arbitrary cluster point $\phi$ of
$(\phi_n,n\geq 1)$ and we have proved that $\ga$ is carried by the
graph of $T=I_W+\nabla\phi$. This implies that $\ga$ and $\phi$ are unique
and  that  the sequence $(\ga_n,n\geq 1)$ has a unique cluster point $\ga$.

Certainly  $(\psi_n,\geq 1)$
converges also in the norm topology of $L^1(\nu)$.
Moreover, from the finite dimensional situation, we have
$\nabla\phi_n(x)+\nabla\psi_n(y)=0$ $\ga_n$-almost everywhere. Hence
$$
E_\nu[|\nabla\psi_n|_H^2]=E[|\nabla\phi_n|_H^2]
$$
this implies the
boundedness of $(\nabla\psi_n,n\geq 1)$ in $L^2(\nu,H)$ (i.e.,
$H$-valued functions). 
To complete the proof
we have to show that, for some measurable, $H$-valued map, say
$\eta$, it holds that 
$x=y+\eta(y)$ $\ga$-almost surely. For this let 
$F$ be a finite dimensional, regular subspace of $H$ and denote by
$\pi_F$ the projection operator onto $F$ which is continuously
extended to $W$, put $\pi_F^\bot=I_W-\pi_F$. We have
$W=F\oplus F^\bot$, with $F^\bot=\ker \pi_F=\pi_F^\bot(W)$. Define the
measures $\nu_F=\pi_F(\nu)$ and $\nu_F^\bot=\pi_F^\bot(\nu)$. From the
construction of $\psi$, we know that, for 
any $v\in F^\bot$, the partial map $u\to \psi(u+v)$ is $1$-convex on
$F$. Let also $A=\{y\in W:\,\psi(y)<\infty\}$, then $A$ is a Borel set
with $\nu(A)=1$ and it is easy to see that, for  $\nu_F^\bot$-almost all
$v\in F^\bot$, one has 
$$
\nu(A|\pi_F^\bot=v)>0\,.
$$
It then follows from Lemma 3.4 of Chapter \ref{ch.convex}, and from the fact that
the regular conditional probability $\nu(\cdot\,|\pi_F^\bot=v)$ is
absolutely continuous with respect to the Lebesgue measure of $F$,
that $u\to \psi(u+v)$ is $\nu(\cdot\,|\pi_F^\bot=v)$-almost everywhere
differentiable on $F$ for $\nu_F^\bot$-almost
all $v\in F^\bot$. It then follows that, $\nu$-almost surely,  $\psi$
is differentiable in the directions of $F$, i.e., there exists
$\nabla_F\psi\in F$ $\nu$-almost surely. Since we also have 
$$
\psi(y+k)-\psi(y)\geq (x-y,k)_H-\frac{1}{2}|k|_H^2\,,
$$
we obtain, $\ga$-almost surely
$$
(\nabla_F\psi(y),k)_H=(x-y,k)_H\,,
$$
for any $k\in F$. Consequently 
$$
\nabla_F\psi(y)=\pi_F(x-y)
$$
$\ga$-almost surely. Let now $(F_n,n\geq 1)$ be a total, increasing
sequence of regular subspaces of $H$, we have a sequence
$(\nabla_n\psi,n\geq 1)$ bounded in $L^2(\nu)$ hence also bounded in
$L^2(\ga)$. Besides $\nabla_n\psi(y)=\pi_nx-\pi_ny$ $\ga$-almost
surely. Since $(\pi_n(x-y),n\geq 1)$ converges in $L^2(\ga,H)$,
$(\nabla_n\psi,n\geq 1)$ converges in the norm topology of
$L^2(\ga,H)$. Let us denote  this limit by  $\eta$, 
then  we have $x=y+\eta(y)$ $\ga$-almost surely. Note that, since
$\pi_n\eta=\nabla_n\psi$, we can even write in a weak sense that
$\eta=\nabla\psi$.  If we define $T^{-1}(y)=y+\eta(y)$, we
see that 
\beaa
1&=&\ga\{(x,y)\in W\times W:T\circ T^{-1}(y)=y\}\\
&=&\ga\{(x,y)\in W\times W:T^{-1}\circ T(x)=x\}\,,  
\eeaa
and this completes the proof of the theorem.
\qed

\begin{remarkk}
\label{nu-closed}
 Assume that the operator $\nabla$ is closable with respect to $\nu$,
then we have $\eta=\nabla\psi$. In particular, if  $\nu$ and $\mu$ are
equivalent, then we have   
$$
T^{-1}=I_W+\nabla\psi\,,
$$
where is $\psi$ is  a $1$-convex function.
\end{remarkk}
\begin{remarkk}
Assume that $L\in \LL_+^1(\mu)$, with $E[L]=1$ and let $(D_k,k\in
\NN)$ be a measurable partition of $W$ such that on each $D_k$, $L$ is
bounded. Define $d\nu=L\,d\mu$ and $\nu_k=\nu(\cdot|D_k)$. It follows
 from Theorem \ref{ineq-thm}, that $d_H(\mu,\nu_k)<\infty$. 
Let then  $T_k$ be the map constructed in Theorem \ref{gaussian-case}
satisfying $T_k\mu=\nu_k$. Define  $n(dk)$ as  the probability
distribution on $\NN$ given by  $n\left(\{k\}\right)=\nu(D_k),\,k\in
\NN$. Then we have  
$$
\int_W f(y)d\nu(y)=\int_{W\times \NN}f(T_k(x))\mu(dx)n(dk)\,.
$$
A  similar result is given in  \cite{Fer2}, the difference with that of
above lies in the fact that we have a more precise information about
the probability space on which $T$ is defined.
\end{remarkk}


\section{Polar factorization of the absolutely continuous
  transformations of the Wiener space}
\label{factorization}
Assume that $V=I_W+v:W\to W$ be an absolutely continuous  transformation and
let $L\in \LL_+^1(\mu)$ be the Radon-Nikodym derivative of $V\mu$ with
respect to $\mu$. Let
$T=I_W+\nabla\phi$ be the 
transport map such that $T\mu=L.\mu$. Then it is easy to see that the
map $s=T^{-1}\circ V$ is a rotation, i.e., $s\mu=\mu$
(cf. \cite{UZ7})  and it can be represented as $s=I_W+\alpha$. In
particular we have  
\begin{equation}
\label{iden-1}
\alpha+\nabla\phi\circ s=v\,.
\end{equation}
Since $\phi$ is a $1$-convex map, we have $h\to
\frac{1}{2}|h|_H^2+\phi(x+h)$ is almost surely convex (cf. Chapter\ref{ch.convex}). Let
$s'=I_W+\alpha'$ be another rotation with $\alpha':W\to H$. By the
$1$-convexity of $\phi$, we have 
$$
\frac{1}{2}|\alpha'|_H^2+\phi\circ s'\geq
\frac{1}{2}|\alpha|_H^2+\phi\circ s
+\Bigl(\alpha+\nabla\phi\circ s,\alpha'-\alpha\Bigr)_H\,,
$$
$\mu$-almost surely. 
Taking the expectation of both sides,  using the fact that $s$ and
$s'$ preserve the Wiener measure $\mu$ and the identity (\ref{iden-1}),
we obtain  
$$
E\left[\half|\alpha|_H^2-(v,\alpha)_H\right]\leq
E\left[\frac{1}{2}|\alpha'|_H^2-(v,\alpha')_H\right] \,.
$$
Hence we have proven the existence part of the following 
\begin{proposition}
Let $\mathcal R_2$ denote the subset of $L^2(\mu,H)$ whose elements
are defined by the property that  $x\to x+\eta(x)$ is a rotation,
i.e., it preserves the Wiener measure. Then $\alpha$ is the unique
element of $\mathcal R_2$ which minimizes the functional
$$
\eta\to M_v(\eta)=E\left[\frac{1}{2}|\eta|_H^2-(v,\eta)_H\right]\,.
$$
\end{proposition}
\proof
To show the uniqueness, assume that $\eta\in \mathcal R_2$ be another
map minimizing $J_v$. Let $\beta$ be the measure on $W\times W$,
defined as 
$$
\int_{W\times W} f(x,y)d\beta(x,y)=\int_W f(x+\eta(x),V(x))d\mu\,.
$$
Then the first marginal of $\beta$ is $\mu$ and the second marginal is
$L.\mu$. Since $\ga=(I_W\times T)\mu$ is the unique solution of the
Monge-Kantorovitch problem, we should have 
$$
\int|x-y|_H^2 d\beta(x,y)>\int
|x-y|_H^2d\ga(x,y)=E[|\nabla\phi|_H^2]\,.
$$
However we have 
\beaa
\int_{W\times W}|x-y|_H^2 d\beta(x,y)&=&E\left[|v-\eta|_H^2\right]\\
&=&E\left[|v|_H^2\right]+2M_v(\eta)\\
&=&E\left[|v|_H^2\right]+2M_v(\alpha)\\
&=&E\left[|v-\alpha|_H^2\right]\\
&=&E\left[|\nabla \phi\circ s|_H^2\right]\\
&=&E\left[|\nabla \phi|_H^2\right]\\
&=&\int_{W\times W}|x-y|_H^2 d\ga(x,y)\\
&=&J(\ga)
\eeaa
and this  gives a contradiction to the uniqueness of $\ga$.
\qed

The following theorem, whose proof is rather easy, gives a better
understanding of  the structure of  absolutely continuous
transformations of the Wiener measure:
\begin{theorem}
\label{pol-fact}
Assume that $U:W\to W$ be a measurable  map and $L\in \LL\log \LL$ a
positive random variable with $E[L]=1$. Assume that the measure
$\nu=L\cdot\mu$ is a Girsanov measure for $U$, i.e., that one has 
$$
E[f\circ U\,L]=E[f]\,,
$$
for any $f\in C_b(W)$. Then there exists a unique map
$T=I_W+\nabla\phi$ with $\phi\in \DD_{2,1}$ is $1$-convex, and a measure
preserving transformation $R:W\to W$ such that $U\circ T=R$
$\mu$-almost surely and $U=R\circ T^{-1}$ $\nu$-almost surely.
\end{theorem}
\proof
By Theorem \ref{gaussian-case} there is a unique map
$T=I_W+\nabla\phi$, with $\phi\in\DD_{2,1}$,  $1$-convex such that $T$
transports $\mu$ to $\nu$. Since $U\nu=\mu$, we have 
\beaa
E[f\circ U\,L]&=&E[f\circ U\circ T]\\
&=&E[f]\,.
\eeaa
Therefore $x\to U\circ T(x)$ preserves the measure $\mu$. The rest is
obvious since $T^{-1}$ exists $\nu$-almost surely.
\qed

Another version of Theorem \ref{pol-fact} can be announced as follows:
\begin{theorem}
\label{pol-fact2}
Assume that $Z:W\to W$ is a measurable map such that $Z\mu\ll\mu$,
with $d_H(Z\mu,\mu)<\infty$. Then $Z$ can be decomposed as 
$$
Z=T\circ s\,,
$$
where $T$ is the unique  transport map of the Monge-Kantorovitch
problem for $\Sigma(\mu,Z\mu)$ and $s$ is a rotation.
\end{theorem}
\proof
Let $L$ be the Radon-Nikodym derivative of $Z\mu$ with respect to 
$\mu$. We have, from Theorem \ref{gaussian-case},  
\beaa
E[f]&=&E[f\circ T^{-1}\circ T]\\
&=&E[f\circ T^{-1}\,L]\\
&=&E[f\circ T^{-1}\circ Z]\,,
\eeaa
for any $f\in C_b(W)$. Hence $T^{-1}\circ Z=s$ is a rotation. Since
$T$ is uniquely defined, $s$ is also uniquely defined.
\qed

Although the following result is a translation of the results of this section,
 it is interesting from the point of view of  stochastic differential
 equations:
\begin{theorem}
\label{weak-sde}
Let  $(W,\mu,H)$ be  the standard  Wiener space on $\reals^d$,
i.e., $W=C(\reals_+,\reals^d)$. Assume that there exists a
probability $P\ll\mu$ which is the weak solution of the stochastic
differential equation  
$$
dy_t=dW_t+b(t,y)dt\,,
$$
such that $d_H(P,\mu)<\infty$. Then there exists a process $(T_t,t\in
\reals_+)$ which is  a pathwise  solution of some stochastic
differential equation whose law is equal to $P$.
\end{theorem}
\proof
Let $T$ be the transport map constructed in Theorem
\ref{gaussian-case} corresponding to $dP/d\mu$. Then it has an inverse
$T^{-1}$ such that $\mu\{T^{-1}\circ T(x)=x\}=1$. Let $\phi$ be the
$1$-convex function  such that $T=I_W+\nabla \phi$ and denote by
$(D_s\phi,s\in \reals_+)$ the representation of $\nabla\phi$ in
$L^2(\reals_+,ds)$. Define $T_t(x)$ as the trajectory $T(x)$ evaluated
at $t\in \reals_+$. Then it is easy to see that $(T_t,t\in \reals_+)$
satifies the stochastic differential equation
$$
T_t(x)=W_t(x)+\int_0^t l(s,T(x))ds\,\,,\,\,t\in \reals_+\,,
$$
where $W_t(x)=x(t)$ and $l(s,x)=D_s\phi\circ T^{-1}(x)$.
\qed


\section{Construction and uniqueness  of the transport map in the
  general case} 
\label{general-case}
In this section we call optimal every probability measure{\footnote{In
  fact the results of this section are essentially true for the
  bounded, positive measures.}} $\ga$ on
$W\times W$ such that $J(\ga)<\infty$ and that $J(\ga)\leq J(\theta)$
for every other probability $\theta$ having the same marginals as
those of $\ga$. We recall that a  finite dimensional subspace $F$  of $W$ is 
called regular if the corresponding  projection is
continuous. Similarly a finite dimensional projection of $H$ is called
regular if it has a continuous extension to $W$.

We begin with the following lemma which answers all kind of questions
of measurability that we may encounter in the sequel:
\begin{lemma}
\label{meas-lemma}
Consider two uncountable Polish spaces $X$ and $T$. Let $t\to \ga_t$
be a Borel family of probabilities on $X$ and let $\calF$ be a
separable sub-$\sigma$-algebra of the Borel $\sigma$-algebra $\calB$
of $X$. Then there exists a Borel kernel 
$$
N_tf(x)=\int_Xf(y)N_t(x,dy)\,,
$$
such that, for any bounded Borel function $f$ on $X$, the following
properties hold true:
\begin{enumerate}
\item[i)]$(t,x)\to N_tf(x)$ is Borel measurable on $T\times X$.
\item[ii)]For any $t\in T$, $N_tf$ is an $\calF$-measurable version of
  the conditional expectation $E_{\ga_t}[f|\calF]$.
\end{enumerate}
\end{lemma}
\proof
Assume first that $\calF$ is finite, hence it is generated by a finite
partition $\{A_1,\ldots,A_k\}$. In this case it suffices to take 
$$
N_tf(x)=\sum_{i=1}^k\frac{1}{\ga_t(A_i)}\left(\int_{A_i}fd\ga_t\right)\,1_{A_i}(x)\,\,\left({\mbox{
    with }}0=\frac{0}{0}\right)\,.
$$
For the general case, take an increasing sequence $(\calF_n,n\geq 1)$
of finite sub-$\sigma$-algebras whose union generates $\calF$. Without
loss of generality we can assume that $(X,\calB)$ is the Cantor set
(Kuratowski Theorem, cf., \cite{D-M}). Then for every clopen set
(i.e., a set which is closed and open at the same time) $G$ and any
$t\in T$, the sequence $(N_t^n1_G,n\geq 1)$ converges $\ga_t$-almost
everywhere. Define
$$
H_G(t,x)=\limsup_{m,n\to\infty}|N^n_t1_G(x)-N^m_t1_G(x)|\,.
$$
$H_G$ is a Borel function on $T\times X$ which vanishes $\ga_t$-almost
all $x\in X$, besides, for any $t\in T$,  $x\to H_G(t,x)$ is
$\calF$-measurable. As there exist only countably many clopen sets in
$X$, the function 
$$
H(t,x)=\sup_GH_G(t,x)
$$
inherits all the  measurability properties. Let $\theta$ be any
probability on $X$, for any clopen $G$, define
$$
\begin{array}{cllrl}
N_t1_G(x)&=&\lim_nN_t^n1_G(x)&{\mbox{ if }}& H(t,x)=0\,,\\
&=&\theta(G)&{\mbox{ if }}& H(t,x)>0\,.
\end{array}
$$
Hence, for any $t\in T$, we get an additive measure on the Boolean
algebra of clopen sets of $X$. Since such a measure is
$\sigma$-additive and extends uniquely as a $\sigma$-additive measure
on $\calB$, the proof is completed.
\qed

\begin{remarkk}
{\rm{
\begin{enumerate}
\item
This result holds in fact for the Lusin spaces since they are
    Borel isomorphic to the Cantor set. Besides it extends easily to
    countable spaces. 
\item The particular case where $T=\calM_1(X)$, i.e., the space of
  probability measures on $X$ under the weak topology and $t\to \ga_t$
  being the identity map, is particularly
  important for the sequel. In this case we obtain a kernel $N$ such
  that  $(x,\ga)\to N_\ga f(x)$ is measurable and $N_\ga f$ is an
  $\calF$-measurable version of $E_\ga[f|\calF]$. 
\end{enumerate}
}}
\end{remarkk}

\begin{lemma}
\label{section-lemma}
Let $\rho$ and $\nu$ be two probability measures on $W$ such that
$$
d_H(\rho,\nu)<\infty
$$
and let $\ga\in\Sigma(\rho,\nu)$ be an optimal measure, i.e.,
$J(\ga)=d_H^2(\rho,\nu)$, where $J$ is given by 
(\ref{J-defn}). Assume that $F$ is a regular 
finite dimensional subspace of $W$ with the 
corresponding projection  $\pi_F$  from $W$ to $F$ and let
$\pi_F^\bot=I_W-\pi_F$ . Define  $p_F$ as  the 
projection from $W\times W$ onto $F$ with $p_F(x,y)=\pi_Fx$ and let
$p_F^\bot(x,y)=\pi_F^\bot x$.  
Consider the Borel disintegration 
\beaa
\ga(\cdot)&=&\int_{F^\bot\times W}\ga(\,\cdot|x^\bot)\ga^\bot(dz^\bot)\\
&=&\int_{F^\bot}\ga(\,\cdot|x^\bot)\rho^\bot(dx^\bot)
\eeaa
along the projection  of $W\times W$ on $F^\bot$, where
$\rho^\bot$ is the measure
$\pi_F^\bot\rho$,  $\ga(\cdot\,|x^\bot)$ denotes the regular 
conditional probability $\ga(\cdot\,|p_F^\bot =x^\bot)$ and $\ga^\bot$
is the measure $p_F^\bot\ga$. Then,
$\rho^\bot$ and $\ga^\bot$-almost surely  
$\ga(\,\cdot|x^\bot)$ is optimal on $(x^\bot+F)\times W$.
\end{lemma}
\proof
Let $p_1,\,p_2$ be the projections of $W\times W$ defined as
$p_1(x,y)=\pi_F(x)$ and $p_2(x,y)=\pi_F(y)$. 
Note first the  following obvious identity:
$$
p_1\ga(\cdot\,|x^\bot)=\rho(\cdot\,|x^\bot)\,,
$$
$\rho^\bot$ and $\ga^\bot$-almost surely.
Define the sets $B\subset F^\bot\times \calM_1(F\times F)$ and $C$ as 
\beaa
B&=&\{(x^\bot,\theta):\,\theta\in\Sigma(p_1\ga(\cdot\,|x^\bot),
p_2\ga(\cdot\,|x^\bot))\}\\
C&=&\{(x^\bot,\theta)\in B:\,J(\theta)<J(\ga(\cdot\,|x^\bot)\}\,,
\eeaa
where $\calM_1(F\times F)$ denotes the set of probability measures on
$F\times F$. Let $K$ be the projection of $C$ on $F^\bot$. Since $B$
and $C$ are Borel measurable, $K$is a
Souslin set, hence it is $\rho^\bot$-measurable.  The selection
theorem (cf. \cite{D-M}) implies the existence of a measurable map
$$
x^\bot\to \theta_{x^\bot}
$$ 
from $K$ to $\calM_1(F\times F)$ such that, $\rho^\bot$-almost
surely,  $(x^\bot,\theta_{x^\bot})\in C$. Define 
$$
\theta(\cdot)=\int_K\theta_{x^\bot}(\cdot)d\rho^\bot(x^\bot)+
\int_{K^c}\ga(\cdot\,|x^\bot)d\rho^\bot(x^\bot)\,.
$$
Then $\theta\in \Sigma(\rho,\nu)$ and we have 
\beaa
J(\theta)&=&\int_K J(\theta_{x^\bot})d\rho^\bot(x^\bot)
+\int_{K^c} J(\ga(\cdot\,|x^\bot))d\rho^\bot(x^\bot)\\
&<&\int_{K} J(\ga(\cdot\,|x^\bot))d\rho^\bot(x^\bot)
+\int_{K^c} J(\ga(\cdot\,|x^\bot))d\rho^\bot(x^\bot)\\
&=&J(\ga)\,,
\eeaa
hence we obtain $J(\theta)<J(\ga)$ which is a contradiction to the
optimality of $\ga$.
\qed

\begin{lemma}
\label{f-dim-lemma}
Assume that  the hypothesis of Lemma \ref{section-lemma} holds  and 
let $F$ be  any regular  finite dimensional  subspace  of $W$. 
Denote by $\pi_F$ the projection operator associated to it and let
$\pi_F^\bot=I_W-\pi_F$. If $\pi_F^\bot\rho$-almost surely,  the
regular conditional probability 
$\rho(\cdot\,|\pi_F^\bot=x^\bot )$ vanishes on the subsets of
$x^\bot+F$ whose  Hausdorff dimension are at most equal to 
${\mbox{\rm  dim}}(F)-1$,  then
there exists a map $T_F:F\times F^\bot\to F$ such that 
$$
\ga\left(\left\{(x,y)\in W\times
  W:\,\pi_F y=T_F(\pi_Fx,\pi_F^\bot x)\right\}\right)=1\,.
$$
\end{lemma}
\proof
Let $C_{x^\bot}$ be the support of the regular conditional probability
$\ga(\cdot\,|x^\bot)$ in $(x^\bot+F)\times W$. 
We know from Lemma \ref{section-lemma}
that the measure  $\ga(\cdot\,|x^\bot)$ is optimal in
$\Sigma(\pi_1\ga(\cdot\,|x^\bot),\pi_2\ga(\cdot\,|x^\bot))$, 
with $J(\ga(\cdot\,|x^\bot))<\infty$ for $\rho^\bot$-almost everywhere
$x^\bot$. From  Theorem 2.3 of \cite{G-Mc} and from \cite{Ab-H}, the
set $C_{x^\bot}$ 
is cyclically monotone, moreover, $C_{x^\bot}$ is a subset of
$(x^\bot+F)\times H$, 
hence the cyclic monotonicity of it implies that the  set
$K_{x^\bot}\subset F\times F$,  defined as 
$$
K_{x^\bot}=\{(u,\pi_Fv)\in F\times F:\,(x^\bot+u,v)\in C_{x^\bot}\}
$$
is cyclically monotone in $F\times F$. Therefore $K_{x^\bot}$ is 
included in the subdifferential of a convex function defined on $F$. 
Since, by hypothesis,  the first marginal of $\ga(\cdot\,|x^\bot)$, i.e.,
$\rho(\cdot\,|x^\bot)$ vanishes on the subsets of $x^\bot+F$ 
of co-dimension one,  the subdifferential
under question, denoted as $U_F(u,x^\bot)$  is
$\rho(\cdot\,|x^\bot)$-almost surely univalent (cf. \cite{An-K,Mc1}). This
implies that  
$$
\ga(\cdot\,|x^\bot)\left(\left\{(u,v)\in
    C_{x^\bot}:\,\pi_Fv=U_F(u,x^\bot)\right\}\right)=1\,, 
$$
$\rho^\bot$-almost surely.
Let 
$$
K_{x^\bot,u}=\left\{v\in W:\,(u,v)\in K_{x^\bot}\right\}\,.
$$
Then $K_{x^\bot,u}$ consists of a single point for
almost all $u$ with respect to $\rho(\cdot\,|x^\bot)$.
Let 
$$
N=\left\{(u,x^\bot)\in F\times
  F^\bot:\,{\mbox{Card}}(K_{x^\bot,u})>1\right\}\,,
$$
note that $N$ is a Souslin set, hence it is universally
measurable. Let $\sigma$ be the measure which is defined as the image
of $\rho$ under the projection $x\to (\pi_Fx,\pi_F^\bot x)$. We then
have 
\beaa
\sigma(N)&=&\int_{F^\bot}\rho^\bot(dx^\bot)\int_F
\won_N(u,x^\bot)\rho(du|x^\bot)\\
&=&0\,.
\eeaa
Hence $(u,x^\bot)\mapsto K_{x^\bot,u}=\{y\}$ is $\rho$ and  $\ga$-almost surely
well-defined and it suffices to denote this map by $T_F$ to achive the
proof.
\qed

\begin{theorem}
\label{monge-general}
Suppose that $\rho$ and $\nu$ are two probability measures on
$W$ such  that
$$
d_H(\rho,\nu)<\infty\,.
$$
Let $(\pi_n,n\geq 1)$ be a total increasing  sequence of regular
projections (of $H$, converging to the identity map of $H$). 
Suppose  that, for any $n\geq 1$, the regular
conditional probabilities $\rho(\cdot\,|\pi_n^\bot=x^\bot)$ vanish
$\pi_n^\bot\rho$-almost surely on 
the subsets of  $(\pi_n^\bot)^{-1}(W)$ with Hausdorff dimension
$n-1$. Then there exists a    
unique solution of the   Monge-Kantorovitch problem, denoted by $\ga\in
\Sigma(\rho,\nu)$ and  $\ga$ is supported by the graph of a Borel
map $T$ which is the solution of the
Monge problem.  $T:W\to W$ is  of the form $T=I_W+\xi$ , where $\xi\in
H$ almost surely. Besides  we have 
\beaa
d_H^2(\rho,\nu)&=&\int_{W\times W}|T(x)-x|_H^2d\ga(x,y)\\
&=&\int_{W}|T(x)-x|_H^2d\rho(x)\,, 
\eeaa
and  
for $\pi_n^\bot\rho$-almost almost all $x_n^\bot$, the map $u\to
u+\xi(u+x_n^\bot)$ is cyclically monotone on
$(\pi_n^\bot)^{-1}\{x_n^\bot\}$, in the sense that 
$$
\sum_{i=1}^N\left(u_i+\xi(x_n^\bot+u_i),u_{i+1}-u_i\right)_H\leq 0
$$
$\pi_n^\bot\rho$-almost surely, for any cyclic sequence
$\{u_1,\ldots,u_N,u_{N+1}=u_1\}$ from $\pi_n(W)$. Finally, if, for any $n\geq
1$, $\pi_n^\bot\nu$-almost surely,  $\nu(\cdot\,|\pi_n^\bot=y^\bot)$
 also vanishes on the $n-1$-Hausdorff dimensional  subsets  
 of $(\pi_n^\bot)^{-1}(W)$, then $T$ is invertible, i.e, there exists
 $S:W\to W$ of the form $S=I_W+\eta$ such that  $\eta\in H$ satisfies
 a similar  cyclic monotononicity property as $\xi$ and that 
\beaa
1&=&\ga\left\{(x,y)\in W\times W: T\circ S(y)=y\right\}\\
&=&\ga\left\{(x,y)\in W\times W: S\circ T(x)=x\right\}\,.
\eeaa
In particular we have 
\beaa
d_H^2(\rho,\nu)&=&\int_{W\times W}|S(y)-y|_H^2d\ga(x,y)\\
&=&\int_{W}|S(y)-y|_H^2d\nu(y)\,. 
\eeaa
\end{theorem}
\begin{remarkk}
In particular, for all  the measures $\rho$ which are  absolutely
continuous with respect to the  Wiener measure $\mu$,  the second
hypothesis is satisfied, i.e., the measure
$\rho(\cdot\,|\pi_n^\bot=x_n^\bot)$ vanishes on the sets of Hausdorff
dimension $n-1$.
\end{remarkk} 
\proof
Let $(F_n,n\geq 1)$ be the  increasing sequence of regular subspaces 
associated to $(\pi_n,n\geq 1)$,  
 whose union is dense in $W$.  From Lemma \ref{f-dim-lemma}, for any
 $F_n$, there exists a map $T_n$, such that
 $\pi_n y=T_n(\pi_n x,\pi_n^\bot x)$ for $\ga$-almost all $(x,y)$, where
 $\pi_n^\bot=I_W-\pi_n$. Write 
 $T_n$ as $I_n+\xi_n$, where $I_n$ denotes the identity map on
 $F_n$. Then we have the following representation:
$$
\pi_ny=\pi_nx+\xi_n(\pi_nx,\pi_n^\bot x)\,,
$$
$\ga$-almost surely.
Since 
\beaa
\pi_n y-\pi_n x&=&\pi_n( y-x)\\
&=&\xi_n(\pi_n x,\pi_n^\bot x)\,
\eeaa 
and since $y-x\in H$ $\ga$-almost surely, $(\pi_ny-\pi_nx,n\geq 1)$
converges $\ga$-almost surely.
Consequently $(\xi_n,n\geq 1)$ converges $\ga$, hence $\rho$ almost
surely to a measurable  $\xi$. Consequently we obtain 
$$
\ga\left(\left\{(x,y)\in W\times W:\, y=x+\xi(x)\right\}\right)=1\,.
$$
Since $J(\ga)<\infty$, $\xi$ takes its values almost surely in the
Cameron-Martin space $H$. The cyclic monotonicity of $\xi$  is
obvious. To prove the uniqueness, assume that we have two 
optimal solutions $\ga_1$ and $\ga_2$ with the same marginals and
$J(\ga_1)=J(\ga_2)$. Since $\beta\to J(\beta)$ is linear,  the measure
defined as  $\ga=\frac{1}{2}(\ga_1+\ga_2)$ is
also optimal and it has also the same marginals $\rho$ and
$\nu$. Consequently, it is also  supported by the graph of a map 
$T$. Note that  $\ga_1$ and $\ga_2$ are absolutely continuous with
respect to $\ga$,  let $L_1(x,y)$ be the Radon-Nikodym density of
$\ga_1$ with respect to $\ga$. For any $f\in C_b(W)$, we then  have 
\beaa
\int_Wfd\rho&=&\int_{W\times W}f(x)d\ga_1(x,y)\\
&=&\int_{W\times W}f(x)L_1(x,y)d\ga(x,y)\\
&=&\int_Wf(x)L_1(x,T(x))d\rho(x)\,. 
\eeaa
Therefore we should have $\rho$-almost surely,  $L_1(x,T(x))=1$, hence
also $L_1=1$ almost everywhere $\ga$ and this implies that
$\ga=\ga_1=\ga_2$. The second part about the invertibility of $T$ is
totally symmetric, hence its proof follows along the same lines as the
proof for $T$.
\qed

\begin{corollary}
\label{c-mono}
Assume that $\rho$ is equivalent to the Wiener measure $\mu$, then for
any $h_1,\ldots,h_N\in H$ and for any permutation $\tau$ of
$\{1,\ldots,N\}$, we have, with the notations of Theorem
\ref{monge-general}, 
$$
\sum_{i=1}^N\left(h_i+\xi(x+h_i),h_{\tau(i)}-h_i\right)_H\leq 0
$$
$\rho$-almost surely.
\end{corollary}
\proof
Again with the notations of the theorem, $\rho_k^\bot$-almost surely,
the graph of the map $x_k\to x_k+\xi_k(x_k,x_k^\bot)$ is cyclically
monotone on $F_k$. Hence, for the case $h_i\in F_n$ for all
$i=1,\ldots,N$ and $n\leq k$, we have 
$$
\sum_{i=1}^N\left(h_i+x_k+\xi_k(x_k+h_i,x_k^\bot),h_{\tau(i)}-h_i\right)_H\leq 0\,.
$$
Since $\sum_i(x_k,h_{\tau(i)}-h_i)_H=0$, we also have 
$$
\sum_{i=1}^N\left(h_i+\xi_k(x_k+h_i,x_k^\bot),h_{\tau(i)}-h_i\right)_H\leq 0\,.
$$
We know that $\xi_k(x_k+h_i,x_k^\bot)$ converges to $\xi(x+h_i)$
$\rho$-almost surely. Moreover $h\to \xi(x+h)$ is continuous from $H$
to $L^0(\rho)$ and the proof follows.
\qed


\section{The Monge-Amp\`ere equation}
\label{equation}
\label{monge-ampere}


Assume   that $W=\R^n$ and take a   density $L\in
\LL\log\LL$. Let $\phi\in \DD_{2,1}$ be the $1$-convex function such
that $T=I+\nabla \phi$ maps $\mu$ to $L\cdot \mu$. Let $S=I+\nabla\psi$
be its inverse with $\psi\in \DD_{2,1}$. 
Let now $\nabla_a^2\phi$ be the
second Alexandrov derivative of 
$\phi$, i.e.,  the Radon-Nikodym derivative of the absolutely
continuous part of the 
vector measure $\nabla^2\phi$ with respect to the Gaussian measure
$\mu$ on $\reals^n$. Since $\phi$ is $1$-convex, it follows that
$\nabla^2\phi\geq -I_{\reals^n}$ in the sense of the distributions,
consequently $\nabla_a^2\phi\geq -I_{\reals^n}$ $\mu$-almost
surely. Define also the Alexandrov version $\calL_a\phi$ of
$\calL\phi$ as the Radon-Nikodym derivative of  the absolutely
continuous part of the distribution 
$\calL\phi$. Since we are in finite dimensional situation, we have the
explicit expression for $\calL_a\phi$ as 
$$
\calL_a\phi(x)=(\nabla\phi(x),x)_{\reals^n}-{\rm
  trace}\left(\nabla_a^2\phi\right)\,.
$$
Let  $\La$ be the Gaussian  Jacobian 
$$
\La=\dett\left(I_{\R^n}+\nabla^2_a\phi\right)\exp\left\{-\calL_a\phi-
\frac{1}{2}|\nabla\phi|_{\R^n}^2\right\}\,.
$$
\begin{remarkk}{\rm
In this expression as well as in the sequel,  the notation
$\dett(I_H+A)$ denotes the modified 
Carleman-Fredholm determinant of the operator $I_H+A$ on a Hilbert
space $H$. If $A$ is an operator of finite rank, then it  is defined  as 
$$
\dett\left(I_H+A\right)=\prod_{i=1}^n(1+l_i)e^{-l_i}\,,
$$
where $(l_i,\,i\leq n)$ denotes the eigenvalues of $A$ counted with
respect to their multiplicity. In fact this determinant has an
analytic extension to the space of Hilbert-Schmidt operators on a
separable Hilbert space, cf.  \cite{D-S} and Appendix A.2 of
\cite{UZ7}. As explained in \cite{UZ7}, the modified determinant
exists for the Hilbert-Schmidt operators while the ordinary
determinant does not, since the latter  requires the existence of the trace of
$A$. Hence the modified Carleman-Fredholm determinant  is particularly
useful when one studies the absolute 
continuity properties of the image of a Gaussian measure under
non-linear transformations in the setting of infinite dimensional
Banach spaces (cf., \cite{UZ7} for further information).  }
\end{remarkk}
It follows from the change of variables formula given in  Corollary
4.3  of \cite{Mc2}, that, for any $f\in C_b(\R^n)$,  
$$
E[f\circ T \,\La]=E\left[f\,1_{\partial\Phi(M)}\right]\,,
$$
where $M$ is the set of non-degeneracy of $I_{\R^n}+\nabla_a^2\phi$,
$$
\Phi(x)=\frac{1}{2}|x|^2+\phi(x)
$$ 
and $\partial \Phi$ denotes  the subdifferential of the convex
function  $\Phi$. Let us note that, in case $L>0$ almost surely, $T$
has  a global inverse $S$, i.e., $S\circ T=T\circ S=I_{\R^n}$
$\mu$-almost surely and $\mu(\partial\Phi(M))=\mu(S^{-1}(M))$. 
Assume now that $\La>0$ almost surely, i.e., that $\mu(M)=1$. Then,
for any $f\in C_b(\R^n)$, we have 
\beaa
E[f\circ T]&=&E\left[f\circ T\,\frac{\La}{\La\circ T^{-1}\circ T}\right]\\
&=&E\left[f\,\frac{1}{\La\circ T^{-1}}1_{\partial\Phi(M)}\right]\\
&=&E[f\,L]\,,
\eeaa
where $T^{-1}$ denotes the left inverse of $T$ whose existence is
guaranteed by Theorem \ref{gaussian-case}. Since $T(x)\in
\partial\Phi(M)$ almost surely, it follows from the above calculations
$$
\frac{1}{\La}=L\circ T\,,
$$
almost surely. Take now  any $t\in [0,1)$, the map
$x\to \frac{1}{2}|x|_H^2+t\phi(x)=\Phi_t(x)$ is strictly convex and a simple
calculation implies that the mapping  $T_t=I+t\nabla\phi$ is
$(1-t)$-monotone (cf. \cite{UZ7}, Chapter 6), consequently it has a
left inverse denoted by $S_t$. Let us  denote by $\Psi_t$ the  Legendre
transformation of $\Phi_t$: 
$$
\Psi_t(y)=\sup_{x\in \R^n}\left\{(x,y)-\Phi_t(x)\right\}\,.
$$
A simple calculation shows that 
\beaa
\Psi_t(y)&=&\sup_x\left[(1-t)\left\{(x,y)-\frac{|x|^2}{2}\right\}+
t\left\{(x,y)-\frac{|x|^2}{2}-\phi(x)\right\}\right]\\
&\leq&(1-t)\frac{|y|^2}{2}+t\Psi_1(y)\,.
\eeaa
Since $\Psi_1$ is the Legendre transformation of
$\Phi_1(x)=|x|^2/2+\phi(x)$ and since $L\in \LL\log \LL$, it is
finite on a convex set of full measure, hence it is finite everywhere.
Consequently $\Psi_t(y)<\infty$ for any $y\in \R^n$. Since a finite,
convex function is almost everywhere differentiable, $\nabla \Psi_t$
exists almost everywhere on  and it is equal almost everywhere on
$T_t(M_t)$ to the left inverse $T_t^{-1}$, where $M_t$ is the set of
non-degeneracy of $I_{\R^n}+t\nabla^2_a\phi$. Note that $\mu(M_t)=1$. 
 The strict convexity implies that $T_t^{-1}$ is Lipschitz with a
Lipschitz constant $\frac{1}{1-t}$. 
Let now $\La_t$ be the Gaussian  Jacobian 
$$
\La_t=\dett\left(I_{\R^n}+t\nabla^2_a\phi\right)\exp\left\{-t\calL_a\phi-
\frac{t^2}{2}|\nabla\phi|_{\R^n}^2\right\}\,.
$$
Since the domain of $\phi$ is the whole space $\R^n$, $\La_t>0$ almost
surely, hence, as we have explained above,   it follows from the
change of variables formula of \cite{Mc2} that $T_t\mu$ is absolutely
continuous with respect to $\mu$ and that 
$$
\frac{1}{\La_t}=L_t\circ T_t\,,
$$
$\mu$-almost surely. 

\noindent
Let us come back to the infinite dimensional case: we first  give an
inequality which may be useful.
\begin{theorem}
\label{monge1-thm}
Assume that $(W,\mu,H)$ is an abstract  Wiener space, assume
that $K,L\in \LL_+^1(\mu)$ with $K>0$ almost surely  and denote by
$T:W\to W$ the transfer map 
$T=I_W+\nabla \phi$, which maps the measure $Kd\mu$ to the measure
$Ld\mu$. Then the following inequality holds:
\begin{equation}
\label{monge1-ineq}
\frac{1}{2}E[|\nabla\phi|_H^2]\leq E[-\log K+\log L\circ T]\,.
\end{equation} 
\end{theorem}
\proof
Let us define $k$ as $k=K\circ T^{-1}$, then for any $f\in C_b(W)$, we
have 
\beaa
\int_Wf(y)L(y)d\mu(y)&=&\int_Wf\circ T(x)K(x)d\mu(x)\\
&=&\int_Wf\circ T(x)k\circ T(x)d\mu(x)\,,
\eeaa
hence 
$$
T\mu=\frac{L}{k}\,.\mu\,.
$$
It then follows from the inequality \ref{tal-ineq} that 
\beaa
\frac{1}{2}E\left[|\nabla\phi|_H^2\right]&\leq&E\left[\frac{L}{k}\log\frac{L}{k}\right]\\
&=&E\left[\log\frac{L\circ T}{k\circ T}\right]\\
&=&E[-\log K+\log L\circ T]\,.
\eeaa
\qed

\noindent
 Suppose that $\phi\in \DD_{2,1}$ is a $1$-convex Wiener
 functional. Let $V_n$ be 
the sigma algebra generated by $\{\delta e_1,\ldots,\delta e_n\}$,
where $(e_n,\,n\geq 1)$ is an orthonormal basis of the Cameron-Martin
space $H$. Then $\phi_n=E[\phi|V_n]$ is again $1$-convex
(cf. Chapter \ref{ch.convex}), hence $\calL\phi_n$ is a measure as
it can be easily verified.  However the sequence
$(\calL\phi_n,\,n\geq 1)$ converges to $\calL\phi$ only  in
$\DD'$. Consequently, there is no reason for the limit 
$\calL\phi$ to be a  measure. In case this happens, we shall denote the
Radon-Nikodym density with respect to $\mu$,  of the 
absolutely continuous part  of this measure  by $\calL_a\phi$. 
\begin{lemma}
\label{app-lemma}
Let $\phi\in \DD_{2,1}$ be $1$-convex and let $V_n$ be defined as
above  and define $F_n=E[\phi|V_n]$. Then
the sequence $(\calL_a F_n,n\geq 1)$ 
is a submartingale, where $\calL_aF_n$ denotes the
$\mu$-absolutely continuous part of the measure $\calL F_n$.
\end{lemma}
\proof
Note that, due to the $1$-convexity, we have $\calL_aF_n\geq
\calL F_n$ for any $n\in \NN$. 
Let $X_n=\calL_aF_n$ and $f\in \DD$ be a positive, $V_n$-measurable
test function. 
Since $\calL E[\phi|V_n]=E[\calL\phi|V_n]$, we have 
\beaa
E[X_{n+1}\,f]&\geq&\langle \calL F_{n+1},f\rangle\\
&=&\langle \calL F_{n},f\rangle\,,
\eeaa
where $\langle\cdot,\cdot\rangle$ denotes the duality bracket for the
dual pair $(\DD',\DD)$. 
Consequently  
$$
E[f\,E[X_{n+1}|V_n]]\geq \langle \calL F_{n},f\rangle\,,
$$ 
for any positive, $V_n$-measurable test function $f$, it follows that  the
absolutely continuous part of $\calL F_n$  is also  dominated by the
same conditional expectation and  this proves the submartingale
property. 
\qed

\begin{lemma}
\label{Lf-lemma}
Assume that $L\in \LL\log \LL$ is a positive random variable whose
expectation is one. Assume further that it  is lower bounded by a  constant
$a>0$. Let $T=I_W+\nabla \phi$ be the transport map such that
$T\mu=L\,.\mu$ and let $T^{-1}=I_W+\nabla \psi$. Then $\calL \psi$ is
a Radon measure on $(W,\calB(W))$. If $L$ is upper bounded by $b>0$,
then $\calL\phi$ is also a Radon measure on $(W,\calB(W))$.
\end{lemma}
\proof
Let $L_n=E[L|V_n]$, then $L_n\geq a$ almost surely. Let
$T_n=I_W+\nabla\phi_n$ be the transport map which satisfies
$T_n\mu=L_n\,.\mu$ and let $T_n^{-1}=I_W+\nabla\psi_n$ be its
inverse. We have 
$$
L_n=\dett\left(I_H+\nabla_a^2\psi_n\right)\exp\left[-\calL_a\psi_n-\frac{1}{2}|\nabla\psi_n|_H^2\right]\,.
$$
By the hypothesis $-\log L_n\leq -\log a$.  Since $\psi_n$ is
$1$-convex, it follows from the finite dimensional results that
$\dett\left(I_H+\nabla_a^2\psi_n\right)\in [0,1]$ almost surely. 
Therefore  we have
$$
\calL_a\psi_n\leq -\log a\,,
$$
besides $\calL\psi_n\leq \calL_a\psi_n$ as distributions, consequently
$$
\calL\psi_n\leq -\log a
$$
as distributions, for any $n\geq 1$. Since $\lim_n\calL\psi_n=\calL
\psi$ in $\DD'$, we obtain $\calL\psi\leq -\log a$, hence $-\log
a-\calL\psi\geq 0$ as a distribution, hence $\calL\psi$  is a Radon measure on
$W$.  This proves the
first claim. Note that whenever $L$ is upperbounded, $\La=1/L\circ T$
is lowerbounded, hence  the proof of the second claim  is similar to
that of the first one.
\qed

\begin{theorem}
\label{girsanov-thm}
Assume that $L$ is a strictly  positive bounded  random variable with
$E[L]=1$. Let $\phi\in \DD_{2,1}$ be
the $1$-convex Wiener functional such that 
$$
T=I_W+\nabla \phi
$$ 
is the transport map realizing the measure $L\,.\mu$ and let
$S=I_W+\nabla\psi$ be 
its inverse. Define $F_n=E[\phi|V_n]$, then the submartingale
$(\calL_aF_n,n\geq 1)$ converges almost surely to $\calL_a\phi$. Let
$\la(\phi)$ be  the random variable defined as 
\beaa
\la(\phi)&=&\lim\inf_{n\to\infty}\La_n\\
&=&\left(\lim\inf_n\dett\left(I_H+\nabla^2_aF_n\right)\right)\exp\left\{-\calL_a\phi-\half|\nabla\phi|_H^2\right\}
\eeaa
where
$$
\La_n=\dett\left(I_H+\nabla_a^2F_n\right)\exp\left\{-\calL_aF_n-\half|\nabla
  F_n|_H^2\right\}\,.
$$
Then it holds true that 
\begin{equation}
\label{sub-gir}
E[f\circ T\,\la(\phi)]\leq E[f]
\end{equation}
for any $f\in C^+_b(W)$, in particular $\la(\phi)\leq \frac{1}{L\circ
  T}$ almost surely. If $E[\la(\phi)]=1$, then the inequality in
(\ref{sub-gir}) becomes an equality and  we also  have  
$$
\la(\phi)=\frac{1}{L\circ T}\,.
$$
\end{theorem}
\proof Let us remark that, due to the $1$-convexity,
$0\leq\dett\left(I_H+\nabla^2_aF_n\right)\leq 1$, hence the $\lim\inf$
exists. Now,  Lemma \ref{Lf-lemma} implies that $\calL\phi$ is a Radon
measure. Let $F_n=E[\phi|V_n]$, then we know from Lemma
\ref{app-lemma} that $(\calL_a F_n,n\geq 1)$ is a submartingale. Let
$\calL^+\phi$ denote the positive part of the measure $\calL\phi$. Since
$\calL^+\phi\geq \calL\phi$, we have also $E[\calL^+\phi|V_n]\geq
E[\calL \phi|V_n]=\calL F_n$. This implies that
$E[\calL^+\phi|V_n]\geq\calL_a^+F_n$. Hence we find that 
$$
\sup_n E[\calL^+_a F_n]<\infty
$$
and this condition implies that the submartingale $(\calL_aF_n,n\geq 1)$
converges almost surely. We shall now identify the limit of this
submartingale. Let $\calL_s G$ be  the singular part of the
measure $\calL G$ for a  Wiener function $G$ such that $\calL G$ is a
measure. We have
\beaa
E[\calL\phi|V_n]&=&E[\calL_a\phi|V_n]+E[\calL_s\phi|V_n]\\
&=&\calL_aF_n+\calL_sF_n\,,
\eeaa
hence 
$$
\calL_aF_n=E[\calL_a\phi|V_n]+E[\calL_s\phi|V_n]_a
$$
almost surely, where $E[\calL_s\phi|V_n]_a$ denotes the absolutely
continuous part of the measure $E[\calL_s\phi|V_n]$. Note that, from
the Theorem of Jessen 
(cf., for example Theorem 1.2.1 of \cite{UZ7}),
$\lim_nE[\calL_s^+\phi|V_n]_a=0$ and
$\lim_nE[\calL_s^-\phi|V_n]_a=0$ almost surely, hence we have 
$$
\lim_n\calL_aF_n=\calL_a\phi\,,
$$
$\mu$-almost surely. To complete the proof, an application of the
Fatou lemma implies that  
\beaa
E[f\circ T\,\la(\phi)]&\leq&E[f]\\
&=&E\left[f\circ T\,\frac{1}{L\circ T}\right]\,,
\eeaa
for any $f\in C_b^+(W)$. Since $T$ is invertible, it follows that  
$$
\la(\phi)\leq \frac{1}{L\circ T}
$$
almost surely. Therefore, in  case $E[\la(\phi)]=1$,  we have
$$
\la(\phi)=\frac{1}{L\circ T}\,,
$$
and this completes the proof.
\qed

\begin{corollary}
\label{monge-amp.eqn}
Assume that $K,L$ are two positive random variables with values in a
bounded interval $[a,b]\subset (0,\infty)$ such that
$E[K]=E[L]=1$. Let
$T=I_W+\nabla\phi$, $\phi\in \DD_{2,1}$,  be the
transport map pushing $Kd\mu$ to $Ld\mu$, i.e, $T(Kd\mu)=Ld\mu$. We
then have
$$
L\circ T\,\la(\phi)\leq K\,,
$$
$\mu$-almost surely. In particular, if $E[\la(\phi)]=1$, then $T$ is
the solution of the Monge-Amp\`ere equation.
\end{corollary}
\proof Since $a>0$, 
$$
\frac{dT\mu}{d\mu}=\frac{L}{K\circ T}\leq \frac{b}{a}\,.
$$
Hence,  Theorem \ref{sub-gir} implies that 
\beaa
E[f\circ T\,L\circ T\,\la(\phi)]&\leq& E[f\,L]\\
&=&E[f\circ T\,K]\,,
\eeaa
consequently
$$
L\circ T\,\la(\phi)\leq K\,,
$$
the rest of the claim is now obvious.
\qed

For later use we give also the folowing result:
\begin{theorem}
\label{interpol-thm}
Assume that $L$ is a positive random variable of class $\LL\log \LL$
such that $E[L]=1$. Let $\phi\in \DD_{2,1}$ be the $1$-convex function
corresponding to the transport map $T=I_W+\nabla\phi$. Define
$T_t=I_W+t\nabla\phi$, where  $t\in [0,1]$. Then,
for any $t\in [0,1]$, $T_t\mu$ is absolutely continuous with respect
to the Wiener measure $\mu$.
\end{theorem}
\proof
Let $\phi_n$ be defined as the transport map corresponding to
$L_n=E[P_{1/n}L_n|V_n]$ and define $T_n$ as $I_W+\nabla\phi_n$. For $t\in
[0,1)$, let $T_{n,t}=I_W+t\nabla\phi_n$.  It follows from the finite
dimensional results which are summarized  in the
beginning of this section, that $T_{n,t}\mu$ is absolutely continuous 
with respect to $\mu$. Let $L_{n,t}$ be the corresponding
Radon-Nikodym density and define $\La_{n,t}$ as 
$$
\La_{n,t}=\dett\left(I_H+t\nabla^2_a\phi_{n}\right)
\exp\left\{-t\calL_a\phi_n-\frac{t^2}{2}|\nabla\phi_n|_H^2\right\}\,.
$$
Besides, for any $t\in [0,1)$, 
\begin{equation}
\label{mon-con}
\left((I_H+t\nabla_a^2\phi_n)h,h\right)_H>0\,,
\end{equation}
$\mu$-almost surely for any $0\neq h\in H$. Since $\phi_n$ is of finite 
rank, \ref{mon-con} implies that $\La_{n,t}>0$ $\mu$-almost surely
and we have shown  at the beginning of this section 
$$
\La_{n,t}=\frac{1}{L_{n,t}\circ T_{n,t}}
$$
$\mu$-almost surely.  
An easy calculation shows that $t\to \log\dett(I+t\nabla_a^2\phi_n)$ is
a non-increasing function. Since $\calL_a\phi_n\geq \calL\phi_n$, we
have $E[\calL_a\phi_n]\geq 0$. Consequently
\beaa
E\left[L_{t,n}\log L_{t,n}\right]&=&E\left[\log L_{n,t}\circ T_{n,t}\right]\\
&=&-E\left[\log\La_{t,n}\right]\\
&=&E\left[-\log\dett\left(I_H+t\nabla^2\phi_n\right)+t\calL_a\phi_n
+\frac{t^2}{2}|\nabla\phi_n|_H^2\right]\\
&\leq&E\left[-\log\dett\left(I_H+\nabla^2\phi_n\right)+\calL_a\phi_n
+\frac{1}{2}|\nabla\phi_n|_H^2\right]\\
&=&E\left[L_n\log L_n\right]\\
&\leq&E[L\log L]\,,
\eeaa
by the  Jensen inequality. Therefore 
$$
\sup_nE[L_{n,t}\log L_{n,t}]<\infty
$$
and this implies that the sequence $(L_{n,t},n\geq 1)$ is uniformly
integrable for any $t\in [0,1]$. Consequently it has a subsequence
which converges weakly in $L^1(\mu)$ to some $L_t$.  Since, from Theorem
\ref{gaussian-case},  $\lim_n\phi_n=\phi$ in $\DD_{2,1}$, where $\phi$
is the transport map associated to $L$, for any $f\in C_b(W)$, we have
\beaa
E[f\circ T_t]&=&\lim_k E\left[f\circ T_{n_k,t}\right]\\
&=&\lim_k E\left[f\,L_{n_k,t}\right]\\
&=&E[f\,L_t]\,,
\eeaa
hence the theorem is proved.
\qed

\subsection{The solution of the Monge-Amp\`ere equation via
  Ito-renormalization} 
\label{sub-equation}
We can interpret the Monge-Amp\`ere equation as follows: given two
probability densities $K$ and $L$, find a map $T:W\to W$ such that 
$$
L\circ T\, J(T)=K
$$
almost surely, where $J(T)$ is a kind of Jacobian to be written in
terms of $T$. In Corollary \ref{monge-amp.eqn}, we have shown the
existence of some $\la(\phi)$ which gives an inequality instead of the
equality. Although in the finite dimensional case there are some
regularity results about the transport map (cf., \cite{Caf}), in the
infinite dimensional case  such techniques do not work.  All these
difficulties  can
be circumvented using the miraculous 
renormalization of the Ito calculus. In fact assume that $K$ and $L$
satisfy the hypothesis of the corollary. First let us indicate that we
can assume $W=C_0([0,1],\reals)$ (cf., \cite{UZ7}, Chapter II, to see
how one can pass from an abstract Wiener space to the standard
one) and in this case the Cameron-Martin space $H$ becomes
$H^1([0,1])$, which is the space of absolutely continuous functions on
$[0,1]$, with a square integrable Sobolev derivative. Let now 
$$
\La=\frac{K}{L\circ T}\,,
$$
where $T$ is as  constructed above. Then $\La.\mu$ is a Girsanov
measure for the map $T$. This means that the law of the stochastic
process $(t,x)\to T_t(x)$ under $\La.\mu$ is equal to the Wiener
measure, where $T_t(x)$ is defined as the evaluation of the trajectory
$T(x)$ at $t\in [0,1]$. In other words the process $(t,x)\to T_t(x)$
is a Brownian motion under the probability $\La.\mu$. Let
$(\calF^T_t,t\in [0,1])$ be its filtration, the invertibility of  $T$
implies  that 
$$
\bigvee_{t\in [0,1]}\calF^T_t=\calB(W)\,.
$$
 $\La$ is upper and lower  bounded $\mu$-almost surely, hence also
 $\La.\mu$-almost surely. The Ito representation theorem implies that
 it can be represented as 
$$
\La=E[\La^2]\exp\left\{-\int_0^1\dot{\alpha}_sdT_s-\half\int_0^1|\dot{\alpha}_s|^2ds\right\}\,,
$$
where $\alpha(\cdot)=\int_0^\cdot\dot{\alpha}_sds$ is an $H$-valued
random variable. In fact $\alpha$ can be calculated explicitly using
the Ito-Clark representation theorem, and it is
given as 
\begin{equation}
\label{alpha}
\dot{\alpha}_t=\frac{E_\La[D_t\La|\calF^T_t]}{E_\La[\La|\calF^T_t]}
\end{equation}
$dt\times \La d\mu$-almost surely, where $E_\La$ denotes the
expectation operator with respect to $\La.\mu$ and $D_t\La$ is the
Lebesgue density of the absolutely continuous map $t\to
\nabla\La(t,x)$. From the relation (\ref{alpha}), it follows that  $\alpha$ is 
a function of $T$, hence we have obtained the strong solution of the
Monge-Amp\`ere equation. Let us announce all this as 
\begin{theorem}
\label{monge-amp.2}
Assume that $K$ and $L$ are upper and lower bounded densities, let $T$
be the transport map constructed in Theorem \ref{monge-general}. Then
$T$ is also the strong solution of the Monge-Amp\`ere equation in the
Ito sense, namely 
$$
E[\La^2]\,L\circ T
\exp\left\{-\int_0^1\dot{\alpha}_sdT_s-\half\int_0^1|\dot{\alpha}_s|^2ds\right\}=K\,,
$$
$\mu$-almost surely, where $\alpha$ is given with (\ref{alpha}).
\end{theorem}

%

\chapter{Stochastic Analysis on Lie Groups}
\label{ch.lie}
\markboth{Stochastic Analysis on Lie Groups}{}
\section*{Introduction}
This chapter  is a partial survey of the construction of  Sobolev-type
analysis on the path space of  a Lie group. The word partial refers to
the fact that we
give some  new results about the quasi-invariance of  anticipative
transformations and the corresponding measure theoretical degree theorems   in
the last section. Almost all the theory has been initiated by S.  Albeverio and
R. H.-Krohn (\cite{A-HK}), L. Gross (\cite{G1,G2}) and M. P. Malliavin and
P. Malliavin (\cite{M-M1}). Although the study of the similar subjects has
already begun in the case of manifolds (cf. \cite{C-M}), we prefer to
understand first the case of the Lie groups because of their relative
simplicity and  this will give a better idea of what is going on in former
situation; since the frame of the Lie group-valued Brownian motion represents
the simplest non-linear and non-trivial  case in which we can construct a
Sobolev type functional analysis on the space of the trajectories.

After some preliminaries in the second section we give the definitions of the
basic tools in the third section, namely the left and right derivatives on the
path space. The fourth section is devoted to the left divergence, in the next
one we study the Ornstein-Uhlenbeck operator, Sobolev spaces and some
applications like the zero-one law. Sixth section is a compilation of the
formulas based essentially on the variation of the constants method of the
ordinary linear differential equations which are to be used in the following
sections. Section seven is devoted to the right derivative which is more
technical and interesting than the left one;  since it contains a
rotation of the path in the sense
of \cite{UZ4}. We also define there the skew-symmetric rotational derivative
and study some of its properties.
 Eighth section is devoted to the quasi-invariance at the left and at the right
 with respect to the multiplicaton of the path with deterministic paths of
 finite variation. Loop space case is also considered there.

Section nine deals with the absolute continuity of the path and loop measures
under the transformation which consists of multiplying from the left
the generic trajectory
with some random, absolutely continuous and anticipative path. We prove a
generalization of the Campbell-Baker-Hausdorff formula which is
fundemental. To prove
this we have been obliged to employ all the recent sophisticated techniques
derived in the flat case. Afterwards, the extension of the Ramer and the degree
theorems are immediate.

In this chapter  we have focuse  our attention to the
probabilistic and functional analytic problems. For the more
general case of Riemannian manifolds cf. \cite{Mall-1} and the
references therein.

\section{Analytic tools on group valued paths}
\markboth{Stochastic Analysis on Lie Groups}{Analytic Tools}
Let $G$ be a finite dimensional, connected, locally compact Lie group and
$\G$ be its Lie algebra of left invariant vector fields which is isomorphic
to the tangent space at identity of $G$, denoted by $T_e(G)$ which  is
supposed to be equipped with an inner product. $C=C_G$ denotes
$C_e([0,1],G)$ (i.e., $p(0)=e$ for $p\in C_G$). $C_\G$ denotes
$C_0([0,1],\G)$. Let
$$
H=H_\G=\left\{h\in C_\G:\int_0^1 |\dot{h}(t)|^2dt=|h|^2<\infty\right\}\,.
$$
Our basic Wiener space is $(C_\G,H,\mu)$. We denote by $p(w)$ the solution
of the following stochastic differential equation:
$$
p_t=e+\int_0^t p_s(w)dW_s(w)
$$
where the integral is in Stratonovitch sense and $W$ is the canonical
Brownian motion on $C_\G$. In general this equation is to be
understood as following: for any smooth function $f$ on $G$, we have
$$
f(p_t)=f(e)+\int_0^t H_if(p_s)dW^i_s\,,
$$
where $(H_i)$ is a basis of $\G$ and $W^i_t=(H_i,W_t)$. Hence
$w\mapsto p(w)$ defines a mapping from $C_\G$ into $C_G$ and we denote
by $\nu$ the image of $\mu$ under this mapping. Similarly, if $h\in H$
then we denote by $e(h)$  the
solution of the following differential equation:
\begin{equation}
\label{e-defn}
e_t(h)=e+\int_0^t e_s(h)\dot{h}_s ds\,.
\end{equation}
\begin{theorem}[Campbell-Baker-Hausdorff Formula]
\label{C-H-formula}
\index[sub]{Campbell-Baker-Hausdorff Formula}
For any $h\in H$ the following identity is valid
almost surely:
\be
\label{C-H-eqn}
p(w+h)=e(\tA p(w)h)p(w)\,,
\ee
where $\tA p(w)h$ is the $H$-valued random variable defined by
$$
\left(\tA p(w)h\right)(t)=\int_0^t \Ad  p_s(w)\dot{h}(s)ds\,.
$$
\end{theorem}
\index[not]{ad@${\tA}p(w)h$}
\remark
In case we work with matrices, $\tA p(w)h$ is defined as
$$
\int_0^t p_s(w)\dot{h}(s)p_s^{-1}(w)ds\,.
$$
\remark This theorem implies in particular that the $C_G$-valued
random variable $w\to p(w)$ has a modification, denoted again by the
same letter $p$, such that $h\to p(w+h)$ is a smooth function of $h\in
H$ for any $w\in C_\G$.
\paragraph{Calculation of $\nabla (f(p_t(w)))$:}
We have $f(p_t(w+\lambda h))=f(e_t(\tA p \lambda h)p_t)$ where
$e_t(h), \,h\in H$ is defined by the equation (\ref{e-defn}) . Let us write
$g=p_t(w)$ and $F(x)=f(xg)$. Then
$$
F(e_t(\lambda\tA ph))=F(e)+\lambda\int_0^t \Ad p_s\dot{h}_s
F(e_s(\lambda\tA ph)) ds\,.
$$
Hence
$$
\frac{d}{d\lambda}F(e_t(\tA p \lambda h))|_{\lambda=0}=
       \int_0^t \Ad p_s\dot{h}(s)F(e)ds\,.
$$
Now if $X$ is a left invariant vector field on $G$, then we have
$XF(x)=X(f(xg))=X(f(gg^{-1}xg))=(\Ad g^{-1}X)f(gx)$ by the left invariance
of $X$. In particular, for $x=e$, we have $XF(e)=(\Ad g^{-1}X)f(g)$. Replacing
$g$ with $p_t(w)$ above, we obtain
\begin{eqnarray}
\nabla_h(f(p_t))&=&\Ad p_t^{-1}\int_0^t\Ad p_s \dot{h}_s f(p_t)ds\\
    &=&\left( \Ad p_t^{-1}\int_0^t\Ad p_s \dot{h}_s ds\right) f(p_t)\,.
\end{eqnarray}
\paragraph{Notation:} In the sequel, we shall denote the map
$h\mapsto \int_0^{\cdot}\Ad p_s\dot{h}(s)ds $ by $\theta_p h$ or
by $\tA ph$ as before, depending on the notational convenience.
\begin{definition}
If $F:C_G\ra \R$ is a cylindrical function, $h\in H$, we define
\begin{eqnarray}
L_hF(p)&=&\frac{d}{d\lambda}F(e(\lambda h)p)|_{\lambda=0}\\
R_hF(p)&=&\frac{d}{d\lambda}F(p e(\lambda h))|_{\lambda=0}\,,
\end{eqnarray}
where $p$ is a generic point of $C_G$. $L$ is called the left derivative and
$R$ is called the right derivative.
\end{definition}

A similar calculation as above gives us
\begin{eqnarray}
L_hf(p_t)&=&\Ad p_t^{-1}h_t f(p_t)\\
R_hf(p_t)&=&h_tf(p_t)\,.
\end{eqnarray}
If $F(p)=f(p_{t_1},\cdots,p_{t_n})$, then
\begin{eqnarray}
L_hF(p)&=&\sum_{i=1}^n \Ad p_{t_i}^{-1}h_{t_i}f(p_{t_1},\cdots,p_{t_n})\\
R_hF(p)&=&\sum_{i=1}^nh_{t_i}f(p_{t_1},\cdots,p_{t_n})\\
\nabla_h(F\circ p(w))&=&\sum_{i=1}^n
      \Ad p_{t_i}^{-1}(w) \theta_{p(w)}h_{t_i}f(p_{t_1},\cdots,p_{t_n})(w)\,.
\end{eqnarray}

\begin{proposition}
$L_h$ is a closable operator on $L^p(\nu)$ for any $p>1$ and $h\in H$.
Moreover, we have
$$
(L_h F)(p(w))=\nabla_{\theta^{-1}_{p(w)}(h)}(F( p(w)))\,.
$$
\end{proposition}
\proof
Suppose that $(F_n)$ is a sequence of cylindrical functions on $C_G$
converging to zero in $L^p(\nu)$ and that $(L_hF_n)$ is Cauchy in
$L^p(\nu)$. Then, from the formulas  (7) and (9), we have
$$
(L_h F_n)(p(w))=\nabla_{\theta^{-1}_{p(w)}(h)}(F_n( p(w)))\,,
$$
since $\nabla$ is a closed operator on $L^p(\mu)$, we have necessarily
$\lim_n L_hF_n=0$ $\nu$-almost surely.
\qed
\begin{remarkk}{\rm
On the cylindrical functions we have the identity
$$
R_hF(p(w))=\nabla_{m(h)}(F(p(w)))
$$
where $m(h)_t=\Ad p_t(w)\int_0^t\Ad p_s^{-1}\dot{h}(s)ds$, but this process
is not absolutely continuous with respect to $t$, consequently, in general,
the right derivative is not a closable operator without further hypothesis on
the structure of $G$, we will come back to this problem later.}
\end{remarkk}
\begin{remarkk}{\rm
\label{devrem}
While working with matrix groups (i.e., the linear case) we can
also define all these in an alternative way (cf. also \cite{G1})
\begin{eqnarray*}
L_hF(p)&=&\frac{d}{d\lambda}F(e^{\lambda h}p)|_{\lambda=0}\\
R_hF(p)&=&\frac{d}{d\lambda}F(p\,e^{\lambda h})|_{\lambda=0}\,,
\end{eqnarray*}
where $e^{h}$ is defined (pointwise) as $e^h(t)=e^{h(t)}$. The advantage of
this definition is that  the right derivative commutes with the right
multiplication (however, as we will see later the corresponding Radon-Nikodym
derivative is more complicated):
$$
\frac{d}{d\lambda}F(p\,e^{\lambda h})=R_hF(p\,e^{\lambda h})\,,
$$
almost surely. Let us also  note the following identity which can be easily
verified on the cylindrical functions:
$$
\frac{d}{d\lambda}F(pe(\lambda h))=R_{\theta_{e(\lambda h)}h}F(p\,e(\lambda
h))\,,
$$
where
$\theta_{e(h)}k\in H$ is defined as
$$
\theta_{e(h)}k(t)=\int_0^t \Ad e_s(h){\dot{k}}(s)ds\,.
$$}
\end{remarkk}
\begin{remarkk}{\rm
On the extended domain of $L$, we have the identity
\begin{eqnarray}
L_hF\circ p(w)&=&\nabla_{\theta^{-1}_{p(w)}(h)}(F\circ p)\\
             &=&(\theta^{-1\star}_{p(w)}\nabla (F\circ p),h)\\
             &=&(\theta_{p(w)}\nabla (F\circ p),h)
\end{eqnarray}
if we assume that the scalar product of $\G$ is invariant with respect to the
inner automorphisms, in which case $G$ becomes of compact type, hence linear,
i.e., a space of matrices and $\theta_p$ becomes an isometry of $H$.}
\end{remarkk}
\begin{proposition}
If $\eta:C_G\ra H$ is a measurable random variable, then we have
\begin{eqnarray*}
(L_{\eta}F)\circ p&=&\nabla_{\theta^{-1}_{p(w)}(\eta\circ p)}(F\circ p)\\
     &=&(\theta_p\nabla(F\circ p),\eta\circ p)\,.
\end{eqnarray*}
\end{proposition}
\proof
By definition, $F\in \Dom(L)$ iff $F\circ p\in \Dom(\nabla)$ and in this case
$h\mapsto L_hF$ induces an $H$-valued random variable, denoted by $LF$. Then,
for any complete orthonormal basis  $(h_i,i\in \N)$  of $H$
\begin{eqnarray*}
L_\eta F\circ p&=&\sum_iL_{h_i}F\circ p(\eta ,h_i)\circ p\\
           &=&\sum_i \nabla_{\theta^{-1}h_i}(F\circ p)(\eta,h_i)_H\circ p\\
         &=&\sum_i \nabla_{\theta^{-1}h_i}(F\circ p)
            (\theta^{-1}\eta\circ p,\theta^{-1}h_i)_H\\
         &=&\nabla_{\theta^{-1}_p(\eta\circ p)}(F\circ p)\\
        &=&\Bigl(\theta_p\nabla(F\circ p),\eta\circ p\Bigr)_H
\end{eqnarray*}
\qed

\section{The left divergence $L^\star$}
\markboth{Stochastic Analysis on Lie Groups}{Left Divergence}
If $\eta:C_G\ra H$ is a cylindrical random variable and if $F$ is a smooth
function  on $C_G$, we have
\begin{eqnarray*}
E_\nu[L_\eta F]&=&E_\mu[(L_\eta F)\circ p]\\
              &=&E_\mu[\nabla_{\theta^{-1}(\eta\circ p)}(F\circ p)]\\
              &=&E_\mu[F\circ p\,\delta(\theta^{-1}(\eta\circ p))]\,.
\end{eqnarray*}
Since $L$ is a closed operator, its adjoint with respect to $\nu$ is
well-defined and we have
\begin{eqnarray*}
E_\nu[L_\eta F]&=&E_\eta[F\,L^\star \eta]\\
              &=&E_\mu[F\circ p\,(L^\star \eta)\circ p]\,.
\end{eqnarray*}
We have
\begin{proposition}\label{ldiv}
The following identity is true:
$$
(L^\star\eta)\circ p=\delta(\theta^{-1}(\eta\circ p))\,.
$$
\end{proposition}
\proof
We have already tested this identity for cylindrical $\eta$ and
$F$. To complete the proof it is sufficient to prove that the
cylindrical $F$ are dense in
$L^p(\nu)$. Then the proof will follow from the closability of
$L$. The density follows from the fact that $(p_t;t\in [0,1])$ and the
Wiener process generate
the same sigma algebra and from the monotone class theorem.
\qed
\begin{lemma}
\label{cond-exp}
Let $(\H_t,t\in [0,1])$ be the filtration (eventually completed) of
the process $(p_t,t\in [0,1])$ and $(\F_t,t\in [0,1])$ be the
filtration of the basic Wiener process. We have
$$
E_\nu[\phi|\H_t]\circ p=E_\mu[\phi\circ p|\F_t]
$$
$\mu$-almost surely.
\end{lemma}
\proof
Let $f$ be a smooth function on $\R^n$. Then
\begin{eqnarray*}
E_\mu[\phi\circ p\,f(p_{t_1}(w),\ldots,p_{t_n}(w))]&=&
                                   E_\nu[\phi\,f(p_{t_1},\ldots,p_{t_n})]\\
 &=&E_\nu[E_\nu[\phi|\H_t]f(p_{t_1},\ldots,p_{t_n})]\\
&=&E_\mu[E_\nu[\phi|\H_t]\circ p f(p_{t_1}(w),\ldots,p_{t_n}(w))]\,,
\end{eqnarray*}
since $E_\nu[\phi|\H_t]\circ p$ is $\F_t$-measurable, the proof follows.
\qed

If $F$ is a nice random variable on $C_G$ and denote by $\pi$ the
optional projection with respect to $(\F_t)$. Using Ito-Clark
representation theorem,
  we have
\begin{eqnarray*}
F\circ p &=& E_\mu[F\circ p]+\delta\left[\pi\nabla(F\circ p)\right]\\
      &=& E_\nu[F]+\delta\left[\theta_p\theta_p^{-1}\pi\nabla(F\circ
        p)\right]\\
    &=& E_\nu[F]+\delta\left[\theta_p^{-1}\pi\theta_p \nabla(F\circ p)\right]\\
    &=& E_\nu[F]+\delta\left[\theta_p^{-1}\pi(LF\circ p)\right]\\
    &=& E_\nu[F]+\delta\left[\theta_p^{-1}({\tilde{\pi}} LF)\circ p)\right]\\
     &=& E_\nu[F]+(L^\star({\tilde{\pi}}LF))\circ p
\end{eqnarray*}
$\mu$-almost surely, where $\tilde{\pi}$ denotes the optional projection
with respect to the filtration $(\H_t)$. Consequently, we have proved the
following
\begin{theorem}
Suppose that $F\in L^p(\nu),\,p>1$ such that $F\circ p\in D_{p,1}$. Then
we have
$$
F=E_\nu[F]+L^\star{\tilde{\pi}}LF
$$
$\nu$-almost surely.
\end{theorem}

\section{Ornstein-Uhlenbeck operator and the Wiener chaos}
\markboth{Stochastic Analysis on Lie Groups}{Wiener Chaos}
Let $F$ be a nice function on $C_G$, then
\begin{eqnarray}
(L^\star L F)\circ p&=&L^\star (LF)\circ p\\
               &=&\delta\left[\theta_p^{-1}(LF\circ p)\right]\\
           &=&\delta\left[\theta^{-1}\theta(\nabla(F\circ p))\right]\\
          &=&\delta\nabla (F\circ p)\\
         &=&\L(F\circ p)\,,
\end{eqnarray}
where $\L=\delta\nabla$ is the Ornstein-Uhlenbeck operator on $W$.
\begin{definition}
We denote by $\K$ the operator $L^\star L$ and call it the Ornstein-Uhlenbeck
operator on $C_G$.
\end{definition}
Let $F$ be a cylindrical function on $G$, for $t\geq 0$, define $Q_tF(p)$ as
$$
Q_tF(p(w))=P_t(F\circ p)(w)\,,
$$
where $P_t$ is the Ornstein-Uhlenbeck semigroup on $C_\G$, i.e.,
$$
P_tf(w)=\int_{C_{\G}} f(e^{-t}w+\sqrt{1-e^{-2t}}y)\mu(dy)\,.
$$
Then it is easy to see that
$$
\frac{d}{dt}Q_tF(p)|_{t=0}=-\K F(p)\,.
$$
Hence we can define the spaces of distributions, verify Meyer inequalities,
etc. ,  as in the flat case (cf. \cite{P-U}): Let $\phi$ be an equivalence
class of random variables on $(C_G,\nu)$ with values in some separable Hilbert
space $X$. For $q>1,\,k\in \N$, we will say that
$\phi$ is in $S_{q,k}(X)$, if there exists a sequence of cylindrical functions
$(\phi_n)$ which converges to $\phi$ in $L^q(\nu,X)$ such that
$(\phi_n\circ p)$ is Cauchy in $\DD_{q,k}(X)$. For $X=\R$, we write simply
$S_{q,k}$ instead of $S_{q,k}(\R)$. We denote by $S(X)$ the projective limit of
the spaces $(S_{q,k};\,q>1,k\in \N)$. Using Meyer inequalities and the fact
that $w\mapsto p(w)$ is smooth in the Sobolev sense, we can show easily that,
for $q>1,\,k\in \Z$
\begin{enumerate}
\item the left derivative $L$ possesses a continuous extension from
  $S_{q,k}(X)$ into $S_{q,k-1}^-(X\otimes H)$, where
$$
S_{q,k}^-(X)=\bigcup_{\epsilon>0}S_{q-\epsilon,k}(X)\,.
$$
\item $L^*$ has a continuous extension as a map from $S_{q,k}(X\otimes H)$ into
  $S_{q,k-1}^-(X)$.
\item Consequently $L$ maps $S(X)$ continuously into $S(X\otimes H)$ and $L^*$
maps $S(X\otimes H)$ continuously into $S(X)$.
\item By duality, $L$ and $L^*$ have continuous extensions, respectively, from
  $S'(X)$ to $S'(X\otimes H)$ and from $S'(X\otimes H)$ to $S'(X)$.
\end{enumerate}
We can now state the $0-1$ law as
a corollary:
\begin{proposition}
Let $A\in \B(C_G)$ such that $A=e(h)A$ $\nu$-almost surely for any $h\in H$,
then $\nu(A)=0$ or $1$.
\end{proposition}
\proof
It is easy to see that $L_h\won_A=0$ (in the sense of the
distributions) for any $h\in H$,  hence,
from Theorem \ref{0-1-law}, we obtain
$$
\won_A=\nu(A)
$$
almost surely.
\qed

Using  the calculations  above we obtain
\begin{proposition}
We have the following identity:
$$
\L^n(F\circ p)=(\K^nF)\circ p
$$
$\mu$-almost surely.
\end{proposition}
\paragraph{Notation:}In the sequel we will denote by $\tau$ the operator
$\theta_p(w)$ whenever $p(w)$  is replaced by the generic trajectory $p$ of
$C_G$.

Let $F$ be a cylindrical function on $C_G$. We know that
$$
F\circ p=E_\mu[F\circ p]+\sum_{i=1}^\infty \frac{1}{n!}
     \delta^n E_\mu[\nabla^n(F\circ p)]\,.
$$
On the other hand
$$
\nabla(F\circ p)=\theta^{-1}(LF\circ p)=(\tau^{-1}LF)\circ p
$$
$\mu$-almost surely. Iterating this identity, we obtain
$$
\nabla^n(F\circ p)=((\tau^{-1}L)^nF)\circ p\,.
$$
Therefore
\begin{eqnarray}
E_\mu[\nabla^n(F\circ p)]&=&E_\mu[((\tau^{-1}L)^nF)\circ p]\\
                       &=&E_\nu[(\tau^{-1}L)^nF]\,.
\end{eqnarray}
On the other hand, for $K$ in $H^{{\hat{\otimes}} n}$ (i.e.,  the symmetric
tensor product), we have
\begin{eqnarray*}
E_\mu[\delta^nK H\circ p]&=&E_\mu\left[(K,\nabla^n(H\circ p))_n\right]\\
    &=&E_\mu\left[(K,((\tau^{-1}L)^nH)\circ p)_n\right]\\
   &=&E_\nu\left[(K,(\tau^{-1} L)^nH)_n\right]\\
  &=&E_\nu\left[(L^\star\tau)^nK \,H\right]\\
  &=&E_\mu\left[((L^\star\tau)^nK)\circ p \, H\circ p\right]\,,
\end{eqnarray*}
for any cylindrical function $H$ on $C_G$, where $(\cdot,\cdot)_n$ denotes
the scalar product in $H^{{\hat{\otimes}} n}$. We have proved the identity
$$
\delta^n K=((L^\star\tau)^n K)\circ p\,,
$$
consequently the following Wiener decomposition holds:
\begin{theorem}
For any $F\in L^2(\nu)$, one has
$$
F=E_\nu[F]+\sum_{n=1}^\infty\frac{1}{n!}(L^\star\tau)^n
\left(E_\nu[(\tau^{-1}L)^n F]\right)\,
$$
where the sum converges in $L^2$.
\end{theorem}

The Ito-Clark representation theorem suggests us a second kind of Wiener chaos
decomposition. First we need the following:
\begin{lemma}
The set
$$
\Psi=\left\{\exp\left(L^\star h-\frac{1}{2}|h|_H^2\right);\,h\in H\right\}
$$
is dense in $L^p(\nu)$ for any $p\geq 1$.
\end{lemma}
\proof
We have
\begin{eqnarray*}
L^\star h\circ p&=&\delta (\theta_{p(w)}^{-1}h)\\
              &=&\int_0^1(\Ad p_s^{-1}{\dot{h}}(s), dW_s)\\
              &=&\int_0^1({\dot{h}}(s),\Ad p_s dW_s)\,.
\end{eqnarray*}
By Paul L\'evy's theorem, $t\mapsto B_t=\int_0^t\Ad p_s dW_s$ defines a
Brownian motion. Hence, to prove the lemma, it suffices to show that
$W$ and $B$ generate the same filtration. To see this, note that the process
$(p_t)$ satisfies the following stochastic differential equation:
$$
df(p_t)=H_if(p_t)\,dW^i_t\,,
$$
($f\in C^\infty(G)$), replacing $dW_t$ by $\Ad p_t dB_t$ we obtain
$$
df(p_t)=\Ad p_t^{-1}H_if(p_t)\,dB^i_t\,.
$$
Since everything is smooth, we see that $p(w)$ is measurable with respect to
the filtration of $B$. But we know that the filtrations of $p$ and $W$ are
equal from the lemma \ref{cond-exp}.
\qed
\begin{remarkk}{\rm
 Using the Brownian motion  $B_t$ defined above we can also represent
the Wiener functionals, this gives another Wiener chaos
decomposition.}
\end{remarkk}

\section{Some useful formulea}
\markboth{Stochastic Analysis on Lie Groups}{Some Useful Formulea}
Let us first recall the variation  of constant method for matrix-valued
equations:
\begin{lemma}\label{varconst}
The solution of the equation
$$
\beta_t( h)=\Phi_t+\int_0^t \beta_s( h){\dot{h}}(s)ds
$$
is given by
$$
\beta_t( h)=\Phi_0+
\left(\int_0^t \frac{d}{ds}\Phi_s e_s( h)^{-1}ds\right)e_t( h)\,.
$$
\end{lemma}
\begin{corollary}
We have
\begin{eqnarray}
\frac{d}{d\lambda}e_t(\lambda h)&=
       &\left(\int_0^t \Ad e_s(\lambda h){\dot{h}}(s)ds\right)e_t(\lambda h)\\
  &=&(\theta_{e(\lambda h)}h)(t)e_t(\lambda h)\,.
\end{eqnarray}
\end{corollary}
\begin{corollary}
We have
$$
\frac{d}{d\lambda}\Ad e_t(\lambda h){\dot{k}}_t=
\left[\int_0^t\Ad e_s(\lambda h){\dot{h}}_sds,\Ad
e_t(\lambda h){\dot{k}}_t\right]\,.
$$
\end{corollary}
\begin{corollary}
We have
$$
\frac{d}{d\lambda}\Ad e^{-1}_t(\lambda h){\dot{k}}_t=-\Ad e_t^{-1}(\lambda h)
\left[\int_0^t\Ad e_s(\lambda h){\dot{h}}_sds,{\dot{k}}_t\right]\,.
$$
\end{corollary}
\proof
Since $\Ad e \Ad e^{-1}=I$, we have
\begin{eqnarray*}
0&=&\frac{d}{d\lambda}\Ad e_t(\lambda h)\Ad e^{-1}_t(\lambda h){\dot{k}}_t\\
&=&\left(\frac{d}{d\lambda}\Ad e_t(\lambda h)\right)
\Ad e^{-1}_t(\lambda h){\dot{k}}_t+\Ad e_t(\lambda
h)\frac{d}{d\lambda}\Ad e^{-1}_t(\lambda h){\dot{k}}_t\,,
\end{eqnarray*}
hence
\begin{eqnarray*}
\frac{d}{d\lambda}\Ad e^{-1}_t(\lambda h){\dot{k}}_t&=&
-\Ad e^{-1}_t(\lambda h)\left(\frac{d}{d\lambda}\Ad e_t(\lambda
  h)\right)\Ad e^{-1}_t(\lambda h){\dot{k}}_t\\
&=&-\Ad e_t^{-1}(\lambda h)\left[\int_0^t\Ad e_s(\lambda h){\dot{h}}_s ds,
\Ad e_t(\lambda h)\Ad e_t^{-1}(\lambda h){\dot{k}}_t\right]\\
&=&-\Ad e_t^{-1}(\lambda h)
\left[\int_0^t\Ad e_s(\lambda h){\dot{h}}_sds,{\dot{k}}_t\right]\,.
\end{eqnarray*}
\qed

In further calculations we shall need to control the terms like
$$
|\Ad e_t^{-1}(v){\dot{h}}_t-\Ad e_t^{-1}(\alpha){\dot{h}}_t|_\G\,.
$$
For this, we have
\begin{eqnarray*}
\Ad e_t^{-1}(v){\dot{h}}_t-\Ad e_t^{-1}(\alpha){\dot{h}}_t&=&
 \int_0^1\frac{d}{d\lambda}\Ad e^{-1}_t(\lambda v+(1-\lambda)\alpha){\dot{h}}_t
   d\lambda\\
  &=&\int_0^1 \frac{d}{d\lambda}\Ad
  e^{-1}_t(\lambda(v-\alpha)+\alpha){\dot{h}}_t d\lambda\,.
\end{eqnarray*}
From the Corollary 6.3, we have
$$
\frac{d}{d\lambda}\Ad e^{-1}_t(\lambda (v-\alpha)+\alpha){\dot{h}}_t=
$$
$$
\qquad -\Ad e^{-1}_t(\lambda(v-\alpha)+\alpha)
\left[\int_0^t\Ad
  e_s(\lambda(v-\alpha)+\alpha)({\dot{v}}_s-{\dot{\alpha}}_s)ds,{\dot{h}}_t\right]\,.
$$
Therefore
$$
|\Ad e_t^{-1}(v){\dot{h}}_t-\Ad e_t^{-1}(\alpha){\dot{h}}_t|_\G \leq
$$
$$
\qquad \int_0^1\left|\left[\int_0^t\Ad e_s(\lambda (v-\alpha)+\alpha)
({\dot{v}}_s-{\dot{\alpha}}_s) ds, {\dot{h}}_t\right]\right|_\G d\lambda\,.
$$
Now we need to control the $\G$-norm of the Lie brackets: for this we introduce
some notations: let $(e_i)$ be a complete, orthonormal basis  of
$\G$. Since $[e_i,e_j]\in \G$ we should have
$$
[e_i,e_j]=\sum_{k=1}^n\gamma_{ij}^ke_k\,.
$$
For $h,k\in \G$,
\begin{eqnarray*}
[h,k]&=&\left[\sum_i h_ie_i,\sum_i k_ie_i\right]\\
    &=&\sum_{i,j}h_i k_j[e_i,e_j]\\
   &=&\sum_{i,j,k}h_i k_j\gamma_{i,j}^k\,.
\end{eqnarray*}
Consequently
\begin{eqnarray*}
|[h,k]|_\G^2&=&\sum_l\left[\sum_{i,j}h_ik_j\gamma_{ij}^l\right]^2\\
       &\leq & \sum_l\left(\sum_{i,j}h_i^2 k_j^2\right)
            \left(\sum_{i,j}(\gamma_{ij}^l)^2\right)\\
      &=&\sum_l |h|_\G^2|k|_\G^2 |\gamma^l|_2^2\\
     &=& |h|_\G^2|k|_\G^2\sum_l |\gamma^l|_2^2\,,
\end{eqnarray*}
where $|\cdot|_2$ refers to the Hilbert-Schmidt norm on $\G$. Although this is
well-known, let us announce the above result as a lemma for later reference:
\begin{lemma}
\label{brac-con}
For any $h,k\in \G$, we have
$$
|[h,k]|_\G\leq |h|_\G|k|_\G\left(\sum_l |\gamma^l|^2_2\right)^{1/2}\,.
$$
\end{lemma}
We have also the immediate consequence
\begin{lemma}
For any $h,k\in H$
$$
\left|\Ad e^{-1}_t(v){\dot{h}}_t-\Ad e^{-1}_t(\alpha){\dot{h}}_t\right|_\G
\leq \|\gamma\|_2 |{\dot{h}}_t|_\G\int_0^t|{\dot{v}}_s-{\dot{\alpha}}_s|_\G ds,
$$
where $\|\gamma\|_2^2=\sum |\gamma^l|_2^2$.
\end{lemma}

\begin{lemma}
We have
$$
\frac{d}{d\lambda}\phi(e(\lambda h)p)=
\left(L\phi(e(\lambda h)p),\tA  e(\lambda h)h\right)_H\,.
$$
\end{lemma}
\proof
We have
$$
e_t(ah)e_t(bh)=e_t(a \tA e^{-1}(bh)h+bh)\,,
$$
hence
\begin{eqnarray*}
e_t(ah+bh)&=&e_t(a\tA e^{-1}(bh)h)e_t(bh)\\
          &=&e_t(b\tA e^{-1}(ah)h)e_t(ah)\,,
\end{eqnarray*}
therefore
$$
e_t((\lambda+\mu)h)=e_t(\mu \tA e(\lambda h)h)e_t(\lambda h)\,,
$$
which gives
$$
\frac{d}{d\lambda}\phi(e(\lambda h)p)=
\left(L\phi(e(\lambda h)p),\tA  e(\lambda h)h\right)_H\,.
$$
\qed

\section{Right derivative}
\markboth{Stochastic Analysis on Lie Groups}{Right Derivative}
Recall that we have defined
$$
R_h\phi(p)=\frac{d}{d\lambda}\phi(p\,e(\lambda h))|_{\lambda=0}\,.
$$
Since $\G$ consists of left invariant vector fields, we have, using
the global notations :
$$
R_hf(p_t)=(h_t f)(p_t)\,,
$$
where $h_tf$ is the function obtained by applying the vector field
$h_t$ to the smooth function $f$. The following is straightforward:
\begin{lemma}
\label{comm1-lemma}
We have
$$
p_t(w)e_t(h)=p_t\left(\int_0^{\cdot}\Ad e_s^{-1}(h) dW_s+h\right)\,,
$$
where $e_t(h)$ for $h\in H$ is defined in (\ref{e-defn}).
\end{lemma}
\begin{lemma}
We have
$$
E_\mu[R_h F\circ p]=E_\mu[F\circ p\,\int_0^1 {\dot{h}}_s dW_s]\,,
$$
for any cylindrical function $F$.
\end{lemma}
\proof From the Lemma \ref{comm1-lemma}, $p_t(w)e_t(\lambda
h)=p_t(\lambda h+\int_0^\cdot \Ad e^{-1}_s(\lambda h)dW_s)$. Since
$\int_0^\cdot \Ad e^{-1}_s(\lambda h)dW_s$ is a Brownian motion, it
follows from the Girsanov theorem that
$$
E\left[F(p(w)e(\lambda h))
\exp\left\{-\lambda
\int_0^1({\dot{h}}_s, \Ad e^{-1}_s(\lambda
h)dW_s)-\frac{\lambda^2}{2}|h|_H^2\right\}\right]
=E[F]\,,
$$
differentiating at $\lambda=0$ gives the result.
\qed
\begin{definition}
For $h\in H$ and $F$ smooth, define
\begin{itemize}
\item
$Q_h F(w)$ by
$$
Q_hF(w)=F\left(\int_0^\cdot \Ad e^{-1}_s(h) dW_s\right)\,,
$$
note that since $\int_0^\cdot \Ad e^{-1}_s(h) dW_s$ is a Brownian motion, the
composition of it with $F$ is well-defined.
\item
And
$$
X_hF(w)=\left.\frac{d}{d\lambda}Q_{\lambda h}F(w)\right|_{\lambda=0}\,.
$$
\end{itemize}
\end{definition}

\begin{example}
{\rm
Let us see how  the derivation operator $X_h$ operates  on  the simple
functional  $F=\exp \delta k$, $k\in H$: we have
\begin{eqnarray*}
Q_{\lambda h}F&=&\exp\int_0^1({\dot{k}}_s,\Ad e^{-1}_s(\lambda h)dW_s)\\
    &=&\exp\int_0^1(\Ad e_s(\lambda h){\dot{k}}_s,dW_s)\,,
\end{eqnarray*}
hence
$$
X_h e^{\delta k}=e^{\delta k}\int_0^1\left([h(s),{\dot{k}}_s],dW_s\right)\,.
$$}
\end{example}
\begin{proposition}
We have the following identity:
$$
(R_hF)\circ p=\nabla_h(F\circ p)+X_h(F\circ p)\,,
$$
for any $F:C_G\ra \R$ smooth. In particular, $R_h$ and $X_h$ are closable
operators.
\end{proposition}
\begin{remarkk}{\rm
From the above definition, we see that
$$
(R_h^2)^\star 1=
\delta^2 h^{\otimes 2}-\int_0^1\left([h_s,{\dot{h}}_s], dW_s\right)\,.
$$
Hence $R^{\star n}$ does not give the pure chaos but mixes them with those
of lower order. Here enters the notion of universal envelopping algebra.}
\end{remarkk}
{\bf{Notation : }}For $h\in H$, we will denote by $\tad h$ the linear  operator
on $H$ defined as
$$
\tad h(k)(t)=\int_0^t[h(s),{\dot{k}}(s)]\,ds\,.
$$
\begin{remarkk}{\rm
Suppose that $R_h\delta k=0$, i.e.,
$$
(h,k)+\int_0^1 [h_s,{\dot{k}}_s]\cdot dW_s=0\,.
$$
Then $(h,k)=0$ and $[h(t),{\dot{k}}(t)]=0$ $dt$-almost surely. Hence this
gives more information than the independence of  $\delta h$ and $\delta k$.}
\end{remarkk}
\begin{remarkk}{\rm
Suppose that $R_hF=0$ a.s. for any $h\in H$. Then we have, denoting
$F=\sum I_n(f_n)$, $R_hF=0$ implies
$$
 n f_n(h)+d\Gamma(\tad h)f_{n-1}=0\,,\,\,k\in H\,.
$$
Since $f_1=0$ (this follows from $E[R_hF]=E[\nabla_hF]=0$), we find that
$f_n(h)=0$ for any $h\in H$, hence $f_n=0$, and $F$ is a constant.}
\end{remarkk}
\begin{remarkk}{\rm
If $X_hF=0$ for any $h\in H$, we find that
$$
 d\Gamma(\tad h)f_{n}=0
$$
for any $h\in H$ and for any $n$. Therefore $f_n$'s take their values in
the tensor spaces constructed from the center of $\G$.}
\end{remarkk}

Recall that in the case of an abstract Wiener space, if $A$ is a deterministic
operator on the Cameron-Martin space $H$, then the operator $d\Gamma(A)$
is defined  on the Fock as
$$
d\Gamma(A)\phi=\frac{d}{dt}\Gamma(e^{tA})\phi|_{t=0}
$$
for any cylindrical Wiener functional $\phi$. We will need the following result
which is well-known in the Quantum Field Theory  folklore:
\begin{lemma}
\label{comm-lemma}
Suppose that $A$ is a skew-symmetric operator on $H$ (i.e., $A+A^*=0$). Then we
have
$$
d\Gamma(A)\phi=\delta A \nabla\phi\,,
$$
for any $\phi\in \cup_{p>1}D_{p,2}$.
\end{lemma}
\proof
By a density argument, it is sufficient to prove the identity for the
functionals $\phi=\exp[\delta h-1/2|h|_H^2]\,,\,h\in H$. In this case we have
\begin{eqnarray*}
\Gamma(e^{tA})\phi&=&\exp\left\{\delta
  e^{tA}h-\frac{1}{2}|e^{tA}h|_H^2\right\}\\
  &=&\exp\left\{\delta e^{tA}h-\frac{1}{2}|h|_H^2\right\}
\end{eqnarray*}
where the last equality follows from the fact that $e^{tA}$ is an isometry of
$H$. Hence, by differentiation, we obtain
$$
d\Gamma(A)\phi=\delta(Ah)\phi\,.
$$
On the other hand
\begin{eqnarray*}
\delta A \nabla \phi&=&\delta\left[Ah\,\,e^{\delta
    h-\frac{1}{2}|h|_H^2}\right]\\
    &=&\left[\delta(Ah)-(Ah,h)_H\right]e^{\delta h-\frac{1}{2}|h|_H^2}\\
    &=&\delta(Ah)e^{\delta h-\frac{1}{2}|h|_H^2}\,,
\end{eqnarray*}
since $(Ah,h)_H=0$.
\qed

As a corollary, we have
\begin{corollary}
For any cylindrical function $F$ on $(C_\G,H,\mu)$, we have the following
commutation relation:
$$
\left[\nabla_h,X_k\right]F=-\nabla_{{\tad} k(h)}F\,,
$$
where $h,\,k\in H$.
\end{corollary}
We have also
\begin{proposition}
Let $\phi$ be a cylindrical function on $(C_\G,H,\mu)$ and $h\in H$. We have
$$
E_\mu[(X_h\phi)^2]\leq \|\gamma\|_2^2|h|_H^2
     E\left\{|\nabla \phi|_H^2+\|\nabla^2 \phi\|_2^2\right\}\,,
$$
where $\gamma$ is the structure constant of $\G$ and $\|\cdot\|_2$ denotes the
Hilbert-Schmidt norm of $H\otimes H$.
\end{proposition}
\proof
From Lemma \ref{comm-lemma}, we have
$
X_h\phi=\delta\left({\tilde{\ad}}h(\nabla \phi)\right)\,.
$
Hence
$$
E\left[(X_h\phi)^2\right]=E[|{\tilde{\ad}}h\nabla\phi|_H^2]+
E\left[{\mbox{trace}}\left(\nabla{\tilde{\ad}}h\nabla \phi\right)^2\right]\,.
$$
From Lemma \ref{brac-con}, we have
$$
\Big|{\tilde{\ad}}h\nabla\phi\Big|_H^2\leq
\|\gamma\|_2^2|h|_H^2|\nabla \phi|_H^2
$$
and
$$
\left|{\mbox{trace}}\left(\nabla{\tilde{\ad}}h\nabla \phi\right)^2\right|
\leq \|\gamma\|_2^2|h|_H^2\|\nabla^2 \phi\|_2^2\,.
$$
\qed

Suppose that $u\in D(H)$ and define $X_u F$, where $F$ is a cylindrical
function on $C_\G$, as $\delta {\tilde{\ad}}u\nabla F$. Then using similar
calculations, we see that
\begin{corollary}
We have the following majoration:
\begin{eqnarray*}
  E[|X_uF|^2]&\leq& \|\gamma\|^2  E\left[|u|_H^2|\nabla F|_H^2\right]\\
 &+&2\|\gamma\|^2E\left[|u|_H^2\|\nabla^2 F\|_2^2+\|\nabla u\|_2^2|\nabla F|_H^2\right]\,.
\end{eqnarray*}
\end{corollary}

\section{Quasi-invariance}
\markboth{Stochastic Analysis on Lie Groups}{Quasi-invariance}
Let $\gamma_t$ be a curve in $G$ such that $t\mapsto \gamma_t$ is absolutely
continuous. We can write it as
\begin{eqnarray*}
d\gamma_t&=&{\dot{\gamma}}_t dt\\
        &=&\gamma_t \gamma^{-1}_t{\dot{\gamma}}_t dt
\end{eqnarray*}
Hence $\gamma_t=e_t(\int_0^\cdot \gamma^{-1}_s{\dot{\gamma}}_s ds)$ provided
$\int_0^1|\gamma^{-1}_t{\dot{\gamma}}_t|^2 dt<\infty$. Under these hypothesis,
we have
$$
\gamma_t\,p_t(w)=p_t\left(w+\int_0^\cdot
{\mbox{Ad}}p_s^{-1}(w)(\gamma_s^{-1}{\dot{\gamma}}_s)ds\right)\,.
$$
For any cylindrical $\phi: G\rightarrow \R$, we have
$$
E_\nu\left[\phi(\gamma\,p) \,J_\gamma\right]=E_\nu[\phi]
$$
where
$$
J_\gamma\circ p(w)=\exp\left\{-\int_0^1
({\mbox{Ad}}p_s^{-1}(\gamma_s^{-1}{\dot{\gamma}}_s),dW_s)-\frac{1}{2}\int_0^1
|\gamma^{-1}_s{\dot{\gamma}}_s|^2 ds\right\}\,.
$$
Similarly
$$
p_t(w)\gamma_t=p_t\left(\int_0^\cdot{\mbox{Ad}}\gamma_s^{-1}dW_s+\int_0^\cdot
\gamma_s^{-1}{\dot{\gamma}}_s ds\right)\,,
$$
hence
$$
E_\nu\left[\phi(p\,\gamma)\,K_\gamma \right]=E_\nu[\phi]
$$
where
\begin{eqnarray}
K_\gamma\circ
p(w)&=&\exp\left
\{-\int_0^1(\gamma^{-1}_s{\dot{\gamma}}_s,{\mbox{Ad}}\gamma^{-1}_s
  dW_s)-\frac{1}{2}\int_0^1|\gamma_s^{-1}{\dot{\gamma}}_s|^2ds\right\}\nonumber\\
&=&\exp\left\{-\int_0^1({\dot{\gamma}}_s\gamma^{-1}_s,dW_s)-\frac{1}{2}\int_0^1|
\gamma_s^{-1}{\dot{\gamma}}_s|^2ds\right\}\,.
\end{eqnarray}
As an application of these results, let us choose $\gamma=e^h$ and denote by
$K_h$ the Radon-Nikodym density defined by
$$
E_\nu[F(pe^{ h})]=E[F\,K_{ h}]\,.
$$
Since $\lambda \mapsto K_{\lambda h}$ is analytic, from Remark
\ref{devrem}, for smooth, cylindrical $F$, we have
\begin{eqnarray*}
E[F(p\,e^{\lambda h})]&=&\sum_{n=0}^\infty\frac{\lambda^n}{n!}E[R_h^nF(p)]\\
              &=&\sum_{n=0}^\infty\frac{\lambda^n}{n!}E[F(p)\,R_h^{n\star}1]\,,
\end{eqnarray*}
hence we have the identity
$$
K_{\lambda h}=\sum_{n=0}^\infty \frac{\lambda^n}{n!}R_h^{n\star}1\,.
$$
Let us now choose $F(p)$ of the form $f(p_1)$, where $f$ is a smooth
function on $G$. Then
\begin{eqnarray*}
E[f(p_1e^{\lambda h(1)})]&=&\sum_{n=0}^\infty
\frac{\lambda^n}{n!}E[R_h^nf(p_1)]\\
    &=&\sum_{n=0}^\infty \frac{\lambda^n}{n!}E[h(1)^nf(p_1)]\,.
\end{eqnarray*}
Let $q(x)dx$ be the law of $p_1$ where $dx$ is the right invariant Haar measure
on $G$. Then
\begin{eqnarray*}
E[h(1)^nf(p_1)]&=&\int_G h(1)^n f(x)\,\,q(x)dx\\
              &=&(-1)^n\int_Gf(x)\,\frac{h(1)^nq(x)}{q(x)}q(x)dx\,.
\end{eqnarray*}
Hence we have proved
\begin{proposition}
We have the following identity:
$$
E_\nu[K_h|p_1=x]=\sum_{n=0}^\infty
\frac{(-1)^n}{n!}\,\,\frac{h(1)^nq(x)}{q(x)}\,,
$$
for all $x\in G$. In particular, if $h(1)=0$ then
$$
E_\nu[K_h|p_1=x]=1
$$
for all $x\in G$.
\end{proposition}
\proof
The only claim to be justified is ``all $x$'' instead of almost all $x$. This
follows from the fact that $x\mapsto E_\nu[K_h|p_1=x]$ is continuous due to the
non-degeneracy of the random variable $p_1$ in the sense of the Malliavin
calculus.
\qed

\noindent
Although the analogue of the following result is obvious in the flat case, in
the case of the Lie groups, the proof requires more work:
\begin{proposition}
The span of $\{K_h;\,h\in H\}$ is dense in $L^r(\nu)$ for any $r>1$.
\end{proposition}
\proof
Let us denote by $\Theta$ the span of the set of the densities. Suppose that
$F\in L^s$ with $E_\nu[F]=0$, where $s$ is the conjugate of $r$, is orthogonal
to $\Theta$. In the sequel we shall denote again by $F$ the random variable
defined as $w\mapsto F\circ p(w)$. From the orthogonality hypothesis, we have
$E[R^n_hF]=0$ for any $h\in H$ and $n\in \N$ ( we have not made any
differentiability hypothesis about $F$ since all  these calculations are
interpreted in the distributional sense). For $n=1$, this gives
\begin{eqnarray*}
0&=&E_\mu[\nabla_hF+X_hF]\\
&=&E_\mu[\nabla_hF]\,,
\end{eqnarray*}
since $X_h+X_h^*=0$.
For $n=2$
\begin{eqnarray*}
0&=&E_\mu[R_h^2F]\\
&=&E_\mu[\nabla^2_hF+X_h\nabla_hF+\nabla_h X_hF+X_h^2 F]\\
&=&E_\mu[\nabla_h^2F]+E_\mu[\nabla_h X_hF]\,.
\end{eqnarray*}
Also we have from the  calculations of the first order
\begin{eqnarray*}
E_\mu[\nabla_h X_hF]&=&E_\mu[X_hF \delta h]\\
              &=&-E_\mu[F\delta(\tad h(h))]\\
             &=&-E_\mu[(\nabla F,\tad h(h))_H]\\
             &=&0\,.
\end{eqnarray*}
By polarization, we deduce that, as
a tensor in $H^\otimes$, $E_\mu[\nabla^2F]=0$. Suppose now that
$E_\mu[\nabla^iF]=0$ for $i\leq n$. We have
$$
E_\mu[R^{n+1}_hF]=E_\mu[\nabla^{n+1}_hF]+{\mbox{ supplementary
    terms}}\,.
$$
Between these supplementary terms, those who begin with $X_h$ or
its powers have automatically zero expectation. We can show via induction
hypothesis that the others are also null. For instance let us take the term
$E_\mu[\nabla_hX_h^nF]$:
\begin{eqnarray*}
 E_\mu\left[\nabla_hX_h^nF\right]&=&E_\mu[X_h^nF\,\delta h]\\
                     &=&(-1)^n E_\mu[F\, \delta((\tad h)^nh)]\\
                     &=&0\,,
\end{eqnarray*}
the other terms can be treated similarly.
\qed

We shall apply these results to the loop measure by choosing a special
form of $\gamma$. Let us first explain the strategy: replace in the
above expressions the random variable
$\phi(p(w))$ by $\phi\circ p(w) f(p_1(w))$. Then we have
\begin{itemize}
\item
$$
E_\mu\Bigl[\phi(\gamma p(w))\,f(\gamma(1)p_1)J_\gamma\circ p(w)\Bigr]=
E_\mu\left[\phi\circ p(w)f(p_1(w))\right]
$$
and
\item
$$
E_\mu\Bigl[\phi( p(w)\gamma)\,f(p_1\gamma(1))K_\gamma\circ p(w)\Bigr]=
E_\mu[\phi\circ p(w)f(p_1(w))]\,.
$$
\end{itemize}
 We shall proceed as follows: let $f:G \rightarrow \R$ be a smooth
cylindrical function. Replace in the above expressions the map $\phi\circ p$ by
$\phi\circ p \,f(p_1(w))$ where $f$ is a smooth function  on $G$. Then we have
 on the  one hand
\begin{equation}
\label{G-1}
E_\nu\left[\phi(\gamma
  p)\,f(\gamma(1)p_1)\,J_\gamma\right]=E_\nu\left[\phi(p)\,f(p_1)\right]
\end{equation}
and on the other hand
\begin{equation}
\label{G-2}
E_\nu\left[\phi(p \gamma
  )\,f(p_1\gamma(1))\,K_\gamma\right]=E_\nu\left[\phi(p)\,f(p_1)\right]\,.
\end{equation}
Choose $\gamma$ such that $\gamma(1)=e$ (i.e., the identity of $G$).
Hence (\ref{G-1}) becomes
$$
E_\nu[\phi(\gamma p)\,f(p_1)\,J_\gamma]=E[\phi(p)\,f(p_1)]\,,
$$
therefore
$$
\int_G E_\nu\left[\phi(\gamma p)J_\gamma(p)\Big|p_1=x\right]f(x) q_1(x) dx =
\int_G E_\nu\left[\phi(p)\Big|p_1=x\right]f(x) q_1(x) dx\,,
$$
where $dx$ is the Haar measure on $G$ and $q_1$ is the density of the law
of $p_1$ with respect to  Haar measure which is smooth and strictly positive.
Consequently we obtain
$$
E_\nu\Bigl[\phi(\gamma
  p)J_\gamma(p)\Big|p_1=x\Bigr]=E_\nu\Bigl[\phi(p)\Big|p_1=x\Bigr]\,.
$$
Since both sides are continuous with respect to $x$, this equality holds
everywhere. We obtain a similar result also for the right perturbation
using the relation (\ref{G-2}).

A natural candidate for $\gamma$ for the loop measure based at $e$,
i.e.,  for  the measure $E_\nu[\,\cdot\,|p_1=e]$ which we will denote by $E_1$,
 would be
$$
\gamma_t(h)=e_t(h)e_1^{-1}(th)\,.
$$
From the calculations of the sixth  section, we have
$$
{\dot{\gamma}}_t(h)=e_t(h)[{\dot{h}}_te^{-1}_1(th)-e_1^{-1}(th)
(\theta_{e(th)}h)(1)]\,.
$$
Hence
\begin{lemma}
For $\gamma_t(h)=e_t(h)e_1^{-1}(th)$, we have
$$
\gamma_t^{-1}(h){\dot{\gamma}}_t(h)=
\Ad e_1(th){\dot{h}}_t-(\theta_{e(th)}h)(1)\,.
$$
\end{lemma}
In this case $J_\gamma$ becomes
\begin{eqnarray*}
J_\gamma\circ p=\exp&-&\int_0^1\left(\Ad p_s^{-1}[\Ad e_1(sh){\dot{h}}_s-
(\theta_{e(sh)}h)(1)],dW_s\right)\\
\exp &-&\frac{1}{2}\int_0^1|\Ad e_1(sh){\dot{h}}_s-(\theta_{e(sh)}h)(1)|^2ds\,.
\end{eqnarray*}

For $K_\gamma$ we have
\begin{eqnarray*}
\Ad \gamma_t (\gamma^{-1}_t{\dot{\gamma}}_t)&=&{\dot{\gamma}}_t\gamma^{-1}_t\\
&=&\Ad e_t(h){\dot{h}}_t-\Ad (e_t(h)e_1^{-1}(th))(\theta_{e(th)}h)(1)\,.
\end{eqnarray*}
Since $|\cdot|$ is $\Ad$-invariant, we have
\begin{eqnarray*}
K_\gamma\circ p=\exp&-&\int_0^1\left(\Ad e_t(h){\dot{h}}_t-\Ad
  (e_t(h)e_1^{-1}(th))
(\theta_{e(th)}h(1)),dW_t\right)\\
\exp&-&\int_0^1|{\dot{h}}_t-\Ad e^{-1}_1(th)(\theta_{e(th)}h(1))|^2dt\,.
\end{eqnarray*}
\begin{remarkk}{\rm
Note that $\gamma$ as chosen above satisfies the following differential
equation:
$$
{\dot{\gamma}}_t=\gamma_t(h)[\Ad e_1(th){\dot{h}}_t-\theta_{e(th)}h(1)]\,.
$$}
\end{remarkk}
\noindent
Let us calculate 
$$
\frac{d}{d\lambda}\phi(p\gamma(\lambda h))|_{\lambda=0} {\mbox{ and }}
\frac{d}{d\lambda}\phi(\gamma(\lambda h)p)|_{\lambda=0}
$$ 
for cylindrical
$\phi$. Denote by $P_0: H\rightarrow H_0$ the orthogonal projection defined by
$$
P_0h(t)=h(t)-t h(1)\,.
$$
Then it is easy to see that
$$
\frac{d}{d\lambda}\phi(\gamma(\lambda h)p)|_{\lambda=0}=L_{P_0h}\phi(p)
$$
and
$$
\frac{d}{d\lambda}\phi(p\gamma(\lambda h))|_{\lambda=0}=R_{P_0h}\phi(p)\,.
$$
Moreover, we have
$$
\frac{d}{d\lambda}J_{\gamma(\lambda h)}(p(w))|_{\lambda=0}=
-\int_0^1\left(\Ad p_s^{-1}(w)({\dot{h}}_s-h(1)),dW_s\right)
$$
and
$$
\frac{d}{d\lambda}K_{\gamma(\lambda h)}(p)|_{\lambda=0}=
-\int_0^1\left({\dot{h}}_s-h(1),dW_s\right)
$$
Consequently we have proven
\begin{theorem}
For any cylindrical function $\phi$ on the loop space of $G$, we have
$$
E_1[L_{P_0h}\phi]=E_1[\phi\,L^{\star}P_0h]
$$
and
$$
E_1[R_{P_0h}\phi]=E_1[\phi\, \delta P_0h]
$$
for any $h\in H$. In particular, the operators $L_{P_0h}$ and $R_{P_0h}$ are
closable on $L^p(\nu_1)$ for any $p>1$.
\end{theorem}
Before closing this section  let us  give a result of L. Gross (cf. \cite{G1}):
\begin{lemma}
For $\alpha<1$ the measure $\nu(\cdot\,|p(1)=e)$ is equivalent to $\nu$ on
$(C_G,\H_\alpha)$ and for any $\H_\alpha$-measurable random variable $F$, we
have
$$
E_\nu[F|p(1)=e]=
 E_\nu\left[F\,\frac{q_{1-\alpha}(p_\alpha,e)}{q_1(e,e)}\right]\,,
$$
where $q_t$ is the density of the law of $p_t$ with respect to the Haar
measure.
\end{lemma}
\proof
Without loss of generality we can suppose that $F$ is a continuous and bounded
function on $C_G$. Let $g$ be a nice function on $G$, from the Markov
property, it follows that
$$
E_\nu[F\,g(p(1))]=E\left[F\int_G q_{1-\alpha}(p_\alpha,y)g(y)dy\right]\,.
$$
On the other hand, from the disintegration of measures, we have
$$
E_\nu[F\,g(p(1))]=\int_G E_\nu[F|p(1)=y]g(y)q_1(e,y)dy\,.
$$
Equating both sides gives
$$
E_\nu[F|p(1)=y]=\frac{1}{q_1(e,y)}E_\nu[F\,q_{1-\alpha}(p_\alpha,y)]
$$
$dy$-almost surely. Since both sides are continuous in $y$ the result follows
if we put $y=e$.
\qed

\begin{remarkk}{\rm
Note that we have the following identity:
$$
E_\nu[F(p)|p(1)=e]=E_\mu[F\circ p(w)|p_1(w)=e]
$$
for any cylindrical function $F$ on $C_G$.}
\end{remarkk}

\section{Anticipative transformations}
\markboth{Stochastic Analysis on Lie Groups}{Anticipative Transformations}
In this section we shall study the absolute continuity of the measures which
are defined as the image of $\nu$ under the mappings which are defined as the
left multiplication of the path $p$ with the exponentials of anticipative
$\G$-valued processes. To be able to use the results of the flat case we need
to extend the Campbell-Baker-Hausdorff formula to this case. We begin
by recalling  the following
\begin{definition}
Let $(W,H,\mu)$ be an abstract Wiener space. A random variable $F$ defined on
this space is said to be of class $R^p_{\alpha,k}$ if $F\in D_{q,r}$ for some
$q>1$, $r\geq 1$ and
$$
\sup_{|h|_H\leq \alpha }|\nabla^k F(w+h)|\in L^p(\mu)\,.
$$
\begin{itemize}
\item If $p=0$, we write $F\in R^0_{\alpha,k}$ ,
\item We write $F\in R^p_{\infty,k}$ if the above condition holds for any
  $\alpha >0$, and $F\in R^p_{\infty,\infty}$ if $F\in R^p_{\infty,k}$ for any
  $k\in \N$.
\item Finally, we say that $F\in R(\infty)$ if $F\in R_{\infty,\infty}^p$ for
  any $p>1$.
\end{itemize}
\end{definition}
\begin{remarkk}{\rm
The importance of this class is easy to realize: suppose that $u$ is an
$H$-valued random variable, and let $F\in R_{\infty,\infty}^0$. If $(u_n)$ is a
sequence of random variables of the form $\sum_{i<\infty} h_i\won_{A_i}$
converging in probability to $u$, with $A_i\cap A_j=\emptyset$ for $i\neq j$,
$h_i\in H$, then we can define $F(w+u_n(w))$ as $\sum
F(w+h_i)\won_{A_i}(w)$ and
evidently the sequence $(F\circ T_n)$ converges in probability, where
$T_n=I_W+u_n$. Furthermore, the limit is independent of the particular choice
of the elements of the equivalence class of $u$. Moreover, if  we  choose a
sequence approximating $F$ as $F_n=E[P_{1/n}F|V_n]$, where $(h_n)$ is a
complete basis of $H$, $V_n$ is the sigma algebra generated by $\delta
h_1,\cdots,\delta h_n$ and $P_{1/n}$ is the Ornstein-Uhlenbeck semigroup at
$t=1/n$, then $\sup_n\sup_{|h|\leq \alpha}|\nabla^k F_n(w+h)|<\infty $ almost
surely  for any
$\alpha>0$, $k\in \N$, and  we can show, using an equicontinuity argument
(cf. \cite{UZ7}) that the limit of $(F\circ T_n)$  is measurable with
respect to the sigma algebra of $T=I_W+u$. }
\end{remarkk}
\begin{lemma}
\label{H-reg-lemma}
For any $t\geq 0$, the random variable $w\mapsto p_t(w)$ belongs to the class
$R(\infty)$. Consequently, for any $H$-valued random variable $u$,
the random variable $w\mapsto p_t(w+u(w))$ is well-defined and it is
independent of the choice of the elements of the equivalence class of $u$.
\end{lemma}
\proof
In fact in \cite{UZ7}, p.175, it has been proven that any diffusion with smooth
coefficients of compact support belongs to $R_{\infty,\infty}^0$. In our
particular case it is easy to see that
$$
\sup_{|h|_H\leq\alpha} \|\nabla^k p_t(w+h)\|\in \cap_p L^p(\mu)
$$
for any $\alpha >0$ and $k,n\in \N$, where $\|\cdot\|$ is the Euclidean  norm
on $M(\R^n)\otimes H^{\otimes k}$ and $M(\R^n)$ denotes the space of linear
operators on $\R^n$.
\qed
\begin{lemma}
\label{div-com-lemma}
Suppose that $\xi\in R^\alpha_{\infty,\infty}(H)\cap D(H)$, and $\delta\xi\in
R^\alpha_{\infty,\infty}$ and that $u\in D(H)$ with $|u|_H\leq a\leq \alpha$
almost surely. Denote by $T$ the mapping $I_W+u$, then we have
$$
(\delta\xi)\circ T=\delta(\xi\circ T)+(\xi\circ
T,u)_H+{\mbox{trace}}(\nabla\xi\circ T\cdot \nabla u)\,,
$$
almost surely.
\end{lemma}
\proof
Let $(e_i)$ be a complete, orthonormal basis in $H$, denote by $V_k$ the sigma
algebra generated by $\{\delta e_1,\cdots,\delta e_k\}$, by $\pi_k$ the
orthogonal projection of $H$ onto the vector space generated by $\{e_1,\cdots,
e_k\}$. Let $u_n$ be defined as $E[\pi_nP_{1/n}u|V_n]$, then $|u_n|_H\leq a$
almost surely again. From the finite dimensional Sobolev injection theorem one
can show that the map $\phi\mapsto \phi\circ T_n$ is continuous from $D$ into
itself  and we have
$$
\nabla(\phi\circ T_n)=(I+\nabla u_n)^*\nabla \phi\circ T_n
$$
(cf. \cite{UZ7}). For $\xi$ as above, it is not difficult to show the
claimed identity, beginning first with a cylindrical $\xi$ then passing to the
limit with the help of the continuity of the map $\phi\mapsto \phi\circ T_n$.
To pass to the limit with respect to $n$, note that we have
$$
|\delta\xi\circ T_n-\delta\xi\circ T|\leq
\sup_{|h|_H\leq \alpha}|\nabla\delta\xi(w+h)|_H|u_n(w)-u(w)|_H\,,
$$
and, from the hypothesis, this sequence converges to zero in all the $L^p$
spaces. For the other terms we proceed similarly.
\qed
\begin{theorem}
\label{iden-thm}
Let $u$ be in $D_{q,1}(H)$ for some $q>1$, then we have
$$
p_t\circ T(w)=e_t(\theta_p u)p_t(w)\,,
$$
where $e_t(\theta_p u)$ is the solution of the ordinary differential equation
given by ${\dot{e}}_t=e_t \Ad p_t {\dot{u}}_t$.
\end{theorem}
\proof Suppose first that $u$ is also bounded.  From  Lemma
\ref{H-reg-lemma}, $p_t$ belongs to $R(\infty)$
hence the same thing is also true for the Stratonovitch  integral
$\int_0^tp_sdW_s$. We can write the Stratonovitch integral as the sum of the
Ito integral of $p_s$ plus $\frac{1}{2}\int_0^tCp_s ds$, where $C$ denotes the
Casimir operator (cf. \cite{Fe}). Since
 $\sup_{r\leq t}|e_r(\theta_p h)|\leq \exp t|h|_H$, $t\mapsto p_t\circ
 T$ is almost surely continuous. Moreover, it is not
difficult to see that $\int_0^t C p_sds$ is in $R(\infty)$. Hence we can
commute the Lebesgue integral with the composition with $T$.  Consequently we
have, using Lemma \ref{div-com-lemma},
\begin{eqnarray*}
\left(\int_0^tp_sdW_s\right)\circ T=& &\int_0^tp_s\circ T\,\delta W_s+\int_0^t
\frac{1}{2}Cp_s \circ T\, ds\\
  &+&\int_0^tp_s\circ T\,{\dot{u}}_sds \\
  &+&\int_0^t\int_0^s( D_rp_s)\circ T \,D_s{\dot{u}}_rdrds
\end{eqnarray*}
where  $\delta W_s$ denotes
the Skorohod integral and $D_s \phi$ is the notation for the  Lebesgue
density of the $H$-valued
random variable $\nabla \phi$.   We can write this expression simply as
$$
\left(\int_0^t p_sdW_s\right)\circ T=\int_0^tp_s\circ T d^\circ
W_s+\int_0^tp_s\circ T{\dot{u}}_sds\,,
$$
where $d^\circ W_s$ represents the anticipative Stratonovitch 
integral, i.e., we add the trace term to the divergence, whenever it is
well-defined. Therefore we obtain the relation
$$
p_t\circ T=e+\int_0^tp_s\circ T d^\circ W_s+\int_0^tp_s\circ T
{\dot{u}}_sds\,. 
$$
 Let us now develop $e_t(\theta_pu)p_t(w)$ using the Ito
formula for anticipative processes (cf. \cite{U5}):
\begin{eqnarray*}
e_t(\theta_pu)p_t(w)&=&e+\int_0^te_s(\theta_pu)p_s(w) d^\circ W_s+
  \int_0^te_s(\theta_pu) \Ad p_s {\dot{u}}_sp_sds\\
&=&e+\int_0^te_s(\theta_pu)p_s(w) d^\circ W_s+
  \int_0^te_s(\theta_pu)  p_s {\dot{u}}_s ds\,.
\end{eqnarray*}
Hence, both $p_t\circ T$ and $e_t(\theta_p u)p_t$ satisfy the same anticipative
stochastic differential equation with the obvious unique solution, therefore
the proof is completed for the case where $u$ is bounded. To get rid of the
boundedness hypothesis, let $(u_n)$ be a sequence in $D_{q,1}(H)$ converging to
$u$ (with respect to $(q,1)$-Sobolev norm) such that $|u_n|_H\leq 2n+1$ and
$u_n=u$ on the set $\{w:\,|u(w)|_H\leq n\}$. Then from the bounded case, we
have $p_t(w+u_n(w))=e_t(\theta_p u_n)(w)p_t(w)$ almost surely. Moreover both
sides of this equality converge in probability respectively to $p_t\circ T$ and
$e_t(\theta_pu)p_t$ and the proof is completed.
\qed

The following results now follow immediately from the flat case and
Theorem \ref{iden-thm}: using the change of variable formula
for the anticipative shifts on the abstract Wiener spaces (cf. \cite{UZ2}),
we can prove
\begin{theorem}
Suppose that $u: C_G\rightarrow H$ be a random variable such that
\begin{enumerate}
\item $\|Lu\|_{L^\infty(\nu,H\otimes H)}<\infty$,
\item $\|\,\|Lu\|_{op}\,\|_{L^\infty(\nu)}\leq c <1$, where $c$ is a fixed
constant.
\end{enumerate}
Then we have
$$
E_\nu\left[F(e(\theta_pu(p))p)\,|J_u|\right]=E_\nu[F]
$$
for any $F\in C_b(C_G)$, where
$$
J_u=\dett(I_H+\theta_p^{-1}Lu)\exp-L^\star(\theta_p u)-
            \frac{1}{2}|u|^2\,.
$$
\end{theorem}
\proof
Let us denote by $u'(w)$ the random variable $u\circ p$ which is defined on
$W=C([0,1],\G)$. From Campbell-Baker-Hausdorff  formula, we have
$$
p(w+u'(w))=e(\theta_{p(w)}u'(w))p(w)
$$
 (in fact here we are dealing with
anticipative processes but the calculations go as if the things were adapted
thanks to the Stratonovitch integral which defines the trajectory $p$). We know
from \cite{UZ2} that
$$
E_\mu[F(p(w+u'(w))\,|\Lambda_{u'}|]=E_\mu[F(p(w))]
$$
where
$$
\Lambda_{u'}=\dett(I_H+\nabla u'(w))\exp-\delta
u'-\frac{1}{2}|u'|^2\,.
$$
To complete the proof it suffices to remark that
\begin{eqnarray*}
\nabla u'(w)&=&\nabla(u\circ p(w))\\
           &=&\theta_p^{-1}\,Lu\circ p(w)\\
\delta u'(w)&=&\delta(u\circ p)(w)\\
           &=&L^\star(\theta_p\,u)\circ p(w)\,.
\end{eqnarray*}
\qed

\noindent
We shall  observe first  the based loop space case.
We need  the following notations: if $\gamma(t)$ is an absolutely
continuous curve with values in $G$,
we  denote by $\ka(\gamma)$ the curve with values in $\G$ defined by
$$
\ka(\gamma)_t=\int_0^t\gamma_s^{-1}{\dot{\gamma}}_sds\,,
$$
where we use, as before, the matrix notation.
\begin{theorem}
Suppose that $\gamma:[0,1]\times C_G\rightarrow G$ be a random variable which
is absolutely continuous with respect to $dt$ and that $\gamma(0)=\gamma(1)=e$,
where $e$ denotes the unit element of $G$. Suppose moreover that
\begin{enumerate}
\item $\|L\theta_p^{-1} \ka(\gamma)\|_{L^\infty(\nu,H\otimes H)}<\infty$,
\item $\|\,\|L\theta_p^{-1} \ka(\gamma)\|_{op}\|_{L^\infty(\nu)}\leq c <1$,
\item $J_\gamma \in S_{r,1}$ for some $r>1$, where $S_{r,1}$ is the Sobolev
  space on $C_G$ which consists of the completion of the cylindrical
  functionals with respect to the norm
  $\|\phi\|_{r,1}=\|\phi\|_{L^r(\nu)}+\|L\phi\|_{L^r(\nu,H)}$.
\end{enumerate}
Then we have
$$
E_1\left[F(\gamma(p)p)\,|J_\gamma|\right]=E_1[F]
$$
for any $F\in C_b(C_G)$, where
$$
J_\gamma=\dett(I_H+\theta_p^{-1}L\theta_p^{-1}\ka(\gamma))
  \exp\left\{-L^\star \ka(\gamma)-\frac{1}{2}|\ka(\gamma)|^2\right\}\,.
$$
\end{theorem}
\proof
It is sufficient to take $u=\theta_p^{-1}\ka(\gamma)$ in the preceding theorem
and then apply the usual conditioning trick to obtain
$$
E_\nu[F(\gamma(p)p)|J_\gamma||p(1)=y]=E_\nu[F|p(1)=y]
$$
$dy$-almost surely.  Note that by the hypothesis, there is some $q>1$ such that
$J_\gamma\circ p$ belongs to the  Sobolev space $\DD_{q,1}$
and $\varepsilon_e\circ p(1)$ ($\varepsilon_e$ denotes the Dirac measure
at $e$) belongs to $\bigcap_s \DD_{s,-1}$ (cf. \cite{W}), hence
both sides of the above equality are continuous with respect to $y$ and the
proof follows.
\qed

\subsection{Degree type results}
\markboth{Stochastic Analysis on Lie Groups}{Degree Type Results}
In this section we will give some straight-forward applications of the measure
theoretic degree theorem on the flat Wiener space to the path and loop spaces
on the Lie group $G$. The following theorem is a direct consequence of
the results of the preceding section  and the degree theory in the
flat case (cf. \cite{UZ5,UZ6},  \cite{UZ7} and Theorem \ref{degree-thm}):
\begin{theorem}
\label{lie-deg-thm}
Let $\gamma:[0,1]\times C_G\rightarrow G$ be a random variable which is
absolutely continuous with respect to $dt$ and that $\gamma(0)=e$. Suppose
moreover that, for some $a>0$,
\begin{enumerate}
\item $J_\gamma\in L^{1+a}(\nu)$,
\item $J_\gamma\left(I_H+\theta_p^{-1}L\theta_p^{-1}\kappa(\gamma)\right)h\in
  L^{1+a}(\nu)$, for any $h\in H$,
\item $\kappa(\gamma)\in S_{r,2}(H)$, for some $r>\frac{1+a}{a}$,
  where $S_{r,2}$ is the Sobolev space of $H$-valued functionals as defined
   before.
\end{enumerate}
Then we have
$$
E_\nu[F(\gamma(p)p)J_\gamma]=E_\nu[F]E_\nu[J_\gamma]\,,
$$
for any $F\in C_b(C_G)$.
\end{theorem}

The following is a consequence of Theorem 3.2 of \cite{UZ6}:
\begin{proposition}
Suppose that $\kappa(\gamma)\in S_{q,1}(H)$ for some $q>1$ and that
$$
\exp\left(-L^*\theta_p^{-1}\kappa(\gamma)+1/2\|L\theta_p^{-1}\kappa(\gamma)\|_2^2\right)\in
L^{1+b}(\nu)\,,
$$
for some $b>1$. Then
$$
E_\nu[J_\gamma]=1\,.
$$
\end{proposition}
Let us look at the loop space case:
\begin{proposition}
Let $\gamma$ be as in Theorem \ref{lie-deg-thm}, with $\gamma(1)=e$
and suppose moreover that 
$J_\gamma\in S_{c,1}$, for some $c>1$. Then
$$
E_1\left[F(\gamma (p)p)J_\gamma\right]=E_1[F]E_\nu[J_\gamma]\,,
$$
for any smooth, cylindrical function $F$.
\end{proposition}
\proof
Let $f$ be a nice function on $G$. From Theorem 9.4, we have
\begin{eqnarray*}
 E_\nu[F(\gamma(p)p)f(p_1)J_\gamma]&=&
 E_\nu[F(\gamma(p)p)f(\gamma_1(p)p_1)J_\gamma]\\
    &=&E_\nu[F(p)f(p_1)]E_\nu[J_\gamma]
\end{eqnarray*}
hence
$$
E_\nu[F(\gamma(p)p)J_\gamma|p_1=y]=E_\nu[F(p)|p_1=y] E_\nu[J_\gamma]
$$
$dy$ almost surely. Since both sides are continuous with respect to $y$, the
equality remains true for every $y\in G$.
\qed
\begin{remarkk}{\rm
Note that the ``degree''  of $\gamma$, namely $E_\nu[J_\gamma]$ remains the
same in both  path and loop spaces.}
\end{remarkk}

\printindex[sub]
\printindex[not]

\end{document}